\documentclass[reqno]{amsart}
\usepackage[margin = 1.5in]{geometry}
\usepackage{amsmath, amssymb, amsthm, fancyhdr, verbatim, graphicx}
\usepackage{enumerate}
\usepackage[all]{xy}
\usepackage[usenames,dvipsnames]{xcolor}
\usepackage{mathrsfs}
\usepackage{tikz-cd}
\usepackage{framed, hyperref}
\usepackage[titletoc]{appendix}
\usepackage{bbm}

\usepackage{mathtools}
\mathtoolsset{showonlyrefs}

\usepackage{lipsum}
\usepackage{adjustbox}

\numberwithin{equation}{section}

\newcommand{\U}{\mathrm{U}}
\newcommand{\rd}{\mathrm{d}}

\newcommand{\iden}{\langle1\rangle}

\DeclareFontFamily{U}{matha}{\hyphenchar\font45}
\DeclareFontShape{U}{matha}{m}{n}{
      <5> <6> <7> <8> <9> <10> gen * matha
      <10.95> matha10 <12> <14.4> <17.28> <20.74> <24.88> matha12
      }{}
\DeclareSymbolFont{matha}{U}{matha}{m}{n}
\DeclareFontFamily{U}{mathx}{\hyphenchar\font45}
\DeclareFontShape{U}{mathx}{m}{n}{
      <5> <6> <7> <8> <9> <10>
      <10.95> <12> <14.4> <17.28> <20.74> <24.88>
      mathx10
      }{}
\DeclareSymbolFont{mathx}{U}{mathx}{m}{n}

\DeclareMathSymbol{\obot}         {2}{matha}{"6B}

\RequirePackage{xspace}

\newcommand{\BA}{\ensuremath{\mathbb{A}}\xspace}

\newcommand{\BC}{\ensuremath{\mathbb{C}}\xspace}

\newcommand{\sL}{\ensuremath{\mathscr{L}}\xspace}

\usepackage{color}

\newcommand\mnote[1]{\marginpar{\tiny #1}}

\newcommand{\F}{\mathbf{F}}

\newcommand{\CC}{\mathbf{C}}
\newcommand{\G}{\mathbf{G}}
\newcommand{\tr}[0]{\operatorname{tr}}
\newcommand{\wt}[1]{\widetilde{#1}}

\newcommand{\pderiv}[2]{\frac{\partial #1}{\partial #2}}
\newcommand{\Q}{\mathbf{Q}}
\newcommand{\Z}{\mathbf{Z}}

\newcommand{\sgn}{\operatorname{sgn}}
\newcommand{\Gal}{\operatorname{Gal}}

\newcommand{\ul}[1]{\underline{#1}}
\newcommand{\ol}[1]{\overline{#1}}
\newcommand{\wh}[1]{\widehat{#1}}

\newcommand{\Cal}[1]{\mathcal{#1}}
\newcommand{\A}{\mathbf{A}}

\newcommand{\mbf}[1]{\mathbf{#1}}

\newcommand{\ft}{{}^{\tau}} 
\newcommand{\co}{\colon}
\newcommand{\mrm}[1]{\mathrm{#1}}

\newcommand{\bs}{\backslash}

\newcommand{\PP}{\mathbf{P}}
\newcommand{\scr}[1]{\mathscr{#1}}

\newcommand{\inj}{\hookrightarrow}

\newcommand{\surj}{\twoheadrightarrow}

\DeclareMathOperator{\GL}{GL}

\DeclareMathOperator{\Frob}{Fr}
\DeclareMathOperator{\gFrob}{Frob}
\DeclareMathOperator{\Fr}{Fr}

\DeclareMathOperator{\coker}{coker}

\DeclareMathOperator{\N}{\mathbf{N}}
\DeclareMathOperator{\Tr}{Tr}

\DeclareMathOperator{\PGL}{PGL}
\DeclareMathOperator{\Sp}{Sp}
\DeclareMathOperator{\SO}{SO}
\DeclareMathOperator{\Hom}{Hom}
\DeclareMathOperator{\Ind}{Ind}

\DeclareMathOperator{\rank}{rank}

\DeclareMathOperator{\Aut}{Aut}
\DeclareMathOperator{\Rep}{Rep}
\DeclareMathOperator{\Nm}{Nm}
\DeclareMathOperator{\Spec}{Spec\,}

\DeclareMathOperator{\End}{End}
\DeclareMathOperator{\Isom}{Isom}

\DeclareMathOperator{\Res}{Res}

\DeclareMathOperator{\Div}{Div}

\DeclareMathOperator{\Stab}{Stab}
\DeclareMathOperator{\Bun}{Bun}
\DeclareMathOperator{\Ext}{Ext}
\DeclareMathOperator{\Pic}{Pic}

\DeclareMathOperator{\Id}{Id}

\DeclareMathOperator{\Gr}{Gr}

\DeclareMathOperator{\Sht}{Sht}

\DeclareMathOperator{\cla}{cl}

\DeclareMathOperator{\Fl}{Fl}

\DeclareMathOperator{\pr}{pr}

\DeclareMathOperator{\Hk}{Hk}

\DeclareMathOperator{\Eis}{Eis}
\DeclareMathOperator{\Ch}{Ch}

\DeclareMathOperator{\Herm}{Herm}
\DeclareMathOperator{\Lagr}{Lagr}
\DeclareMathOperator{\Coh}{Coh}
\DeclareMathOperator{\supp}{supp}
\DeclareMathOperator{\Den}{Den}

\DeclareMathOperator{\St}{St}

\DeclareMathOperator{\Prym}{Prym}
\DeclareMathOperator{\all}{all}

\newcommand{\KR}[1]{\Cal{Z}_{\Cal{#1}}}

\DeclareMathOperator{\Spr}{Spr}
\DeclareMathOperator{\HSpr}{HSpr}

\DeclareMathOperator{\ev}{ev}
\DeclareMathOperator{\Int}{Int}

\newcommand{\SH}{\textup{SkHm}}
\newcommand{\Ug}{\textup{U}}
\newcommand\Sig{\Sigma}
\newcommand{\bt}{\boxtimes}
\newcommand{\bi}{\mathbf{i}}

\newcommand{\incl}{\hookrightarrow}
\newcommand{\isom}{\stackrel{\sim}{\to}}
\newcommand{\bij}{\leftrightarrow}

\newcommand{\Qlbar}{\overline{\Q}_\ell}
\newcommand{\twtimes}[1]{\stackrel{#1}{\times}}

\renewcommand{\j}[1]{\langle{#1}\rangle}
\newcommand\mat[4]{\left(\begin{array}{cc} #1 & #2 \\ #3 & #4 \end{array}\right)}  
\newcommand\un{\underline}
\newcommand{\bu}{\bullet}
\newcommand{\ov}{\overline}

\newcommand\sss{\subsubsection}
\newcommand\xr{\xrightarrow}
\newcommand\op{\oplus}
\newcommand\ot{\otimes}
\newcommand\bigot{\bigotimes}
\newcommand\one{\mathbf{1}}

\newcommand{\oll}[1]{\overleftarrow{#1}}
\newcommand{\orr}[1]{\overrightarrow{#1}}

\renewcommand\c{\circ}
\newcommand\vn{\varnothing}

\newcommand{\cohog}[2]{\textup{H}^{#1}({#2})}     
\newcommand{\hBM}[2]{\textup{H}^{\textup{BM}}_{#1}({#2})}  

\newcommand\gl{\mathfrak{gl}}
\renewcommand\sp{\mathfrak{sp}}

\newcommand\Og{\mathrm{O}}

\renewcommand\a\alpha
\renewcommand\b\beta
\newcommand\g\gamma
\renewcommand\d\delta
\newcommand\D\Delta
\newcommand{\e}{\varepsilon}
\newcommand{\io}{\iota}
\renewcommand{\th}{\theta}

\newcommand{\ph}{\varphi}
\newcommand{\s}{\sigma}
\renewcommand{\t}{\tau}

\newcommand{\y}{\eta}
\newcommand{\z}{\zeta}
\newcommand{\ep}{\varepsilon}
\newcommand{\vp}{\varpi}
\renewcommand{\l}{\lambda}
\renewcommand{\L}{\Lambda}
\newcommand{\om}{\omega}

\renewcommand{\r}{\rho}

\newcommand\dm{\diamondsuit}
\newcommand\hs{\heartsuit}

\newcommand\sh{\sharp}
\newcommand\da{\dagger}

\newcommand\cA{\mathcal{A}}
\newcommand\cB{\mathcal{B}}

\newcommand\cE{\mathcal{E}}
\newcommand\cF{\mathcal{F}}
\newcommand\cG{\mathcal{G}}
\newcommand\cH{\mathcal{H}}
\newcommand\cI{\mathcal{I}}

\newcommand\cK{\mathcal{K}}
\newcommand\cL{\mathcal{L}}
\newcommand\cM{\mathcal{M}}
\newcommand\cN{\mathcal{N}}
\newcommand\cO{\mathcal{O}}
\newcommand\cP{\mathcal{P}}
\newcommand\cQ{\mathcal{Q}}
\newcommand\cR{\mathcal{R}}

\newcommand\cT{\mathcal{T}}

\newcommand\cV{\mathcal{V}}

\newcommand\cY{\mathcal{Y}}
\newcommand\cZ{\mathcal{Z}}

\newcommand\frD{\mathfrak{D}}

\newcommand\frI{\mathfrak{I}}
\newcommand\frJ{\mathfrak{J}}

\newcommand\frP{\mathfrak{P}}

\newcommand\frR{\mathfrak{R}}

\newcommand\frY{\mathfrak{Y}}

\newcommand\frg{\mathfrak{g}}

\newcommand\fm{\mathfrak{m}}

\newcommand\fro{\mathfrak{o}}
\newcommand\frp{\mathfrak{p}}

\newcommand\cMa{\overline{\mathcal{M}}}
\newcommand\dds{\frac{\mrm{d}}{\mrm{d}s}}

\newtheorem{thm}{Theorem}[section]
\newtheorem{lemma}[thm]{Lemma}
\newtheorem{prop}[thm]{Proposition}
\newtheorem{cor}[thm]{Corollary}

\newtheorem{defn-prop}[thm]{Definition-Proposition}
\newtheorem{claim}[thm]{Claim}

\theoremstyle{remark}
\newtheorem{remark}[thm]{Remark} 
\newtheorem{defn}[thm]{Definition}

\makeatletter
\def\th@remark{%
  \thm@headfont{\bfseries}%
  \normalfont 
  \thm@preskip \thm@preskip 
  \thm@postskip\thm@preskip
}
\def\imod#1{\allowbreak\mkern5mu({\operator@font mod}\,\,#1)}
\makeatother

\numberwithin{equation}{section}

\title[Higher Siegel--Weil: non-singular terms]{Higher Siegel--Weil formula
for unitary groups: \\ the non-singular terms}

\author{Tony Feng, Zhiwei Yun, Wei Zhang}

\begin{document}

\begin{abstract} We construct special cycles on the moduli stack of hermitian shtukas. We prove an identity between (1) the $r^{\rm th}$ central derivative of non-singular Fourier coefficients of a normalized Siegel--Eisenstein series, and  (2) the degree of special cycles of ``virtual dimension 0'' on the moduli stack of hermitian shtukas with $r$ legs. This may be viewed as a function-field analogue of the Kudla-Rapoport Conjecture, that has the additional feature of encompassing all higher derivatives of the Eisenstein series. 
\end{abstract}

\maketitle

\tableofcontents

\section{Introduction} 

The classical Siegel--Weil formula (\cite{Siegel1951, Weil1965}) relates the special values  of Siegel--Eisenstein series on the symplectic group (resp. the unitary group) to theta functions, which are generating series of representation numbers of quadratic (resp. Hermitian) forms over number fields. In particular, by exploiting the factorization of the non-singular Fourier coefficients into a product of local terms, one arrives at Siegel's formula for representation numbers of global quadratic or Hermitian forms in terms of local representation densities. 

In \cite{Kudla1997} Kudla began to study an arithmetic version of the Siegel--Weil formula and he discovered a relation between an ``arithmetic theta function'' --- a generating series of arithmetic cycles on an integral model of a Shimura curve---and the  first central derivative of a Siegel--Eisenstein series on $\Sp_4$. In a series of papers, Kudla and Rapoport developed this paradigm by defining the non-singular terms of a generating series of special cycles on suitable integral models of Shimura varieties for $\SO(n-1,2)$ with $n\leq 4$ and for all $\U(n-1,1)$. Of particular relevance to our paper, in \cite{KRI, KRII} Kudla and Rapoport defined the sought-after special cycles on integral models of unitary Shimura varieties, now known as \emph{Kudla--Rapoport cycles}, and conjectured a relationship to the non-singular Fourier coefficients  of the central derivative of the Siegel--Eisenstein series.  Their conjecture has been recently proved by Li and one of us \cite{LZ1}; we also refer to the introduction of \cite{LZ1} for a more detailed account of recent advances in some other related directions (see also \cite{LZ2} for the orthogonal analog). With the Kudla--Rapoport  conjecture proved in \cite{LZ1} and its archimedean counterpart proved by Liu \cite{Liu11}  and independently by Garcia and Sankaran \cite{GS19} as some of the key ingredients, Li and Liu \cite{LL1,LL2} have recently proved an arithmetic Rallis inner product formula relating the height pairing of the generating series to the first derivative of L-functions for unitary groups, from which they deduced cases of Beilinson--Bloch conjecture for certain high rank motives.

In this paper we study a function field analogue of the arithmetic Siegel--Weil formula, for unitary groups. 
In particular, we will construct special cycles on the moduli space of hermitian shtukas with arbitrary number of legs. Then we formulate and prove the analogue of the Kudla-Rapoport conjecture for derivatives of {\em arbitrary order} at the center of the Siegel--Eisenstein series, relating the {\em non-singular} Fourier coefficients of such higher derivatives to the degrees of special cycles. We remark that the proofs here follow a completely different strategy than in \cite{LZ1}.

In the sequel \cite{FYZ2}, we will construct the complete generating series of special cycles (including singular terms) and give evidence for their modularity.

\subsection{Statement of main result}

To formulate the result, let $X$ be a smooth, proper and geometrically connected curve over $k=\F_q$ of characteristic $p\ne2$, and $\nu \co X' \rightarrow X$ be an \'{e}tale double cover, with the non-trivial automorphism denoted $\sigma \in \Aut(X'/X)$. Let $F$ be the function field of $X$ and let $F'$ be the ring of rational functions on $X'$ (we allow $X'=X\coprod X$). In \S \ref{sec: unitary shtukas} we recall the definition of the moduli stack $\Sht_{U(n)}^r$ parametrizing rank $n$ ``Hermitian shtukas'' with $r$ legs.  Roughly speaking it classifies chains of vector bundles with Hermitian structures
\begin{equation}\label{intro Herm Sht}
\cF_{0}\dashrightarrow\cF_{1}\dashrightarrow\cdots \dashrightarrow \cF_{r}\cong {}^{\t}\cF_{0}
\end{equation}
related by elementary modifications. It admits a fibration $\Sht_{U(n)}^r \rightarrow (X')^r$, and will play the role of Shimura varieties in the function field context. 

\subsubsection{Special cycles}
Drawing inspiration from the construction of Kudla-Rapoport cycles on unitary Shimura varieties \cite{KRII}, we introduce in \S \ref{s:int prob} certain stacks $\Cal{Z}_{\cE}^r(a)$ over $\Sht_{U(n)}^{r}$ indexed by $\cE$, a vector bundle of rank $m$ with $1\leq m\leq n$ on $X'$, and a Hermitian map\footnote{A map of vector bundles of the form $a \co \cE \rightarrow \sigma^* \cE^\vee$ is \emph{Hermitian} if $\sigma^* a^\vee = a$.} $a \co \cE \rightarrow \sigma^* \cE^{\vee}$ where $\cE^{\vee}:= \un\Hom(\cE, \omega_{X'})$ is the Serre dual of $\cE$. They classify Hermitian shtukas as in \eqref{intro Herm Sht} together with compatible maps $\cE\to \cF_{i}$ such that $a$ is induced from the Hermitian form on $\cF_{0}$. 

If $\cE$ is a line bundle on $X'$,  $\Cal{Z}_{\cE}^r(a)$ is an analogue of Kudla-Rapoport divisors although they have dimension $r$ less than $\Sht_{U(n)}^{r}$. In general, $\cZ_{\cE}^{r}(a)$ are analogs of special cycles for function fields.

We will be particularly interested in the case $m=n$ and $a \co \cE \rightarrow \sigma^* \cE^{\vee}$ is injective (by this we shall always mean as a map of coherent sheaves). In this case, the ``virtual dimension'' of $\Cal{Z}_{\cE}^r(a)$ is $0$. However, as is already seen in the number field context \cite{KRII}, the literal dimension of $\Cal{Z}_{\cE}^r(a)$ is often significantly larger; this problem is exacerbated as $r$ increases. Nevertheless, under the assumption that $a \co \cE \rightarrow \sigma^* \cE^{\vee}$ is injective (as a map of coherent sheaves), we are able to construct an appropriate ``virtual fundamental cycle'' $[\Cal{Z}_{\cE}^r(a) ] \in \Ch_0(\Cal{Z}_{\cE}^r(a))_{\Q}$. Interestingly, it turns out that there are some new difficulties present in this construction that do \emph{not} appear in the Shimura variety setting. For $a$ injective, it turns out that $\Cal{Z}_{\cE}^r(a)$ is proper over $\F_{q}$, so that $[\Cal{Z}_{\cE}^r(a) ] $ has a well-defined degree $\deg[\Cal{Z}_{\cE}^r(a) ] \in \Q$.

\subsubsection{The main result}
Let $E(g,s, \Phi)$ be the Siegel--Eisenstein series for the standard split $F'/F$-skew-Hermitian space of dimension $2n$, with respect to the unramified standard section $\Phi$. For a rank $n$ vector bundle $\cE$ on $X'$ as above, $E(g,s, \Phi)$ admits a Fourier expansion with respect to $\cE$ indexed by Hermitian maps $a \co \cE \rightarrow \sigma^* \cE^{\vee}$. We let $ \wt{E}_{a}(m(\cE),s,\Phi)$ be the $a^{\mrm{th}}$ Fourier coefficient multiplied by certain normalization factors, explained precisely in \eqref{eq: normalized Eisenstein coefficients}. 


In our normalization, $s=0$ is the center of the functional equation for $\wt{E}_a(m(\cE), s, \Phi)$. Our main theorem relates the Taylor expansion at this central point to the degrees of special cycle classes. 

\begin{thm}\label{th:intro} Let $n \geq 1$ and $r \geq 0$. Let $\cE$ be a rank $n$ vector bundle on $X'$ and $a \co \cE \rightarrow \sigma^*  \cE^{\vee}$ be an injective Hermitian map. Then we have
\begin{equation}\label{eq: intro higher derivatives}
\frac{1}{(\log q)^r}  \left(\dds \right)^r\Big|_{s=0} \left( q^{ds} \wt{E}_{a}(m(\cE),s,\Phi) \right)  =  \deg [\cZ_{\cE}^r(a)],
\end{equation}
where $d=-\deg(\cE)+n\deg\omega_X =-\chi(X',\cE)$.
\end{thm}

\sss{Initial comments on the proof} We stress that \eqref{eq: intro higher derivatives} holds for \emph{all} $r$, regardless of the order of vanishing of $\wt{E}_{a}(m(\cE),s,\Phi)$ at $s=0$. This is a distinguishing novelty of Theorem \ref{th:intro} compared to all other works on the Seigel-Weil or arithmetic Siegel-Weil formula. The first results of this nature, giving motivic interpretations of Taylor coefficients of automorphic $L$-functions even ``beyond the leading term'', were proved in \cite{YZ, YZII} for $\PGL_2$. Our results here are the first higher derivative formulas to be proved for groups of \emph{arbitrary} rank. Our proof shares some common ingredients with these earlier works, but also has a number of interesting new ones. For example, a key discovery for us was a connection between the Fourier coefficients of Siegel--Eisenstein series and certain perverse sheaves arising from Springer theory. Another key realization was that the special cycles are governed by certain variants of the Hitchin fibration, whose geometry can also be described in terms of Springer theory. In particular, the geometry behind Theorem \ref{th:intro} is much more complicated than that in \cite{YZ, YZII} as soon as $n>2$. An overview of the proof will be given in \S \ref{ssec: overview of proof}.

Another feature of the proof of Theorem \ref{th:intro} is that it is completely uniform in $r$, and in particular unites the ``Siegel--Weil formula'' and ``arithmetic Siegel--Weil formula'' in the same framework. For this reason, we propose to call \eqref{eq: intro higher derivatives} a {\em higher Siegel--Weil formula}. This formula will serve as the first step to establish a higher order derivative version of the aforementioned recent work of Li and Liu \cite{LL1,LL2} over function fields, which would give a geometric interpretation of higher derivatives of Langlands $L$-functions.

\begin{remark}When $r=0$, the coarse moduli space of $\Sht_{U(n)}^r$ is just the discrete set of points which form the domain of everywhere unramified automorphic forms for $U(n)$. In that case, Theorem \ref{th:intro} specializes to (the non-singular Fourier coefficients of) the classical Siegel--Weil formula, which can be found in \cite{Weil1965}.

One should imagine that when $r=1$, $\Sht_{U(n)}^r \rightarrow X'$ is analogous to (the integral model of) a unitary Shimura variety. Now, under the technical assumptions of the present paper (namely the everywhere unramifiedness assumptions) this space is always empty, corresponding to the fact that the sign of the functional equation for the Siegel--Eisenstein series is $+1$ (so that all odd order derivatives vanish). However, with a mild modification of the setup, the same methods may be used to prove variants of Theorem \ref{th:intro} in which the sign of the functional equation is $-1$. More precisely, in this paper we consider rank $n$ vector bundles on $X'$ with a Hermitian pairing valued in the canonical bundle $\omega_{X'}\cong \nu^* \omega_X$; if we replace $\omega_{X} $ here by a line bundle on $X$ which is not a norm from $X'$, then the sign of the functional equation is $-1$ when $n$ is odd. The precise formulation is in \cite[\S9]{FYZ2}. We mention also the work \cite{Wei19} over function fields, which should be thought of as being similar to the special case of Theorem \ref{th:intro} for $r=1$ and $n=1$.


When $r>1$, no analogue of the spaces $\Sht_{U(n)}^r$ is presently known in the number field setting. Consequently, we do not know how to formulate an analogue of the main result for number fields.

\end{remark}


\subsubsection{Construction of virtual fundamental cycles}
For a vector bundle $\cE$ of rank $m$ on $X'$ and a Hermitian map $a\co \cE\to \s^{*}\cE^{\vee}$,  the dimension of $\Cal{Z}_{\cE}^r(a)$ differs from its ``virtual dimension'', which is $r(n-m)$. The situation gets worse if $a$ is {\em singular} (i.e., not injective, in analogy to the terminology of \cite{KRII}). For example, when $a=0$, $\cZ^{r}_{\cE}(0)$ contains $\Sht^{r}_{U(n)}$ as a substack. It is a nontrivial task to define a cycle class $[\Cal{Z}_{\cE}^r(a)]$ in the expected dimension $r(n-m)$. 

Our companion paper \cite{FYZ2} proposes two solutions to this problem, one using classical intersection theory and the other using derived algebraic geometry. There, we construct cycle classes $[\cZ^{r}_{\cE}(a)]$ for all $\cE$ of rank $\le n$ and \emph{possibly singular} $a \co \cE \rightarrow \sigma^* \cE^{\vee}$. Moreover,  we assemble them into generating series valued in the Chow groups of $\Sht_{U(n)}^r$ and conjecture it to be automorphic, in analogy to known results over number fields \cite{BHKRY1}, which fall under the umbrella of the Kudla program.

In this paper, we use a more elementary method to define the $0$-cycle $[\Cal{Z}_{\cE}^r(a)]$ in the case $m=n$ and $a$ injective. First, we prove that when $\cL$ is a line bundle and $a \co \cL \rightarrow \sigma^* \cL^{\vee}$ is an injective Hermitian map, $\Cal{Z}_{\cL}^r(a)$ has the expected dimension (cf. Proposition \ref{prop: dim Z_L} and Remark \ref{rem: actual dimension}). Next, when $\cE=\oplus_{i=1}^n \cL_i$ is a direct sum of line bundles, the class $[\Cal{Z}_{\cE}^r(a) ] \in \Ch_0(\Cal{Z}_{\cE}^r(a))_{\Q}$ can be defined as (the restriction to $\Cal{Z}_{\cE}^r(a)$ of) the intersection product of $\Cal{Z}_{\cL_i}^r(a_{ii})$ for the diagonal entries $a_{ii}$ of $a$; this is similar to the number field case.  However, compared to the number field case, a new difficulty arises since $\cE$ is not necessarily a direct sum of line bundles. We overcome this difficulty in \S\ref{ssec: intersection problem} by introducing the notion of a \emph{good framing} for $\cE$ to reduce to the case of a sum of line bundles. A nontrivial task is to verify that the cycle class $[\Cal{Z}_{\cE}^r(a) ]$ is independent of the choice of the good framing, which occupies much of the sections \S\ref{sec: hitchin spaces}--\S\ref{s:compare cycles}.


%

\subsection{Method of proof}\label{ssec: overview of proof} To summarize, we prove Theorem \ref{th:intro} by constructing two perverse sheaves that encode the two sides of \eqref{eq: intro higher derivatives} in the sense of sheaf-function correspondence, and then identifying these two perverse sheaves using a Hermitian variant of Springer theory, which labels these perverse sheaves by representations of the appropriate Weyl group. In this way, Theorem \ref{th:intro} is eventually unraveled into an elementary identity between representations of the Weyl group for type B/C.



On the geometric side, the connection between special cycles and Springer theory comes via the geometry of a moduli stack that resembles the Hitchin moduli space. On the other side, the connection between the Fourier coefficients of Siegel-Eisenstein series and Springer theory goes through local density formulas of Cho-Yamauchi.

Let us briefly explain the connection between the higher Siegel-Weil formula and the Hitchin moduli stack and Hermitian Springer theory, and defer details to the later paragraphs.  The degree of the special cycle that appear on the right side of \eqref{eq: intro higher derivatives} is essentially an intersection number of cycles on $\Sht^{r}_{U(n)}$. The ambient space $\Sht^{r}_{U(n)}$ can itself be realized an intersection of a Hecke correspondence with the graph of a Frobenius endomorphism. We use this to ``unfold'' all the intersections, and then redo them in a different order, performing the \emph{linear} intersections (i.e., those not involving the Frobenius map) first, and leaving the Frobenius semi-linear intersection till the last step (cf. \eqref{9 term} --  \eqref{row Hk}). In this process, a Hitchin-type moduli stack $\cM_{d}$ appears naturally as we perform linear intersections (cf. \eqref{Mbar d}). The degree of the special cycle $[\cZ^{r}_{\cE}(a)]$ can be expressed as a weighted counting of $k$-points on the fiber of a map $f_{d}: \cM_{d}\to \cA_{d}$ (analogue of Hitchin fibration) over the point $(\cE,a)\in \cA_d(k)$, where $(\cE,a)$ are as in the statement of Theorem \ref{th:intro}. 


The cokernel $\cQ=\coker(a)$ is a torsion sheaf on $X'$ with a Hermitian structure inherited from $a$. This motivates the introduction of the moduli stack $\Herm_{2d}(X'/X)$ that parametrizes torsion coherent sheaves on $X'$ of length $2d$ together with a Hermitian structure, so that $\cQ$ is a $k$-point of $\Herm_{2d}(X'/X)$ (where $2d=\dim_{k}\Gamma(X',\cQ)$). We show that the fiber of $f_{d}:\cM_{d}\to \cA_{d}$ over $(\cE,a)$ depends only on $\cQ=\coker(a)$, therefore the degree of $[\cZ^{r}_{\cE}(a)]$ depends only on the $k$-point $\cQ$ of $\Herm_{2d}(X'/X)$.



On the other hand, the Eisenstein series side of \eqref{eq: intro higher derivatives} can be written as a product of local terms -- representation density functions for Hermitian lattices. These density functions again only depend on the torsion sheaf $\cQ$ together with its Hermitian structure, i.e., a $k$-point in $\Herm_{2d}(X'/X)$.

Therefore we reduce to proving that two quantities attached to a $k$-point in $\Herm_{2d}(X'/X)$ are equal. A key realization is that both quantities are of motivic nature: they come by the sheaf-to-function correspondence from two   (graded, virtual) perverse sheaves on $\Herm_{2d}(X'/X)$. This is where Hermitian Springer theory enters.  Classically, starting with a reductive Lie algebra $\frg$, Springer theory outputs a perverse sheaf $\Spr_{\frg}$ on $\frg$, defined as the direct image complex of the Grothendieck-Springer resolution $\pi_{\frg}: \wt\frg\to \frg$,  together with an action of the Weyl group $W$. In our setting, $\Herm_{2d}(X'/X)$ will play the role of $\frg$.  In  \S\ref{sec: Herm}, we construct a perverse sheaf $\Spr^{\Herm}_{2d}$ on $\Herm_{2d}(X'/X)$ together with an action of $W_{d}=(\Z/2\Z)^{d}\rtimes S_{d}$ analogous to the Springer sheaf.  If $\Herm_{2d}(X'/X)$ is replaced by $\Coh_{d}(X)$, the moduli of torsion coherent sheaves on $X$ of length $d$, such a Springer sheaf was constructed by Laumon \cite{Lau87}. The Springer sheaf on $\Coh_{d}(X)$ (resp. $\Herm_{2d}(X'/X)$) can be viewed as a global version of the Springer sheaf for $\gl_{d}$ (resp. $\mathfrak{o}_{2d}$). The perverse sheaves on $\Herm_{2d}(X'/X)$ that govern both sides of \eqref{eq: intro higher derivatives} will be constructed from direct summands of the Hermitian Springer sheaf $\Spr^{\Herm}_{2d}$.

Thus, the proof of Theorem \ref{th:intro} is completed in three steps:
\begin{enumerate}
\item Construct a graded perverse sheaf on $\Herm_{2d}(X'/X)$
\begin{equation*}
\cK^{\Eis}_{d}=\bigoplus_{i=0}^{d}\cK^{\Eis}_{d,i}
\end{equation*}
whose Frobenius trace at $\cQ$ is related to the LHS of \eqref{eq: intro higher derivatives}. More precisely, 
\begin{equation}\label{intro KEis}
\wt{E}_{a}(m(\cE),s, \Phi)= \sum_{i=0}^{d}\Tr(\Frob_{\cQ}, (\cK^{\Eis}_{d,i})_{\cQ})q^{-2is}.
\end{equation}
\item Construct a graded perverse sheaf on $\Herm_{2d}(X'/X)$
\begin{equation*}
\cK^{\Int}_{d}=\bigoplus_{i=0}^{d}\cK^{\Int}_{d,i}
\end{equation*}
whose Frobenius trace at $\cQ$ is relate to the RHS of \eqref{eq: intro higher derivatives}. More precisely,
\begin{equation}\label{intro KInt}
\deg [\cZ_{\cE}^r(a)]=\sum_{i=0}^{d}\Tr(\Frob_{\cQ}, (\cK^{\Int}_{d,i})_{\cQ})\cdot (d-2i)^{r}.
\end{equation}
\item Prove that 
\begin{equation}\label{KK}
\cK^{\Eis}_{d}\cong \cK^{\Int}_{d}
\end{equation}
as graded perverse sheaves on $\Herm_{2d}(X'/X)$.
\end{enumerate}

These three steps correspond to the three parts of the paper. We elaborate on the main ideas involved in each step. 

\subsubsection{Step (1)} After a standard procedure expressing the nonsingular Fourier coefficients of Eisenstein series in terms of local density of Hermitian lattices, we use the formula of Cho and Yamauchi \cite{CY} for these densities (more precisely, the unitary variant developed in \cite{LZ1}). We also need an extension of their formula in the split case (Theorem \ref{th:CY}). The formula of Cho and Yamauchi depends only on the Hermitian torsion sheaf $\cQ=\coker(a)$, which gives the hope that the local density, as a function on the set of Hermitian torsion sheaves, comes from a sheaf on $\Herm_{2d}(X'/X)$ via Grothendieck's sheaf-to-function dictionary. We do this by developing an analog of Springer theory over  $\Herm_{2d}(X'/X)$ (\S\ref{sec: springer}-\S\ref{sec: Herm}).


The key observation here is that the term in the Cho--Yamauchi formula resembles the Frobenius trace function for a certain linear combination of Springer sheaves for $\gl_{d}$ or $\Coh_{d}(X)$, except for some signs. To match the signs exactly we consider an analogous linear combination of Springer sheaves on $\Herm_{2d}(X'/X)$, and we compare the Frobenius actions on the cohomology of Springer fibers over $\Coh_{d}(X)$ and over $\Herm_{2d}(X'/X)$, see \S\ref{ss:HSpr stalk} and \S\ref{ss:trace comp Spr}.

\subsubsection{Step (2)} This step consists of three substeps. 
\begin{itemize}
\item First, we define special cycles for nonsingular $a$ (\S\ref{sec: unitary shtukas}-\S\ref{s:int prob}). When $\cE$ is a direct sum of line bundles $\cL_{i}$, we define, following Kudla and Rapoport, $[\cZ_{\cE}^{r}(a)]$ as the intersection of cycle classes $[\cZ_{\cL_{i}}^{r}(a_{ii})]$, which, despite not being divisors in our setting, always have the ``expected'' dimension (more precisely, codimension $r$ in $\Sht^{r}_{U(n)}$). The definition of $[\cZ^{r}_{\cE}(a)]$ for general vector bundles $\cE$ requires choosing a ``good framing'' on $\cE$, i.e., an injective map from a direct sum of line bundles $\cE'=\op_{i=1}^{n}\cL_{i}\inj \cE$ satisfying certain conditions.  In any case, the RHS of \eqref{eq: intro higher derivatives} is an intersection number of cycles on $\Sht^{r}_{U(n)}$.  
\item The well-definedness of $[\cZ^{r}_{\cE}(a)]$ is proved in the second substep (\S\ref{sec: hitchin spaces}-\S\ref{s:compare cycles}), which also gives a different definition of these cycle classes without any choices. The idea is similar to the one used in \cite{YZ}, namely by exchanging the order of intersection, we perform ``linear intersections'' first to form Hitchin-type moduli stacks (denoted $\cM_{d}$, making sense over any base field), and in the last step we perform a shtuka-type construction by intersecting with the graph of Frobenius. 
\item In the last substep (\S\ref{s:int trace}) we use the Lefschetz trace formula to express the degree of $[\cZ^{r}_{\cE}(a)]$, formulated using the Hitchin-type moduli stack $\cM_{d}$, as the trace of Frobenius composed with the $r^{\mrm{th}}$ power of an endomorphism $C$ on the direct image complex $Rf_{*}\Qlbar$ of the Hitchin map  $f:\cM_{d}\to \cA_{d}$. Now, the ``Hitchin base'' $\cA_d$ has a canonical smooth map to $\Herm_{2d}(X'/X)$, and it turns out that $Rf_{*}\Qlbar$ together with its endomorphism $C$ descends through this map to a perverse sheaf $\cK^{\Int}_{d}$ on $\Herm_{2d}(X'/X)$ with an endomorphism $\ov C$. The decomposition of $\cK^{\Int}_{d}$ into graded pieces $\cK^{\Int}_{d,i}$ is according to the eigenvalues of the $\ov C$-action, which are of the form $(d-2i)$. Combining these facts we get \eqref{intro KInt}.
\end{itemize}

\subsubsection{Step (3)} Both $\cK^{\Eis}_{d}$ and $\cK^{\Int}_{d}$ are linear combinations of isotypical summands of $\Spr^{\Herm}_{2d}$ under the action of $W_{d}$. The isomorphism \eqref{KK} then comes from an isomorphism of two graded virtual representations of $W_{d}$, which we verify directly.

\subsection*{Acknowledgment}

We thank Michael Harris, Chao Li, Alan Peng, Zhiyu Zhang, and the anonymous referee for comments on a draft. We thank the referee for helpful suggestions and corrections. TF was supported by NSF Postdoctoral Fellowship DMS-1902927, NSF Grant DMS-2302520, and the Friends of the Institute for Advanced Study.  ZY was partially supported by the Packard Fellowship, and the Simons Investigator grant. WZ is partially supported by the NSF grant DMS \#1901642  and the Simons Investigator grant.

\subsection{Notation}\label{ssec: notation}

Throughout this paper, $k=\F_{q}$ is a finite field of odd characteristic $p$. Let $\ell\ne p$ be a prime.   Let $\psi_{0}: k\to\Qlbar^{\times}$ be a nontrivial character. For a stack $\cY$ over $k$, we write $\Fr$ or $\Fr_{\cY}$ for its $q$-power Frobenius endomorphism. We will use $\gFrob$ or $\gFrob_{y}$ for geometric Frobenius at an $\F_{q}$-point $y$. 

Let $X$ denote a smooth curve over $k$. With the exception of \S\ref{sec: springer} and \S\ref{sec: Herm}, $X$ is assumed to be projective and geometrically connected. Let $\om_{X}$ be the line bundle of $1$-forms on $X$. 

Let $F=k(X)$ denote the function field of $X$.  Let $|X|$ be the set of closed points of $X$. For $v\in |X|$, let $\cO_{v}$ be the completed local ring of $X$ at $v$ with fraction field $F_{v}$ and residue field $k_{v}$.   Let $\BA=\BA_{F}$ denote the ring of ad\`eles of $F$, and $\wh\cO=\prod_{v\in |X|}\cO_{v}$. Let $\deg(v)=[k_{v}:k]$, and $q_{v}=q^{\deg(v)}=\#k_{v}$. A uniformizer of $\cO_{v}$ is typically denoted $\vp_{v}$. Let $|\cdot|_{v}: F^{\times}_{v}\to q^{\Z}_{v}$ be the absolute value such that  $|\vp_{v}|_{v}=q^{-1}_{v}$. Let $|\cdot|_{F}: \BA_{F}^{\times}\to q^{\Z}$ be the absolute value that is $|\cdot|_{v}$ on $F_{v}^{\times}$.

Let $X'$ be another smooth curve over $k$ and $\nu:X'\to X$ be a finite map of degree $2$ that is generically \'etale. We denote by $\sigma$ the non-trivial automorphism of $X'$ over $X$. With the exception of \S\ref{ss:Herm} and \S\ref{ss:Herm Spr}, $\nu$ is assumed to be \'etale. We emphasize that the case $X'=X\coprod X$ is allowed.  Let $F'$ be the ring of rational functions on $X'$, which is either a quadratic extension of $F$ or $F\times F$. We let $k'$ be the ring of constants in $F'$. The notations $\om_{X'}, |X'|, F'_{v'}, \cO_{v'}, k_{v'}, \BA_{F'}, |\cdot|_{v'}, |\cdot|_{F'}, q_{v'}$ and $\deg(v')$ (for $v'\in |X'|$) are defined similarly as their counterparts for $X$. Additionally, for $v\in |X|$, we use $\cO'_{v}$ to denote the completion of $\cO_{X'}$ along $\nu^{-1}(v)$, and define $F'_{v}$ to be its total ring of fractions.

For a vector bundle $\cE$ on $X'$, let $\cE^{\vee}=\un\Hom(\cE, \om_{X'})$ be its Serre dual.  For a torsion sheaf $\cT$ on $X'$, let $\cT^{\vee}=\un\Ext^{1}(\cT, \om_{X'})$.  

When $X$ (hence $X'$) is projective, let $\Bun_{\GL_{n}}$ (resp. $\Bun_{\GL_{n}'}$) be the moduli stack of rank $n$ vector bundles over $X$ (resp. $X'$). Let $g$ be the genus of $X$ and $g'$ be the arithmetic genus of $X'$. Note that whenever $\nu$ is \'{e}tale, we have $g' = 2g-1$. 

For an algebraic stack $\cY$, $\Ch(\cY)$ denotes its \emph{rationalized} Chow group and $D^b(\cY,\Qlbar)$  its bounded derived category of constructible $\Qlbar$-sheaves.

%


\part{The analytic side}
\section{Fourier coefficients of Eisenstein series}\label{s:Eis}

In this section we will define the Siegel-Eisenstein series featuring into our main theorem, and explain how to express their non-singular Fourier coefficients in terms of local density polynomials, which will be geometrized in later sections.

\subsection{Siegel--Eisenstein series}\label{ss:Eis} 
For any one-dimensional $F$-vector space $L$, let $\Herm_{n}(F,L)$ be the $F$-vector space of $F'/F$-Hermitian forms $h: F'^{n}\times F'^{n}\to L\ot_{F}F'$ (with respect to the involution $1\ot\s$ on $L\ot_{F}F'$). For any $F$-algebra $R$, $\Herm_{n}(R, L):=\Herm_{n}(F,L)\ot_{F}R$ is the set of $L\ot_{F}R'$-valued $R'/R$-Hermitian forms on $R'^{n}$, where $R'=R\ot_{F}F'$.  When $L=F$ we write $\Herm_{n}(F)=\Herm_{n}(F,F)$ and $\Herm_{n}(R)=\Herm_{n}(F)\ot_F R$ for any $F$-algebra $R$.

Let $W$ be the standard split $F'/F$-skew-Hermitian space of dimension $2n$. 
Let $H_n=U(W)$. Write $\BA := \BA_F$ for the ring of adeles of $F$. Let $P_n(\BA)=M_n(\BA)N_n(\BA)$ be the standard Siegel parabolic subgroup of $H_n(\BA)$, where
\begin{align*}
  M_n(\BA)&=\left\{m(\alpha)=\begin{pmatrix}\alpha & 0\\0 &{}^t\bar \alpha^{-1}\end{pmatrix}: \alpha\in \GL_n(\BA_{F'})\right\},\\
  N_n(\BA) &= \left\{n(\beta)=\begin{pmatrix} 1_n & \beta \\0 & 1_n\end{pmatrix}: \beta\in\Herm_n(\BA_F)\right\}.
\end{align*}

Let $\eta: \BA_{F}^\times/F^\times\rightarrow \mathbb{C}^\times$ be the quadratic character associated to $F'/F$ by class field theory. Fix $\chi: \BA_{F'}^\times/F'^\times\rightarrow \mathbb{C}^\times$ a character such that $\chi|_{\BA_F^\times}=\eta^n$. We may view $\chi$ as a character on $M_n(\BA)$ by $\chi(m(\alpha))=\chi(\det(\alpha))$ and extend it to $P_n(\BA)$ trivially on $N_n(\BA)$. Define the \emph{degenerate principal series} to be the unnormalized smooth induction 
$$
I_n(s,\chi)=\Ind_{P_n(\BA)}^{H_n(\BA)}(\chi\cdot |\cdot|_{F'}^{s+n/2}),\quad s\in \mathbb{C}.
$$

For a standard section $\Phi(-, s)\in I_n(s,\chi)$, define the associated \emph{Siegel--Eisenstein series} 
$$E(g,s, \Phi)= \sum_{\gamma\in P_n(F)\backslash H_n(F)}\Phi(\gamma g, s),\quad g\in H_n(\BA),
$$ which converges for $\Re(s)\gg 0$ and admits meromorphic continuation to $s\in \mathbb{C}$. Notice that $E(g,s,\Phi)$ depends on the choice of $\chi$.

\begin{remark}\label{rem: chi exists}
In this paper, we will choose $\chi$ to be unramified everywhere. To see that such $\chi$ exists, observe that since $\CC^{\times}$ is injective (in the category of abelian groups), it suffices to check that $\eta^n$ is trivial on $\ker (\Pic(X) \rightarrow \Pic(X') )$. If $X'/X$ is the trivial double cover or the double cover corresponding to $\F_{q^2}/\F_q$, then then this kernel is trivial so the result is immediate. Otherwise, the cover is geometrically non-trivial. Since $\mrm{char}(k) \neq 2$, the kernel consists of the 2-torsion line bundle whose class in $\cohog{1}{X, \mu_2}$ agrees with $\eta \in H^1(X, \Z/2\Z)$ under the isomorphism $\mu_2 \cong \Z/2\Z$. If $n$ is even then there is nothing to check; if $n$ is odd then the desired vanishing property amounts (when $\mrm{char}(k) \neq 2$) to the alternating property of the cup product pairing $H^1(X_{\ol{\F}_q},\Z/2\Z)\times \cohog{1}{X_{\ol{\F}_q},\Z/2\Z}\to \Z/2$, which follows from the graded commutativity of the cup product and the fact that the geometric $\Z_2$-cohomology of curves in characteristic $\neq 2$ is torsion-free.
\end{remark}

As justified by Remark \ref{rem: chi exists}, we may choose $\chi$ to be everywhere unramified. Then  $I_n(s,\chi)$ is unramified and we fix $\Phi(-, s)\in I_n(s,\chi)$  as the unique $K=H_n(\wh\cO)$-invariant section normalized by
$$
\Phi(1_{2n}, s)=1.
$$
Similarly we normalize $\Phi_v \in I_n(s,\chi_v)$ for every $v\in|X|$ and we then have a factorization $\Phi=\bigotimes_{v\in|X|}\Phi_v$.

\subsection{Fourier expansion} 
Let $\om_{F}$ be the generic fiber of the canonical bundle of $X$, and $\BA_{\om_{F}}=\BA\ot_{F}\om_{F}$. The sum of the residues induces a pairing $\BA_{\om_{F}}\times \BA_{F}\to k$ induces a pairing  
$$\j{\cdot,\cdot}: \Herm_{n}(\BA, \om_{F})\times \Herm_{n}(\BA)\to k$$
given by $\j{T,b}=\Res(-\Tr(Tb))$. Composing this pairing with the fixed nontrivial additive character $\psi_{0}: k\to \BC^\times$  exhibits $\Herm_{n}(\BA,\om_{F})$ as the Pontryagin dual of $\Herm_{n}(\BA)$. Moreover, it exhibits $\Herm_{n}(F,\om_{F})$ as the Pontryagin dual of $\Herm_{n}(F)\bs \Herm_{n}(\BA)=N_{n}(F)\bs N_{n}(\BA)$. The global residue pairing is the sum of local residue pairings $\j{\cdot,\cdot}_{v}: \Herm_{n}(F_{v}, \om_{F_v})\times \Herm_{n}(F_{v})\to k$ defined by $\j{T,b}_{v}=\tr_{k_v/k}\Res_{v}(-\Tr(Tb))$.


We have a Fourier expansion 
$$E(g,s,\Phi)=\sum_{T\in\Herm_n(F,\om_{F})}E_T(g,s,\Phi),$$ 
where
 $$E_T(g,s,\Phi)=\int_{N_n(F)\backslash N_n(\mathbb{A})} E(n(b)g,s,\Phi)\psi_{0}(\j{T,b})\,\rd n(b),$$ 
and the Haar measure $\rd n(b)$ is normalized such that $N_{n}(F)\bs N_{n}(\BA)$ has volume $1$. 
For any $\a\in M_{n}(\BA)$ we have
\begin{equation}\label{Eis ch base}
E_{T}(m(\a)g,s,\Phi)=\chi(\det(\bar\alpha))^{-1}|\det(\alpha)|^{-s+n/2}_{F'}E_{{}^{t}\ov\a T \a}(g,s,\Phi).
\end{equation}

Suppose $T$ is \emph{nonsingular}, meaning that for one (equivalently, any) choice of trivialization of $\omega_F$ it has non-vanishing determinant, for a factorizable $\Phi=\bigotimes_{v\in|X|}\Phi_v$ we have a factorization of the Fourier coefficient into a product (cf. \cite[\S 4]{Kudla1997})
\begin{equation}\label{Eis factor Wh}
E_T(g,s,\Phi)=|\om_{X}|_{F}^{-n^{2}/2}\prod_v W_{T,v}(g_v, s, \Phi_v),
\end{equation}
where the \emph{local (generalized) Whittaker function} is defined by 
$$W_{T,v}(g_v, s, \Phi_v)=\int_{N_n(F_{v})}\Phi_v(w_n^{-1}n(b)g_v,s) \psi_{0}(\j{T,b}_{v})
\, \rd_{v} n(b),\quad w_n=
\begin{pmatrix}
0  & 1_n\\
  -1_n & 0\\
\end{pmatrix}$$ 
and has analytic continuation to $s\in \mathbb{C}$. Here the local Haar measure $\rd_{v} n(b)$ is the one such that the volume of $N_{n}(\cO_{v})$ is $1$. The factor $|\om_{X}|_{F}^{-n^{2}/2}$ is the ratio between the global measure $\rd n$ and the product of the local measures $\prod_{v}\rd_{v}n$.

Note that for $\a\in M_{n}(F_{v})$,
\begin{equation}\label{eq:ch base}
W_{T,v}(m(\alpha), s, \Phi_v)=\chi(\det(\bar\alpha))^{-1}|\det(\alpha)|^{-s+n/2}_{F'_v}W_{^t\bar\alpha\,T\alpha,v}(1, s, \Phi_v).
\end{equation}

We define the {\em regular} part of the Eisenstein series to be
\begin{equation}\label{eq:def reg}
E^{\rm reg}(g,s,\Phi)=\sum_{T\in\Herm_n(F,\om_{F})\atop \rank T=n}E_T(g,s,\Phi).
\end{equation}

\subsection{Local densities for Hermitian lattices}\label{ss:loc den}
The local density for Hermitian lattices in the non-split case has been studied in \cite[\S3]{LZ1} following the strategy of Cho--Yamauchi \cite{CY}. Here we recall the result of \cite{LZ1} and extend the results to the split case.
 
From now on until \S\ref{sec:relation-with-local}, let $F$ be a non-archimedean local field of characteristic not equal to $2$ (but possibly with residue characteristic $2$). Let $F'$ be either an unramified quadratic field extension or  the split quadratic $F$-algebra $F'=F\times F$. Denote by $\cO_F$ (resp. $\cO_{F'}$) the ring of integers in $F$ (resp. $F'$). In the split case we have $\cO_{F'}=\cO_F\times \cO_F$.  Let $\eta=\eta_{F'/F}: F^\times\to\{\pm1\}$ be the quadratic character attached to $F'/F$ by class field theory. Let $\varpi$ be a uniformizer of $F$, $k$ the residue field, $q=\#k$. 

Let $L, M$ be two Hermitian $\cO_{F'}$-lattices. In the split case,  the datum of a \emph{Hermitian $\cO_{F'}$-lattice $L$} is a pair $(L_1,L_2)$ of $\cO_{F}$-lattices together with an $\cO_{F}$-bilinear pairing
\begin{equation*}
(\cdot,\cdot): L_{1}\times L_{2} \to \cO_{F}
\end{equation*}
that is perfect after base change to $F$.  We will define $L^\vee=(L_1^\vee, L_2^\vee)$ where $L_1^\vee= \{x\in L_{1}\otimes_{\cO_F} F: (x,L_2)\subset \cO_F  \}$ and similarly for $L_2^\vee$.

Let $\Rep_{M,L}$ be the \emph{scheme of integral representations of $M$ by $L$}, an $\cO_F$-scheme such that for any $\cO_F$-algebra $R$, 
\begin{align*}\label{def: Rep}
\Rep_{M,L}(R)=\Herm(L \otimes_{\cO_F}R, M \otimes_{\cO_F}R),
\end{align*} 
where $\Herm$ denotes the set of Hermitian $R$-module homomorphisms. In the split case, if we write $L$ and $M$ in terms of pairs $(L_1,L_2)$ and $(M_1,M_2)$ with their $\cO_{F}$-bilinear pairings, then a  Hermitian module homomorphism consists of a pair of $R$-linear maps $\phi_i:L_i \otimes_{\cO_F}R\to M_i \otimes_{\cO_F}R $ preserving the base change to $R$ of the $\cO_{F}$-bilinear pairings.

The \emph{local density} of integral representations of $M$ by $L$ is defined to be $$\Den(M,L)\colon= \lim_{N\rightarrow +\infty}\frac{\#\Rep_{M,L}(\cO_F/\varpi^N)}{q^{N\cdot\dim (\Rep_{M,L})_{F}}}.$$ Note that if $L, M$ have $\cO_{F'}$-rank $n, m$ respectively and the generic fiber $(\Rep_{M,L})_{F}\ne\varnothing$, then $n\le m$ and
\begin{equation*}
\dim (\Rep_{M,L})_{F}=\dim \U_m-\dim \U_{m-n}= n\cdot (2m-n).  
\end{equation*}

\subsection{Cho--Yamauchi formula for local density}\label{ssec: Cho-Yamauchi}

%
%

\begin{defn}\label{def:maT}
For $a\in \Z_{\ge0}$ we define  a polynomial of degree $a$
\begin{equation*}
\fm(a; T) := \prod_{i=0}^{a-1}(1-(\y(\varpi)q)^i T) \in \Z[T].
\end{equation*}
\end{defn}
Note that $\fm(a; T)$ depends on $F'/F$. 

In both the non-split and the split cases, for a finite torsion $\cO_{F}$-module $\cT$ we define 
\begin{eqnarray*}
\label{def ell}\ell(\cT)&:=&\mbox{ length of $\cT$ as an $\cO_{F}$-module};\\
\label{def t} t(\cT)&:=&\dim_{k}(\cT\ot_{\cO_{F}}k).
\end{eqnarray*}
For an $\cO_{F'}$-Hermitian lattice $L$, we define its type 
$$
t(L):=t(L^\vee/L)
$$
where we view the finite torsion $\cO_{F'}$-module $L^\vee/L$ as an $\cO_{F}$-module.

When $F'/F$ is non-split, for a finite torsion $\cO_{F'}$-module $\cT$ we define 
\begin{eqnarray*}
\label{def ell'}\ell'(\cT)&:=&\mbox{length of $\cT$ as an $\cO_{F'}$-module};\\
\label{def t'}t'(\cT)&:=&\dim_{k'}(\cT\ot_{\cO_{F'}}k').
\end{eqnarray*}
 Then we have 
\begin{equation}\label{eq:ell&t}
\ell(\cT)=2\ell'(\cT),\quad t(\cT)=2t'(\cT).
\end{equation}
When $F'=F\times F$ is split, for a finite torsion $\cO_{F'}$-module $\cT$ we may define $\ell'(\cT)$ and $t'(\cT)$ by \eqref{eq:ell&t}. Moreover, for  $\cO_{F'}$-Hermitian lattices $L=(L_1,L_2)$ and $L'=(L'_1,L'_2)$ such that $L\subset L'$ (meaning that $L_1 \subset L_1'$ and $L_2 \subset L_2'$), we have 
$$
\ell(L'/L)=\ell(L'_1/L_1)+\ell(L'_2/L_2)
$$
and
$$
t'(L^\vee/L)=t(L_2^{\vee}/L_1)=t(L_1^{\vee}/L_2).
$$
In both the split and non-split case, we define 
$$
t'(L)=t'(L^\vee/L).
$$


We have the following analog of Cho--Yamauchi formula \cite{CY}.
\begin{thm}\label{th:CY} 
Let $j\ge0$ be an integer. Let $\iden^j$ be the self-dual Hermitian $\cO_{F'}$-lattice of rank $j$ with Hermitian form given the identity matrix $\mathbf{1}_j$. Let $L$ be a Hermitian $\cO_{F'}$-lattice of rank $n$. 
  \begin{enumerate}
  \item We have
\begin{equation}
  \label{eq: iden}
\Den(\iden^{n+j}, \iden^n)=\prod_{i=1}^n(1-(\eta(\varpi)q)^{-i}T)\bigg|_{T= (\eta(\varpi)q)^{-j}}.
\end{equation}
\item
There is  a (unique) polynomial $\Den(T,L)\in \mathbb{Z}[T]$, called  (normalized) \emph{local Siegel series} of $L$,  such that for all $j\geq 0$,
$$\Den((\eta(\varpi) q)^{-j},L)=\frac{\Den(\iden^{n+j}, L)}{\Den(\iden^{n+j}, \iden^n)}.$$ 
  \item We have
\begin{equation}\label{eq:CY}
\Den(T,L) = \sum_{L \subset L '\subset L'^{\vee} \subset L^{\vee}} T^{2\ell'(L'/L)}\fm(  t'(L');T). 
\end{equation}
Here the sum is over $\cO_{F'}$-lattices $L'$ (in the $F'$-Hermitian space spanned by $L$) containing $L$ on which the Hermitian form is integral.
\end{enumerate}
\end{thm}

\begin{proof}
The non-split case is proved in \cite[Thm.~3.5.1]{LZ1} and here we indicate the necessary change in the split case.  Now suppose $F'=F\times F$ and hence $k'=k\times k$. Let $L_k=L\otimes_{\cO_F} k$ and $\iden^m_{k}=\iden^m\otimes_{\cO_F} k$, which are free $k'$-modules with the induced $k'/k$-Hermitian forms. In particular, $\iden^m_{k}$ is non-degenerate and the radical of $L_k=L\otimes_{\cO_F} k$  has $k'$-rank equal to $t'(L)=t(L_1^\vee/L_1)=t(L_2^\vee/L_2)$. Let $\Isom_{\iden^m_{k}, L_k}$ be the $k$-scheme of ``isometric embeddings" from $L_k$ to  $\iden^m_{k}$, i.e., {\em injective} $k'$-linear maps from $L_k$ to  $\iden^m_{k}$ preserving the Hermitian forms.

Similar to the orthogonal case \cite[\S3.3]{CY},  we have
  $$
  \Den(\iden^m,L)= q^{-\dim \Rep(\iden^m,L)_{F}}\sum_{L\subset L'\subset L'^\vee} \# (L'/L)^{-(m-n)} \# \Isom_{\iden^m_{k}, L_k}(k),
  $$
  where $\dim \Rep(\iden^m,L)_{F}=m^2-(m-n)^2=2mn-n^2$.

It remains to show that 
  \begin{align}\label{eq:k pts}
   \#\Isom_{\iden^m_{k}, L_k}(k)= q^{m^2-(m-n)^2}\cdot \prod_{i=0}^{n+a-1}(1-q^{i-m})
  \end{align}
  where $a=t'(L)$ is the $k'$-rank of the radical of $L_k$. Note that up-to-isomorphism, $L_k$ is determined by its rank and the rank of its radical. Let $U_{n-a,a}$ be a $k'/k$-Hermitian space of rank $n$ with radical of rank $a$. Let $V_m=U_{m,0}$ be a (non-degenerate) $k'/k$-Hermitian space of dimension $m\ge n$. Then it is easy to see that $U_{n-a,a}\simeq U_{n-a,0}\oplus U_{0,a}$ and
  \begin{equation}\label{eq:k pts1}
 \# \Isom_{V_m, U_{n-a,a}}(k)= \# \Isom_{V_m, U_{n-a,0}}(k)\cdot\#  \Isom_{V_{m-(n-a)}, U_{0,a}}(k).
  \end{equation}
By \eqref{eq:k pts1} (note that $   \#\Isom_{\iden^m_{k}, L_k}(k)= \# \Isom_{V_m, U_{n-a,a}}(k)$), it suffices to show \eqref{eq:k pts} in the two extreme cases: $a=0$ and $a=n$.
 
 First we consider the case $a=n$. Then, to give an isometric embedding from $U=U_{0,n}=k'^n$ to $V=U_{m,0}=k'^m$ is equivalent to give an injective $k$-linear map $\phi:k^n\to k^m$ and then an injective $k$-linear map $\varphi:k^n\to {\rm Im}(\phi)^\perp\subset k^m$. Therefore, denoting  by $\Hom^*_k(k^n,k^m)$ the set of  injective $k$-linear maps $\phi:k^n\to k^m$, we have
  \begin{align*}
 \# \Isom_{V_{m}, U_{0,n}}(k)= &\#\Hom^*_k(k^n,k^m)\cdot  \#\Hom^*_k(k^n,k^{m-n})\\
 =&q^{mn}\prod_{i=0}^{n-1}(1-q^{i-m})\cdot  q^{(m-n)n}\prod_{i=0}^{n-1}(1-q^{i-m+n})
 \\=&  q^{2mn-n^2}\prod_{i=0}^{2n-1}(1-q^{i-m}).
  \end{align*}
  
It remains to consider the case $a=0$. Then  a similar argument shows
   \begin{align*}
 \# \Isom_{V_{m}, U_{n,0}}(k)= &\#\Hom^*_k(k^n,k^m)\cdot  \#\Hom_k(k^n,k^{m-n})\\
 =&q^{mn}\prod_{i=0}^{n-1}(1-q^{i-m})\cdot  q^{(m-n)n}
 \\=&  q^{2mn-n^2}\prod_{i=0}^{n-1}(1-q^{i-m}).
  \end{align*}
  This completes the proof.
  
 \end{proof}
 
 \begin{remark}\label{r:Den tor} By Theorem \ref{th:CY}, the polynomial $\Den(T,L)$ depends only on the induced Hermitian form on the torsion module $L^\vee/L$. Indeed, for a Hermitian torsion module\footnote{By this we mean a torsion $\cO_{F'}$-module with an $\cO_{F'}/\cO_F$-Hermitian form.} $Q$ we define $\Den(T,Q)$ by the formula
\[
\Den(T,Q) = \sum_{Q'\subset (Q')^{\vee} \subset Q} T^{2\ell'(Q')}\fm(  t'(Q');T). 
\]
Then by \eqref{eq:CY}, we have $\Den(T, L) = \Den(T, L^\vee/L)$. 
\end{remark}

\begin{remark}
In the split case, write $L=(L_1,L_2)$ and $L'=(L'_1,L'_2)$. Then the formula reads
\begin{eqnarray*}
\Den(T,L) = \sum_{L_1 \subset L'_1 \subset L_2^{'\vee} \subset L_2^{\vee}  } T^{\ell(L_1'/L_1)+\ell(L_2'/L_2)}\fm(t({L_2'}^{\vee}/L'_1);T).
\end{eqnarray*}
\end{remark}

\begin{remark}
The local Siegel series satisfies a functional equation 
\begin{equation}
  \label{eq:functionalequation}
  \Den(T,L)=(\eta(\varpi) T)^{\ell'(L^\vee/L)}\cdot \Den\left(\frac{1}{T},L\right).
\end{equation}A proof in the inert case can be found in \cite[Theorem 5.3]{Hironaka2012}. By Theorem, \ref{th:CY} the constant term of $  \Den(T,L)$ is $1$. It follows that the degree of the polynomial $\Den(T,L)$ is equal to $\ell'(L^\vee/L)$. We will not use this fact in this paper.
See Corollary \ref{cor:FE Int} for the (global) geometric analog.
\end{remark}

\subsection{Relation with local Whittaker functions}\label{sec:relation-with-local} We continue to let $F$ be a local field. 
Define the local L-function
\begin{equation*}
\sL_{n,F'/F}(s) :=\prod_{i=1}^nL(i+2s, \eta^i)=\prod_{i=1}^n\frac{1}{1-\eta^i (\varpi)q^{-i-2s}}.
\end{equation*}

\begin{lemma}\label{lem:Whi}
Let $L$ be a Hermitian $\cO_{F'}$-lattice of rank $n$. Let $T=((x_i, x_j))_{1\le i,j\le n}$ be the fundamental matrix of an $\cO_{F'}$-basis $\{x_1,\ldots, x_n\}$ of $L$, an $n\times n$ Hermitian matrix over $F$. Let $\th$ be a generator of $\om_{\cO_{F}}$ so that $T\th\in \Herm_{n}(F,\om_{F})$. Then
  \begin{align*}
W_{T\th}(1, s, \Phi)=\sL_{n,F'/F}(s)^{-1} \Den(q^{-2s},  L).
  \end{align*}
  Here $\Phi$ is the local unramified section normalized by $\Phi(1_{2n},s)=1$.
\end{lemma}
\begin{proof}
Note that by Theorem \ref{th:CY} 
$$\sL_{n,F'/F}( j)=\Den(\iden^{n+2j},  \iden^{n})^{-1}. $$
It is known that $W_{T\th}(1, s, \Phi)$ is a rational function in $q^{s}$.
Therefore the formula is equivalent to
  \begin{align*}
W_{T\th}(1, j, \Phi)= \Den(\iden^{n+2j},  L).
  \end{align*}
  for all integer $j\geq 0$. 
In the non-split case this is essentially \cite[Prop.~10.1]{KRII} (cf. \cite[\S3.3]{LZ1}), which can be easily modified to the split case. We note that $W_{T\th}$ is the same as $W_{T}$ in {\em loc. cit.}.
\end{proof}
%
%
%
%

%


\subsection{Fourier coefficients revisited}\label{ss:FC} Now we return to the global situation. 
We need the following global L-function to normalize the Eisenstein series
\begin{equation}\label{eq: L-factor}
\sL_n(s)=\prod_{i=1}^nL(i+2s, \eta^i).
\end{equation}


We now consider the restriction of the regular part  $E^{\rm reg}(\cdot,s,\Phi)$ (as a function in $g\in H_n(\BA)$, cf. \eqref{eq:def reg}) to the Levi subgroup $M_n(\BA)$. Since the restriction is left $M_n(F)$-invariant and right $K$-invariant, it descends to a function on
$$M_n(F)\bs M_n(\BA)/M_n(\wh\cO)\simeq \Bun_{M_n}(k)\simeq \Bun_{\GL_n'}(k),
$$
via the canonical identifications.
From now on we will freely switch between $g=m(\alpha)\in M_n(\BA)$ and the corresponding element $\cE\in  \Bun_{\GL_n'}(k)$ and we will write 
$$
E^{\rm reg}(m(\cE),s,\Phi)=E^{\rm reg}(m(\alpha),s,\Phi).
$$
Note that the absolute value on $\BA_{F'}^\times$ is normalized such that $|\det(\alpha)|_{F'}=q^{\deg(\cE)}$. 
By abuse of notation we also view $\chi$ as a function on $\Bun_{\GL_1'}(k)$. 

Recall that $\cE^\vee=\ul{\Hom}_{\cO_{X'}}(\cE,\omega_{X'})$ denotes the Serre dual of $\cE$. Consider a rational Hermitian map $a:\cE\dashrightarrow \s^{*}\cE^{\vee}$ (i.e., $a$ is a map defined over the generic point of $X'$, such that $\s^* a = a$). Given a pair $(\cE,a)$ as above, we shall define the Fourier coefficient
\begin{equation*}
E_{a}(m(\cE), s,\Phi)
\end{equation*}
as follows. For any generic trivialization $\t: \cE_{F'}\isom (F')^{n}$, the pair $(\cE,\t)$ gives a point $\a=\a(\cE,\t)\in M_n(\BA)/M_n(\wh\cO)$ such that $\cE$ is glued from $(F')^{n}$ and the lattices $\a_{v}\cO^{n}_{F'_{v}}$. Under $\t$, the restriction of $a$ at the generic point gives an $\om_{F}$-valued Hermitian form on $(F')^{n}$ which we denote by $T=T(a,\t)$. Then we define
\begin{equation}\label{eq:def Ea cE}
E_{a}(m(\cE), s,\Phi):=E_{T(a,\t)}(m(\a(\cE,\t)), s,\Phi).
\end{equation}
If we change $\t$ to $\g\t$ for some $\g\in M_{n}(F)=\GL_{n}(F')$, then $\a(\cE,\g\t)=\g\a(\cE,\t)$ and $T(a,\g\t)={}^{t}\ov\g^{-1}T(a,\t)\g^{-1}$. By  \eqref{Eis ch base}, we have
\begin{equation*}
E_{T(a,\g\t)}(m(\a(\cE,\g\t)), s,\Phi)=E_{T(a,\t)}(m(\a(\cE,\t)), s,\Phi)
\end{equation*}
for all $\g\in M_{n}(F)$. Therefore $E_{a}(m(\cE),s,\Phi)$ is well-defined.

Now suppose $a: \cE\inj \s^{*}\cE^\vee$ is an injective Hermitian map.  Let $(\cE_v,a)$ denote the Hermitian $\cO_{v}'$-lattice (valued in $\om_{\cO_{F_v}}:=\om_{X}\otimes_{\cO_{X}} \cO_{F_v}$) induced by $a$ at $v\in|X|$. Choosing a generator of the free $\cO_{F_v}$-module $\om_{\cO_{F_v}}$ of rank one, we obtain a Hermitian lattice $\cE_v$ (valued in  $\cO_{F_v}$) and hence  the density polynomial $\Den(T, \cE_v)$ defined by  \eqref{eq:CY}  relative to $F'_{v}/F_{v}$. We define   
 the density polynomial for $(\cE_v,a)$ as 
 \begin{equation}
 \label{def Den Ev}
 \Den_{v}(T, (\cE_v,a)):=\Den(T, \cE_v).
\end{equation} It is easy to see that the result is independent of the choice of the generator of  $\om_{\cO_{F_v}}$. 
We then define 
 $$
\Den(q^{-2s}, (\cE,a))=\prod_{v\in|X|} \Den_{v}(q_v^{-2s}, (\cE_v,a)),
$$
Note that the degree of $\Den(q^{-2s}, (\cE,a))$ (as a polynomial of $q^{-s}$) is $$\deg( \s^{*}\cE^\vee )-\deg(\cE)=-2\deg(\cE)+2n\deg\omega_X.$$

\begin{thm}\label{th:Eis Den} 
Let $\cE$ be a vector bundle over $X'$ of rank $n$.  
Then 
\begin{equation}\label{eq:reg}
E^{\rm reg}(m(\cE),s,\Phi)=\sum_{a: \cE\inj \s^{*}\cE^\vee}E_{a}(m(\cE), s,\Phi)
\end{equation}
where the sum runs over all injective Hermitian maps  $a: \cE\to \s^{*}\cE^\vee $. Moreover, we have
\begin{equation}\label{Ea Den}
E_{a}(m(\cE),s,\Phi)=\chi(\det(\cE))q^{-\deg(\cE)( s-n/2)-\frac{1}{2}n^2\deg\omega_X }\sL_n(s) ^{-1} \Den(q^{-2s}, (\cE,a)).
\end{equation}
\end{thm}
\begin{proof} 
From the definitions it is clear that
\begin{equation*}
E^{\rm reg}(m(\cE),s,\Phi)=\sum_{a:\cE\dashrightarrow\s^{*}\cE^{\vee}}E_{a}(m(\cE),s,\Phi)
\end{equation*}
where $a$ runs over {\em rational} Hermitian maps $\cE\dashrightarrow \s^* \cE^{\vee}$ that are generically nonsingular.

Now let $a:\cE\dashrightarrow\s^{*}\cE^{\vee}$ be such a rational nonsingular Hermitian map. We continue with the convention defining  \eqref{eq:def Ea cE} 
\begin{equation}\label{EaET}
E_{a}(m(\cE),s,\Phi)=E_{T}(m(\a), s,\Phi).
\end{equation}
By \eqref{Eis factor Wh} and \eqref{eq:ch base}, and noting that the character $\chi$ is trivial on the norm of $\BA_{F'}^\times$, we have 
\begin{equation}\label{ETWh}
E_{T}(m(\alpha),s,\Phi)=\chi(\det(\alpha))|\alpha|_{F'}^{-s+n/2} |\om_{X}|^{-\frac{1}{2}n^2} \prod_{v\in |X|}W_{T_{v}}(1,s,\Phi_{v}).
\end{equation}
If the ($\omega_{F_v}$-valued) Hermitian form $T_{v}$ does not have integral entries, then $W_{T_{v}}(1,s,\Phi_{v})=0$ (since $\Phi$ is invariant under $N_{n}(\cO_{v})$). Therefore $E_{T}(m(\alpha),s,\Phi)$ is nonzero only when  $T_{v}$ is integral for all $v$, i.e., $a$ is an everywhere regular Hermitian map $\cE\inj\s^{*}\cE^{\vee}$. This proves \eqref{eq:reg}.

For such $a: \cE\inj\s^{*}\cE^{\vee}$, by Lemma \ref{lem:Whi}, the right side of \eqref{ETWh} is
\begin{equation}\label{WhDen}
\chi(\det(\alpha))|\det (\alpha)|_{F'}^{-s+n/2} |\om_{X}|^{-\frac{1}{2}n^2} \prod_{v\in |X|}\frac{1}{\sL_{n, F_v'/F_v}(s) }  \Den(q_v^{-2s},  \cE_{v}),
\end{equation}
where, by choosing a generator of $\om_{\cO_{F_v}}$, the $\cO_{F_v'}$-module $\cE_v$ is endowed with the Hermitian form (valued in $\cO_{F_v}$) induced by the $\om_{X}$-valued Hermitian form $a$.

By $\chi(\det(\alpha))=\chi(\det(\cE))$, $|\det (\alpha)|_{F'}=q^{\deg(\cE)}$, and \eqref{def Den Ev}, we obtain
\begin{equation*}
\Den(q^{-2s}, (\cE,a))=\prod_{v\in |X|}\Den(q^{-2s}_{v}, (\cE_{v},a_{v}))=\prod_{v\in |X|}\Den(q_v^{-2s},  \cE_{v}).
\end{equation*}
Combining these facts with \eqref{EaET}, \eqref{ETWh} and \eqref{WhDen}, we get \eqref{Ea Den}.

\end{proof}


\section{Springer theory for torsion coherent sheaves}\label{sec: springer}
In this section we review the construction of the Springer sheaf on the moduli stack of torsion coherent sheaves on a curve following Laumon \cite{Lau87}. We also compute the Frobenius trace function of a particular summand of the Springer sheaf called the Steinberg sheaf. 

In this section let $X$ be any smooth (not necessarily projective or connected) curve over $k=\F_{q}$.  For $d\in \N$, let $X_{d}$ be the $d^{\mrm{th}}$ symmetric power of $X$.

\subsection{Local geometry of $\Coh_{d}$}\label{ss:Coh loc}
Let $\Coh_{d}=\Coh_{d}(X)$ be the moduli stack of torsion coherent sheaves on $X$ of length $d$. For any $k$-scheme $S$, $\Coh_d(S)$ is the groupoid of coherent sheaves on $X \times S$ whose pushforward to $S$ is locally free of rank $d$. 

Let $s^{\Coh}_{d}: \Coh_{d}\to X_{d}$ be the support map. Recall that for any $k$-scheme $S$, $[\gl_{d}/\GL_{d}](S)$ is the groupoid of $(V, T)$ where $V$ is a vector bundle of rank $d$ on $S$ and $T$ is an endomorphism of $V$. When $X=\A^{1}$, we have a canonical isomorphism
\begin{equation*}
\Coh_{d}(\A^{1})\cong [\gl_{d}/\GL_{d}]
\end{equation*}
given as follows. For $\cQ\in \Coh_{d}(\A^{1})(S)$, $\Gamma(\A^{1}_{S},\cQ)$ is a locally free rank $d$ $\cO_{S}$-module equipped with an endomorphism given by the affine coordinate $t$ for $\A^{1}$, giving an $S$-point of $[\gl_{d}/\GL_{d}]$; conversely, given an object $(V,T) \in [\gl_{d}/\GL_{d}]$ we may view $V$ as an $\cO_{S}[t]$-module $\cQ$ (viewed as a coherent sheaf on $\A^{1}_{S}$) with $t$ acting as $T$.  

Let $U\subset X_{\ov k}$ be open and $f:U\to \A^{1}_{\ov k}$ be an \'etale map. Such a pair $(U,f)$ is called an {\em \'etale chart} for $X_{\ov k}$. It induces a map $f^{\Coh}_{d}: \Coh_{d}(U)\to \Coh_{d}(\A^{1}_{\ov k})$ sending $\cQ$ to $f_{*}\cQ$ which is compatible with the symmetric power $f_{d}: U_{d}\to (\A^{1}_{\ov k})_{d}$ under $s^{\Coh}_{d}$. Let $\frD_{d,\A^{1}}\subset (\A^{1}_{\ov k})_{d}$ and $\frD_{d,U}\subset U_{d}$  be the discriminant divisors, i.e., they parametrize divisors with  multiplicities.  Clearly $\frD_{d,U}\subset f^{-1}_{d}(\frD_{d,\A^{1}})$, therefore we may write $f^{-1}_{d}(\frD_{d,\A^{1}})=\frD_{d,U}+\frR_{d,f}$ as Cartier divisors on $U_{d}$.  

\begin{lemma}\label{l:not in R}
Let $D$ be an effective divisor of degree $d$ on $U$. Then $D\in U_{d}\bs \frR_{d,f}$ if and only if for all pairs of distinct points $x,y$ in the support of $D$, $f(x)\ne f(y)$. 
\end{lemma}
\begin{proof}
Let $\pi_{d}: U^{d}\to U_{d}$ be the quotient map by the symmetric group $S_{d}$. We compute the divisor $\pi_{d}^{-1}(\frR_{d,f})$. Consider $U\times_{\A^{1}_{\ov k}} U$. Since $U$ is \'etale over $\A^{1}_{\ov k}$, $\D(U)\subset U\times_{\A^{1}} U$ is open and closed, hence we can write $U\times_{\A^{1}_{\ov k}}U=\D(U)\coprod \frR$. For geometric points $x,y\in U$,  $(x,y)\in \frR$ if and only if $x\ne y$ and $f(x)=f(y)$. 

For $1\le i<j\le d$, let $p_{ij}: U^{d}\to U\times U$ be the projection to the $i$th and $j$th coordinates. Let 
$$\wt\D_{ij}=p_{ij}^{-1}(U\times_{\A^{1}}U)=\D_{ij}\coprod \frR_{ij}, \mbox{ where }\D_{ij}=p_{ij}^{-1}(\D(U)), \frR_{ij}=p_{ij}^{-1}(\frR).$$ 
By definition, $\pi_{d}^{-1}(\frD_{d,\A^{1}})=\sum_{1\le i< j\le d}\wt\D_{ij}$, $\pi_{d}^{-1}(\frD_{d,U})=\sum_{1\le i< j\le d}\D_{ij}$ as divisors on $U^{d}$. Therefore
\begin{equation}\label{Rdf sum}
\pi_{d}^{-1}(\frR_{d,f})=\sum_{1\le i< j\le d}\frR_{ij}.
\end{equation}
From this we see, if $D=x_{1}+x_{2}+\cdots+x_{d}$, where $x_{i}\in U(\ov k)$, then $D\notin \frR_{d,f}$ if and only if $(x_{1},\cdots, x_{d})\notin \pi_{d}^{-1}(\frR_{d,f})$. By \eqref{Rdf sum}, the latter happens if and only if for all $1\le i<j\le d$, $(x_{i},x_{j})\notin \frR$, i.e., either $x_{i}=x_{j}$ or $f(x_{i})\ne f(x_{j})$. 
\end{proof}

Let $\Coh_{d}(U)^{f}\subset \Coh_{d}(U)$ be the preimage of $U_{d}\bs \frR_{d,f}$.  Then $\Coh_{d}(U)^{f}$ is an open substack of $\Coh_{d}(X)_{\ov k}=\Coh_{d}(X_{\ov k})$.

The following lemma shows that $\Coh_{d}(X)$ is \'etale locally isomorphic to $\Coh_{d}(\A^{1})\cong[\gl_{d}/\GL_{d}]$.
\begin{lemma}\label{l:Coh local}
\begin{enumerate}
\item For any \'etale chart $(U,f)$ of $X_{\ov k}$, the map $f^{\Coh}_{d}:\Coh_{d}(U)\to \Coh_{d}(\A^{1})_{\ov k}$  is \'etale when restricted to $\Coh_{d}(U)^{f}$.
\item The stack $\Coh_{d}(X)_{\ov k}$ is covered by the substacks $\Coh_{d}(U)^{f}$ for various \'etale charts $(U,f)$ of $X_{\ov k}$. 
\end{enumerate}
\end{lemma}
\begin{proof}
(1) For any $\cQ\in \Coh_{d}(U)^{f}(\ov k)$, the tangent map of $f_{d}^{\Coh}$  at $\cQ$ is $\Ext^{*}_{U}(\cQ,\cQ)\to \Ext^{*}_{\A^{1}}(f_{*}\cQ, f_{*}\cQ)$. By Lemma \ref{l:not in R}, different points in the support of $\cQ$ map to different points in $\A^{1}$, the above map is the direct sum of $\t_{z}: \Ext^{*}_{\cO_{U,z}}(\cQ_{z},\cQ_{z})\to \Ext^{*}_{\cO_{\A^{1},f(z)}}(\cQ_{z}, \cQ_{z})$ over $z\in \supp(\cQ)$. Since $f$ is \'etale at each such $z$, $\t_{z}$  are isomorphisms, and hence $f_{d}^{\Coh}$ is \'etale at $\cQ$ by the Jacobian criterion. 

(2) For every point $\cQ\in \Coh_{d}(X)(\ov k)$ we will construct an \'etale chart $(U,f)$ such that $\cQ\in \Coh_{d}(U)^{f}(\ov k)$.  Let $Z\subset X(\ov k)$ be the support of $\cQ$.  For $z\in Z$, let $\cO_{z}$ be the completed local ring of $X_{\ov k}$ at $z$ with a uniformizer $\vp_{z}$. The map of sheaves $r: \cO_{X_{\ov k}}\to \op_{x\in Z}\cO_{z}/\vp_{z}^{2}$ is surjective. Let $c:Z\to \ov k$ be any injective map of sets. Then there exists an open neighborhood $U_{1}$ of $Z$ and $f\in \cO(U_{1})$ such that $r(f)=(c(z)+\varpi_{z})_{z\in Z}$. Viewing $f$ as a map $f: U_{1}\to \A^{1}_{\ov k}$, it is then \'etale at $Z$, hence \'etale in an open neighborhood $U\subset U_{1}$ of $Z$, i.e., $(U,f)$ is an \'etale chart. Since  $\{f(z)=c(z)\}_{z\in Z}$ are distinct points in $\A^{1}_{\ov k}$, we see that $\cQ\in \Coh_{d}(U)^{f}_{\ov k}$ by Lemma \ref{l:not in R}.
 \end{proof}

\subsection{Springer theory for $\Coh_{d}$}
Let $\wt\Coh_{d}(X)$ be the moduli stack classifying a full flag of torsion sheaves  on $X$
\begin{equation*}
0\subset \cQ_{1}\subset \cQ_{2}\subset \cdots\subset \cQ_{d}=\cQ
\end{equation*}
where $\cQ_{j}$ has length $j$. Let 
\begin{equation*}
\pi^{\Coh}_{d,X}: \wt\Coh_{d}(X)\to \Coh_{d}(X)
\end{equation*}
be the forgetful map recording only $\cQ=\cQ_{d}$.  

\begin{lemma}[Laumon {\cite[Theorem 3.3.1]{Lau87}}]\label{lem: Laumon small} The stacks $\wt\Coh_{d}(X)$ and $\Coh_{d}(X)$ are smooth of dimension zero, and the map $\pi^{\Coh}_{d,X}$ is proper and small.
\end{lemma}
\begin{proof} It is enough to check the same statements after base change to $\ov k$. We give a quick alternative proof using Lemma \ref{l:Coh local}: for an \'etale chart $(U,f)$ (over $\ov k$), we have a diagram in which both squares are Cartesian:
\begin{equation*}
\xymatrix{     \wt\Coh_{d}(X)_{\ov k}   \ar[d]^{\pi^{\Coh}_{d,X}} &  \wt\Coh_{d}(U)^{f}   \ar[d]^{\pi^{\Coh}_{d,U}} \ar[r] \ar@{_{(}->}[l] & \wt\Coh_{d}(\A^{1})_{\ov k}   \ar[d]^{\pi^{\Coh}_{d,\A^{1}}} \\
\Coh_{d}(X)_{\ov k}   & \Coh_{d}(U)^{f} \ar[r]^{f^{\Coh}_{d}} \ar@{_{(}->}[l]  & \Coh_{d}(\A^{1})_{\ov k}}
\end{equation*}
Here $\wt\Coh_{d}(U)^{f}$ is the preimage of $\Coh_{d}(U)^{f}$ in $\wt\Coh_{d}(U)$. Since the horizontal maps are \'etale and the $\Coh_{d}(U)^{f}$ cover $\Coh_{d}(X)_{\ov k}$ by Lemma \ref{l:Coh local}, the desired properties of $\pi^{\Coh}_{d,X}$ follow from the same properties of $\pi^{\Coh}_{d,\A^{1}}$, which is the Grothendieck alteration $\pi_{\gl_{d}}: [\wt\gl_{d}/\GL_{d}]\to [\gl_{d}/\GL_{d}]$.
\end{proof}

Let $X_{d}^{\c}\subset X_{d}$ be the open subset of multiplicity-free divisors (i.e., the complement of $\frD_{d,X}$), and let $\Coh_{d}(X)^{\c}$ (resp. $\wt\Coh_{d}(X)^{\c}$) be its  preimage under $s^{\Coh}_{d}$ (resp. under $s^{\Coh}_{d}\c\pi_{d}^{\Coh}$). Then $\wt\Coh_{d}(X)^{\c}\to \Coh_{d}(X)^{\c}$ is an $S_{d}$-torsor. 

\begin{cor}[Laumon {\cite[p.320]{Lau87}}]\label{c:Coh Spr} The complex 
\begin{equation*}
\Spr_{d}:=R\pi^{\Coh}_{d,X*}\Qlbar\in D^{b}(\Coh_{d}(X), \Qlbar)
\end{equation*}
is a perverse sheaf on $\Coh_{d}(X)$ that is the middle extension from its restriction to $\Coh_{d}(X)^{\c}$. In particular, the natural $S_{d}$-action on $\Spr_d|_{\Coh_{d}(X)^{\c}}$ extends to the whole $\Spr_{d}$.
\end{cor}

\subsection{Springer fibers}\label{ss:Spr fiber}
Let $\cQ\in \Coh_{d}(X)(\ov k)$ with image $D$ in $X_{d}(\ov k)$, an effective divisor of degree $d$.  Let $Z=(\supp D)(\ov k)$. Let $\Sigma(Z)$ be the set of maps $y:\{1,2,\cdots, d\}\to Z$ such that $\sum_{i=1}^{d}y(i)=D$. Let $\cB_{\cQ}$  be the fiber of $\pi^{\Coh}_{d}$ over $\cQ$. Then $\cB_{\cQ}$ classifies complete flags of subsheaves $0\subset \cQ_{1}\subset \cQ_{2}\subset\cdots\subset \cQ_{d-1}\subset \cQ$. We write $\cohog{*}{-} := \cohog{*}{-; \Q_\ell}$ for $\ell$-adic cohomology (regarded as a graded $\Q_{\ell}$-vector space). By Corollary \ref{c:Coh Spr}, $\cohog{*}{\cB_{\cQ}}=(\Spr_{d})_{\cQ}$ carries an action of $S_{d}$.

For $y\in \Sigma(Z)$, let $\cB_{\cQ}(y)$ be the open and closed subscheme of $\cB_{\cQ}$ defined by the condition $\supp\cQ_{i}/\cQ_{i-1}=y(i)$. Then $\cB_{\cQ}$ is the disjoint union of $\cB_{\cQ}(y)$ for $y\in \Sigma(Z)$. Hence
\begin{equation*}
\cohog{*}{\cB_{\cQ}}\cong \bigoplus_{y\in \Sigma(Z)}\cohog{*}{\cB_{\cQ}(y)}.
\end{equation*}
There is an action of $S_{d}$ on $\Sigma(Z)$ by precomposing. 

\begin{lemma}\label{l:perm Coh}
The action of $w\in S_{d}$ on $\cohog{*}{\cB_{\cQ}}$ sends $\cohog{*}{\cB_{\cQ}(y)}$ to $\cohog{*}{\cB_{\cQ}(y\c w^{-1})}$, for all $y\in \Sigma(Z)$.
\end{lemma}
\begin{proof}
It suffices to check the statement for each simple reflection $s_{i}$ switching $i$ and $i+1$ ($1\le i\le d-1$). Let $\wt\Coh^{i}_{d}(X)$ be the moduli stack classifying chains of torsion coherent sheaves $0\subset\cQ_{1}\subset\cdots\subset\cQ_{i-1}\subset\cQ_{i+1}\subset\cdots \subset\cQ_{d}$ with $\cQ_{i}$ missing. Then we have a factorization
\begin{equation*}
\pi^{\Coh}_{d}: \wt\Coh_{d}(X)\xr{\r_{i}}\wt\Coh^{i}_{d}(X)\xr{\pi_{i}}\Coh_{d}(X).
\end{equation*}
The map $\r_{i}$ is an \'etale double cover over the open dense locus  $\wt\Coh^{i,\hs}_{d}(X)$  where $\cQ_{i+1}/\cQ_{i-1}$ (which has length $2$) is supported at two distinct points. The map $\r_{i}$ is small by Lemma \ref{lem: Laumon small}, and $R\r_{i*}\Qlbar$ carries an involution $\wt s_{i}$, which induces an involution $\wt s_{i}$ on $R\pi_{i*}R\r_{i*}\Qlbar\cong \Spr_{d}$. This action coincides with the action of $s_{i}$ over $\Coh_{d}(X)^{\c}$, hence coincides with $s_{i}$ everywhere. 

Let $\cB_{\cQ}^{i}=\pi_{i}^{-1}(\cQ)$. By considering the support of the successive quotients, we have a decomposition of $\cB_{\cQ}^{i}$ by the orbit set $\Sigma(Z)/\j{s_{i}}$. When $y\in\Sigma(Z)$ satisfies $y\ne y\c s_{i}$, the $s_{i}$-orbit $\y=\{y,y\c s_{i}\}$ gives an open and closed substack $\cB^{i}_{\cQ}(\y)\subset \cB_{\cQ}^{i}$, such that $\r^{-1}_{i}(\cB^{i}_{\cQ}(\y))=\cB_{\cQ}(y)\coprod\cB_{\cQ}(y\c s_{i})$, and $\cB^{i}_{\cQ}(\y)\subset \wt\Coh^{i,\hs}_{2d}$. Therefore in this case the action of $\wt s_{i}$ on $\cohog{*}{\r^{-1}_{i}(\cB^{i}_{\cQ}(\y))}$ comes from the involution on $\cB_{\cQ}(y)\coprod\cB_{\cQ}(y\c s_{i})$ that interchanges the two components. Since $\wt s_{i}=s_{i}$, this proves the statement for $s_{i}$ and $y$ such that $y\ne y\c s_{i}$. For $y=y\c s_{i}$ the statement is vacuous. This finishes the proof. 
\end{proof}

Let $\cQ_{x}$ be the direct summand of $\cQ$ supported at $x\in Z$. Let $d_{x}=\dim_{\ov k}\cQ_{x}$. Then for any $y\in \Sigma(Z)$, there is a canonical isomorphism over $\ov k$
\begin{equation}\label{BQy}
\b_{y}: \cB_{\cQ}(y)\cong \prod_{x\in Z}\cB_{\cQ_{x}}
\end{equation}
sending $(\cQ_{i})\in \cB_{\cQ}(y)$ to the full flag of $\cQ_{x}$ given by taking the summands of $\cQ_{i}$ supported at $x$.

The proof above implies the following statement that we record for future reference.
\begin{lemma}\label{l:min w} Let $y,y'\in \Sigma(Z)$ and let $w\in S_{d}$ be such that $y\c w^{-1}=y'$. Assume that $w$ has minimal length (in terms of the simple reflections $s_{1},\cdots, s_{d-1}$) among such elements (such $w$ is unique). Then the Springer action $w: \cohog{*}{\cB_{\cQ}(y)}\to \cohog{*}{\cB_{\cQ}(y')}$ is induced by the composition of the canonical isomorphisms $\b_{y,y'}:=\b_{y'}^{-1}\c\b_{y}: \cB_{\cQ}(y)\isom \cB_{\cQ}(y')$. In particular, $w$ sends the fundamental class of $\cB_{\cQ}(y)$ to the fundamental class of $\cB_{\cQ}(y')$.
\end{lemma}
\begin{proof}
Let $w^{-1}=s_{i_{1}}\cdots s_{i_{N}}$ be a reduced word for $w^{-1}$. Let $y_{j}=ys_{i_{1}}\cdots s_{i_{j}}$, $1\le j\le N$. Let $y_{0}=y$, and $y'=y_{N}$. Since $w$ has minimal length among $w'\in S_{d}$ such that $y\c w'^{-1}=y'$, for each $1\le j\le N$, $y_{j-1}\ne y_{j}$ for otherwise one could delete $s_{i_{j}}$ to shorten $w$. Since $y_{j}=y_{j-1}\c s_{i_{j}}\ne y_{j-1}$, the proof of Lemma \ref{l:perm Coh} shows that the Springer action of $s_{i_{j}}: \cohog{*}{\cB_{\cQ}(y_{j-1})}\to \cohog{*}{\cB_{\cQ}(y_{j})}$ is induced by the canonical isomorphism $\s_{j}=\b^{-1}_{y_{j}}\c\b_{y_{j-1}}: \cB_{\cQ}(y_{j-1})\isom \cB_{\cQ}(y_{j})$. The action $w: \cohog{*}{\cB_{\cQ}(y)}\to \cohog{*}{\cB_{\cQ}(y')}$, being the composition $\s_{N}\c\cdots\c\s_{1}$, is then equal to $\b_{y'}^{-1}\c\b_{y}:\cB_{\cQ}(y)\isom \cB_{\cQ}(y')$. 
\end{proof}


\begin{cor}\label{c:Sind} Let $y\in \Sigma(Z)$ and $S_{y}\cong \prod_{x\in Z}S_{d_{x}}$ be the stabilizer of $y$ under $S_{d}$.  There is an isomorphism of graded $S_{d}$-representations
\begin{equation*}
\cohog{*}{\cB_{\cQ}}\cong \Ind^{S_{d}}_{S_{y}}\cohog{*}{\cB_{\cQ}(y)}\cong \Ind^{S_{d}}_{S_{y}}\left(\bigotimes_{x\in Z}\cohog{*}{\cB_{\cQ_{x}}}\right).
\end{equation*}
Here on the right side, each factor $S_{d_{x}}$ of $S_{y}$ acts on the tensor factor indexed by $x$ (for $x\in Z$) via the Springer action in Corollary \ref{c:Coh Spr} on $(\Spr_{d_{x}})_{\cQ_{x}}$.
\end{cor}
\begin{proof} By Lemma \ref{l:perm Coh},  $\cohog{*}{\cB_{\cQ}(y\c w^{-1})}=w\cohog{*}{\cB_{\cQ}(y)}$ for $w\in S_{d}$. In particular, $\cohog{*}{\cB_{\cQ}(y)}$ is stable under $S_{y}$, and  $\cohog{*}{\cB_{\cQ}}\cong \Ind^{S_{d}}_{S_{y}}\cohog{*}{\cB_{\cQ}(y)}$. By \eqref{BQy} and the K\"unneth formula, we have $\cohog{*}{\cB_{\cQ}(y)}\cong \ot_{x\in Z}\cohog{*}{\cB_{\cQ_{x}}}$. 

It remains to check that the action of $S_{y}$ on $\cohog{*}{\cB_{\cQ}(y)}$ (as the restriction of the $S_{d}$-action on $\cohog{*}{\cB_{\cQ}}$) is the same as the tensor product of the Springer action of $S_{d_{x}}$ on  $\cohog{*}{\cB_{\cQ_{x}}}$. Since the action of $S_{d}$ on $\Sigma(Z)$ is transitive, it suffices to check this statement for a particular $y\in \Sigma(Z)$. 

Order points in $Z$ as $x_{1},\cdots, x_{r}$.  Let $y_{0}\in \Sigma(Z)$ be the unique increasing function, i.e. such that if $i< j$ then the index of $y_0(i)$ is less than or equal to the index of $y_0(j)$. Let $d_{i}=d_{x_{i}}$. Let $\d=(\d_{i})_{1\le i\le r}$ be the increasing sequence $\d_{i}=d_{1}+\cdots+d_{i}$. Let $\Coh_{\d}(X)$ be the moduli stack of partial chains of torsion coherent sheaves $0\subset \cQ_{\d_{1}}\subset \cdots\subset \cQ_{\d_{r-1}}\subset \cQ_{\d_{r}}=\cQ$ such that $\cQ_{\d_{i}}$ has length $\d_{i}$. The map $\pi^{\Coh}_{d}$ then factorizes as
\begin{equation*}
\wt\Coh_{d}(X)\xr{\pi_{\d}}\Coh_{\d}(X)\xr{\nu_{\d}}\Coh_{d}(X).
\end{equation*}
We have a Cartesian diagram
\begin{equation}\label{partial Spr Coh}
\xymatrix{\wt\Coh_{d}(X)\ar[d]^{\pi_{\d}}\ar[r]^-{\wt c} & \prod_{i=1}^{r} \wt\Coh_{d_{i}}\ar[d]^{\prod\pi^{\Coh}_{d_{i}}} \\
\Coh_{\d}(X)\ar[r]^-{c} & \prod_{i=1}^{r}\Coh_{d_{i}}(X)}
\end{equation}
where $c$ sends $(\cQ_{\d_{i}})$ to $(\cQ_{\d_{i}}/\cQ_{\d_{i-1}})$. By proper base change we have $R\pi_{\d*}\Qlbar\cong c^{*}(\boxtimes_{i=1}^{r}\Spr_{d_{i}})$, and the latter carries the Springer action of $S_{d_{1}}\times\cdots\times S_{d_{r}}=S_{y_{0}}$ (pulled back along $c$). Pushing forward along $\nu_{\d}$, this induces an action of $S_{y_{0}}$ on $R\nu_{\d*}R\pi_{\d*}\Qlbar=\Spr_{d}$. This action coincides with the restriction of the action of $S_{d}$ because both actions come from deck transformations over $\Coh_{d}(X)^{\c}$. 

Now $\nu_{\d}^{-1}(\cQ)$ contains the point $\cQ^{\da}\in \Coh_{\d}(X)$ where $\supp\cQ_{\d_{i}}/\cQ_{\d_{i-1}}=\{x_{i}\}$ for $1\le i\le r$. This is an isolated point in $\nu^{-1}_{\d}(\cQ)$, and $\cB_{\cQ}(y_{0})=\pi_{\d}^{-1}(\cQ^{\da})$.   Moreover, the isomorphism \eqref{BQy} is the one given by taking the Cartesian diagram \eqref{partial Spr Coh} and restricting to $\cQ^{\da}\in \Coh_{\d}(X)$. The above discussion shows that the action of $S_{y_{0}}\subset S_{d}$ on $\cohog{*}{\cB_{\cQ}(y_{0})}\subset \cohog{*}{\cB_{\cQ}}$ is the same as the Springer action of $\prod_{i}S_{d_{i}}$ on $\ot_{i}\cohog{*}{\cB_{\cQ_{x_{i}}}}$ via the isomorphism \eqref{BQy}.
\end{proof}

\subsection{The Steinberg sheaf} Let $\St_{d}\in D^{b}(\Coh_{d}(X),\Qlbar)$ be the direct summand of $\Spr_{d}$ where $S_{d}$ acts through the sign representation. We will describe its Frobenius trace function below.   The result is well-known but we include a self-contained proof.

We call $\cQ\in \Coh_{d}(X)(\ov k)$ {\em semisimple} if it is a direct sum of skyscraper sheaves at closed points.

\begin{prop}\label{p:St} 
\begin{enumerate}
\item If $\cQ\in \Coh_{d}(X)(\ov k)$ is not semisimple, then the stalk of $\St_{d}$ at $\cQ$ is zero.
\item Let $\cQ=\op_{v\in |X|}k_v^{\op d_{v}}\in \Coh_{d}(X)(k)$ be semisimple. Then the stalk of  $\St_{d}$ at $\cQ$ is $1$-dimensional, and the geometric Frobenius $\gFrob$ acts on the stalk $\St_{d,\cQ}$ by the scalar
\begin{equation*}
\e(\cQ)\prod_{v\in \supp\cQ}q_{v}^{d_{v}(d_{v}-1)/2}
\end{equation*}
where $\e(\cQ)\in\{\pm1\}$ is the sign of Frobenius permuting the geometric points in the support of $\cQ$ counted with multiplicities (as a multi-set of cardinality $d$).
\end{enumerate}
\end{prop}
\begin{proof}
Let $\cQ\in \Coh_{d}(X)(\ov k)$. Let $Z\subset X(\ov k)$ be the geometric points in the support of $\cQ$ and $y\in \Sigma(Z)$. By Corollary \ref{c:Sind} and Frobenius reciprocity,
\begin{align}\label{St factor}
\St_{d,\cQ} &\cong \Hom_{S_{d}}(\sgn, \Ind^{S_{d}}_{S_{y}}(\ot \cohog{*}{\cB_{\cQ_{x}}}))\cong \Hom_{S_{y}}(\sgn, \ot \cohog{*}{\cB_{\cQ_{x}}}) \nonumber \\
&\cong  \ot_{x\in Z}\Hom_{S_{d_{x}}}(\sgn,\cohog{*}{\cB_{\cQ_{x}}})\cong \ot_{x\in Z}\St_{d_{x}, \cQ_{x}}.
\end{align}

(1) By the above factorization of $\St_{d,\cQ}$, it suffices to show that if $\cQ_{x}$ is not semisimple, then $\St_{d_{x}, \cQ_{x}}=0$. By Lemma \ref{l:Coh local} we may reduce to the case $X=\A^{1}$ and $\cQ$ is concentrated at $x=0$. In this case $\Spr$ is the usual Springer sheaf on $[\gl_{d}/\GL_{d}]$, and $\cQ$ corresponds to a nilpotent element $e\in \scr{N}_{d}\subset \gl_{d}$ (here $\scr{N}_d$ is the nilpotent cone in $\gl_d$). It is well-known that $\St_{d}|_{\scr{N}_{d}}\cong \d_{0}[-d(d-1)]$ where $\d_{0}$ is the skyscraper sheaf at $0\in \scr{N}_{d}$. Indeed, by \cite[\S3.4, Corollary (b)]{BM}, for any nonzero nilpotent element $e\in \scr{N}_{d}$, the sign representation of $S_{d}$ does not appear in $\cohog{*}{\cB_{e}}$ ($\cB_{e}$ is the Springer fiber for $e$). For $e=0$, the sign representation of $S_{d}$ only appears in the top degree $\cohog{d(d-1)}{\cB_{e}}$, which is one-dimensional. This implies that $\St_{d}|_{\scr{N}}\cong \d_{0}[-d(d-1)]$.
In particular, $\St_{d,e}=0$ for all nilpotent $e\ne 0$.

(2) Let $\cQ\in \Coh_{d}(X)(k)$ be semisimple. Let $|Z|$ be the set of closed points in the support of $\cQ$. The above discussion shows that $\St_{d_{x},\cQ_{x}}\cong \cohog{top}{\cB_{\cQ_{x}}} =\cohog{d_{x}(d_{x}-1)}{\Fl_{d_{x}}}$ where $\Fl_{d_{x}}$ is the flag variety for $\GL_{d_x}$. By \eqref{St factor}, $\St_{d,\cQ}$ is $1$-dimensional and is in the top degree cohomology of $\cohog{*}{\cB_{\cQ}}$. Let 
\[
N=\dim \cB_{\cQ}=\sum_{x\in Z}d_{x}(d_{x}-1)/2=\sum_{v\in|Z|}\deg(v)d_{v}(d_{v}-1)/2
\]
(here $d_{v}=d_{x}$ for any $x|v$). Let $0\ne \xi\in \St_{d,\cQ}\subset \op_{y\in\Sigma(Z)}\cohog{2N}{\cB_{\cQ}(y)}$. Let $\Frob:\cB_{\cQ}\to\cB_{\cQ}$  be the Frobenius morphism. We need to show that $\Frob^{*}\xi=\e(\cQ)q^{N}\xi$. 

For $y\in \Sigma(Z)$, let $\y_{y}\in \cohog{2N}{\cB_{\cQ}(y)}$ be the fundamental class of $\cB_{\cQ}(y)$. Then $\Frob$ sends $\cB_{\cQ}(y)$ onto $\cB_{\cQ}(\Frob(y))$ (here $\Frob(y)$ means post-composing $y$ with the Frobenius permutation on $Z$), and hence $\Frob^{*}\y_{\Frob(y)}=q^{N}\y_{y}$. On the other hand, let $w\in S_{d}$ be the minimal length element such that $\Frob(y)=y\c w^{-1}$. By Lemma \ref{l:min w}, the Springer action of $w$ satisfies $w\y_{y}=\y_{\Frob(y)}$. Write  $\xi=(\xi_{y})_{y\in \Sigma(Z)}$ where $\xi_{y}=c_{y}\y_{y}$ for some $c_{y}\in \Qlbar^{\times}$. Since $w\xi=\sgn(w)\xi$, we see that $w\xi_{y}=\sgn(w)\xi_{\Frob(y)}$. Since  $w\y_{y}=\y_{\Frob(y)}$, we have $c_{y}=\sgn(w)c_{\Frob(y)}$. Therefore
\begin{equation*}
(\Frob^{*}\xi)_{y}=\Frob^{*}(\xi_{\Frob(y)})=c_{\Frob(y)}\Frob^{*}\y_{\Frob(y)}=q^{N}c_{\Frob(y)}\y_{y}=\sgn(w)q^{N}c_{y}\y_{y}=\sgn(w)q^{N}\xi_{y}.
\end{equation*}Note that, for any choice of $y$ and $w$ above, $\sgn(w)$ is equal to the sign of the Frobenius permutation of the multiset $\{y(i)\}_{1\le i\le d}$, which is $\e(\cQ)$. This implies $\Frob^{*}\xi=\e(\cQ)q^{N}\xi$ as desired.
\end{proof}

\section{Springer theory for Hermitian torsion sheaves}\label{sec: Herm}

In this section we extend the construction in \S\ref{sec: springer} to the case of Hermitian torsion sheaves. The main output is a perverse sheaf $\Spr^{\Herm}_{2d}$ on the moduli stack of Hermitian torsion sheaves with an action of $W_{d}:=(\Z/2\Z)^{d}\rtimes S_{d}$. We will compare the stalks and Frobenius trace functions of $\Spr^{\Herm}_{2d}$ with those of $\Spr_{d}$.

As in \S\ref{sec: springer}, $X$ is a smooth curve over $k$ (not necessarily projective or connected). Recall from \S \ref{ssec: notation} that $\nu: X'\to X$ is a finite map of degree $2$ that is assumed to be {\em generically  \'etale} (and $X'$ is smooth over $k$). We develop the Hermitian  Springer theory in this generality.  Starting from \S\ref{ss:Wd act} we will assume $\nu$ to be \'etale, which is the case needed for proving the main theorem. Let $\s\in\Gal(X'/X)$ be the nontrivial involution. 

\subsection{Local geometry of $\Herm_{d}$}\label{ss:Herm} 
We say that a map of torsion coherent sheaves $a \co \cQ \rightarrow \sigma^* \cQ^\vee$ is \emph{Hermitian} if $\sigma^* a^\vee= a$. 
Let $d\in \N$. Let
\begin{equation*}
\Herm_{d}(X'/X), \textup{ or simply } \Herm_{d}
\end{equation*}
be the moduli stack of pairs $(\cQ,h)$ where $\cQ$ is a torsion coherent sheaf on $X'$ of length $d$, and $h$ is a Hermitian isomorphism $\cQ\isom \s^{*}\cQ^{\vee}:=\s^{*}\un\Ext^{1}(\cQ,\om_{X'})$. 


We offer two other ways to think about a Hermitian torsion sheaf $(\cQ,h)$. For a torsion sheaf $\cQ$ on $X'$ of length $d$, the datum of a Hermitian isomorphism $h \co \cQ \xrightarrow{\sim} \sigma^* \cQ^\vee$ is equivalent to either:
  \begin{enumerate}
\item a symmetric $k$-bilinear nondegenerate pairing 
$$(\cdot,\cdot): V\times V\to k
$$ on $V=\Gamma(X', \cQ)$ satisfying $(fv_{1},v_{2})=(v_{1},\s^{*}(f)v_{2})$ for any  function $f$ on $X'$ regular near the support of $\cQ$, or 
\item an $\cO_{X'}$-sesquilinear nondegenerate pairing 
$$
\j{\cdot,\cdot}: \cQ\times \cQ\to \om_{F'}/\om_{X'}
$$ satisfying $\j{v_{1},v_{2}}=\s^{*}\j{v_{2},v_{1}}$. Here $\om_{F'}$ is the constant (and quasi-coherent) sheaf on $X'$ whose local sections are the rational $1$-forms on $X'$.
\end{enumerate}
For example, to pass from (2) to (1), form cohomology and apply the trace map $\cohog{0}{\cQ} \rightarrow k$. To pass from $h$ to the pairing in (1), observe that $\ul{\Ext}^1(Q, \omega_{X'}) \cong \ul{\Hom}(Q, \omega_{F'}/\omega_{X'})$ by the long exact sequence associated to $\omega_{X'} \rightarrow  \omega_{F'} \rightarrow \omega_{F'}/\omega_{X'}$. Therefore, $h$ is equivalent to a sesquilinear pairing $\cQ \times \cQ \rightarrow \omega_{F'}/\omega_{X'}$, which upon taking global sections and applying the residue map $\cohog{0}{\omega_{F'}/\omega_{X'}} \rightarrow k$ gives the pairing in (1). 

We refer to $h$, or any of the above equivalent data, as a \emph{Hermitian structure} on $\cQ$. 
  
We have the support map
\begin{equation*}
s^{\Herm}_{d}: \Herm_{d}(X'/X)\to (X'_{d})^{\s}.
\end{equation*}
Note that we have an isomorphism $ (X'_{2d})^{\s} \cong X_d$, sending a $\sigma$-invariant divisor on $X'$ to its descent on $X$. so we will also allow ourselves to view the support map as $s^{\Herm}_{2d}: \Herm_{2d}(X'/X)\to X_d$. 

\begin{remark} When $\nu$ is \'etale and $d$ is odd, $(X'_{d})^{\s}=\vn$ hence $\Herm_{d}(X'/X)=\vn$. 

In general, when $\nu$ is ramified over the points $R\subset X(\ov k)$, $(X'_{d})^{\s}$ has a decomposition into open and closed subschemes according to the parity of the multiplicities of the divisor at each point $x\in R$.
\end{remark}

Let $\A^{1}_{\sqrt{t}}\to \A^{1}_{t}$ be the square map of affine lines.
\begin{lemma}\label{l:Herm A1} There is a canonical isomorphism
\begin{eqnarray*}
\Herm_{d}(\A^{1}_{\sqrt{t}}/\A^{1}_{t})\cong [\fro_{d}/\Og_{d}].
\end{eqnarray*}
Here $\Og_{d}$ denotes the orthogonal group on a $d$-dimensional nondegenerate quadratic space over $k$ and $\fro_{d}$ is its Lie algebra (the stack $[\fro_{d}/\Og_{d}]$ is independent of the quadratic form).
\end{lemma}
\begin{proof} We give the map $\Herm_{d}(\A^{1}_{\sqrt{t}}/\A^{1}_{t})\to [\fro_{d}/\Og_{d}]$ on $S$-points. For an $S$-point $(\cQ,h)$ of $\Herm_{d}(\A^{1}_{\sqrt{t}}/\A^{1}_{t})$, $V=\Gamma(\A^{1}_{S}, \cQ)$ is a locally free $\cO_{S}$-module of rank $d$ with a nondegenerate symmetric self-duality $(\cdot,\cdot)$, i.e., an $\Og_{d}$-torsor over $S$. Moreover the action of $\sqrt{t}$ on $V$ satisfies $(\sqrt{t}v_{1},v_{2})=-(v_{1},\sqrt{t}v_{2})$ since $\s^{*}\sqrt{t}=-\sqrt{t}$. Therefore $\sqrt{t}$ gives a section of the adjoint bundle of $V$. It is easy to check this map is an equivalence of groupoids $\Herm_{d}(\A^{1}_{\sqrt{t}}/\A^{1}_{t})(S)\isom [\fro_{d}/\Og_{d}](S)$.  
\end{proof}

An {\em $\s$-equivariant \'etale chart} of $X'_{\ov k}$ is  a pair $(U,f)$, where $U\subset X_{\ov k}$ is an open subset (with preimage $U'\subset X'_{\ov k}$) and a regular function $f: U'\to \A^{1}_{\sqrt{t}, \ov k}$ that is an \'etale map satisfying $\s^{*}f=-f$. Note that if $\nu$ is \'etale, the image of $f$ has to lie in $\A^{1}_{\sqrt{t}, \ov k}\bs\{0\}$.  

A  $\s$-equivariant \'etale chart $(U,f)$ of $X'_{\ov k}$ induces a map 
\[
f^{\Herm}_{d}: \Herm_{d}(U'/U)\to \Herm_{d}(\A^{1}_{\sqrt{t}}/\A^{1}_{t})_{\ov k}
\]
by sending $\cQ$ to  $f_{*}\cQ$. Let $\Herm_{d}(U'/U)^{f}$  be the preimages of $(U'_{d})^{\s}\bs \frR^{\s}_{d,f}$ under the support maps (here $\frR_{d,f}\subset U'_{d}$ is defined using the map $f:U'\to \A^{1}_{\sqrt{t},\ov k}$; see \S\ref{ss:Coh loc}).

We have an analog of Lemma \ref{l:Coh local} in the Hermitian setting.

\begin{lemma}\label{l:Herm local}
\begin{enumerate}
\item Let $(U,f)$ be a $\s$-equivariant \'etale chart for $X'_{\ov k}$. Then the map $f^{\Herm}_{d}$ is \'etale when restricted to $\Herm_{d}(U'/U)^{f}$. 
\item Assume $\nu$ is ramified at at most one point (over $\ov k$). Then the stack $\Herm_{d}(X'/X)_{\ov k}$  is covered by $\Herm_{d}(U'/U)^{f}$ for various $\s$-equivariant \'etale charts $(U,f)$ of $X'_{\ov k}$. In particular,  $\Herm_{d}(X'/X)$ is \'etale locally isomorphic to  $[\fro_{d}/\Og_{d}]$.
\item In general, $\Herm_{d}(X'/X)$ is smooth of dimension $0$.
\end{enumerate}
\end{lemma}
\begin{proof}
(1) is similar to that of Lemma \ref{l:Coh local}(1). 

(2) We only need to construct for $(\cQ,h)\in \Herm_{d}(X'/X)(\ov k)$, a $\s$-equivariant \'etale chart $(U,f)$ such that $(\cQ,h)\in \Herm_{d}(U'/U)^{f}$. Let $Z'$ be its support  in $X'$; since the Hermitian property of $(\cQ, h)$ implies that $Z'$ is stable under $\sigma$, hence is the preimage of some $Z\subset X(\ov k)$ under $\nu$. Let $\cL=(\nu_{*}\cO_{X'_{\ov k}})^{\s=-1}$, a line bundle over $X$. Then the map $r=(r_{z})_{z\in Z}: \cL\to \op_{z\in Z}\cL_{z}/\vp_{z}^{2}$ is surjective. Let $Z_{0}\subset Z$ be the points over which $\nu$ is \'etale (so $Z-Z_{0}$ is empty or has one point). For each $z\in Z_{0}$, upon choosing $z'\in Z'$ over $z$, we may identify $\cL_{z}$ with $\cO_{z}=\ov k[[\vp_{z}]]$; changing $z'$ to $\s(z')$ changes the identification by a sign.  If $z\in Z-Z_{0}$, then $\cL_{z}\cong \sqrt{\vp_{z}}\ov k[[\vp_{z}]]$. Choose a map $c:Z_{0}\to \ov k^{\times}$ such that $c(z)^{2}$ are distinct for $z\in Z$. Let $f$ be a section of $\cL$ over some open neighborhood $U_{1}\subset X_{\ov k}$ of $Z$ such that $r_{z}(f)=c(z)+\vp_{z}$ for $z\in Z_{0}$ under one of the two identifications $\cL_{z}\cong \cO_{z}$, and $r_{z}(f)\equiv\sqrt{\vp_{z}}\mod\vp_{z}$ for $z\in Z-Z_{0}$.  Then $f$ restricts to an \'etale map $U'=\nu^{-1}(U)\to \A^{1}_{\sqrt{t},\ov k}$ for some open neighborhood $U$ of $Z$ in $U_{1}$. The definition of $\cL$ implies $\s^{*}(f)=-f$. Now $\{f(z')|z'\in Z'\}$ is the union of $\{c(z),-c(z)|z\in Z_{0}\}$ and possibly $\{0\}$ if $Z-Z_{0}$ is nonempty, which are all distinct points in $\A^{1}_{\sqrt{t},\ov k}$ by construction. We conclude that $(\cQ,h)\in \Herm_{d}(U'/U)^{f}$. 

(3) Let $R\subset X_{\ov k}$ be the ramification locus of $\nu$. The case $|R|\le 1$ is treated in (2), so we may assume $|R|\ge2$. For $x\in R$, let $Y_{x}=X\bs (R\bs \{x\})$ and let $Y'_{x}=\nu^{-1}(Y_{x})$.   For any function $\d:R\to \Z_{\ge0}$ such that $\sum_{x\in R} \d(x)=d$ we have a map $\frY_{\d}:=\prod_{x\in R}(Y'_{x,d(x)})^{\s}\to (X'_{d})^{\s}$ by adding divisors. Let $\frY_{\d}^{\hs}\subset \frY_{\d}$ be the open locus where the divisors indexed by different $x\in R$ are disjoint. It is clear that $\frY^{\hs}_{\d}\to (X'_{d})^{\s}$ is \'etale and for varying $\d$ their images cover $(X'_{d})^{\s}$. To prove the statement it suffices to show that the base change $\Herm_{d}(X'/X)|_{\frY_{\d}^{\hs}}$ is smooth of dimension $0$ for each $\d$. Observe that $\Herm_{d}(X'/X)|_{\frY_{\d}^{\hs}}$ is isomorphic to the restriction of the product $\prod_{x\in R}\Herm_{\d(x)}(Y'_{x}/Y_{x})$ to $\frY_{\d}^{\hs}$. Since $\nu|_{Y_{x}}: Y'_{x}\to Y_{x}$ is ramified at one point, by (2) $\Herm_{\d(x)}(Y'_{x}/Y_{x})$ is smooth of dimension $0$. Therefore $\Herm_{d}(X'/X)|_{\frY_{\d}^{\hs}}\cong \prod_{x\in R}\Herm_{\d(x)}(Y'_{x}/Y_{x})|_{\frY_{\d}^{\hs}}$ is smooth of dimension $0$.  
\end{proof}


\begin{remark}
There is an obvious notion of skew-Hermitian torsion sheaves. Let $\SH_{d}(X'/X)$ be the moduli stack of skew-Hermitian torsion sheaves on $(X',\s)$ of length $d$. Then $d$ is even if $\SH_{d}(X'/X)\ne\vn$. The skew-Hermitian analog of Lemma \ref{l:Herm local} says that $\SH_{d}(X'/X)$ is \'etale locally isomorphic to $[\sp_{d}/\Sp_{d}]$, at least  when $\nu$ is \'etale.
\end{remark}

\subsection{The Hermitian Springer sheaf}\label{ss:Herm Spr}  
Let $\wt\Herm_{d}(X'/X)$ be the moduli stack classifying $(\cQ,h)\in \Herm_{d}(X'/X)$ together with a full flag
\begin{equation*}
0\subset \cQ_{1}\subset \cdots \subset \cQ_{i}\subset \cdots \subset \cQ_{d-1}=\cQ_{1}^{\bot}\subset \cQ_{d}=\cQ,
\end{equation*} 
where $\cQ_{i}$ has length $i$ and $\cQ_{d-i}=\cQ_{i}^{\bot}$ (the orthogonal of $\cQ_{i}$ under the Hermitian pairing $\cQ\times \cQ\to \om_{F'}/\om_{X'}$). Let
\begin{equation*}
\pi^{\Herm}_{d}: \wt\Herm_{d}(X'/X)\to \Herm_{d}(X'/X)
\end{equation*}
be the forgetful map.  Let $\Herm_{d}(X'/X)^{\c}\subset \Herm_{d}(X'/X)$ be the preimage of the multiplicity-free part  $X'^{\c}_{d}$ under the support map $s^{\Herm}_{d}$.


We recall the Grothendieck alteration for the full orthogonal group $\Og(V,Q)$ for some vector space $V$ of dimension $d$ over $k$ and a nondegenerate quadratic form $Q$ on $V$. Let $\Fl(V,Q)$ be the flag variety that parametrizes full isotropic flags $V_{\bu}=(V_{1}\subset \cdots\subset V_{d}=V)$ in $V$. Note that when $d$ is even, this is different from the flag variety of $\SO(V,Q)$ but rather a double cover of it because there are two choices for $V_{d/2}$ given the rest of members of a flag. Let $\fro(V,Q)$ be the Lie algebra of $\Og(V,Q)$. Let $\wt\fro(V,Q)$ be the moduli space of pairs $(A, V_{\bu})\in \fro(V,Q)\times\Fl(V,Q)$ such that $AV_{i}\subset V_{i}$ for all $i$. The Grothendieck alteration for $\Og(V,Q)$ is the $\Og(V,Q)$-equivariant map $\wt\fro(V,Q)\to \fro(V,Q)$ forgetting the flag. The quotient stacks $[\fro(V,Q)/\Og(V,Q)]$ and $[\wt\fro(V,Q)/\Og(V,Q)]$ are canonical; they are independent of the quadratic form $Q$ and only depend on $d=\dim V$. Therefore we also write the Grothendieck alteration as $\pi_{\Og_{d}}: [\wt\fro_{d}/\Og_{d}]\to [\fro_{d}/\Og_{d}]$.

\begin{prop}\label{p:SH Spr}
\begin{enumerate}
\item If $\nu$ is ramified at at most one point, then the map $\pi^{\Herm}_{d}$ is \'etale locally isomorphic to the Grothendieck alteration $\pi_{\Og_{d}}:[\wt\fro_{d}/\Og_{d}]\to [\fro_{d}/\Og_{d}]$. 
\item In general, $\wt\Herm_{d}(X'/X)$ is smooth of dimension $0$ and $\pi^{\Herm}_{d}$ is a small map. In particular, the complex
\begin{equation*}
\Spr^{\Herm}_{d}:=R\pi^{\Herm}_{d,*}\Qlbar
\end{equation*}
is the middle extension perverse sheaf of its restriction to $\Herm_{d}(X'/X)^{\c}$.
\end{enumerate}
\end{prop}
\begin{proof} 
(1) The proof is similar to that of Corollary \ref{c:Coh Spr}. For a $\s$-equivariant \'etale chart $(U,f)$ for $X'_{\ov k}$ we have a diagram with Cartesian squares and \'etale horizontal maps by Lemma \ref{l:Herm local}
\begin{equation}\label{SH chart}
\xymatrix{    \wt\Herm_{d}(X'/X)_{\ov k} \ar[d]^{\pi^{\Herm}_{d,X'/X}}& \wt\Herm_{d}(U'/U)^{f} \ar[r]\ar@{_{(}->}[l] \ar[d]^{\pi^{\Herm}_{d,U'/U}}& \wt\Herm_{d}(\A^{1}_{\sqrt{t}}/\A^{1}_{t}\ar[d]^{\pi^{\Herm}_{d, \A^{1}_{\sqrt{t}}/\A^{1}_{t}}})_{\ov k}        \\
\Herm_{d}(X'/X)_{\ov k}    & \Herm_{d}(U'/U)^{f}\ar[r]^{f^{\Herm}_{d}}\ar@{_{(}->}[l] & \Herm_{d}(\A^{1}_{\sqrt{t}}/\A^{1}_{t})_{\ov k}     }
\end{equation}
Using the isomorphism in Lemma \ref{l:Herm A1}, we identify $\pi^{\Herm}_{d, \A^{1}_{\sqrt{t}}/\A^{1}_{t}}$ with the Grothendieck alteration $\pi_{\Og_{d}}$.  Since $\Herm_{d}(U'/U)^{f}$ cover $\Herm_{d}(X'/X)$ by Lemma \ref{l:Herm local}(2), $\pi^{\Herm}_{d,X'/X}$ is 
\'etale locally  isomorphic to $\pi^{\Herm}_{d, \A^{1}_{\sqrt{t}}/\A^{1}_{t}}=\pi_{\Og_{d}}$.

(2) We use the notation from the proof of Lemma \ref{l:Herm local}(3). We may assume $|R|\ge2$. For each function $\d: R\to \Z_{\ge0}$ satisfying $\sum_{x\in R}\d(x)=d$, the base change of $\pi^{\Herm}_{d}$ along $\frY_{\d}^{\hs}\to (X'_{d})^{\s}$ is a disjoint union of the restriction of
\begin{equation*}
\prod_{x}\pi^{\Herm}_{\d(x)}: \prod_{x\in R}\wt\Herm_{\d(x)}(Y'_{x}/Y_{x})\to \prod_{x\in R}\Herm_{\d(x)}(Y'_{x}/Y_{x})
\end{equation*}
to $\frY_{\d}^{\hs}$. The disjoint union comes from different ways to distribute $\supp(\cQ_{i}/\cQ_{i-1})$ among various factors in the product $\prod_{x}(Y'_{x,\d(x)})^{\s}$. By (1), $\pi^{\Herm}_{\d(x)}: \wt\Herm_{\d(x)}(Y'_{x}/Y_{x})\to \Herm_{\d(x)}(Y'_{x}/Y_{x})$ has smooth $0$-dimensional source and is small for each $x\in R$, the same holds true for the base change of $\pi^{\Herm}_{d}$ to $\frY_{\d}^{\hs}$. Since $\{\frY_{\d}^{\hs}\}_{\d}$ form an \'etale covering of $(X'_{d})^{\s}$, the same is true for $\pi^{\Herm}_{d}$.
\end{proof}

\subsection{The action of $W_{d}$} \label{ss:Wd act} 
From now on we assume that $\nu:X'\to X$ is an \'etale double cover.  In this case, $(X'_{2d})^{\s}$ can be identified with $X_{d}$ via $\nu^{-1}(D)\bij D$. Let
\begin{equation*}
W_{d}:=(\Z/2\Z)^{d}\rtimes S_{d}
\end{equation*}
be the Weyl group for $\Og_{2d}$. Then $\pi^{\Herm}_{2d}$ is a $W_{d}$-torsor over $\Herm_{2d}(X'/X)^{\c}$.  

\begin{cor}[of Proposition \ref{p:SH Spr}(1)]\label{c:Wd action} If $\nu:X'\to X$ is an \'etale double cover, then there is a canonical action of $W_{d}$ on $\Spr^{\Herm}_{2d}$ extending the geometric action on its restriction to $\Herm_{2d}(X'/X)^{\c}$.
\end{cor}

\begin{defn}\label{d:Spr isotypic} 
\begin{enumerate}
\item For any representation $\r$ of $W_{d}$, we define $\Spr^{\Herm}_{2d}[\r]$ to be the perverse sheaf  on $\Herm_{2d}(X'/X)$:
\begin{equation*}
\Spr^{\Herm}_{2d}[\r]=(\r^{\vee}\ot \Spr^{\Herm}_{2d})^{W_{d}}\in D^b(\Herm_{2d}(X'/X),\Qlbar).
\end{equation*}
\item We define the Hermitian  analog of the Springer sheaf $\Spr$ as
\begin{eqnarray*}
\HSpr_{d}:=(\Spr^{\Herm}_{2d})^{(\Z/2\Z)^{d}}\in D^b(\Herm_{2d}(X'/X),\Qlbar).
\end{eqnarray*}
\end{enumerate}
\end{defn}
Note that the notation shifts from the subscript $2d$ to $d$. By Corollary \ref{c:Wd action}, $\HSpr_{d}$ carries a canonical $S_{d}$-action.

\begin{remark} In the case $\nu$ is ramified with  ramification locus $R\subset X(\ov k)$, the stack $\Herm_{2d}(X'/X)_{\ov k}$ decomposes into the disjoint union of open and closed substacks $\Herm^{\ep}_{2d}(X'/X)_{\ov k}$ indexed by $\e: R\to \{0,1\}$ where the length of $\cQ_{x}$ has parity $\e(x)$ for all $x\in R$. Then $\Spr^{\Herm}_{2d}|_{\Herm^{\ep}_{2d}(X'/X)_{\ov k}}$ carries a canonical action of $W_{d'}$ where $d'=(d-\sum_{x\in R}\ep(x))/2$.  
\end{remark}


\subsection{The Springer fibers over $\Herm_{2d}$}\label{ss:Herm fiber} 
Let $(\cQ,h)\in \Herm_{2d}(X'/X)(\ov k)$ and consider its Hermitian Springer fiber
\begin{equation*}
\cB^{\Herm}_{\cQ}:=\pi^{\Herm,-1}_{2d}(\cQ,h).
\end{equation*}
This is a proper scheme over $\ov k$. In this subsection we prove the Hermitian analogs of results in \S\ref{ss:Spr fiber}.

Let $Z'=\supp\cQ\subset X'(\ov k)$.  Let $D'=s^{\Herm}_{2d}(\cQ)\in (X'_{2d})^{\s}(\ov k)$, which is of the form $D'=\nu^{-1}(D)$ for some $D\in X_{d}(\ov k)$. Let $Z=\nu(Z')$, the support of $D$. Write $D'=\sum_{z\in Z'}d_{z}z$. 

Let $\Sigma(Z')$ be the set of maps $y': \{1,2,\cdots,2d\}\to Z'$ satisfying $y'(2d+1-i)=\s(y'(i))$ for all $i$ and $\sum_{i=1}^{2d}y'(i)=D'$. Identifying $W_{d}$ with permutations of $\{1,2,\cdots,2d\}$ commuting with the involution $i\mapsto 2d+1-i$, we get an action of $W_{d}$ on $\Sigma(Z')$ by $w: y'\mapsto y'\c w^{-1}$.

Similarly let $\Sigma(Z)$ be the set of maps $y:\{1,\cdots, d\}\to Z$ such that $\sum_{i=1}^{d}y(i)=D$. Then the natural map $\Sigma(Z')\to \Sigma(Z)$ (sending $y'$ to $y$ defined by $y(i)=\nu(y'(i))$) is a $(\Z/2\Z)^{d}$-torsor. 

For $y'\in\Sigma(Z')$, let $\cB^{\Herm}_{\cQ}(y')$ be the subscheme of $\cB^{\Herm}_{\cQ}$ consisting of isotropic flags $\cQ_{\bu}$ such that $\supp(\cQ_{i}/\cQ_{i-1})=y'(i)$ for all $1\le i\le 2d$. Then we have a decomposition into open and closed subschemes
\begin{equation*}
\cB^{\Herm}_{\cQ}=\coprod_{y'\in \Sigma(Z')}\cB^{\Herm}_{\cQ}(y').
\end{equation*}
Accordingly we get a decomposition of cohomology
\begin{equation*}
\cohog{*}{\cB^{\Herm}_{\cQ}}=\bigoplus_{y'\in \Sigma(Z')}\cohog{*}{\cB^{\Herm}_{\cQ}(y')}.
\end{equation*}

\begin{lemma}\label{l:permute comp}
The action of $w\in W_{d}$ on $\cohog{*}{\cB^{\Herm}_{\cQ}}$ sends the direct summand $\cohog{*}{\cB^{\Herm}_{\cQ}(y')}$ to the direct summand $\cohog{*}{\cB^{\Herm}_{\cQ}(y'\c w^{-1})}$.
\end{lemma}
\begin{proof}
It suffices to check the statement for each simple reflection $s_{i}$, $i=1,\cdots, d$. Here, for $1\le i\le d-1$,  $s_{i}=(i,i+1)(2d-i,2d+1-i)$; for $i=d$, $s_{d}=(d,d+1)$. For $1\le i\le d$, let $\wt\Herm^{i}_{2d}$ be the moduli stack classifying isotropic flags that only misses the terms $\cQ_{i}$ and $\cQ_{2d-i}$ (for $i=d$ only misses $\cQ_{d}$). Then we have a factorization
\begin{equation*}
\pi^{\Herm}_{2d}: \wt\Herm_{2d}\xr{\r_{i}}\wt\Herm^{i}_{2d}\xr{\pi_{i}}\Herm_{2d}.
\end{equation*}
The map $\r_{i}$ is an \'etale double cover over the open dense locus  $\wt\Herm^{i,\hs}_{2d}$  where $\cQ_{i+1}/\cQ_{i-1}$ (which has length $2$) is supported at two distinct points. The map $\r_{i}$ is small, and $R\r_{i*}\Qlbar$ carries an involution $\wt s_{i}$, which induces an involution $\wt s_{i}$ on $R\pi_{i*}R\r_{i*}\Qlbar\cong R\pi^{\Herm}_{2d*}\Qlbar$. This action coincides with the action of $s_{i}$ over $\Herm_{2d}^{\c}$, hence coincides with $s_{i}$ everywhere. 

Let $(\cQ,h)\in\Herm_{2d}(\ov k)$, and $\cB_{\cQ}^{i}=\pi_{i}^{-1}(\cQ,h)$. We have a decomposition of $\cB_{\cQ}^{i}$ by the orbit set $\Sigma(Z')/\j{s_{i}}$. When $y'\in\Sigma(Z')$ satisfies $y'\ne y'\c s_{i}$, the $s_{i}$-orbit of $\y'=\{y',y'\c s_{i}\}$ gives an open closed substack $\cB^{i}_{\cQ}(\y')\subset \cB_{\cQ}^{i}$, such that $\r^{-1}_{i}(\cB^{i}_{\cQ}(\y'))=\cB_{\cQ}(y')\coprod\cB_{\cQ}(y'\c s_{i})$, and $\cB^{i}_{\cQ}(\y')\subset \wt\Herm^{i,\hs}_{2d}$. Therefore in this case the action of $\wt s_{i}$ on $\cohog{*}{\r^{-1}_{i}(\cB^{i}_{\cQ}(\y'))}$ comes from the involution on $\cB_{\cQ}(y')\coprod\cB_{\cQ}(y'\c s_{i})$ that interchanges the two components. Since $\wt s_{i}=s_{i}$, this proves the statement for $s_{i}$ and $y'$ such that $y'\ne y'\c s_{i}$. For $y'=y'\c s_{i}$ the statement is vacuous. This finishes the proof. 
\end{proof}

Choose $Z^{\sh}\subset Z'$ such that $Z^{\sh}\coprod\s(Z^{\sh})=Z'$. Then for each $x\in Z$ there is a unique $x^{\sh}\in Z^{\sh}$ above $x$. For  $y'\in \Sig(Z')$ with image $y\in \Sig(Z)$, we have an isomorphism
\begin{equation}\label{gZ}
\g_{Z^{\sh}, y'}: \cB^{\Herm}_{\cQ}(y')\isom\cB_{\nu_{*}(\cQ|_{Z^{\sh}})}(y)\cong \prod_{x\in Z}\cB_{\cQ_{x^{\sh}}}
\end{equation}
mapping  $(\cQ_{i})_{1\le i\le 2d}$ to the (non-strictly increasing) flag  $(\cQ_{i,x^{\sh}})$ of $\cQ_{x^{\sh}}$.

If $y',y''\in \Sig(Z')$, the composition
\begin{equation*}
\g_{y',y''}:=\g_{Z^{\sh},y''}^{-1}\c\g_{Z^{\sh}, y'}:\cB^{\Herm}_{\cQ}(y')\isom\cB^{\Herm}_{\cQ}(y'')
\end{equation*}
is independent of the choice of $Z^{\sh}$.

\begin{lemma}\label{l:min w Herm}
Let $y',y''\in \Sig(Z')$ and let $w\in W_{d}$ be a minimal length element such that $y''=y'\c w^{-1}$. Then the Springer action $w: \cohog{*}{\cB^{\Herm}_{\cQ}(y')}\to  \cohog{*}{\cB^{\Herm}_{\cQ}(y'')}$ is induced by the isomorphism $\g_{y', y''}$.
\end{lemma}
\begin{proof}
Similar to the proof of Lemma \ref{l:min w}.
\end{proof}

\subsection{Comparing stalks of $\HSpr_{d}$ and $\Spr_{d}$}\label{ss:HSpr stalk}
In this subsection we abbreviate $\Herm_{2d}(X'/X)$ by $\Herm_{2d}$. Consider the stack $\Lagr_{2d}$ classifying pairs $(\cL\subset\cQ)$ where $\cQ\in\Herm_{2d}$ and $\cL\subset \cQ$ is a Lagrangian subsheaf, i.e., $\cL$ has length $d$ and the composition $\cL\inj \cQ\xr{h}\s^{*}\cQ^{\vee}\to \s^{*}\cL^{\vee}$ is zero. We have natural maps
\begin{equation*}
\xymatrix{\Herm_{2d} & \Lagr_{2d} \ar[l]_{\upsilon_{2d}}\ar[r]^{\e'_{d}}& \Coh_{d}(X')\ar[r]^{\nu_{*}} & \Coh_{d}(X)}
\end{equation*}
where $\upsilon_{2d}(\cL\subset\cQ)=\cQ$ and $\e'_{d}(\cL\subset\cQ)=\cL$.  Let $\e_{d}=\nu_{*}\c \e'_{d}: \Lagr_{2d}\to \Coh_{d}(X)$.
 
Let $(X'_{d})^{\dm}\subset X'_{d}$ be the open subscheme parametrizing $D\in X'_{d}$ such that $D\cap\s(D)=\vn$.  
Let $\Lagr^{\dm}_{2d}\subset \Lagr_{2d}$ be the preimage of $(X'_{d})^{\dm}$ under the map $\Lagr_{2d}\xr{\e'_{d}}\Coh_{d}(X')\xr{s^{\Coh}_{d,X'}} X'_{d}$.  It is easy to see that $\e'_{d}$ restricts to an isomorphism
\begin{equation*}
\Lagr^{\dm}_{2d}\cong \Coh_{d}(X')^{\dm}
\end{equation*}
whose inverse is given by $\cL\mapsto (\cL\subset \cQ=\cL\op \s^{*}\cL^{\vee})$. 

Let $\upsilon^{\dm}_{2d}$ and $\e^{\dm}_{d}$ be the restrictions of $\upsilon_{2d}$ and $\e_{d}$ to $\Lagr^{\dm}_{2d}$. Thus we view $\Lagr^{\dm}_{2d}\cong \Coh_{d}(X')^{\dm}$ as a correspondence between $\Herm_{2d}$ and $\Coh_{d}(X)$
\begin{equation*}
\xymatrix{ \Herm_{2d} &  \Lagr^{\dm}_{2d}\ar[l]_{\upsilon^{\dm}_{2d}}\ar[r]^{\e^{\dm}_{d}} &\Coh_{d}(X) }
\end{equation*}
Note that both $\upsilon^{\dm}_{2d}$ and $\e^{\dm}_{d}$ are surjective. Now both $\Herm_{2d}$ and $\Coh_{d}(X)$ carry Springer sheaves $\HSpr_{d}$ and $\Spr_{d}$ with $S_{d}$-actions. The following proposition says that they become isomorphic after pullback to $\Lagr^{\dm}_{2d}$.

\begin{prop}\label{p:Spr HSpr}
There is a canonical $S_{d}$-equivariant isomorphism of perverse sheaves on $\Lagr^{\dm}_{2d}$
\begin{equation*}
\upsilon^{\dm,*}_{2d}\HSpr_{d}\cong \e^{\dm,*}_{d}\Spr_{d}.
\end{equation*}
\end{prop}
\begin{proof} The map $\pi^{\Herm}_{2d}$ factors as
\begin{equation*}
\pi^{\Herm}_{2d}: \wt\Herm_{2d}\xr{\l_{2d}} \Lagr_{2d}\xr{\upsilon_{2d}} \Herm_{2d}.
\end{equation*}
Let $\l^{\dm}_{2d}: \wt\Herm^{\dm}_{2d}\to \Lagr^{\dm}_{2d}$ be the restriction of $\l_{2d}$ to $\Lagr^{\dm}_{2d}$.
We have a commutative diagram
\begin{equation}\label{3sq}
\xymatrix{\Herm_{2d}\times_{X_{d}}X^{d}\ar[d] & \wt\Herm^{\dm}_{2d}\ar[l]\ar[dl]_{\pi^{\Herm,\dm}_{2d}}\ar[d]^{\l^{\dm}_{2d}}\ar[r]^-{\wt \e'^{\dm}_{d}} & \wt\Coh_{d}(X')^{\dm}\ar[d]^{\pi^{\Coh,\dm}_{X',d}}\ar[r]^{\nu_{*}} &\wt\Coh_{d}(X)\ar[d]^{\pi^{\Coh}_{d}} \\
\Herm_{2d} & \Lagr^{\dm}_{2d}\ar@/_2pc/[rr]^-{\e^{\dm}_{d}}\ar[l]_{\upsilon^{\dm}_{2d}}\ar[r]^-{\e'^{\dm}_{d}} & \Coh_{d}(X')^{\dm}\ar[r]^{\nu_{*}} &   \Coh_{d}(X)
}
\end{equation}
Here $\Coh_{d}(X')^{\dm}$ and $\wt\Coh_{d}(X')^{\dm}$ are the preimages of $(X'_{d})^{\dm}$ under the support map.  We have:
\begin{itemize}
\item The middle square is Cartesian. This is true even before restricting to the $\dm$ locus. 
\item Since $\e'^{\dm}_{d}$  is an isomorphism, so is $\wt \e'^{\dm}_{d}$.
\item The rightmost square is Cartesian. 
\end{itemize}
From these properties we get maps
\begin{eqnarray*}
\a: \upsilon^{\dm *}_{2d}\HSpr_{d} &\to& \upsilon^{\dm*}_{2d}\Spr^{\Herm}_{2d}=\upsilon^{\dm*}_{2d}R\upsilon_{2d*}R\l_{2d*}\Qlbar\to \upsilon^{\dm*}_{2d}R\upsilon^{\dm}_{2d*}R\l^{\dm}_{2d*}\Qlbar \\
&\to&  R\l^{\dm}_{2d*}\Qlbar\cong \e^{\dm*}_{d}\pi^{\Coh}_{d*}\Qlbar=\e^{\dm*}_{d}\Spr_{d}.
\end{eqnarray*}

To check $\a$ is an isomorphism, it suffices to check on geometric stalks. Let $\cL\in \Coh_{d}(X')^{\dm}(\ov k)$ with support $Z^{\sh}\subset X'(\ov k)$. Let $\cQ=\cL\op \s^{*}\cL^{\vee}\in \Herm_{2d}(\ov k)$. The support of $\cQ$ is $Z'=Z^{\sh}\coprod \s(Z^{\sh})$, with image $Z\subset X(\ov k)$. We have $(\cL\subset \cQ)\in \Lagr^{\dm}_{2d}(\ov k)$, with image $\nu_{*}\cL\in\Coh_{d}(X)(\ov k)$. The stalk of $\a$ at $(\cL\subset\cQ)$ is
\begin{equation}\label{astalk}
\a_{(\cL\subset\cQ)}: \cohog{*}{\cB^{\Herm}_{\cQ}}^{(\Z/2\Z)^{d}}\to\cohog{*}{\l^{-1}_{2d}(\cL\subset\cQ)}\isom \cohog{*}{\cB_{\nu_{*}\cL}}
\end{equation}
Recall $\cB^{\Herm}_{\cQ}=\coprod_{y'\in\Sig(Z')}\cB^{\Herm}_{\cQ}(y')$. Let $\Sig(Z^{\sh})\subset \Sig(Z')$ be the set of $y'$ such that $y'(i)\in Z^{\sh}$ for $1\le i\le d$. Then we have a natural bijection $\Sig(Z)\bij\Sig(Z^{\sh}), y\mapsto y^{\sh}$. The fiber $\l^{-1}_{2d}(\cL\subset\cQ)$ is the disjoint union $\coprod_{y\in\Sig(Z)}\cB^{\Herm}_{\cQ}(y^{\sh})$. Recall the isomorphism $\g_{Z^{\sh}, y^{\sh}}:  \cB^{\Herm}_{\cQ}(y^{\sh})\isom  \cB_{\nu_{*}\cL}(y)$ from \eqref{gZ}. Using these descriptions we may rewrite \eqref{astalk} as
\begin{equation}
\cohog{*}{\cB^{\Herm}_{\cQ}}^{(\Z/2\Z)^{d}}  \surj \op_{y'\in\Sig(Z^{\sh})}\cohog{*}{\cB^{\Herm}_{\cQ}(y')}\cong \op_{y\in \Sig(Z)}\cohog{*}{\cB_{\nu_{*}\cL}(y)}=\cohog{*}{\cB_{\nu_{*}\cL}}.
\end{equation}
It remains to show that the first map above is an isomorphism. But this follows from the fact that  $(\Z/2\Z)^{d}$ acts freely on $\Sig(Z')$ with orbit representatives $\Sig(Z^{\sh})$, and Lemma \ref{l:permute comp}. This shows that $\a$ is an isomorphism.

Finally we show that $\a$ is $S_{d}$-equivariant. By Proposition  \ref{p:SH Spr}, $\pi^{\Herm}_{2d}$ and hence $\l_{2d}$ is small,  $R\l_{2d*}\Qlbar$ is the middle extension from a dense open substack of $\Lagr_{2d}$. Therefore the same is true for $R\l^{\dm}_{2d*}\Qlbar$. Since $\a$ is an isomorphism, both $\upsilon^{\dm,*}_{2d}\HSpr_{d}$ and $\e^{\dm*}_{d}\Spr_{d}$ are middle extension perverse sheaves from a dense open substack of $\Lagr^{\dm}_{2d}$. To check that $\a$ is $S_{d}$-equivariant it suffices to check it over the dense open substack which is the preimage of the multiplicity-free locus $X_{d}^{\c}$. Over $X_{d}^{\c}$, all squares in \eqref{3sq} are Cartesian,  and all vertical maps are $S_{d}$-torsors. The $S_{d}$-actions on $\HSpr_{d}|_{\Herm_{2d}^{\c}}$ and $\Spr_{d}|_{\Coh_{d}(X)^{\c}}$ come from the vertical $S_{d}$-torsors in the diagram, so $\a$ is $S_{d}$-equivariant when restricted over $X_{d}^{\c}$. This finishes the proof.
\end{proof}

\subsection{Comparing Frobenius traces of $\HSpr_{d}$ and $\Spr_{d}$}\label{ss:trace comp Spr}
In this subsection we will prove a relationship between Frobenius trace functions for $\HSpr_{d}$ and for $\Spr_{d}$. Since these sheaves live on different stacks, to make sense of the comparison we first need to identify the isomorphism classes of the $k$-points of these stacks.

For a groupoid $\cG$, let $|\cG|$ denote its set of isomorphism classes.
\begin{lemma}\label{l:bij}
There is a canonical bijection of sets
\begin{equation*}
|\Herm_{2d}(X'/X)(k)|\cong |\Coh_{d}(X)(k)|
\end{equation*}
respecting the support maps to $X_{d}(k)$. 
\end{lemma}
\begin{proof} Let $\cP(d)$ be the set of partitions of $d\in\Z_{\ge0}$, and $\cP=\coprod_{d\ge0}\cP(d)$. Let $\cP_{|X|}$ be the set of functions $\l: |X|\to \cP$ such that $\l(v)$ is the zero partition for almost all $v\in |X|$. For $\l\in \cP_{|X|}$, let $|\l|=\sum_{v}|\l(v)|\deg(v)$. Let $\cP_{|X|}(d)$ be the subset of those $\l\in \cP_{|X|}$ with $|\l|=d$. Let $s_{d}: \cP_{|X|}(d)\to X_{d}(k)$ be the map sending  $\l$ to the divisor $\sum_{v}|\l(v)|v$. 

By taking the Jordan type of a torsion sheaf at each closed point, we get a canonical bijection $\L^{\Coh}_{d}: |\Coh_{d}(X)(k)|\isom\cP_{|X|}(d)$. The map $s_{d}^{\Coh}:|\Coh_d(X)(k)|\to X_{d}(k)$ corresponds to $s_{d}$ under this bijection.

We define a map $\L^{\Herm}_{d}: |\Herm_{2d}(X'/X)(k)|\to \cP_{|X|}(d)$ as follows. For $(\cQ,h)\in \Herm_{2d}(X'/X)(k)$ and $v\in |X|$,  let $\l(v)$ be the Jordan type of $\cQ_{v'}$ (the summand $\cQ_{v'}$ supported at $v'$) for any $v'\in |X'|$ above $v$. When $v$ is split in $X'$, the two choices of $v'$ give the same Jordan type. The support map $s_{2d}^{\Herm}: |\Herm_{2d}(X'/X)(k)|\to X_{d}(k)$ is the composition $s_{d}\c \L^{\Herm}_{d}$. 

We claim that $\L^{\Herm}_{d}$ is a bijection. Then $\L^{\Coh,-1}_{d}\c\L^{\Herm}_{d}: |\Herm_{2d}(X'/X)(k)|\isom |\Coh_{d}(X)(k)|$ is the desired bijection. 

To prove the claim, since any torsion coherent sheaf splits as a direct sum over the finitely many points in its support, it suffices to fix a closed point $v\in |X|$ and show that the set of isomorphism classes $|\Herm_{v}|$ of Hermitian torsion sheaves supported above $v$ maps bijectively (by the restriction of $\L^{\Herm}_{d}$) to $\cP$ which are supported at $v$. 

If $v$ splits into $v'$ and $v''$ in $|X'|$, then any $(\cQ,h)\in |\Herm_{v}|$ has the form $\cQ_{v'}\op \s^{*}\cQ_{v'}$ equipped with the canonical Hermitian structure. In this case we see that the isomorphism class of $(\cQ,h)$ is determined by the Jordan type of $\cQ_{v'}$, and conversely each Jordan type arises from some $(\cQ, h)$. 

If $v$ is inert with preimage  $v'\in |X'|$, let $(\cQ,h)\in |\Herm_{v}|$ be of length $d$ over $k_{v'}$. Then $V=\Gamma(X',\cQ)$ is a $d$-dimensional Hermitian $k_{v'}$-vector space with a self-adjoint nilpotent endomorphism $e$ given by the action of a uniformizer $\vp\in \cO_{v}$. Fix a $d$-dimensional Hermitian space $(V_{d},h)$ over $k_{v'}$ (unique up to isomorphism), then the isomorphism classes of $(\cQ,h)\in |\Herm_{v}|$ with length $d$ over $k_{v'}$ is in bijection with the adjoint orbits of the unitary group $\Ug(V_{d},h)(k_v)$ acting on the nilpotent cone $\cN(V_{d},h)(k_v)$ of self-adjoint nilpotent endomorphisms of $V_{d}$. Being a Galois twisted version of the usual nilpotent orbits under $\GL_{d}$, the orbits $\cN(V_{d},h)(k_v)/\Ug(V_{d},h)(k_v)$ are again classified by partitions of $d$ according the Jordan types of $e\in \cN(V_{d},h)(k_v)$ (here we use that the centralizer $C_{\GL_{d}}(e)$ is connected, and Lang's theorem implies $\cohog{1}{k_v, C_{\GL_{d}}(e)}=\{1\}$). Therefore the isomorphism class of $(\cQ,h)\in |\Herm_{v}|$ is determined by the Jordan type of $\cQ$, and conversely each Jordan type arises from a $(\cQ, h)$. This shows that $\L^{\Herm}_{d}$ is a bijection.  
\end{proof}

\subsubsection{Further notations}\label{sss:choose points}
Now let $(\cQ,h)\in \Herm_{2d}(k)$. We write $\cQ_{\ov k}$ for the base change of $\cQ$ over $X'_{\ov k}$, and adapt the notations $Z'\subset X'(\ov k),Z\subset X(\ov k),\Sig(Z'),\Sig(Z)$ from \S\ref{ss:Herm fiber}. Let $|Z'|$ and $|Z|$ be the set of closed points contained in $Z'$ and $Z$. We have a decomposition
\begin{equation*}
|Z|=|Z|_{s}\coprod |Z|_{i}
\end{equation*}
into split and inert places. For each closed point $v\in |Z|$ we choose a geometric point $x'_{v}\in Z'$ above $v$ and denote its image in $Z$ by $x_{v}$.

Let $\Frob:X'\to X'$ be the Frobenius morphism.  Let $Z^{\sh}$ be the following subset of $Z'$
\begin{equation*}
Z^{\sh}=\left\{F^{i}(x'_{v}): v\in |Z|, 0\le i<\deg(v)\right\}.
\end{equation*}
When $v$ splits into $v',v''$ in $|X'|$, with $x'_{v}|v'$, then $Z^{\sh}$ contains all geometric points above $v'$ and not any above $v''$. When $v$ is inert with preimage $v'\in|X'|$,  $Z^{\sh}$ contains half of the geometric points above $v'$ which form a chain under the Frobenius, starting with $x'_{v}$. Therefore $Z'=Z^{\sh}\coprod\s(Z^{\sh})$. For $x\in Z$ let $x^{\sh}\in Z^{\sh}$ be the unique element above $x$. This induces a section $\Sig(Z)\isom \Sig(Z^{\sh})\subset \Sig(Z')$ which we denote $y\mapsto y^{\sh}$.

Let $\cQ^{\flat}\in \Coh_{d}(X)(k)$ be the point corresponding to the isomorphism class of $(\cQ,h)$ under the bijection in  Lemma \ref{l:bij}. Then $\cQ^{\flat}_{\ov k}\cong\nu_{*}(\cQ|_{Z^{\sh}})$. Recall the isomorphism \eqref{BQy} for each $y\in \Sig(Z)$
\begin{equation*}
\b_{y}: \cB_{\cQ^{\flat}}(y)\isom \prod_{x\in Z}\cB_{\cQ^{\flat}_{x}}.
\end{equation*}
From this and the K\"unneth formula we get an identification
\begin{equation*}
\cohog{*}{\cB_{\cQ^{\flat}}(y)}=\bigotimes_{x\in Z} \cohog{*}{\cB_{\cQ^{\flat}_{x}}}.
\end{equation*}
By \cite[Corollary 2.3(1)]{HS77}, $\cohog{*}{\cB_{\cQ^{\flat}_{x}}}$ is concentrated in even degrees. Let $\cohog{*}{\cB_{\cQ^{\flat}}(y)}^{+}$ be the direct sum of all $\ot_{x\in Z} \cohog{2i_{x}}{\cB_{\cQ^{\flat}_{x}}}$ (for varying $(i_{x})_{x\in Z}$) such that 
\begin{equation}\label{sum ixv}
\mbox{$\sum_{v\in|Z|_{i}}i_{x_{v}}$ is even.}
\end{equation}
Similarly, let $\cohog{*}{\cB_{\cQ^{\flat}}(y)}^{-}$ be the direct sum of all $\ot_{x\in Z} \cohog{2i_{x}}{\cB_{\cQ^{\flat}_{x}}}$ such that the  quantity in \eqref{sum ixv} is odd. We have
\begin{equation}\label{HBy pm}
\cohog{*}{\cB_{\cQ^{\flat}}(y)}=\cohog{*}{\cB_{\cQ^{\flat}}(y)}^{+}\op\cohog{*}{\cB_{\cQ^{\flat}}(y)}^{-}.
\end{equation}
Taking direct sum over all $y\in \Sig(Z)$ we get a decomposition
\begin{equation}\label{HBQ pm}
\cohog{*}{\cB_{\cQ^{\flat}}}=\cohog{*}{\cB_{\cQ^{\flat}}}^{+}\op\cohog{*}{\cB_{\cQ^{\flat}}}^{-}.
\end{equation}
Note that this decomposition depends on the choice of a geometric point $x_{v}$ over each inert $v$. By Corollary \ref{c:Sind}, the action of $S_{y}\subset S_{d}$ on $\cohog{*}{\cB_{\cQ^{\flat}}(y)}$ preserves the decomposition \eqref{HBy pm} since it is the same as the tensor product of the Springer actions on each factor $\cohog{*}{\cB_{\cQ^{\flat}_{x}}}$. Therefore the decomposition \eqref{HBQ pm} is stable under the $S_{d}$-action. 

Now $(\cQ|_{Z^{\sh}}\subset \cQ)$ gives a geometric point of $\Lagr^{\dm}_{2d}$, which is not defined over $k$ if $|Z|_{i}\ne\vn$.  Using this geometric point in $\Lagr^{\dm}_{2d}$, Proposition \ref{p:Spr HSpr} gives an isomorphism $\a^{\sh}:=\a_{(\cQ|_{Z^{\sh}}\subset\cQ)}$ on the level of stalks (see \eqref{astalk}):
\begin{equation*}
\a^{\sh}: \cohog{*}{\cB_{\cQ}^{\Herm}}^{(\Z/2\Z)^{d}}\cong \cohog{*}{\cB_{\cQ^{\flat}}}.
\end{equation*}
This isomorphism is $S_{d}$-equivariant. Both sides now carry geometric Frobenius actions which we denote by $\gFrob_{\cQ}$ and $\gFrob_{\cQ^{\flat}}$, which are not necessarily intertwined under $\a^{\sh}$ because the point $(\cQ|_{Z^{\sh}}\subset \cQ)$ is not necessarily defined over $k$. The next result gives the relation between the two Frobenius actions.

\begin{prop}\label{p:trace relation}  Let $\th$ be the involution on $\cohog{*}{\cB_{\cQ^{\flat}}}$ which is $1$ on $\cohog{*}{\cB_{\cQ^{\flat}}}^{+}$ and $-1$ on $\cohog{*}{\cB_{\cQ^{\flat}}}^{-}$. Then under the isomorphism $\a^{\sh}$, $\gFrob_{\cQ}$ corresponds to $\gFrob_{\cQ^{\flat}}\c\th$.
\end{prop}
\begin{proof} Recall from the proof of Proposition \ref{p:Spr HSpr} that $\a^{\sh}$ is the composition
\begin{equation*}
\cohog{*}{\cB_{\cQ}^{\Herm}}^{(\Z/2\Z)^{d}}\cong \bigoplus_{y\in \Sig(Z)}\cohog{*}{\cB^{\Herm}_{\cQ}(y^{\sh})}\xr{\op\g_{Z^{\sh}, y^{\sh}}}\bigoplus_{y\in \Sig(Z)}\cohog{*}{\cB_{\cQ^{\flat}}(y)}.
\end{equation*}
Here $\g_{Z^{\sh},y^{\sh}}$ is defined in \eqref{gZ}.

Then $\Fr^{*}$ on $\cohog{*}{\cB^{\Herm}_{\cQ}}$ maps $\cohog{*}{\cB^{\Herm}_{\cQ}(\Fr(y^{\sh}))}$ to $\cohog{*}{\cB^{\Herm}_{\cQ}(y^{\sh})}$. Note that $\Fr(y^{\sh})$ and $(\Fr y)^{\sh}$ are in general different: if $y(i)=\Fr^{-1}(x_{v})$ for some inert  $v$, then $\Fr(y^{\sh})(i)=\s(x'_{v})$ while $(\Fr y)^{\sh}(i)=x'_{v}$. In other words, the only difference between $\Fr (y^{\sh})$ and $(\Fr y)^{\sh}$ is the switch of all $x'_{v}$ and $\s(x'_{v})$ for all inert $v$. Therefore there is a unique element $\t_{y}\in(\Z/2\Z)^{d}$ such that $(\Fr y)^{\sh}=\Fr (y^{\sh})\c \t_{y}$.

Identifying $\cohog{*}{\cB_{\cQ}^{\Herm}}^{(\Z/2\Z)^{d}}$ with $\op_{y\in \Sig(Z)}\cohog{*}{\cB_{\cQ}^{\Herm}(y^{\sh})}$, the geometric Frobenius endomorphism $\gFrob_{\cQ}$ on $\cohog{*}{\cB_{\cQ}^{\Herm}}^{(\Z/2\Z)^{d}}$ is the direct sum of the following compositions
\begin{equation*}
\cohog{*}{\cB_{\cQ}^{\Herm}((\Fr y)^{\sh})}\xr{\t_{y}}\cohog{*}{\cB_{\cQ}^{\Herm}(\Fr (y^{\sh}))} \xr{\Fr^{*}}\cohog{*}{\cB_{\cQ}^{\Herm}(y^{\sh})}
\end{equation*}
where the first map is the Springer action of $\t_{y}\in W_{d}$ on $(\Spr^{\Herm}_{2d})_{\cQ}$.

On the other hand, let $w_{y}\in W_{d}$ be the {\em minimal length} element such that $(\Fr y)^{\sh}= \Fr(y^{\sh})\c w_{y}$. Write
\begin{equation*}
\t_{y}=w_{y}u_{y}, \mbox{ for a unique } u_{y}\in \Stab_{W_{d}}((\Fr y)^{\sh})=\Stab_{S_{d}}(\Fr y)\subset S_{d}.
\end{equation*}
Note that $\Stab_{S_{d}}(Fy)=\prod_{x\in Z}S_{I_{x}}$ where $I_{x}\subset \{1,2,\cdots, d\}$ is the preimage of $x$ under $y$. An easy calculation shows that $u_{y}=(u_{y,x})_{x\in Z}$ where  $u_{y,x}\in S_{I_{x}}$ is 
\begin{equation*}
u_{y,x}=\begin{cases} w_{I_{x}}, & \mbox{ if }x=x_{v}, v\in |Z|_{i},\\ 1, & \mbox{ otherwise.}\end{cases}
\end{equation*}
Here $w_{I_{x}}\in S_{I_{x}}$ is the involution that reverses the order of $I_{x}$.

We use abbreviated notation
\begin{align*}
H(y') &:=\cohog{*}{\cB^{\Herm}_{\cQ}(y')}, \mbox{ for $y'\in\Sig(Z')$},\\
C(y) &:=\cohog{*}{\cB_{\cQ^{\flat}}(y)}, \mbox{ for $y\in\Sig(Z)$}.
\end{align*}
For each $y\in \Sig(Z)$, consider the following diagram
\begin{equation}\label{HC}
\xymatrix{H((\Fr y)^{\sh})\ar[r]^{u_{y}}\ar[d]^{\g_{Z^{\sh}, (\Fr y)^{\sh}}} & H((\Fr y)^{\sh})\ar[rr]^-{w_{y}}\ar[d]^{\g_{Z^{\sh}, (\Fr y)^{\sh}}} & & H(\Fr (y^{\sh})) \ar[dll]^{\g_{Z^{\sh}, \Fr (y^{\sh})}}\ar[r]^{\Fr ^{*}}\ar[d]^{\g_{\Fr (Z^{\sh}), \Fr (y^{\sh})}} & H(y^{\sh})\ar[d]^{\g_{Z^{\sh}, y^{\sh}}}\\
C(Fy)\ar[r]^{u_{y}} & C(Fy)\ar[rr]^{\d^{*}} && C(Fy)\ar[r]^{\Fr^{*}} & C(y)}
\end{equation}
The left square is commutative by the $S_{d}$-equivariance of $\a^{\sh}$ proved in Proposition \ref{p:Spr HSpr} (here $u_{y}\in \Stab_{W_{d}}((\Fr y)^{\sh})\subset S_{d}$). The middle upper triangle is commutative by Lemma \ref{l:min w Herm}. The map $\d^{*}$ is defined to make the lower middle triangle commutative. The right square is clearly commutative.  The composite of the upper row is the restriction of $\Frob_{\cQ}$ to $H((\Fr y)^{\sh})$. Let us compute the composite of the lower row. 

The map $\d^{*}$ is the pullback along the automorphism $\d$ of $\cB_{\cQ^{\flat}}(\Fr y)$ that makes the following diagram commutative
\begin{equation*}
\xymatrix{& \cB^{\Herm}_{\cQ}(\Fr (y^{\sh})) \ar[d]^{\g_{\Fr(Z^{\sh}), \Fr(y^{\sh})}}\ar[dl]_{\g_{Z^{\sh}, \Fr(y^{\sh})}}\\
\cB_{\cQ^{\flat}}(\Fr y) & \cB_{\cQ^{\flat}}(\Fr y) \ar[l]_{\d}}
\end{equation*}
Under the isomorphism $\b_{\Fr y}:\cB_{\cQ^{\flat}}(\Fr y)\isom \prod_{x\in Z}\cB_{\cQ^{\flat}_{x}}$, $\d$ is the product of automorphisms $\d_{x}$ for each $\cB_{\cQ^{\flat}_{x}}$. If $x$ is not of the form $x=x_{v}$ for $v\in |Z|_{i}$, 
$\d_{x}$ is the identity. If $x=x_{v}$ for some $v\in |Z|_{i}$, then $\cQ^{\flat}_{x}=\cQ_{x'_{v}}$ and the Hermitian structure on $\cQ$ gives an isomorphism $\io: \cQ_{\s(x'_{v})}\cong\cQ_{x'_{v}}^{\vee}$. On the other hand, $\Fr^{\deg(v)}$ gives an isomorphism $\phi: \cQ_{x'_{v}}\cong \cQ_{\s(x'_{v})}$ since $\s(x'_{v})=\Fr^{\deg(v)}(x'_{v})$.  Combining $\io$ and $\phi$ we get a perfect symmetric pairing $(\cdot,\cdot)_{x'_{v}}$ on $\cQ_{x'_{v}}$ itself. Then $\d_{x}$ sends a full flag $\cR_{\bu}$ of $\cQ^{\flat}_{x_{v}}=\cQ_{x'_{v}}$ to $\cR^{\bot}_{\bu}$ under the pairing $(\cdot,\cdot)_{x'_{v}}$.

By the description of $u_{y}$ and $\d$ above, under the isomorphism $\b_{Fy}$, the composition $\d^{*}\c u_{y}$ takes the form
\begin{equation*}
\ot(\d_{x}^{*}\c u_{y,x}): \bigotimes_{x\in Z}\cohog{*}{\cB_{\cQ^{\flat}_{x}}}\to\bigotimes_{x\in Z}\cohog{*}{\cB_{\cQ^{\flat}_{x}}}.
\end{equation*}
The automorphisms $\d_{x}^{*}\c u_{y,x}$ are the identity maps except when $x=x_{v}$ for some $v\in |Z|_{i}$, in which case $u_{y,x}=w_{I_{x}}$. Let us compute $\d_{x}^{*}\c w_{I_{x}}$ on $\cohog{*}{\cB_{\cQ^{\flat}_{x}}}$ for $x=x_{v}$ and $v\in |Z|_{i}$.  For this we switch to the following notation. Let $V=\cQ_{x'_{v}}=\cQ^{\flat}_{x_{v}}$, a vector space of dimension $m$ over $\ov k$. We have argued that $V$ carries a symmetric self-duality $(\cdot,\cdot)$; the action of a uniformizer at $x_{v}$ gives a nilpotent element $e\in \End_{\ov k}(V)$, which is self-adjoint under $(\cdot,\cdot)$. Let $\cB$ be the flag variety of $\GL(V)$ and $\cB_{e}$ be the Springer fiber of $e$. Then $S_{m}$ acts on $\cohog{*}{\cB_{e}}$. Let $w_{0}$ be the longest element in $S_{m}$. Let $\d: \cB_{e}\isom \cB_{e}$ be the map sending a flag $V_{\bu}$ to $V_{\bu}^{\bot}$. We claim that $\d^{*}\c w_{0}$ acts on $\cohog{2i}{\cB_{e}}$ by $(-1)^{i}$. Indeed, by \cite[Corollary 2.3(2)]{HS77}, the restriction map $\cohog{*}{\cB}\to \cohog{*}{\cB_{e}}$ is surjective and is clearly equivariant under $\d^{*}\c w_{0}$, so it suffices to show that $\d^{*}\c w_{0}$ acts by $(-1)^{i}$ on $\cohog{2i}{\cB}$. Since $\d^{*}\c w_{0}$ preserves the cup product on $\cohog{*}{\cB}$, it suffices to show it acts by $-1$ on $\cohog{2}{\cB}$ (which generates all of $\cohog{2*}{\cB}$ under the cup product). For $1\le j\le m$, let $\xi_{j}$ be  the Chern class of the tautological line bundle on $\cB$ whose fiber at $V_{\bu}$ is  $V_j/V_{j-1}$. Then $\cohog{2}{\cB}$ is spanned by $\xi_j$ for $1\le j\le m$. Now we have $w_{0}(\xi_{j})=\xi_{m+1-j}$ since $\cohog{2}{\cB}$ is the reflection representation of $S_{m}$, and $\d^{*}\xi_{j}=-\xi_{m+1-j}$ by the definition of $\d$. Therefore $\d^{*}\c w_{0}(\xi_{j})=-\xi_{j}$ for all $1\le j\le m$, which proves that $\d^{*}\c w_{0}$ acts by $-1$ on $\cohog{2}{\cB}$, and hence acts by $(-1)^{i}$ on $\cohog{2i}{\cB}$ and on $\cohog{2i}{\cB_{e}}$. 

The above argument shows that $\d^{*}\c u_{y}$ acts by $1$ on $\cohog{*}{\cB_{\cQ^{\flat}}(\Fr y)}^{+}$ and acts by $-1$ on $\cohog{*}{\cB_{\cQ^{\flat}}(\Fr y)}^{-}$. Therefore the bottom row of \eqref{HC} is $\Frob_{\cQ^{\flat}}\c\th$. By the commutativity of \eqref{HC}, $\Frob_{\cQ^{\flat}}\c\th$ corresponds under $\a^{\sh}$ to the composite of the top row, which is $\Frob_{\cQ}$. This finishes the proof.
\end{proof}

\section{Geometrization of local densities}\label{s:geom density}
The goal of this section is to give a sheaf-theoretic interpretation of the formula of Cho-Yamauchi on the representation density of Hermitian lattices, see Theorem \ref{th:Den Wd}. This will complete the geometrization of the analytic side of our proposed Siegel-Weil formula, at least for non-singular Fourier coefficients. The technical part of the proof of the theorem is a Frobenius trace calculation that uses properties of the Hermitian Springer action proved in \S\ref{ss:trace comp Spr}.

\subsection{Density function for torsion sheaves}\label{ssec: density function for torsion sheaves}

Following Remark \ref{r:Den tor}, for any Hermitian torsion sheaf $(\cQ,h)\in \Herm_{2d}(X'/X)(k)$, we may define the density  polynomial $\Den(T,\cQ)$ using the Cho-Yamauchi formula as follows. Let $\cQ_{v}$ be the summand of $\cQ$ supported over $v\in |X|$, we define
\begin{equation*}
\Den(T,\cQ):=\prod_{v\in|X|}\Den_{v}(T^{\deg(v)}, \cQ_{v})
\end{equation*}
where
\begin{equation*}
\Den_{v}(T, \cQ_{v})=\sum_{0\subset\cI\subset \cI^{\bot}\subset \cQ_{v}} T^{2\ell'_{v}(\cI)}\fm_{v}(t'_{v}(\cI^{\bot}/\cI);T). 
\end{equation*}
Here $I^\perp$ is the orthogonal of $I$ under the Hermitian form on $\cQ_v$, so in other words the sum is over all subsheaves of $\cQ_{v}$ that are isotropic under $h_{v}=h|_{\cQ_{v}}$, and we write $\fm_{v}(-)$ to emphasize the dependence of $\fm(a;T)$ on $F'_{v}/F_{v}$ (see Definition \ref{def:maT}). The functions $\ell'_{v}(-)$ and $t'_{v}(-)$ are the functions $\ell'(-)$ and $t'(-)$ defined in \S \ref{ssec: Cho-Yamauchi} for $F'_{v}/F_{v}$. 

Expanding the product into a summation, we see
\begin{equation}\label{Den sum}
\Den(T, \cQ)=\sum_{0 \subset \cI \subset \cI^{\bot} \subset \cQ} T^{\dim_{k}\cI}\prod_{v\in |X|}\fm_{v}(t'_v(\cI^{\bot}/\cI);T^{\deg(v)}).
\end{equation}

Given an injective Hermitian map $a \co \cE \inj \s^* \cE^\vee$, we have a Hermitian structure on the torsion sheaf $\coker(a)$ as follows. Applying $\s^{*}R\ul{\Hom}(-,\omega_{X'})$ to the short exact sequence
\begin{equation}\label{EQ}
0\to \cE\xr{a} \s^{*}\cE^{\vee}\to \cQ\to 0
\end{equation}
yields a short exact sequence
\begin{equation}\label{EQ'}
0\to \cE\xr{\s^{*}a^{\vee}} \s^{*}\cE^{\vee}\to \s^{*}\ul{\Ext}^{1}(\cQ,\omega_{X'})\to 0.
\end{equation}
Since $\s^*a^{\vee}=a$, we may identify \eqref{EQ} and \eqref{EQ'} and get an isomorphism $h_{\cQ}: \cQ\isom \s^{*}\cQ^{\vee}$. Using this, we may restate Theorem \ref{th:Eis Den} as follows.
\begin{thm}\label{th:Eis Den torsion} Let $\cE$ be a rank $n$ vector bundle over $X'$ and $a:\cE\inj\s^{*}\cE^{\vee}$ be an injective Hermitian map. Then
\[
E_{a}(m(\cE),s,\Phi)=\chi(\det(\cE))q^{-\deg(\cE)( s-n/2)-\frac{1}{2}n^2\deg(\omega_X) }
\times \sL_n(s) ^{-1} \Den(q^{-2s}, \coker(a)).
\]
\end{thm}

\subsection{Density sheaves}
We will define a graded perverse sheaf on $\Herm_{2d}(X'/X)$ whose Frobenius trace at $\cQ$ recovers $\Den(T,\cQ)$.  We will suppress $X'/X$ from the notations.

For $0 \leq i \leq d$, let $\Herm_{i,2d}$ be the stack classifying 
\[
\{ (\cI, (\cQ,h)) \in \Coh_i(X') \times \Herm_{2d} \co \text{$\cI\subset \cQ$ and is  isotropic under $h$}\}.
\]
We have the following maps
\begin{equation*}
\xymatrix{& \Herm_{i,2d}\ar[dl]_{\oll{f}_{i}}\ar[dr]^{\orr{f_{i}}=(\orr{f}'_{i},\orr{f}''_{i})} \\
\Herm_{2d}& & \Coh_i(X') \times  \Herm_{2(d-i)}
}
\end{equation*}
Here $\oll{f}_{i}$ takes $(\cI,(\cQ,h))$ to $(\cQ,h)$, $\orr{f}'_{i}$ takes it to $\cI$ and $\orr{f}''_{i}$ takes it to $\cI^{\bot}/\cI$ with the Hermitian structure induced from $h$.

Recall the perverse sheaf $\HSpr_{d}$ on $\Herm_{2d}$ from Definition \ref{d:Spr isotypic}(2). It is obtained from the Springer sheaf on $\Herm_{2d}$ by taking $(\Z/2\Z)^{d}$-invariants, and $\HSpr_{2d}$ carries an action of $S_{d}$.

\begin{defn}\label{def: K_Eis} We define the following graded virtual perverse sheaves on $\Herm_{2d}$. (In the notation below, the degree of the formal variable $T$ encodes the grading.) 
\begin{enumerate}
\item $\frP_{d}(T)=\bigoplus_{j=0}^{d}(-1)^{j}(\HSpr_{d})^{(S_{j}\times S_{d-j}, \sgn_{j}\bt \one)}T^{j}$.
\item $\cK^{\Eis}_{d}(T)=\bigoplus_{i=0}^{d} R\oll{f}_{i,!}R\orr{f}^{*}_{i}(\Qlbar\cdot T^{i}\bt \frP_{d-i}(T)).$
\end{enumerate}
\end{defn}

\begin{thm}\label{th:Den Wd} For any $\cQ\in \Herm_{2d}(k)$, we have
\begin{equation*}
\Den(T,\cQ)=\Tr(\Frob, \cK^{\Eis}_{d}(T)_{\cQ}).
\end{equation*}
\end{thm}
\begin{proof}
By the Grothendieck-Lefschetz trace formula, we have
\begin{equation*}
\Tr(\Frob, \cK^{\Eis}_{d}(T)_{\cQ})= \sum_{i=0}^d \sum_{(\cI,\cQ)\in \Herm_{i,2d}(k)}T^{\dim_{k}\cI}\Tr(\Frob, \frP_{d-i}(T)_{\cI^{\bot}/\cI}).
\end{equation*}
Comparing with the expansion \eqref{Den sum} for $\Den(T,\cQ)$, it suffices to prove, changing notation $d-i$ to $d$ and $\cI^{\bot}/\cI$ to $\cQ$, that
\begin{equation}\label{Pm}
\Tr(\Frob, \frP_{d}(T)_{\cQ})=\prod_{v\in |X|}\fm_{v}(t'_{v}(\cQ);T^{\deg(v)}).
\end{equation}
This will be proved in Proposition \ref{p:PQ trace}.
\end{proof}

The rest of the section is devoted to the proof of \eqref{Pm}. The idea is to relate the Frobenius trace to a similar Frobenius trace of a graded perverse sheaf on $\Coh_{d}(X)$ using results from \S\ref{ss:trace comp Spr}, and then calculate the latter explicitly.


\subsection{Comparison of two graded Frobenius modules}
For $(\cQ,h)\in \Herm_{2d}(k)$, write
\[
\frP_{\cQ}(T)=\frP_{d}(T)_{\cQ}=\bigoplus_{j=0}^{d} (-1)^{j}(\HSpr_{d})_{\cQ}^{(S_{j}\times S_{d-j}, \sgn_{j}\bt1)}T^{j} = \bigoplus_{j=0}^{d} (-1)^{j}\cohog{*}{\cB^{\Herm}_{\cQ}}^{(W_{j}\times W_{d-j}, \ov\sgn_{j}\bt1)}T^{j}.
\]
Here $\ov\sgn_{j}$ is the inflation of the sign representation of $S_{j}$ under $W_{j}\surj S_{j}$. We view $\frP_{\cQ}(T)$  as a $\Z$-graded virtual Frobenius module, with the $\Z$-grading indicated by the power of $T$.  Let $\cQ^{\flat}\in \Coh_{d}(X)(k)$ be in the isomorphism class that corresponds to $(\cQ,h)$ under the bijection in Lemma \ref{l:bij}. Define
\[
\frP_{\cQ^{\flat}}(T)=\bigoplus_{j=0}^{d} (-1)^{j}(\Spr_{d})_{\cQ^{\flat}}^{(S_{j}\times S_{d-j}, \sgn_{j}\bt1)}T^{j} =\bigoplus_{j=0}^{d} (-1)^{j}\cohog{*}{\cB_{\cQ^{\flat}}}^{(S_{j}\times S_{d-j}, \sgn_{j}\bt1)}T^{j}.
\]
Define the Frobenius traces 
\begin{equation*}
P_{\cQ}(T):=\Tr(\Frob,\frP_{\cQ}(T)), \quad  P_{\cQ^{\flat}}(T):=\Tr(\Frob,\frP_{\cQ^{\flat}}(T)) \in \Qlbar[T].
\end{equation*}
The goal is to get a relationship between $P_{\cQ}(T)$ and $P_{\cQ^{\flat}}(T)$. Note that $\frP_{\cQ^{\flat}}(T)$ is a special case of $\frP_{\cQ}(T)$ when the double cover $X'=X\sqcup X$ (and $\cQ$ is the direct sum of $\cQ^{\flat}$ on one copy of $X$ and $\cQ^{\flat,\vee}$ on the other). We shall apply Proposition \ref{p:trace relation} to express $P_{\cQ}(T)$ in terms of $P_{\cQ^{\flat}}(T)$. For this we need to calculate the decomposition of $\frP_{\cQ^{\flat}}(T)$ given in \eqref{HBQ pm}.

Recall notations $Z, Z', \Sig(Z)$ and $\Sig(Z')$ from \S\ref{ss:Herm fiber}.   First we show that $\frP_{\cQ^{\flat}}(T)$ factorizes according to the support of $\cQ$ in $X$. Let $|Z|$ be the set of closed points in $Z$.   For $v\in |Z|$, let $\cQ^{\flat}_{v}$ denote the direct summand of $\cQ^{\flat}$ whose support is over $v$. 

\begin{lemma}\label{l:Spr factor}
We have a natural isomorphism of graded $\Frob$-modules
\begin{equation*}
\frP_{\cQ^{\flat}}(T)\cong\bigot_{v\in |Z|}\frP_{\cQ^{\flat}_{v}}(T).
\end{equation*}
In particular,
\begin{equation*}
P_{\cQ^{\flat}}(T)=\prod_{v\in |Z|}P_{\cQ^{\flat}_{v}}(T).
\end{equation*}
\end{lemma}
\begin{proof}
Choose any $y\in \Sig(Z)$ and let $|y|$ be the resulting map $\{1,2,\cdots, d\}\to |Z|$. Let $S_{|y|}=\Stab_{S_{d}}(|y|)$, then $S_{|y|}\cong \prod_{v\in|Z|}S_{I_{v}}$, where $I_{v}=|y|^{-1}(v)$. Applying  Corollary \ref{c:Sind} to $\cQ^{\flat}$ and to each $\cQ^{\flat}_{v}$ gives an isomorphism of $(S_{d},\Frob)$-modules 
\begin{equation}\label{HBQflat ind}
\cohog{*}{\cB_{\cQ^{\flat}}}\cong \Ind^{S_{d}}_{S_{|y|}}\left(\bigot_{v\in|Z|}\cohog{*}{\cB_{\cQ^{\flat}_{v}}}\right).
\end{equation}
Here the factor $S_{I_{v}}$ of $S_{|y|}$ acts on $\cohog{*}{\cB_{\cQ^{\flat}_{v}}}$ by the Springer action. 

Write $H_{v}=\cohog{*}{\cB_{\cQ^{\flat}_{v}}}$ as a $(S_{I_{v}},\Frob)$-module. Let $d_{v}=\# I_{v}$. By Mackey theory, the $(S_{j}\times S_{d-j}, \sgn_{j}\bt1)$-isotypical part of the right side of \eqref{HBQflat ind} is a direct sum over double cosets $(S_{j}\times S_{d-j})\bs S_{d}/S_{|y|}$, which can be identified with the set of functions $i: |Z|\to \Z_{\ge0}$, $v\mapsto i_{v}\le d_{v}$, such that $\sum_{v}i_{v}=j$. The stabilizer of the $S_{j}\times S_{d-j}$-action on the orbit indexed by $i$ is isomorphic to $\prod_{v\in |Z|}S_{i_{v}}\times S_{d_{v}-i_{v}}$ (where the $S_{i_{v}}$ factor lies in $S_{j}$, $S_{d_{v}-i_{v}}$ lies in $S_{d-j}$, and $S_{i_{v}}\times S_{d_{v}-i_{v}}$ is naturally a subgroup of $S_{I_{v}}\cong S_{d_{v}}$). The contribution of the summand indexed by $i$ is 
\begin{equation*}
\bigot_{v\in |Z|}H_{v}^{(S_{i_{v}}\times S_{d_{v}-i_{v}}, \sgn_{i_{v}}\bt 1)}.
\end{equation*}
This implies that
\begin{eqnarray*}
\frP_{\cQ^{\flat}}(T)&\cong &\sum_{i:|Z|\to \Z_{\ge0},  i_{v}\le d_{v}}(-1)^{\sum i_{v}}\left(\bigot_{v\in |Z|} H_{v}^{(S_{i_{v}}\times S_{d_{v}-i_{v}}, \sgn_{i_{v}}\bt 1)}\right)T^{\sum i_{v}}\\
&\cong&\bigot_{v\in |Z|}\left(\sum_{i_{v}=0}^{d_{v}} (-1)^{i_{v}}H_{v}^{(S_{i_{v}}\times S_{d_{v}-i_{v}}, \sgn_{i_{v}}\bt 1)}T^{i_{v}}\right) = \bigot_{v\in |Z|}\frP_{\cQ^{\flat}_{v}}(T).
\end{eqnarray*}
\end{proof}

Let $\cQ_{v}$ be the direct summand of $\cQ$ supported over $v$.  Then $\cQ_{v}$ and $\cQ^{\flat}_{v}$ correspond under the bijection in Lemma \ref{l:bij}. 

For any Frobenius module $M$ with integral weights,  let $\Gr^{W}_{i}M$ be the pure of weight $i$ part of $M$. This notation applies also to graded Frobenius modules by taking $\Gr^{W}_{i}$ on each graded piece.

\begin{prop}\label{l:HSpr factor} Let $\cQ\in \Herm_{2d}(k)$.
\begin{enumerate}
\item We have
\begin{equation}\label{HSpr factor}
\frP_{\cQ}(T)\cong\bigot_{v\in |Z|}\frP_{\cQ_{v}}(T), \quad
P_{\cQ}(T)=\prod_{v\in |Z|}P_{\cQ_{v}}(T).
\end{equation}
\item If $v\in |Z|$ is split in $X'$, then
\begin{equation}\label{PQ split}
P_{\cQ_{v}}(T)=P_{\cQ^{\flat}_{v}}(T).
\end{equation}
\item If $v\in |Z|$ is inert in $X'$, then
\begin{equation}\label{PQ inert}
P_{\cQ_{v}}(T)=\sum_{i}(-1)^{i}\Tr(\Frob,\Gr^{W}_{2i\deg(v)}\frP_{\cQ^{\flat}_{v}}(T)).
\end{equation}
Moreover, $\Tr(\Frob,\Gr^{W}_{i}\frP_{\cQ^{\flat}_{v}}(T))=0$ if $i$ is not a multiple of $2\deg(v)$.
\end{enumerate}
\end{prop}
\begin{proof} We will use the notations from \S\ref{sss:choose points}. For each $y\in \Sig(Z)$, the summand $\cohog{*}{\cB_{\cQ^{\flat}}(y)}\cong \ot_{x\in Z}\cohog{*}{\cB_{\cQ^{\flat}_{x}}}$ is graded by multidegree $\bi: Z\to \Z_{\ge0}$
\begin{equation*}
\cohog{2\bi}{\cB_{\cQ^{\flat}}(y)}=\bigot_{x\in Z}\cohog{2\bi(x)}{\cB_{\cQ^{\flat}_{x}}}.
\end{equation*}
We define $\cohog{2\bi}{\cB_{\cQ^{\flat}}}$ as the direct sum of  $\cohog{2\bi}{\cB_{\cQ^{\flat}}(y)}$ over all $y\in \Sig(Z)$. Then each $\cohog{2\bi}{\cB_{\cQ^{\flat}}}$ is stable under $S_{d}$. Accordingly, $\frP_{\cQ^{\flat}}(T)$ decomposes into the direct sum of $\frP^{2\bi}_{\cQ^{\flat}}(T)$, which is by definition $\op_{j=0}^{d}(-1)^{j}\cohog{2\bi}{\cB_{\cQ^{\flat}}}^{(S_{j}\times S_{d-j}, \sgn_{j}\bt1)}T^{j}$. Let $\bi_{v}$ be the restriction of $\bi$ to those $x|v$, then under the factorization isomorphism in Lemma \ref{l:Spr factor} we have
\begin{equation}\label{frPi}
\frP^{2\bi}_{\cQ^{\flat}}(T)\cong \bigot_{v\in |Z|} \frP^{2\bi_{v}}_{\cQ^{\flat}_{v}}(T).
\end{equation}

(1) Recall the involution $\th$  on $\cohog{*}{\cB_{\cQ^{\flat}}}$ in Proposition \ref{p:trace relation}. Using the above notation,  we see that $\th$ acts on $\cohog{2\bi}{\cB_{\cQ^{\flat}}}$ by $\prod_{v\in |
Z|_{i}}(-1)^{\bi(x_{v})}$ (where $x_{v}\in Z$ is a chosen geometric point over $v$, as in \S\ref{sss:choose points}). Because of \eqref{frPi}, the action of $\th$ on $\frP_{\cQ^{\flat}}(T)$ factorizes as the tensor product of the similarly-defined $\th_{v}$ on each $\frP_{\cQ^{\flat}_{v}}(T)$. By Proposition \ref{p:trace relation}, $\frP_{\cQ}(T)$ is the Frobenius module obtained by modifying the Frobenius action on $\frP_{\cQ^{\flat}}(T)$ by composing with $\th$. By Lemma \ref{l:Spr factor}, this modified Frobenius structure on  $\frP_{\cQ^{\flat}}(T)$ is the tensor product of the similarly modified Frobenius modules $\frP_{\cQ^{\flat}_{v}}(T)$, which in turn are isomorphic to $\frP_{\cQ_{v}}(T)$ by Proposition \ref{p:trace relation}. This implies \eqref{HSpr factor}.

(2) From the definition we see that if $v$ is split, then $\th_{v}$ is the identity on $\frP_{\cQ^{\flat}_{v}}(T)$. Hence $\frP_{\cQ_{v}}(T)=\frP_{\cQ^{\flat}_{v}}(T)$ and $P_{\cQ_{v}}(T)=P_{\cQ^{\flat}_{v}}(T)$.

(3) Since the cohomology of $\cB_{\cQ^{\flat}_{v}}$ is pure by \cite[Corollary 2.3(3)]{HS77}, we see that $\Gr^{W}_{i}\frP_{\cQ^{\flat}_{v}}(T)$ is the sum of $\frP^{2\bi_{v}}_{\cQ^{\flat}_{v}}(T)$ where $\bi_{v}: \{x\in Z: x|v\}\to \Z_{\ge0}$ satisfies $\sum_{x|v}2\bi_{v}(x)=i$. The action of Frobenius sends $\frP^{2\bi_{v}}_{\cQ^{\flat}_{v}}(T)$ to $\frP^{2F^{*}\bi_{v}}_{\cQ^{\flat}_{v}}(T)$, where $F^{*}\bi_{v}$ means precomposing $\bi_{v}$ with Frobenius. Therefore only constant multi-degrees (i.e., constant functions $\bi_{v}$) contribute to the Frobenius trace of $\frP_{\cQ^{\flat}_{v}}(T)$. This implies $\Tr(\Frob, \Gr^{W}_{i}\frP_{\cQ^{\flat}_{v}}(T))=0$ unless $i$ is a multiple of $2\deg(v)$.

By the discussion above,
\begin{equation*}
P_{\cQ_{v}}(T)=\Tr(\Frob\c\th_{v},  \frP_{\cQ^{\flat}_{v}}(T))=\sum_{i\ge0}\Tr(\Frob\c\th_{v}, \frP^{(2i,2i,\cdots,2i)}_{\cQ^{\flat}_{v}}(T)).
\end{equation*}
Since $v$ is inert, $\th_{v}$ acts by $(-1)^{i}$ on $\frP^{(2i,2i,\cdots,2i)}_{\cQ^{\flat}_{v}}(T)$, hence
\begin{equation*}
\Tr(\Frob\c\th_{v}, \frP^{(2i,2i,\cdots,2i)}_{\cQ^{\flat}_{v}}(T))=(-1)^{i}\Tr(\Frob, \frP^{(2i,2i,\cdots,2i)}_{\cQ^{\flat}_{v}}(T)).
\end{equation*}
On the other hand,  $\Tr(\Frob, \Gr^{W}_{2i\deg(v)}\frP_{\cQ^{\flat}_{v}}(T))$ is the sum of $\Tr(\Frob, \frP^{2\bi_{v}}_{\cQ^{\flat}_{v}}(T))$ with total degree $\sum_{x|v}2\bi_{v}(x)=2i\deg(v)$. Since only constant multi-degrees contribute to the trace, we again conclude
\begin{equation*}
\Tr(\Frob, \Gr^{W}_{2i\deg(v)}\frP_{\cQ^{\flat}_{v}}(T))=\Tr(\Frob, \frP^{(2i,2i,\cdots,2i)}_{\cQ^{\flat}_{v}}(T)).
\end{equation*}
Combining the above identities we get \eqref{PQ inert}.
\end{proof}


\subsection{Calculation of $P_{\cQ^{\flat}}(T)$ and $P_{\cQ}(T)$}

Let $\cQ^{\flat}\in \Coh_{d}(X)(k)$ with support $Z\subset X(\ov k)$. For each $v\in |Z|$ recall $t_{v}(\cQ^{\flat})$ from \eqref{def t} for the local field $F_{v}$. 

\begin{prop}\label{p:PQ Coh trace} For $\cQ^{\flat}\in \Coh_{d}(X)(k)$, we have
\begin{equation*}
P_{\cQ^{\flat}}(T)=\prod_{v\in |Z|}(1-T^{\deg(v)})(1-q_{v}T^{\deg(v)})\cdots (1-q_{v}^{t_{v}(\cQ^{\flat})-1}T^{\deg(v)}).
\end{equation*}
\end{prop}
\begin{proof} We write $\Coh_{d}(X)$ simply as $\Coh_{d}$ in the proof. For $0\le j\le d$, consider the correspondence
\begin{equation*}
\xymatrix{ & \Coh_{j,d}\ar[dl]_{p}\ar[dr]^{r}\\
\Coh_{d} && \Coh_{j}\times \Coh_{d-j}.}
\end{equation*}
Here $\Coh_{j,d}$ classifies pairs $(\cQ_{j}\subset \cQ)$ of torsion sheaves of length $j$ and $d$ respectively, and the map $r$ sends $(\cQ_{j}\subset \cQ)$ to $(\cQ_{j}, \cQ/\cQ_j)$. We claim that
\begin{equation}\label{Sj ind}
\Spr_{d}^{(S_{j}\times S_{d-j}, \sgn_{j}\bt1)}\cong Rp_{*}Rr^{*}(\St_{j}\bt\Qlbar).
\end{equation}
Indeed, consider the following diagram with Cartesian parallelogram
\begin{equation*}
\xymatrix{ & & \wt\Coh_{d}\ar[dl]_{\pi_{j,d}}\ar[dr]^{\wt r}\\
& \Coh_{j,d}\ar[dl]_{p}\ar[dr]^{r} & & \wt\Coh_{j}\times\wt\Coh_{d-j}\ar[dl]_{\pi_{j}\times\pi_{d-j}}\\
\Coh_{d} && \Coh_{j}\times \Coh_{d-j}
}
\end{equation*}
Here $\pi_{j}=\pi_{j}^{\Coh}$, etc. The composition $p\c \pi_{j,d}=\pi_{d}=\pi^{\Coh}_{d}$. By the proper base change theorem, we get
\begin{equation*}
\Spr_{d}=Rp_{*}R\pi_{j,d*}\Qlbar\cong Rp_{*}Rr^{*}R(\pi_{j}\times\pi_{d-j})_{*}\Qlbar=Rp_{*}Rr^{*}(\Spr_{j}\bt\Spr_{d-j}).
\end{equation*}
This isomorphism is $S_{j}\times S_{d-j}$-equivariant by checking easily over the open substack $\Coh_{d}^{\c}$. Taking $(S_{j}\times S_{d-j},\sgn_{j}\bt1)$-isotypical parts of both sides we get \eqref{Sj ind}.

By Lemma \ref{l:Spr factor} it remains to compute $P_{\cQ^{\flat}_{v}}(T)$ for each $v$. Therefore we may assume $\cQ^{\flat}=\cQ^{\flat}_{v}$.  By \eqref{Sj ind} and the Grothendieck-Lefschetz trace formula we have
\begin{eqnarray*}
P_{\cQ^{\flat}_{v}}(T)&=& \sum_{j}(-1)^{j}\Tr(\Frob, (\Spr_{d})^{(S_{j}\times S_{d-j}, \sgn_{j}\bt1)}_{\cQ^{\flat}_{v}})T^{j}\\
 &= &\sum_{j}(-1)^{j}\Tr(\Frob, Rp_{*}Rr^{*}(\St_{j}\bt\Qlbar)_{\cQ_{v}^{\flat}})T^j\\
&=& \sum_{\substack{\cR\subset \cQ^{\flat}_{v}, \\ \dim_{k_v} (\cR)=j}}(-1)^{\deg(v)j}\Tr(\Frob, (\St_{\deg(v)j})_{\cR})T^{\deg(v)j}.
\end{eqnarray*}
By Proposition \ref{p:St}, $(\St_{\deg(v)j})_{\cR}$ is zero unless $\cR\cong k_{v}^{\op j}$, in which case the Frobenius trace is $(-1)^{j(\deg(v)-1)}q_{v}^{j(j-1)/2}$. Let $\cQ_{v}[\vp]$ be the kernel of the action of a uniformizer $\vp$ at $v$. Note $V:=\cQ^{\flat}_{v}[\vp]$ has dimension $t=t_{v}(\cQ^{\flat})$ over $k_{v}$.  Then we only need to sum over $k_{v}$-subspaces $\cR$ of $V$. The above sum becomes
\begin{equation*}
\sum_{j=0}^{t}(-1)^{j}q_{v}^{j(j-1)/2}T^{\deg(v)j}\#\Gr(j,V)(k_{v}).
\end{equation*}
Recall that the ``$q$-binomial theorem'' says that $q_v^{j(j-1)/2} \# \Gr(j, V)(k_{v})$ is the coefficient of $x^j$ in $(1+x)(1+q_vx)\ldots(1+q_v^{t-1}x)$. Making the change of variables $x=-T^{\deg(v)}$, we get
\[
\sum_{j=0}^{t}(-1)^{j}q_{v}^{j(j-1)/2}T^{\deg(v)j}\#\Gr(j,V)(k_{v}) =(1-T^{\deg(v)})(1-q_{v}T^{\deg(v)})\ldots(1-q_v^{t-1}T^{\deg(v)})
\]
as desired.
\end{proof}


Now we are ready to prove \eqref{Pm}.


\begin{prop}\label{p:PQ trace} For $\cQ\in \Herm_{2d}(X'/X)(k)$ with support $Z$, we have
\[
P_{\cQ}(T)=\prod_{v\in |X|}\fm_{v}(t'_{v}(\cQ); T^{\deg(v)}) = \prod_{v\in |Z|}\prod_{j=0}^{t'_{v}(\cQ)-1}(1-(\y(\vp_{v})q_{v})^{j}T^{\deg(v)}).
\]
\end{prop}
\begin{proof}
By \eqref{HSpr factor} it suffices to treat the case $\cQ$ is supported over a single place $v$. Let $\cQ^{\flat}\in \Coh_{d}(X)(k)$ be the corresponding point. If $v$ is split, we have $P_{\cQ}(T)=P_{\cQ^{\flat}}(T)$ by \eqref{PQ split}, and the formula follows from Proposition \ref{p:PQ Coh trace}. 

If $v$ is inert, let $t=t'_{v}(\cQ)=t_{v}(\cQ^{\flat})$. From the form of $P_{\cQ^{\flat}}(T)$ computed in Proposition \ref{p:PQ Coh trace}, which is valid for any extension of $k$,  the trace of the pure weight pieces of $\frP_{\cQ^{\flat}}(T)$ are separated by different powers of $q_{v}$, i.e., $q_{v}^{-i}\Tr(\Frob, \Gr^{W}_{2i\deg(v)}\frP_{\cQ^{\flat}}(T))$ is   the coefficient of $q_{v}^{i}$ in $\prod_{j=0}^{t-1}(1-q_{v}^{j}T^{\deg(v)})$.  By \eqref{PQ inert},
\begin{equation*}
P_{\cQ}(T)=\sum_{i}(-1)^{i}\Tr(\Frob, \Gr^{W}_{2i\deg(v)}\frP_{\cQ^{\flat}}(T)) =\prod_{j=0}^{t-1}(1-(-q_{v})^{j}T^{\deg(v)})
\end{equation*}
which is what we want because $\y(\vp_{v})=-1$ in this case.
\end{proof}

\part{The geometric side}

\section{Moduli of hermitian shtukas}\label{sec: unitary shtukas}

In this section we introduce some of the fundamental geometric objects in our story, in particular the moduli stacks of unitary (also called Hermitian) shtukas, which play an analogous role to that of unitary Shimura varieties in the work of Kudla-Rapoport. 

\subsection{Hermitian bundles}


We adopt the notation of \S \ref{ssec: notation}, and in particular for the remainder of the paper enforce the assumption that $X$ is \emph{proper}, and $\nu \co X' \rightarrow X$ is a finite \'{e}tale double cover (possibly trivial). 

\begin{defn} A rank $n$ \emph{Hermitian} (also called \emph{unitary}) bundle on $X \times S$ with respect to $\nu:X'\to X$ is a vector bundle $\Cal{F}$ of rank $n$ on $X' \times S$, equipped with an isomorphism $h \co \Cal{F} \xrightarrow{\sim} \sigma^* \Cal{F}^{\vee}$ such that $\sigma^* h^{\vee} = h$. We refer to $h$ as the \emph{Hermitian structure} on $\cF$. 
\end{defn}

We denote by $\Bun_{U(n)}$  the moduli stack of rank $n$ unitary bundles on $X$, which sends a test scheme $S$ to the groupoid of rank $n$ unitary bundles on $X \times S$.  The notation is justified by the following remark.

\begin{remark}\label{r:Un} There is an equivalence of categories between the groupoid of Hermitian bundles on $X \times S$, and the groupoid of $G$-torsors for the group scheme $G = U(n)$ over $X$ defined as 
\[
\{ g \in \Res_{X'/X}\GL_n \co \sigma({}^tg^{-1}) = g \}.
\]
Indeed, we choose a square root $\om^{1/2}_{X}$ of $\om_{X}$ (which exists over $k = \F_{q}$ by \cite[p.291, Theorem 13]{Weil}). Then $\cF_{1}:=\nu^{*}\om^{1/2}_{X'}$ is equipped with the canonical Hermitian structure $h_{1}: \cF_{1}\cong \s^{*}\cF_{1}\cong  \s^{*}\cF^{\vee}_{1}$, and $(\cF_{n},h_{n}):=(\cF_{1},h_{1})^{\op n}$ is a rank $n$ Hermitian bundle on $X$ whose automorphism group scheme is $U(n)$.  To a Hermitian bundle $(\cF,h)$ on $X \times S$, $\un\Isom_{X\times S}((\cF_{n}\boxtimes\cO_{S},h_{n}\boxtimes\Id), (\cF,h))$ (the scheme of unitary isometries)   is a right torsor for $U(n)$ over $X\times S$. Conversely, for a right $U(n)$-torsor $\cG$ over $X\times S$, the contracted product $\cG\twtimes{U(n)}\cF_{n}$ is a Hermitian bundle on $X \times S$. 
\end{remark}

\subsection{Hecke stacks} We now define some particular Hecke correspondences for $\Bun_{U(n)}$.

\begin{defn}\label{defn: Hk U(n)}
Let $r \geq 0$ be an integer. The \emph{Hecke stack} $\Hk_{U(n)}^{r}$ has as $S$-points the groupoid of the following data: 
\begin{enumerate}
\item $x_i' \in X'(S)$ for $i = 1, \ldots, r$, with graphs denoted by $\Gamma_{x'_i} \subset X' \times S$.
\item A sequence of vector bundles $\Cal{F}_0, \ldots, \Cal{F}_r$ of rank $n$ on $X' \times S$, each equipped with Hermitian structure $h_i \co \Cal{F}_i \xrightarrow{\sim} \sigma^* \Cal{F}_i^{\vee}$. 
\item Isomorphisms $f_i \co \Cal{F}_{i-1}|_{X' \times S - \Gamma_{x'_i} - \Gamma_{\sigma(x'_i)}} \xrightarrow{\sim} \Cal{F}_i|_{X' \times S - \Gamma_{x'_i}-\Gamma_{\sigma(x'_i)}}$, for $1 \leq i \leq r$, compatible with the Hermitian structures, with the following property: there exists a rank $n$ vector bundle $\Cal{F}^{\flat}_{i-1/2}$ and a diagram of vector bundles
\begin{equation}\label{eq: modification diagram 1}
\begin{tikzcd}
&\Cal{F}_{i-1/2}^{\flat} \ar[dl, "f_i^{\leftarrow}"', hook']  \ar[dr, "f_i^{\rightarrow}", hook ]  &  \\
\Cal{F}_{i-1} & & \Cal{F}_i
\end{tikzcd}
\end{equation}
such that $\coker(f_i^{\leftarrow})$ is locally free of rank $1$ over $\Gamma_{x'_i}$, and $\coker(f_i^{\rightarrow})$ is locally free of rank $1$ over $\Gamma_{\sigma(x'_i)}$. In particular, $f_i^{\leftarrow}$ and $f_i^{\rightarrow}$ are invertible upon restriction to $X' \times S - \Gamma_{x'_i}-\Gamma_{\sigma(x'_i)}$, and the composition 
\[
\Cal{F}_{i-1}|_{X' \times S - \Gamma_{x'_i}-\Gamma_{\sigma(x'_i)}} \xrightarrow{(f_i^{\leftarrow})^{-1}} \Cal{F}_{i-1/2}^{\flat}|_{X' \times S - \Gamma_{x'_i}-\Gamma_{\sigma(x'_i)}} \xrightarrow{f_i^{\rightarrow}} \Cal{F}_{i} |_{X' \times S - \Gamma_{x'_i}-\Gamma_{\sigma(x'_i)}} 
\]
agrees with $f_i$. 
\end{enumerate}

\end{defn}

\begin{remark}\label{rem: formulation of upper/lower}
Condition (3) above is equivalent to asking for the existence of a diagram 
\[
\begin{tikzcd}
\Cal{F}_{i-1}  \ar[dr, "h_i^{\leftarrow}"', hook] & & \Cal{F}_i  \ar[dl, "h_i^{\rightarrow}", hook']  \\
& \Cal{F}_{i-1/2}^{\sharp} 
\end{tikzcd}
\]
such that $\coker(h_i^{\leftarrow})$ is flat of length $1$ over $\Gamma_{\sigma(x'_i)}$, and $\coker(h_i^{\rightarrow})$ is flat of length $1$ over $\Gamma_{x'_i}$. In particular, $h_i^{\leftarrow}$ and $h_i^{\rightarrow}$ are invertible upon restriction to $X' \times S - \Gamma_{x'_i}-\Gamma_{\sigma(x'_i)}$, and the composition 
\[
\Cal{F}_{i-1}|_{X' \times S - \Gamma_{x'_i}-\Gamma_{\sigma(x'_i)}} \xrightarrow{h_i^{\leftarrow}} \Cal{F}_{i-1/2}^{\sharp}|_{X' \times S - \Gamma_{x'_i}-\Gamma_{\sigma(x'_i)}} \xrightarrow{(h_i^{\rightarrow})^{-1}} \Cal{F}_{i} |_{X' \times S - \Gamma_{x'_i}-\Gamma_{\sigma(x'_i)}} 
\]
agrees with $f_i$. 
\end{remark}

\begin{defn}[Terminology for modifications of vector bundles]
Given two vector bundles $\Cal{F}$ and $\Cal{F}'$ on $X' \times S$, we will refer to an isomorphism between $\Cal{F}$ and $\Cal{F}'$ on the complement of a relative Cartier divisor $D\subset X'\times S$ 
as a ``modification'' between $\Cal{F}$ and $\Cal{F}'$, and denote such a modification by $\Cal{F} \dashrightarrow \Cal{F}'$. Given $x,y \in X'(S)$, we say that the modification is \emph{``lower'' of length $1$ at $x$ and ``upper'' of length $1$ at $y$} if it is as in Definition \ref{defn: Hk U(n)} (3), i.e. if there exists a diagram 
\begin{equation}
\begin{tikzcd}
&\Cal{F}^{\flat} \ar[dl, "f^{\leftarrow}"', hook']  \ar[dr, "f^{\rightarrow}", hook ]  &  \\
\Cal{F} & & \Cal{F}'
\end{tikzcd}
\end{equation}
such that $\coker(f^{\leftarrow})$ is flat of length $1$ over $\Gamma_{x}$, and $\coker(f^{\rightarrow})$ is flat of length $1$ over $\Gamma_{y}$, and $\Cal{F} \dashrightarrow \Cal{F}'$ agrees with the composition 
\[
\Cal{F}|_{X' \times S - \Gamma_{x}-\Gamma_{y}} \xrightarrow{(f^{\leftarrow})^{-1}} \Cal{F}^{\flat}|_{X' \times S - \Gamma_x-\Gamma_y} \xrightarrow{f^{\rightarrow}} \Cal{F} |_{X' \times S - \Gamma_x-\Gamma_y} .
\]
The condition admits a reformulation as in Remark \ref{rem: formulation of upper/lower}.
\end{defn}

\subsection{Hermitian shtukas}
For a vector bundle $\Cal{F}$ on $X' \times S$, we denote by $\ft \Cal{F} := (\Id_{X'} \times \Frob_S)^* \Cal{F}$. If $\Cal{F}$ has a Hermitian structure $h \co \Cal{F} \xrightarrow{\sim} \sigma^* \Cal{F}^{\vee}$, then $\ft \Cal{F}$ is equipped with the Hermitian structure $\ft h$; we may suppress this notation when we speak of the ``Hermitian bundle'' $\ft \Cal{F}$. 

\begin{defn}
Let $r \geq 0$ be an integer. We define $\Sht_{U(n)}^{r}$ by the Cartesian diagram
\[
\begin{tikzcd}
\Sht_{U(n)}^{r} \ar[r] \ar[d] & \Hk_{U(n)}^{r} \ar[d] \\
\Bun_{U(n)} \ar[r, "{(\Id, \Frob})"]  & \Bun_{U(n)} \times \Bun_{U(n)}
\end{tikzcd}
\]
A point of $\Sht_{U(n)}^{r}$ will be called a ``$U(n)$-shtuka''. 

Concretely, the $S$-points of $\Sht_{U(n)}^{r}$ are given by the groupoid of the following data: 
\begin{enumerate}
\item $x'_i \in X'(S)$ for $i = 1, \ldots, r$, with graphs denoted $\Gamma_{x'_i} \subset X \times S$. These are called the \emph{legs} of the shtuka. 
\item A sequence of vector bundles $\Cal{F}_0, \ldots, \Cal{F}_n$ of rank $n$ on $X' \times S$, each equipped with a Hermitian structure $h_i \co \Cal{F}_i \xrightarrow{\sim} \sigma^* \Cal{F}_i^{\vee}$. 
\item Isomorphisms $f_i \co \Cal{F}_{i-1}|_{X' \times S - \Gamma_{x'_i} - \Gamma_{\s(x'_i)}} \xrightarrow{\sim} \Cal{F}_i|_{X' \times S - \Gamma_{x'_i}-\Gamma_{\s(x'_i)}}$ compatible with the Hermitian structure, which as modifications of the underlying vector bundles on $X' \times S$ are lower of length $1$ at $x'_i$ and upper of length $1$ at $\sigma(x'_i)$. 

\item An isomorphism $\varphi \co \Cal{F}_r \cong \ft \Cal{F}_0$ compatible with the Hermitian structure. 
\end{enumerate}
\end{defn}

\begin{lemma}\label{lem: shtuka empty parity} The stack $\Sht_{U(n)}^r$ is empty if and only if $r$ is odd. 
\end{lemma}

\begin{proof} We first treat the case $n=1$. Let $\Nm_{X'/X}: \Pic_{X'}\to \Pic_{X}$ be the norm map. Then $\Bun_{U(1)}\cong \Nm^{-1}_{X'/X}(\om_{X})$, hence it is a torsor under $\Prym(X'/X)=\ker(\Nm_{X'/X})$. Moreover, $\Sht_{U(1)}^r$ fits into a Cartesian square
\[
\begin{tikzcd}
\Sht_{U(1)}^r \ar[r] \ar[d] & \Bun_{U(1)} \ar[d, "\text{Lang}"] & \cF\ar[d, maps to] \\
(X')^r \ar[r] & \Prym(X'/X)  & \cF^{-1}\ot {}^{\t}\cF 
\end{tikzcd}
\]
with the bottom horizontal map sending $D \mapsto  \Cal{O}(D - \sigma D )$. If $X'$ is geometrically connected, then the stack $\Prym(X'/X)$ has two connected components, and by a result of Wirtinger, explained in \cite[\S 2]{Mum71}, the bottom horizontal map lands in the identity component if and only if $r$ is even. If $X'$ is geometrically disconnected (i.e. it is either $X \coprod X$ or $X_{k'}$), then we have $\pi_0(\Prym(X'/X)_{\ol{k}}) \cong \Z$, the Lang map lands in (therefore surjects onto, by Lang's Theorem) the identity component, and the bottom horizontal map hits the identity component if and only if $r$ is even. This shows that, in all cases, $\Sht_{U(1)}^r$ is empty if and only if $r$ is odd. 

For general $n$, taking determinant of a hermitian shtuka gives a map $\Sht^{r}_{U(n)}\to \Sht^{r}_{U(1)}$. From this we see that if $r$ is odd, then $\Sht_{U(n)}^r$ is empty for any $n$ since $\Sht^{r}_{U(1)}$ is empty.

On the other hand if $r$ is even, then $\Sht_{U(1)}^r$ is non-empty. If $n>1$, from an $S$ point of $\Sht_{U(1)}^r$, we can produce an $S$-point of $\Sht_{U(n)}^r$ by formation of direct sum with (the base change to $X\times S$ of) a unitary bundle of rank $n-1$ on $X$ (e.g. we can take $(\cF_{n-1},h_{n-1})$ from Remark \ref{r:Un}). 
\end{proof}

\subsection{Geometric properties}

\begin{lemma}\label{lem: smooth Bun_G}
The stack $\Bun_{U(n)}$ is smooth and equidimensional.  
\end{lemma}

\begin{proof}
The standard tangent complex argument, cf. \cite[Prop. 1]{Hei10}.
\end{proof}

\begin{lemma}\label{lem: shtuka geometry} 
(1) The projection map $(\pr_{X}, \pr_{r}):\Hk_{U(n)}^r \rightarrow (X')^r \times \Bun_{U(n)}$ recording $\{x_{i}\}$ and $(\cF_{r},h_{r})$ is smooth of relative dimension $r(n-1)$. 

(2) $\Sht_{U(n)}^r$ is a Deligne-Mumford stack locally of finite type. The map $\Sht_{U(n)}^r \rightarrow (X')^r$ is smooth, separated, equidimensional of relative dimension $r(n-1)$. 
\end{lemma}

\begin{proof}
The statements about $\Sht_{U(n)}^r$ being locally finite type and separated are well-known properties of moduli of $G$-shtukas for general $G$ \cite[Proposition 2.16 and Theorem 2.20]{V04}.\footnote{See also \cite[paragraph after Theorem 5.4]{YZ} for a sketch of the separatedness in a similar situation, which readily adapts here.}

Part (2) follows from (1) by \cite[Lemma 2.13]{Laff18}. 

So it suffices to check (1).  As a self-correspondence of $\Bun_{U(n)}$, $\Hk^{r}_{U(n)}$ is the $r$-fold composition of $\Hk^{1}_{U(n)}$. This allows us to reduce to the case $r=1$.  In this case, the map $(\pr_{X},\pr_{1}): \Hk^{1}_{U(n)}\to X'\times\Bun_{U(n)}$ exhibits $\Hk^{1}_{U(n)}$ as a $\PP^{n-1}$-bundle whose fiber over $(x', \cF_{1},h_{1})$ classifies hyperplanes in $\cF_{1,\s(x')}$. Indeed, a hyperplane in $\cF_{1, \s(x')}$ determines a lower modification at $\s(x')$, and the upper modification at $x'$ is then determined from the lower modification by the Hermitian structure. This shows that $(\pr_{X},\pr_{1})$ is smooth, separated and equidimensional of relative dimension $(n-1)$ in the case $r=1$, and the general case follows.

\end{proof}

\section{Special cycles: basic properties}\label{s:int prob}

In this section we define special cycles over the moduli stacks of hermitian shtukas, and construct corresponding cycle classes. The latter task is rather subtle, as the cycles are in most cases of a highly ``derived'' nature, with their ``virtual dimension'' differing significantly from their actual dimension. 

\subsection{Special cycles}

\begin{defn}\label{def: Z} Let $\Cal{E}$ be a rank $m$ vector bundle on $X'$.

 We define the stack $\Cal{Z}_{\Cal{E}}^{r}$ whose $S$-points are the groupoid of the following data: 
 \begin{itemize}
 \item A $U(n)$-shtuka with $(\{x'_1, \ldots, x'_r\}, \{\Cal{F}_0, \ldots, \Cal{F}_r\}, \{f_1, \ldots, f_r\}, \varphi) \in  \Sht_{U(n)}^{r}(S)$.
 \item Maps of coherent sheaves  $t_i \co \Cal{E} \boxtimes \cO_{S} \rightarrow \Cal{F}_i$ on $X'\times S$ such that the isomorphism $\varphi \co \Cal{F}_r \cong \ft \Cal{F}_0$ intertwines $t_r$ with $\ft t_0$, and the maps $t_{i-1}, t_{i}$ are intertwined by the modification $f_i \co  \Cal{F}_{i-1} \dashrightarrow \Cal{F}_{i}$ for each $i = 1, \ldots, r$, i.e. the diagram below commutes. 
\[
\begin{tikzcd}
\Cal{E}\boxtimes\cO_{S} \ar[d, "t_0"] \ar[r, equals] &  \Cal{E}\boxtimes\cO_{S} \ar[r, equals]  \ar[d, "t_1"]& \ldots \ar[d] \ar[r, equals] & \Cal{E}\boxtimes\cO_{S} \ar[r, "\sim"] \ar[d, "t_{r}"] & \ft (\Cal{E} \boxtimes\cO_{S})\ar[d, "\ft t_0"] \\
\Cal{F}_0 \ar[r, dashed, "f_0"] & \Cal{F}_1 \ar[r, dashed, "f_1"]  &  \ldots \ar[r, dashed, "f_r"] &  \Cal{F}_{r} \ar[r, "\sim"] & \ft \Cal{F}_0  
\end{tikzcd}
\]
\end{itemize}
In the sequel, when writing such diagrams we will usually just omit the ``$\boxtimes \cO_S$'' factor from the notation.

We will call the $\Cal{Z}_{\Cal{E}}^r$ (or their connected components) {\em special cycles of corank $m$} (with $r$ legs).
\end{defn}

There is an evident map $\Cal{Z}_{\Cal{E}}^{r} \rightarrow \Sht_{U(n)}^r$ projecting to the data in the first bullet point. When $\rank \Cal{E}=1$, the $\Cal{Z}_{\Cal{E}}^{r} $ are function field analogues (with multiple legs) of the \emph{Kudla-Rapoport divisors} introduced in \cite{KRI, KRII}. 

\subsection{Indexing via Hermitian maps}\label{ssec: sht index by a} 
\begin{defn}\label{def: A(k)}Let $\Cal{A}^{\rm all}_{\Cal{E}}(k)$ be the $k$-vector space of Hermitian maps\footnote{We will later in \S\ref{ss:Hit base} introduce a space $\Cal{A}^{\rm all}$ over $\Bun_{\GL'_m}$ for which $\Cal{A}^{\rm all}_{\Cal{E}}(k)$ is the $k$-rational points of the fiber over $\cE\in \Bun_{\GL'_m}(k)$, justifying the notation.} $a \co \Cal{E} \rightarrow \sigma^* \Cal{E}^{\vee}$ such that $\sigma(a)^{\vee} =a$. Let $\Cal{A}_{\Cal{E}}(k)\subset \Cal{A}^{\rm all}_{\Cal{E}}(k)$ be the subset where the map $a \co \Cal{E} \rightarrow \sigma^* \Cal{E}^{\vee}$ is injective (as a map of coherent sheaves). 
\end{defn}

Let $(\{x'_i\}, \{\Cal{F}_i\}, \{f_i\}, \varphi, \{t_i\}) \in \Cal{Z}_{\Cal{E}}^r(S)$. By the compatibilities between the $t_i$ in the definition of $\Cal{Z}_{\Cal{E}}^r$, the compositions 
\begin{equation}\label{eq: hitchin map formula}
\Cal{E} \boxtimes \cO_S \xrightarrow{t_i} \Cal{F}_i \xrightarrow{h_i} \sigma^* \Cal{F}_i^{\vee} \xrightarrow{\sigma^* t^{\vee}_{i}} \sigma^* \Cal{E}^{\vee} \boxtimes \cO_S
\end{equation}
agree for each $i$, and \eqref{eq: hitchin map formula} for $i=r$ also agrees with the Frobenius twist of \eqref{eq: hitchin map formula} for $i=0$. Hence \eqref{eq: hitchin map formula} for every $i$ gives the same map $\Cal{E} \boxtimes \cO_S  \rightarrow \Cal{E}^{\vee} \boxtimes \cO_S$, which moreover must come by pullback from $\Cal{A}^{\rm all}_{\Cal{E}}(k)$. This defines a map $\Cal{Z}_{\Cal{E}}^{r} \rightarrow \Cal{A}^{\rm all}_{\Cal{E}}(k)$. For $a \in \Cal{A}^{\rm all}_{\Cal{E}}(k)$, we denote by $\Cal{Z}^{r}_{\Cal{E}}(a)$ the fiber of $\Cal{Z}^{r}_{\Cal{E}}$ over $a$. We have
\begin{equation*}
\cZ_{\cE}^{r}=\coprod_{a\in \cA^{\rm all}_{\cE}(k)} \cZ^{r}_{\cE}(a).
\end{equation*}

\begin{defn}\label{def:Da}
For $a\in \cA_{\cE}(k)$, let $D_{a}$ be the effective divisor on $X$ such that $\nu^{-1}(D_{a})$ is the divisor of the Hermitian map $\det(a): \det(\cE)\to \s^{*}\det(\cE)^{\vee}$.
\end{defn}

\begin{defn}\label{def:Z(0)}
For any $a \in \cA_{\cE}^{\all}(k)$, we define: 
\begin{itemize}
\item $\cZ_{\cE}^{r}(a)^{\circ}\subset\cZ^{r}_{\cE}(a)$ to be the open substack classifying $(\{x'_{i}\}, \{\cE\xr{t_{i}}\cF_{i}\})$ with the additional condition that all $t_{i}$ are \emph{injective} when restricted to $X'_{\ov s}$ for any geometric point $\ov s$ of the test scheme $S$. 
\item $\cZ_{\cE}^{r}(a)^{*}\subset\cZ^{r}_{\cE}(a)$ to be the open substack classifying $(\{x'_{i}\}, \{\cE\xr{t_{i}}\cF_{i}\})$ with the additional condition that all $t_{i}$ are \emph{non-zero} when restricted to $X'_{\ov s}$ for any geometric point $\ov s$ of the test scheme $S$.
\end{itemize}
Note that if $\rank(\cE)=1$, then the inclusion $\cZ_{\cE}^{r}(a)^{\circ} \inj \cZ_{\cE}^{r}(a)^{*}$ is an isomorphism, and both include isomorphically into $\cZ_{\cE}^r(a)$ unless $a=0$.
\end{defn}

\subsection{Finiteness properties} 

We next establish that the projection map $\Cal{Z}_{\Cal{E}}^{r}(a) \rightarrow \Sht_{U(n)}^{r}$ is finite, which will eventually allow us to construct cycle classes on $ \Sht_{U(n)}^{r}$ associated to $\Cal{Z}_{\Cal{E}}^{r}(a) $. 

\begin{prop}\label{prop: Z finite}
Let $\cE$ be any vector bundle of rank $m$ on $X'$ and let  $a \in \Cal{A}^{\rm all}_{\Cal{E}}(k)$. Then the projection maps $\Cal{Z}_{\Cal{E}}^{r}(a) \rightarrow \Sht_{U(n)}^{r}$ and $\Cal{Z}_{\Cal{E}}^{r}(a)^* \rightarrow \Sht_{U(n)}^{r}$ are both finite.
\end{prop}

\begin{proof}
Note that $\Cal{Z}_{\Cal{E}}^{r}(a)$ has a closed substack where the map $t \co \cE \rightarrow \cF$ is $0$, which projects isomorphically to $\Sht_{U(n)}^r$. The complement of this closed substack is $\Cal{Z}_{\Cal{E}}^{r}(a)^*$, so it suffices to show the finiteness of $\Cal{Z}_{\Cal{E}}^{r}(a)^* \rightarrow \Sht_{U(n)}^{r}$. We will show that it is proper and quasi-finite. First we establish the properness. It suffices to show this locally on the target, so we pick a Harder-Narasimhan polygon $P$ for $\Bun_{U(n)}$ and consider the truncation $\Bun_{U(n)}^{\leq P}$. Define $\Sht_{U(n)}^{r, \leq P} $ to be the open substack of $\Sht_{U(n)}^r$ obtained as the pullback of $\Bun_{U(n)}^{\leq P} \inj \Bun_{U(n)}$ via the tautological projection $\pr_{0}: \Sht_{U(n)}^r \rightarrow \Bun_{U(n)}$ recording $\cF_{0}$, and $\Cal{Z}_{\Cal{E}}^{r, \leq P}(a) \inj \Cal{Z}_{\Cal{E}}^{r}(a)$ the analogous pullback.  

We can then pick a sufficiently anti-ample vector bundle $\cE'$ of rank $m$ on $X'$ and an injection $\iota \co \Cal{E}' \inj \Cal{E} $ so that the stack $\ul{\Hom}(\Cal{E}', -)^{\leq P}$ parametrizing $\{(\Cal{F} \in \Bun_{U(n)}^{\leq P}, t \in \Hom(\Cal{E}', \Cal{F}))\}$ forms a vector bundle over $\Bun_{U(n)}^{\leq P}$, with respect to the obvious projection map. Let $a' := (\sigma^{*} \iota^{\vee}) \circ a \circ \iota \co \Cal{E}' \rightarrow \Cal{E}'^\vee$.  Then we have a closed embedding $\Cal{Z}_{\Cal{E}}^{r, \leq P}(a)  \inj \cZ_{\Cal{E}'}^{r, \leq P}(a')$ cut out by the condition that the map $t \co \Cal{E}' \rightarrow \Cal{F}$ factors through $\iota$, which fits into a commutative diagram 
\[
\begin{tikzcd}
\Cal{Z}_{\Cal{E}}^{r, \leq P}(a)  \ar[r, hook] \ar[dr] & \cZ_{\Cal{E}'}^{r, \leq P}(a') \ar[d] \\
& \Sht_{U(n)}^{r, \leq P}.
\end{tikzcd}
\]
It suffices to show the open substack $\cZ_{\Cal{E}'}^{r, \leq P}(a')^*$ defined by the condition that  $t_{0}\ne0$ fiberwise over the test scheme is proper over $ \Sht_{U(n)}^{r, \leq P} $. We can factorize this map as the composition of two maps in the diagram below: 
\[
\begin{tikzcd} 
\cZ_{\Cal{E}'}^{r, \leq P}(a')^* \ar[dr] \ar[r, "j"] & \PP(\ul{\Hom}(\Cal{E}', -)^{\leq P}) \times_{\Bun_{U(n)}^{\leq P}} \Sht_{U(n)}^{r, \leq P}\ar[d, "\pr_2"] \\
& \Sht_{U(n)}^{r, \leq P}
\end{tikzcd}
\]
where $j$ is determined by the map $\cZ_{\Cal{E}'}^{r, \leq P}(a')^* \rightarrow \PP(\ul{\Hom}(\Cal{E}', -)^{\leq P})$ sending 
\[
(\Cal{F}_0 \dashrightarrow \ldots \dashrightarrow \Cal{F}_r \cong \ft \Cal{F}_0, (t_i)_{i=0}^r) \mapsto (\Cal{F}_0 , t_0 \co \Cal{E}' \rightarrow \Cal{F}_0).
\]
The map $\pr_2$ is a projective bundle by design, so in order to establish that $\pr_2 \circ j$ is proper it suffices to show that $j$ is finite. Indeed, since data of all the $t_i$ is determined by $t_0$, the analogous map $\cZ_{\Cal{E}'}^{r, \leq P}(a')^*  \rightarrow \ul{\Hom}(\Cal{E}', -)^{\leq P} \times_{\Bun_{U(n)}^{\leq P}} \Sht_{U(n)}^{r, \leq P}$ is a closed embedding. The requirement $t_0 = \ft t_0$ in the definition of $\cZ_{\Cal{E}'}^{r}$ therefore implies that the map $j$ is a $k^{\times}$-torsor onto its image, which is a closed substack of $\PP(\ul{\Hom}(\Cal{E}', -)^{\leq P})$. This completes the proof of properness. 

It remains to show that $\Cal{Z}_{\Cal{E}}^{r}(a) \rightarrow \Sht_{U(n)}^{r}$ is quasi-finite. Since the map is already established to be proper, it suffices by \cite[Tag 01TC]{stacks-project} to check that the fibers over field-valued points are finite. Let 
\[
(\{x'_{i}\}_{1\le i\le r}, (\Cal{F}_0 ,h_0) \dashrightarrow (\Cal{F}_1, h_1) \dashrightarrow \ldots (\Cal{F}_{r}, h_{r}) \xrightarrow{\sim} (\ft\Cal{F}_0, \ft h_0)) \in \Sht_{U(n)}^{r}(\kappa)
\]
be such a point valued in a field $\kappa$. Its fiber in $\Cal{Z}_{\Cal{E}}^{r}(a)(\kappa)$ consists of $\{t_i \co \cE \rightarrow \cF_i \}_{0\le i\le r}$ fitting into commutative diagrams 
\[
\begin{tikzcd}
\Cal{E}_{\kappa} \ar[d, "t_0"] \ar[r, equals] &  \Cal{E}_{\kappa} \ar[r, equals]  \ar[d, "t_1"]& \ldots \ar[r, equals] & \Cal{E}_{\kappa} \ar[r, "\sim"] \ar[d, "t_{r-1}"] & \ft \Cal{E}_{\kappa} \ar[d, "\ft t_0"] \\
\Cal{F}_0 \ar[r, dashed] & \Cal{F}_1 \ar[r, dashed] \ar[r] &  \ldots \ar[r, dashed] &  \Cal{F}_{r-1} \ar[r, "\sim"] & \ft \Cal{F}_0  
\end{tikzcd}
\]
such that $\sigma^* t_i^{\vee} \c h_{i}\circ t_i = a \in \Cal{A}_{\Cal{E}}^{\all}(\kappa)$ for each $i = 0, \ldots, r$. We want to show that there are finitely many possibilities for such $t_i \in \cohog{0}{X'_\kappa, \Cal{E}^{\vee}_{\kappa} \otimes \Cal{F}_i}$. 

The situation can be abstracted to the following semi-linear algebra problem.  

\begin{lemma}\label{lem: semlinear alg}
Suppose that $\kappa$ is any field over $k$, and we have finite-dimensional $\kappa$-vector spaces $V_1, V_2 \subset V$ with an injective $\Frob$-semi-linear map $\tau \co V_1 \hookrightarrow V_2$. 
\[
\begin{tikzcd}
V_1 \ar[rr,"\tau"] \ar[dr, hook] & & V_2 \ar[dl, hook'] \\
& V 
\end{tikzcd}
\]
Then the set $\{ x \in V_1 \co \tau(x) = x  \in V\}$ is finite. 
\end{lemma}

We assume Lemma \ref{lem: semlinear alg} for the moment and use it to conclude the proof of Proposition \ref{prop: Z finite}. We apply it to the situation above with $V_1 := \Hom_{X'_{\kappa}} (\Cal{E}_{\kappa}, \Cal{F}_{0})$, $V_2 := \Hom_{X'_{\kappa}}( \ft \Cal{E}_{\kappa},\ft \Cal{F}_0)$, which are both viewed as subspaces of 
\[
V := \Hom_{X'_{\kappa}}\left( \Cal{E}_{\kappa},  \Cal{F}_0(\sum_{j=1}^r (x'_j + \sigma(x'_j))) \right)
\]
by the obvious inclusion. The map $V_1 \rightarrow V_2$ is the twist by $\tau$. Then Lemma \ref{lem: semlinear alg} shows that there are finitely many possibilities for $t_0$ since $\t(t_{0})=t_{0}$. The other $t_i$ are determined by $t_0$ (if they exist) because the $t_i$ as well as the modifications $\cF_i \dashrightarrow \cF_{i+1}$ are all isomorphisms over an open subset of $X'_{\kappa}$. 
\end{proof}

\begin{proof}[Proof of Lemma \ref{lem: semlinear alg}] By replacing $\kappa$ with an algebraic closure, it suffices to consider the case when $\kappa$ is algebraically closed. Let us call a subspace $V_1' \subset V_1$ ``$\tau$-fixed'' if $\tau(V_1') = V_1' \subset V$. Since a sum of $\tau$-fixed subspaces is evidently $\tau$-fixed, there is a well-defined largest $\tau$-fixed subspace $V^{\circ}_1 \subset V_1$. It is a sub $\kappa$-vector space of $V_1$, hence necessarily finite-dimensional. Since $\tau:V_1\to V_2$ is injective, the restriction of  $\tau$ to $V_1^{\circ}$ is a $\Frob$-semi-linear bijection. The set $\{ x \in V_1 \co \tau(x) = x  \in V\}$ is evidently contained in $(V_1^{\circ})^{\tau}$, which is an $k$-form of $V_1^{\circ}$ (because $\kappa$ is algebraically closed) and therefore finite-dimensional over $k$. 
\end{proof}

\subsection{Variation with $\cE$} 
Let $\cE',\cE$ be two vector bundles (with possibly distinct ranks) on $X'$ and $s \co \Cal{E}' \to \Cal{E}$ be a map of coherent sheaves. Given $a \co \Cal{E} \to \sigma^* \Cal{E}^{\vee}$ in $\cA^{\rm all}_{\cE}(k)$, let $a'  = (\sigma^{*} s^{\vee} ) \circ a \circ s \co \Cal{E}' \to \sigma^* (\Cal{E}')^{\vee}$ be the corresponding element in $\cA^{\rm all}_{\cE'}(k)$. Therefore, composing with $s$ defines a map 
\begin{equation}\label{eq: Z_E to Z_E'}
z_{s}: \Cal{Z}_{\Cal{E}}^r(a) \to \Cal{Z}_{\Cal{E}'}^r(a')
\end{equation}
sending
\[
\begin{tikzcd}
\Cal{E} \ar[d] \ar[r , equals ] &  \ldots \ar[r, dashed]  & \Cal{E} \ar[r, equals] \ar[d]  & \ft \Cal{E}   \ar[d] \\
 \Cal{F}_0 \ar[r, dashed]  & \ldots \ar[r, dashed] & \Cal{F}_r \ar[r, "\sim"] & \ft \Cal{F}_0
\end{tikzcd} \mapsto 
\begin{tikzcd}
\Cal{E}' \ar[d] \ar[r , equals ] &  \ldots \ar[r, dashed]  & \Cal{E}' \ar[r, equals]  \ar[d]  & \ft \Cal{E}'   \ar[d] \\
 \Cal{F}_0 \ar[r, dashed]  & \ldots \ar[r, dashed] & \Cal{F}_r \ar[r, "\sim"] & \ft \Cal{F}_0
\end{tikzcd}
\]

The following lemma follows directly from definitions.

\begin{lemma}\label{l:KR direct sum}
If $\cE=\cE_{1}\op \cE_{2}$, and $a_{i}\in \cA^{\rm all}_{\cE_{i}}(k)$ for $i=1,2$, then there is a canonical isomorphism
\begin{equation}
\cZ^{r}_{\cE_{1}}(a_{1})\times_{\Sht^{r}_{U(n)}}\cZ^{r}_{\cE_{2}}(a_{2})\cong \coprod_{a=\mat{a_{1}}{*}{*}{a_{2}}}\cZ^{r}_{\cE}(a)
\end{equation}
where the union runs over all  Hermitian maps $a:\cE\to \s^{*}\cE^{\vee}$ whose restriction to $\cE_{i}$ is $a_{i}$ (for $i=1,2$). The map from the right side to the left is given by $(z_{\io_{1}},z_{\io_{2}})$, where $\io_{i}:\cE_{i}\incl \cE$ is the inclusion.
\end{lemma}

\begin{lemma}\label{l:change E} Under the  notations of the beginning of this subsection,
\begin{enumerate}
\item If $s:\cE'\to \cE$ is generically surjective, then $z_{s}:\cZ^{r}_{\cE}(a)\to \cZ^{r}_{\cE'}(a')$ is a closed embedding.
\item Suppose that $s$ is generically an isomorphism (in particular $\cE$ and $\cE'$ have the same rank). Let $D_{s}\subset X'$ be the divisor of $\det(s)$. Then the restriction of $z_{s}$ over $(X\bs \nu(D_{s}))^{r}$ 
\begin{equation}
\cZ^{r}_{\cE}(a)|_{(X\bs \nu(D_{s}))^{r}}\subset \cZ^{r}_{\cE'}(a')|_{(X\bs \nu(D_{s}))^{r}}
\end{equation}
is open and closed. Here we write $\cZ^{r}_{\cE}(a)|_{(X\bs \nu(D_{s}))^{r}}$ for the preimage of $(X\bs \nu(D_{s}))^{r}$ under the leg map $\cZ^{r}_{\cE}(a)\to \Sht^{r}_{U(n)}\to X'^{r}\xr{\nu^{r}} X^{r}$.
\end{enumerate}
\end{lemma}
\begin{proof}
(1) Let $\ov\cE'=\cE'/\ker(s)\xr{\ov s}\cE$, equipped with the induced Hermitian map $\ov a' \co \ov\cE' \rightarrow \sigma^* (\ov\cE')^{\vee}$. Then $s^{*}$ factors as $\cZ^{r}_{\cE}(a)\xr{z_{\ov s}}\cZ_{\ov\cE'}^{r}(\ov a')\subset \cZ^{r}_{\cE'}(a')$, the latter being evidently a closed embedding. Therefore it suffices to show $z_{\ov s}$ is a closed embedding. We thus reduce to the case $s$ is generically an isomorphism. 

Let $D$ be an effective divisor on $X'$ such that $\cE(-D)\incl \cE'\xr{s}\cE$. Let $\cF^{\rm univ}$ be the universal Hermitian bundle over $X'\times \Bun_{U(n)}$, and $\cF^{\rm univ}_{D}$ its restriction to $D\times \Bun_{U(n)}$. Let $\cV_{D}=\pr_{2*}\un\Hom(\pr_{1}^{*}\cE(-D)|_{D},\cF^{\rm univ}_{D})$, where $\pr_{1},\pr_{2}$ are the projections of $D\times \Bun_{U(n)}$ to the two factors. Then $\cV_{D}$ is a vector  bundle of rank equal to $n\rank(\cE)\deg(D)$ over $\Bun_{U(n)}$. Let $\cV_{i,D}$ be the pullback of $\cV_{D}$ over $\cZ^{r}_{\cE'}(a')$ via the map $\pr_{i}: \cZ^{r}_{\cE'}(a')\to \Bun_{U(n)}$. Then $\cV_{i,D}$ has a section $v_{i}$ whose value at $(\{x'_{j}\}, \{\cF_{j}\}, \{t_{j}: \cE'\to \cF_{j}\})\in \cZ^{r}_{\cE'}(a')$ is the restriction of $t_{i}$ to $\cE(-D)|_{D}\to \cF_{i}|_{D}$.  Then $t_{i}$ extends to $\cE$ if and only if $v_{i}$ vanishes at the point $(\{x'_{j}\}, \{\cF_{j}\}, \{t_{j}\})$. This identifies $\cZ^{r}_{\cE}(a)$ as the common zero locus of the sections $(v_{i})_{0\le i\le r-1}$ of the vector bundles $\cV_{i,D}$ over $\cZ^{r}_{\cE'}(a')$, hence closed in $\cZ^{r}_{\cE'}(a')$. 

(2) By (1),  it remains to show the openness of $z_s$ when restricted to $(X\bs \nu(D_{s}))^{r}$. Let $U'=X'\bs \supp(D_{s}+\s D_{s})$. Let $\Sht^{r}_{U(n),D_{s}}$ be the moduli stack of Hermitian shtukas $(\{x'_{i}\}, \{\cF_{i}\})$ of rank $n$ with legs in $U'^{r}$, and trivializations of $\cF_{i}|_{D_{s}}$ (as a vector bundle over $D_{s}$ of rank $n$) compatible with the shtuka structures. Then $\l: \Sht^{r}_{U(n),D_{s}}\to \Sht^{r}_{U(n)}|_{U'^{r}}$ is a $\GL_{n}(\cO_{D_{s}})$-torsor. Let $\cZ^{r}_{\cE}(a)_{D_{s}}$ and $\cZ^{r}_{\cE'}(a')_{D_{s}}$ be the base changes of $\cZ^{r}_{\cE}(a)$ and $\cZ^{r}_{\cE'}(a')$ along $\l$. Since $\cZ^{r}_{\cE'}(a')_{D_{s}}\to \cZ^{r}_{\cE'}(a')|_{U'^{r}}$ is finite \'etale surjective, it suffices to show that the inclusion $\cZ^{r}_{\cE}(a)_{D_{s}}\inj \cZ^{r}_{\cE'}(a')_{D_{s}}$ is open. Using the trivializations of $\cF_{i}|_{D_{s}}$, we get an evaluation map
\begin{equation}
\ev_{D_{s}}: \cZ^{r}_{\cE'}(a')_{D_{s}}\to \Hom_{D_{s}}(\cE'|_{D_{s}}, \cO_{D_{s}}^{\op n})
\end{equation}
where the target is a discrete set. Then $\cZ^{r}_{\cE}(a)_{D_{s}}$ is the preimage of the image of
\begin{equation}
\Hom_{D_{s}}(\cE|_{D_{s}}, \cO_{D_{s}}^{\op n}) \xr{(-)\c s} \Hom_{D_{s}}(\cE'|_{D_{s}}, \cO_{D_{s}}^{\op n})
\end{equation}
under $\ev_{D_{s}}$. Indeed, a map $\cE'\to \cF_{i}$ extends to $\cE\to \cF_{i}$ if and only if $\cE'|_{D_{s}}\to \cF_{i}|_{D_{s}}$ vanishes on $\ker(\cE'|_{D_{s}}\to \cE|_{D_{s}})$ (this can be checked locally using elementary divisors). Since the target of $\ev_{D_{s}}$ is discrete, $\cZ^{r}_{\cE}(a)_{D_{s}}\subset \cZ^{r}_{\cE'}(a')_{D_{s}}$ is open and closed.
\end{proof}

\subsection{Corank $1$ special cycles} 

A special role is played by the case $m=1$, i.e. where $\Cal{E}$ is a line bundle on $X'$, because it is only in this case that we can appropriately control the dimension of the cycles $\Cal{Z}_{\Cal{E}}^{r}$. We will write $\Cal{L} := \Cal{E}$ to emphasize that it is a line bundle. 

Note that in this case $a \in \cA_{\cL}(k)$ if and only if $a \neq 0$. 

\begin{prop}\label{prop: dim Z_L}
We have $\dim \cZ_{\cL}^{r}(a)^{\circ} \leq r(n-1)$. 
\end{prop}

This is established later, in Proposition \ref{p:ShM dim} (for $a\ne0$) and Proposition \ref{p:dim ZL0} (for $a=0$),  as a consequence of a more refined study of the geometry of $\cZ_{\cL}^{r}(a)^{\circ}$. 

\begin{remark}\label{rem: actual dimension}
One can show that when $a\ne0$, in fact $\cZ_{\cL}^{r}(a)^{\circ}$ is LCI of pure dimension $r(n-1)$. This will appear in a future paper; it relies on some ideas from \cite{FYZ2}. This fact will not be used in the present paper, but it may be psychologically helpful. 
\end{remark}

\begin{defn}\label{def: m=1 cycle class} 
For $a \in \cA^{\all}_{\cL}(k)$, we define $[ \cZ_{\cL}^{r}(a)^{\circ} ] \in \Ch_{r(n-1)}( \cZ_{\cL}^{r}(a)^{\circ})$ to be the cycle class of the union of the irreducible components of $\cZ_{\cL}^{r}(a)^{\circ} $ with dimension $r(n-1)$, throwing away the irreducible components of dimension $<r(n-1)$. (According to Remark \ref{rem: actual dimension}, there are no such components to be thrown away at least when $a\ne0$, but we neither prove nor use this in the present paper.)  
\end{defn}

\subsection{Corank $n$ special cycles}
In this paper we are mainly concerned with the case where the rank of $\cE$ is $m=n$.   The following proposition contains basic geometric information about $\cZ^{r}_{\cE}(a)$. 

\begin{lemma}\label{l:KR cycle legs}
Let $\cE$ be a vector bundle on $X'$ of rank $n$, and $a\in \cA_{\cE}(k)$. Then the map $\cZ^{r}_{\cE}(a)\to X'^{r}$ recording  the legs has image in $(\supp\nu^{-1}(D_{a}))^{r}$.
\end{lemma}
\begin{proof} Let $(\{x'_{i}\}, \{\cF_{i}\}, \{t_{i}\})$ be a geometric point of $\cZ^{r}_{\cE}(a)$. For each $1\le i\le r$,  the Hermitian map $a$ factorizes as $\cE\incl \cF^{\flat}_{i-1/2}\incl \cF_{i-1}\incl \s^{*}\cE^{\vee}$, we see that $x'_{i}$ (the support of $\cF_{i-1}/\cF^{\flat}_{i-1/2}$) is in the support of $\sigma^*\cE^{\vee}/a(\cE)$, i.e., $x'_{i}\in \supp \nu^{-1}(D_{a})$.
\end{proof}

\begin{prop}\label{prop: KR cycle proper} 
Let $\cE$ be a vector bundle on $X'$ of rank $n$, and $a\in \cA_{\cE}(k)$. Then $\cZ^{r}_{\cE}(a)$ is a proper scheme over $k$ that depends only on the torsion sheaf $\cQ=\coker(a)=\s^{*}\cE^{\vee}/\cE$ together with the Hermitian structure $\ov a$ on $\cQ$ induced from $a$ (see \S\ref{ss:Herm} for the notion of Hermitian structure on a torsion sheaf).
\end{prop}
The proof involves a few ideas not yet introduced, and will be given later in \S \ref{sssec: complete proof}.

The next goal is to equip the proper scheme $\cZ_{\cE}^{r}(a)$ with a $0$-cycle class in its Chow group. The ``virtual dimension'' of $\cZ_{\cE}^{r}(a)$ is at most zero, for if $\cE$ is a direct sum of line bundles $\cL_{1}\op\cdots\op \cL_{n}$, then $\cZ_{\cE}^{r}(a)$ is contained in the intersection of $Z^{r}_{\cL_{i}}(a_{ii})$ for $1\le i\le r$, each having codimension at least $r$ in $\Sht^{r}_{U(n)}$ by Proposition \ref{prop: dim Z_L} (which can be shown to be an equality, cf. Remark \ref{rem: actual dimension}).  However the actual dimension of $\cZ_{\cE}^{r}(a)$ can be strictly positive. Our task is to find the correct virtual fundamental class of $\cZ_{\cE}^{r}(a)$.

\subsection{Intersection theory on stacks}\label{ss:int stacks}
Recall the discussion of intersection theory on Deligne-Mumford stacks from \cite[Appendix A]{YZ}. Let $Y$ be a smooth, separated, locally finite type  Deligne-Mumford stack over $k$ of pure dimension $d$. Let $Y_{1}, \cdots, Y_{n}$ be Deligne-Mumford stacks with maps $f_{i}: Y_{i}\to Y$. Then there is an intersection product
\[
(-) \cdot_{Y} (-)  \cdot_Y \cdots \cdot_Y (-) \co \Ch_{i_{1}}(Y_{1}) \times \Ch_{i_{2}}(Y_{2})\times\cdots\times\Ch_{i_{n}}(Y_{n}) \rightarrow \Ch_{i_{1}+\cdots+i_{n}-d(n-1)}(Y_{1}\times_{Y}\cdots\times_{Y}Y_{n}).
\]
For $\z_{i}\in \Ch_{*}(Y_{i})$, the intersection product $\z_{1}\cdot_{Y}\cdots\cdot_{Y}\z_{n}$ is defined as the Gysin pullback of the external product $\z_{1}\times\cdots\times \z_{n}\in \Ch_{*}(Y_{1}\times\cdots\times Y_{n})$ along the diagonal map $\Delta: Y\to Y^{n}$, which is a regular embedding of codimension $d(n-1)$.

\subsection{Intersection problem: the case of a direct sum of line bundles}\label{ssec: KR cycle classes}
We now formulate the cycle classes which enter into our intersection problem.  We first consider the case $\Cal{E}$ a direct sum of $m$ line bundles on $X'$,
\begin{equation}
\Cal{E} \cong \Cal{L}_1 \oplus \ldots \oplus \Cal{L}_m.
\end{equation}
Let $a \in \Cal{A}_{\Cal{E}}(k)$. We write $a$ as an $m\times m$-matrix with entries $a_{ij}\in \Hom(\cL_{j}, \s^{*}\cL^{\vee}_{i})$. 

Let 
\begin{equation}
\cZ^{r}_{\cL_{1},\cdots, \cL_{m}}(a_{11},\cdots, a_{mm})^{\circ}:=\Cal{Z}^{r}_{\Cal{L}_1}(a_{11})^{\circ} \times_{\Sht^{r}_{U(n)}} \ldots \times_{\Sht^{r}_{U(n)}} \Cal{Z}^{r}_{\Cal{L}_m}(a_{mm})^{\circ}.
\end{equation}
In Definition \ref{def: m=1 cycle class} we defined a fundamental class $[\Cal{Z}_{\Cal{L}}^{r}(a)^{\circ}] \in \Ch_{r(n-1)}(\Cal{Z}_{\Cal{L}}^{r}(a)^{\circ})$. Applying the intersection product construction in \S\ref{ss:int stacks} for $Y=\Sht^{r}_{U(n)}$ (the hypotheses apply by Lemma \ref{lem: shtuka geometry} (2)), we obtain a class  
\begin{equation}\label{int m Z}
[\Cal{Z}^r_{\Cal{L}_1}(a_{11})^{\circ}] \cdot_{\Sht_{U(n)}^r} \ldots \cdot_{\Sht_{U(n)}^r} [\Cal{Z}^{r}_{\Cal{L}_m}(a_{mm})^{\circ}]\in\Ch_{r(n-m)}(\Cal{Z}^{r}_{\Cal{L}_1,\cdots, \cL_{m}}(a_{11},\cdots, a_{mm})^{\circ}).
\end{equation}
Let $\cA_{\cE}^{\all}(a_{11},\cdots, a_{mm})(k)$ be the finite set of Hermitian maps $a:\cE\to \s^{*}\cE^{\vee}$ (not assumed to be injective) such that its restriction to $\cL_{i}$ is $a_{ii}$ for $i=1,\cdots, m$. By Lemma \ref{l:KR direct sum}, there is a map
\begin{equation}
\Cal{Z}_{\Cal{L}_1}^r(a_{11}) \times_{\Sht^{r}_{U(n)}} \ldots \times_{\Sht^{r}_{U(n)}} \Cal{Z}^r_{\Cal{L}_m}(a_{mm})\to \cA_{\cE}^{\all}(a_{11},\cdots, a_{mm})(k)
\end{equation}
such that the fiber over $a \in \cA_{\cE}(k)$ identifies with $\cZ^{r}_{\cE}(a)$. Since $a$ is injective, the image of $\cZ^{r}_{\cE}(a)\to \cZ^{r}_{\cL_{i}}(a_{ii})$ lies in $\cZ^{r}_{\cL_{i}}(a_{ii})^{\circ}$. In particular, 
\begin{equation}
\cZ^{r}_{\cE}(a)\subset \Cal{Z}^{r}_{\Cal{L}_1,\cdots,\Cal{L}_m}(a_{11}, \cdots, a_{mm})^{\circ}
\end{equation}
is open and closed. Restricting \eqref{int m Z} to $\cZ^{r}_{\cE}(a)$ gives a cycle class
\[
\z^{r}_{\cL_{1},\cdots, \cL_{m}}(a):=  \left([\Cal{Z}^r_{\Cal{L}_1}(a_{11})^{\circ}] \cdot_{\Sht_{U(n)}^r} \ldots \cdot_{\Sht_{U(n)}^r} [\Cal{Z}^{r}_{\Cal{L}_m}(a_{mm})^{\circ}]\right)  |_{\cZ^{r}_{\cE}(a)}
 \in  \Ch_{r(n-m)}(\cZ^{r}_{\cE}(a)).
\]

\begin{remark} Our notation suggests that $\z^{r}_{\cL_{1},\cdots, \cL_{m}}(a)$ (as a cycle class on $\cZ_{\cE}^{r}(a)$) depends, at least a priori, on the decomposition of $\cE$ into a direct sum of line bundles $\cL_{1},\cdots, \cL_{m}$.  However, we will show later in Theorem \ref{th:compare KR}  that, at least when $m=n$, it only depends on $\cE$, and is equal to the cycle class $[\cZ^{r}_{\cE}(a)]$ that we will define for general rank $n$ bundle $\cE$.
\end{remark}



\subsection{Intersection problem: $m=n$ and $\cE$ arbitrary}\label{ssec: intersection problem}

To define a $0$-cycle $[\KR{E}^r(a)] $ for general rank $n$ vector bundle $\Cal{E}$ on $X'$, we need to make some auxiliary choice first; eventually we will show that the definition is independent of the choice.

\begin{defn}\label{def: auxiliary E'}  Let $\cE$ be a rank $n$ vector bundle over $X'$ and $a\in \cA_{\cE}(k)$. A {\em good framing} of $(\cE,a)$ is an $n$-tuple $(s_{i}: \cL_{i}\to \cE)_{1\le i\le n}$ of $\cO_{X'}$-linear maps from line bundles $\cL_{i}\in \Pic(X')$ satisfying:
\begin{enumerate}
\item The map $s=\op s_{i}: \cE':=\op_{i=1}^{n}\cL_{i}\to \cE$ is injective.
\item Let $D_{s}$ be the divisor of the nonzero map $\det(s): \ot_{i=1}^{n}\cL_{i}\to \det\cE$. Then $\nu(D_{s})$ (image in $X$) is  disjoint from $D_{a}$ (see Definition \ref{def:Da}). 
\end{enumerate}
\end{defn}

\begin{lemma} For any rank $n$ bundle $\cE$ on $X'$ and $a\in \cA_{\cE}(k)$, there exists a good framing for $(\cE,a)$ in the sense of Definition \ref{def: auxiliary E'}.
\end{lemma}

\begin{proof} For notational convenience we give the argument for $X'$ connected; the case $X'=X\coprod X$ can be proved with obvious changes. 

We strengthen condition (2) on $s:\op_{i=1}^{n}\cL_{i}\to \cE$ slightly by asking $\nu(D_{s})$ to avoid a prescribed divisor $D_{0}$ on $X$, instead of $D_{a}$.  We prove the existence of $s$ satisfying this stronger condition by induction on $n$. 

The base case $n=1$ is trivial: take $\cL_{1}=\cE$. 

For the inductive step, start by picking any saturated line bundle $\Cal{L}_{1}\inj \cE$. 
Then $\Cal{E}_{n-1} := \Cal{E}/\Cal{L}_{1}$ is a vector bundle of rank $n-1$. By induction hypothesis  we may pick $\ov s: \op_{i=2}^{n}\cL'_{i} \inj \Cal{E}_{n-1}$ satisfying the conditions of Definition \ref{def: auxiliary E'} and such that $\nu(D_{\ov s})$ avoids the given divisor $D_{0}$. Let $D_{2},\cdots, D_{n}$ be effective divisors on $X'$ such that
\begin{enumerate}
\item $\nu(D_{2}),\cdots, \nu(D_{n})$ are disjoint from $\nu(D_{0})$, and
\item $\deg \cL'_{i}-D_{i}+2g'-2<  \deg \cL_{1}$ for $i=2,\cdots, n$.
\end{enumerate}
Let $\cL_{i}=\cL'_{i}(-D_{i})$. By the inequality above we see that $\Ext^1(\Cal{L}_i, \Cal{L}_1) = 0$, so the map $\ov s_{i}: \cL_{i} \inj \Cal{E}_{n-1}=\cE/\cL_1$ lifts to a map $s_{i}: \cL_{i} \inj \Cal{E}$, $i=2,\cdots, n$. 

Now we have an injection $s: \op_{i=1}^{n}\cL_{i}\inj \cE$ whose divisor $D_{s}$ satisfies $D_{s}=D_{\ov s}+D_{2}+\cdots+D_{n}$. Since $\nu(D_{2}),\cdots, \nu(D_{n}), \nu(D_{\ov s})$ are disjoint from $D_{0}$ by construction, the same is true for $\nu(D_{s})$. 
\end{proof}

\begin{cor}\label{c:ZE open closed} If $s:\cE'=\op_{i=1}^{n}\cL_{i}\inj \cE$ is a good framing, then the map \eqref{eq: Z_E to Z_E'} realizes $\Cal{Z}_{\Cal{E}}(a)$ as an open and closed subscheme of $\Cal{Z}_{\Cal{E}'}(a')$. 
\end{cor}
\begin{proof} Closedness is proved in Lemma \ref{l:change E}(1). Only the openness requires an argument.
By the definition of a good framing, $\nu(D_{s})$ is disjoint from $D_{a}$, and therefore disjoint from all legs of all points of $\cZ^{r}_{\cE}(a)$ by Lemma \ref{l:KR cycle legs}. Let $U'=X'\bs \supp(D_{s}+\s D_{s})$, then $\cZ^{r}_{\cE}(a)=\cZ^{r}_{\cE}(a)|_{U'^{r}}$. By Lemma \ref{l:change E}(2), the inclusion
\begin{equation}
\cZ^{r}_{\cE}(a)=\cZ^{r}_{\cE}(a)|_{U'^{r}}\inj \Cal{Z}_{\Cal{E}'}(a')|_{U'^{r}}
\end{equation}
is open, hence the inclusion $\cZ^{r}_{\cE}(a)\inj \cZ^{r}_{\cE'}(a')$ is open.
\end{proof}

\begin{defn-prop}\label{defn:general KR} Let $\cE$ be a vector bundle of rank $n$ over $X'$ and $a\in \cA_{\cE}(k)$. Let $s:\op_{i=1}^{n}\cL_{i}\inj \cE$ be a good framing of $(\cE,a)$. Let 
\begin{equation}
[\Cal{Z}^{r}_{\Cal{E}}(a)]:=\z^r_{\cL_{1},\cdots, \cL_{n}}(a')|_{\cZ^{r}_{\cE}(a)}\in \Ch_{0}(\cZ_{\cE}^{r}(a)).
\end{equation}
Here we are using Corollary \ref{c:ZE open closed} to make sense of the restriction, as it implies that $\Cal{Z}^{r}_{\Cal{E}}(a)$ is a union of connected components of $\Cal{Z}^{r}_{\Cal{E}'}(a')$. Then the cycle class $[\Cal{Z}^{r}_{\Cal{E}}(a)]$ thus defined is independent of the good framing $s:\cE'=\op_{i=1}^{n}\cL_{i}\inj \cE$. 
\end{defn-prop}
The independence of good framing  will be proved in Theorem \ref{th:compare KR} after some preparation in \S\ref{sec: hitchin spaces}. The idea is to construct another $0$-cycle class on $\cZ^{r}_{\cE}(a)$ without making auxiliary choices (which is done by introducing Hitchin shtukas), and show that the two constructions agree. 


%



By Proposition \ref{prop: KR cycle proper}, $\cZ^{r}_{\cE}(a)$ is proper over $k$, 
therefore the {\em degree} of the $0$-cycle of $[\cZ^{r}_{\cE}(a)]\in \Ch_{0}(\cZ^r_{\cE}(a))$ is a well-defined  number in $\Q$.  The main problem we are concerned with in this Part is to determine $\deg [\KR{E}^r(a)] \in \Q$.


\section{Hitchin-type moduli spaces}\label{sec: hitchin spaces}

In this section we introduce certain ``Hitchin-type moduli stacks'' which will help to analyze the special cycles. In particular, we will be able to use these to give an alternative construction of the cycle classes associated to special divisors, that is manifestly independent of auxiliary choices. 

\subsection{Hitchin stacks}\label{ssec: Hitchin stacks}
Until \S\ref{ss:vb to tor}, we fix an arbitrary positive integer $m$.

\begin{defn} The \emph{Hitchin stack} $\Cal{M}^{\rm all}(m,n)$ (sometimes denoted $\Cal{M}^{\rm all}$ when $m,n$ are understood) has $S$-points the groupoid consisting of the following data. 
\begin{itemize}
\item $\Cal{E}$ a rank $m$ vector bundle on $X' \times S$.
\item $\Cal{F}$ a rank $n$ vector bundle on $X' \times S$, equipped with a Hermitian map $h \co \Cal{F} \xrightarrow{\sim} \sigma^* \Cal{F}^{\vee}$. 
\item A map of underlying coherent sheaves $t \co \Cal{E} \rightarrow \Cal{F}$ over $X'\times S$.\end{itemize}
We define  $\cM(m,n)\subset\cM^{\rm all}$ (sometimes denoted $\Cal{M}$ when $m,n$ are understood) to be the open substack where the map $t$ base changes to an injective map on $X'_{\ol{s}}$ for each geometric point $\ol{s} \rightarrow S$. 
\end{defn}
Let us emphasize that \emph{both} $\Cal{E}$ and $(\Cal{F},h)$ are varying in this definition. We will usually suppress the dependence on $m,n$ from the notation.

\subsection{Hitchin base} \label{ss:Hit base}
\begin{defn} We define the following two versions of the \emph{Hitchin base}. 
\begin{enumerate}
\item $\cA^{\mrm{all}}(m)$ (sometimes denoted $\cA^{\mrm{\all}}$ when $m$ is understood) to be the stack whose $S$-points is the groupoid of the following data: 
\begin{itemize}
\item $\Cal{E}$ a rank $m$ vector bundle on $X'\times S$;
\item $a \co \Cal{E} \rightarrow \sigma^* \Cal{E}^{\vee}$ is a map of coherent sheaves on $X'\times S$ such that $\sigma(a^{\vee}) = a$.
\end{itemize}

\item We define $\cA \subset \cA^{\mrm{all}}$ to be the open substack where $a\co \Cal{E} \rightarrow \sigma^* \Cal{E}^{\vee}$ is injective after base change to $X'_{\ol{s}}$ for every geometric point $\ol{s} \rightarrow S$.
\end{enumerate}
\end{defn}

\begin{defn} For integers $1 \leq m \leq n$, we define the \emph{Hitchin fibration} for $\cM^{\rm all} = \cM^{\rm all}(m,n)$ to be the map $f \co \cM^{\rm all} \rightarrow \Cal{A}^{\mrm{all}}$ sending $(\Cal{E}, (\Cal{F},h), t)$ to
the composition 
\[
a\co \Cal{E} \xrightarrow{t} \Cal{F} \xrightarrow{h} \sigma^* \Cal{F}^{\vee} \xrightarrow{\sigma^* t^{\vee}} \sigma^* \Cal{E}^{\vee}.
\]	
\end{defn}

\begin{remark}\label{r:M to A}In general the Hitchin fibration does not send $\cM(m,n)$ to $\cA(m)$ even when $m\le n$. However, in the special case $m=n$,  the Hitchin map does send $\cM(n,n)$ to $\cA(n)$ because when $t:\cE\to \cF$ is generically injective, it is generically an isomorphism for rank reasons, hence the induced Hermitian map $a$ on $\cE$ is generically non-degenerate.
\end{remark}


\subsection{Hitchin shtukas}\label{ssec: Hitchin shtukas}
We now discuss a notion of shtukas for Hitchin stacks. Throughout, $\cM = \cM(m,n)$. 

\begin{defn}[Hecke stacks for Hitchin spaces]\label{defn: hitchin hecke}
For $r \geq 0$, we define $\Hk_{\Cal{M}^{\all}}^{r}$ to be the stack whose $S$-points are given by the groupoid of the following data: 	
\begin{enumerate}
\item $(\{x'_i\}_{1\le i\le r},\{(\cF_{i},h_{i})\}_{0\le i\le r})\in \Hk^{r}_{U(n)}(S)$. 
\item A vector bundle $\Cal{E}$ of rank $m$ on $X' \times S$. 
\item Maps $t_i \co \Cal{E} \rightarrow \Cal{F}_i$ fitting into the commutative diagram
\end{enumerate}
\[
\begin{tikzcd}
\Cal{E} \ar[r, equals] \ar[d, "t_0"] & \Cal{E} \ar[r, equals]  \ar[d, "t_1"] & \ldots \ar[r, equals] & \Cal{E} \ar[d, "t_r"]  \\
\Cal{F}_{0}  \ar[r, dashed] & \Cal{F}_1 \ar[r, dashed] & \ldots \ar[r, dashed] &  \Cal{F}_r
\end{tikzcd}
\]
We define the open substack $\Hk_{\Cal{M}}^{r} \subset \Hk_{\Cal{M}^{\all}}^{r}$ by the condition that $t_0 $ base changes to an injective map along every geometric point $\ol{s} \rightarrow S$ (equivalently, every $t_i$ has this property). Let $\pr_{i}: \Hk^{r}_{\cM}\to \cM$ (resp. $\pr_i^{\all} \co \Hk^r_{\cM^{\all}} \to \cM^{\all}$) be the map recording $(\cE, \cF_{i}, h_{i}, t_{i})$, for $0\le i\le r$.

\end{defn}

\begin{defn}[Shtukas for Hitchin stacks]\label{defn: hitchin shtuka}
For $r \geq 0$, we define $\Sht_{\Cal{M}^{\all}}^{r}$ as the fibered product 
\begin{equation}\label{eq: shtuka all as intersection}
\begin{tikzcd}
\Sht_{\Cal{M}^{\all}}^{r}  \ar[r] \ar[d] & \Hk_{\Cal{M}^{\all}}^{r} \ar[d, "{(\pr_0^{\all}, \pr_r^{\all})}"] \\
\Cal{M}^{\all} \ar[r, "{(\Id, \Frob)}"] & \Cal{M}^{\all}  \times \Cal{M}^{\all}
\end{tikzcd}
\end{equation}
and the open substack $\Sht_{\Cal{M}}^{r} \subset \Sht_{\Cal{M}^{\all}}^{r}$ as the fibered product 
\begin{equation}\label{eq: shtuka as intersection}
\begin{tikzcd}
\Sht_{\Cal{M}}^{r}  \ar[r] \ar[d] & \Hk_{\Cal{M}}^{r} \ar[d, "{(\pr_0, \pr_r)}"] \\
\Cal{M} \ar[r, "{(\Id, \Frob)}"] & \Cal{M}  \times \Cal{M}
\end{tikzcd}
\end{equation}
Explicitly, the stack $\Sht_{\Cal{M}^{\all}}^{r}$ parametrizes diagrams of the form below, with notation as in Definition \ref{defn: hitchin hecke},
\begin{equation}\label{eq: sht_M point}
\begin{tikzcd}
\Cal{E} \ar[d, "t_0"] \ar[r, equals] &  \Cal{E} \ar[r, equals]  \ar[d, "t_1"]& \ldots \ar[r, equals] & \Cal{E} \ar[r, "\sim"] \ar[d, "t_{r}"] & \ft \Cal{E} \ar[d, "\ft t_0"] \\
\Cal{F}_0 \ar[r, dashed] & \Cal{F}_1 \ar[r, dashed]  &  \ldots \ar[r, dashed] &  \Cal{F}_{r} \ar[r, "\sim"] & \ft \Cal{F}_0  
\end{tikzcd}
\end{equation}
and $\Sht_{\Cal{M}}^{r}$ is the open substack where $t_0$ base changes to an injective map along every geometric point $\ol{s} \rightarrow S$ (equivalently, the same property holds for every $t_i$). 

\end{defn}


In particular, $\Cal{E} \xrightarrow{\sim} \ft \Cal{E}$ is a shtuka with no legs, exhibiting $\Cal{E}$ as arising from a rank $m$ vector bundle on $X'$, i.e. coming from $\Bun_{\GL'_m}(k)$. Note that in \eqref{eq: sht_M point}, $\cE$ is not fixed, so the automorphisms of $\cE$ are present in the functor of points of $\Sht_{\Cal{M}^{\all}}^{r}$. Therefore, if we define
\[
\ol{\Cal{Z}}_{\Cal{E}}^{r} :=  [\Cal{Z}_{\Cal{E}}^{r} /  (\Aut(\cE)(k)],
\]
then $\Sht_{\Cal{M}^{\all}}^{r}$ decomposes as a disjoint union of special cycles
\begin{equation}\label{eq: hitchin shtuka decomposition E}
\Sht_{\Cal{M}^{\all}}^{r} = \coprod_{\Cal{E} \in \Bun_{\GL_{m}'}(k)} \ol{\Cal{Z}}_{\Cal{E}}^{r}.
\end{equation}
This decomposition can be refined.
The compositions $f \circ \pr_i \co \Hk_{\Cal{M}^{\all}}^{r} \rightarrow \Cal{A}^{\mrm{all}}$ all coincide, and they induce a map 
\begin{equation}\label{eq: hitchin shtuka decomposition A}
\Sht_{\Cal{M}^{\all}}^{r} \rightarrow \Cal{A}^{\mrm{all}}(k).
\end{equation}
Here $\cA^{\mrm{all}}(k)$ is the groupoid of pairs $\{\cE, \cE \xrightarrow{a} \sigma^* \cE^\vee\}$. Let us write $\Aut(a) := \Aut(\cE \xrightarrow{a} \sigma^* \cE^\vee)\subset \Aut(\cE)$. Then each $\cE \xrightarrow{a} \sigma^* \cE^\vee$ defines a map $a: \Spec k \to B \Aut(a)(k)\incl \cA^{\all}(k)$, the latter map being an open-closed inclusion, and we have a Cartesian square
\begin{equation}\label{eq: hitchin shtuka decomposition E and A}
\begin{tikzcd}
\cZ_{\cE}^r(a) \ar[r] \ar[d] & \Sht_{\Cal{M}^{\all}}^{r}  \ar[d, "f"] \\
\Spec k \ar[r, "a"] & \cA^{\mrm{all}}(k) 
\end{tikzcd}
\end{equation}
In particular, $\cZ_{\cE}^r(a)$ is a finite \'{e}tale cover of an open-closed substack of $\Sht_{\Cal{M}^{\all}}^{r}$ isomorphic to $[\cZ_{\cE}^r(a)/\Aut(a)(k)]$. 
Similarly, restricting to the injective locus (see Definition \ref{def:Z(0)}) we have a 
Cartesian square
\begin{equation}\label{eq: hitchin shtuka decomposition E and A inj}
\begin{tikzcd}
\cZ_{\cE}^r(a)^{\circ} \ar[r] \ar[d] & \Sht_{\Cal{M}}^{r}  \ar[d, "f"] \\
\Spec k \ar[r, "a"] & \cA^{\mrm{all}}(k) 
\end{tikzcd}
\end{equation}
If $m=n$ and $a\in \cA(k)$, then we have $\cZ_{\cE}^r(a)^{\circ}=\cZ_{\cE}^r(a)$ and we can replace $ \cA^{\mrm{all}}(k)$ by $ \cA(k)$ in the above diagram (see Remark \ref{r:M to A}).


\subsection{From vector bundles to torsion sheaves}\label{ss:vb to tor}
For the rest of the section, we concentrate on the case $m=n$. In this case, we will relate $\cM = \cM(n,n)$ to the moduli stack of Hermitian torsion sheaves introduced in \S\ref{ss:Herm}.  We introduce the following abbreviated notations.

\begin{defn} Let $d\in \Z_{\ge0}$.
\begin{enumerate}
\item Let $\cM_{d} = \cM(n,n)_d$ be the open-closed substack of $\cM = \cM(n,n)$ consisting of $(\cE\xr{t}\cF)$ where $d=\frac{\deg\cE^\vee-\deg(\cE)}{2}=-\chi(X',\cE)$.  
\item Let $\cA_{d}  = \cA(n)_d$ be the open-closed substack of $\cA = \cA(n)$ consisting of $(\cE,a)$ where $d=\frac{\deg\cE^\vee-\deg(\cE)}{2}=-\chi(X',\cE)$.
\end{enumerate}
\end{defn}

By Remark \ref{r:M to A}, the Hitchin map for $\cM ^{\rm all}= \cM^{\rm all}(n,n)$ restricts to a map 
\begin{equation}
f_{d}: \cM_{d}\to \cA_{d}.
\end{equation}
When $d$ is understood, we abbreviate $f$ for $f_{d}$.

Recall that  $\Herm_{2d}=\Herm_{2d}(X'/X)$ is the moduli stack of length $2d$ torsion coherent sheaves $\cQ$ on $X'$ equipped with a Hermitian structure $h_{\cQ} \co \cQ \xrightarrow{\sim} \s^{*}\cQ^{\vee}$, where $\cQ^{\vee} := \ul{\Ext}^1(\cQ, \omega_{X'})$ such that $\sigma^* h_{\cQ}^{\vee}=h_{\cQ}$.  Alternatively we may think of $h_{\cQ}$ as the datum of a perfect pairing 
\begin{equation}
h'_{\cQ}: \cQ\ot_{\cO_{X'}}\s^{*}\cQ\to \omega_{X',F'}/\omega_{X'}
\end{equation}
where $\omega_{X',F'}$ is the constant Zariski sheaf of rational differential form on $X'$. The Hermitian condition is equivalent to $h'_{\cQ}(u,v)=\s^{*}h'_{\cQ}(v,u)$ for local sections $u,v$ of $\cQ$.
 
In \S\ref{ss:HSpr stalk} we have also introduced the moduli stack $\Lagr_{2d}=\Lagr_{2d}(X'/X)$ classifying $(\cQ,h_{\cQ}, \cL)$ where $(\cQ,h_{\cQ})\in \Herm_{2d}$ and $\cL\subset \cQ$ is a Lagrangian subsheaf.

There is a canonical map $g: \Cal{A}_d \rightarrow \Herm_{2d}$ sending $(\Cal{E},a)$ to the torsion sheaf $\cQ=\sigma^* \Cal{E}^{\vee}/\Cal{E}$ together with the Hermitian structure $h_{\cQ}$ defined in \S \ref{ssec: density function for torsion sheaves}. 


We have a map $g_{\cM}: \Cal{M}_d \rightarrow \Lagr_{2d}$ sending $(\Cal{E} \xr{t} \Cal{F}) \in \Cal{M}_d$ to the torsion Hermitian sheaf $(\cQ=\sigma^* \Cal{E}^{\vee}/\Cal{E}, h_{\cQ})$ constructed above together with the Lagrangian  $\cL=\Cal{F}/\cE$.

\begin{lemma}\label{lem: hitchin cartesian} The maps defined above fit into a Cartesian diagram 
\begin{equation}\label{eq: hitchin cartesian}
\begin{tikzcd}
\Cal{M}_d \ar[d, "f"] \ar[r,"g_{\cM}"]  & \Lagr_{2d} \ar[d, "\upsilon_{2d}"] \\
\Cal{A}_d \ar[r, "g"] & \Herm_{2d}
\end{tikzcd}
\end{equation}
\end{lemma}

\begin{proof} Given $a: \Cal{E} \rightarrow \sigma^* \Cal{E}^{\vee}$ that is injective, the datum of a subsheaf $\cL\subset \sigma^* \Cal{E}^{\vee}/\Cal{E}$ is the same as a coherent sheaf $\cF$ such that $\Cal{E} \subset \Cal{F} \subset \sigma^* \Cal{E}^{\vee}$. It is easy to see that $\cL$ is Lagrangian if and only if $\cF$ is self-dual under the Hermitian map $a$.
\end{proof}

\begin{cor}\label{c: hitchin proper}
The Hitchin fibration $f_{d} \co \Cal{M}_{d} \rightarrow \Cal{A}_{d}$ is proper.
\end{cor}
\begin{proof}Apply Lemma \ref{lem: hitchin cartesian} and the fact that $\upsilon_{2d}$ is proper.
\end{proof}

\begin{lemma}\label{l:Lag small}
The map $\upsilon_{2d} \co \Lagr_{2d} \rightarrow \Herm_{2d}$ is small.
\end{lemma}
\begin{proof}
The map $\pi^{\Herm}_{2d}$ from \S\ref{ss:Herm Spr} factors as 
\begin{equation}\label{factor pi Herm}
\pi^{\Herm}_{2d}: \wt\Herm_{2d}\xr{\l_{2d}}\Lagr_{2d}\xr{\upsilon_{2d}}\Herm_{2d}.
\end{equation}
Since $\l_{2d}$ is surjective and $\pi^{\Herm}_{2d}$ is small by Proposition \ref{p:SH Spr}, we get the desired statement.
\end{proof}


\subsubsection{Proof of Proposition \ref{prop: KR cycle proper}}\label{sssec: complete proof}
Let $\Lagr(\cQ)$ be the moduli space of Lagrangian subsheaves of $\cQ$. Let $\Hk^{r}_{\Lagr(\cQ)}$ be its Hecke version, classifying points $\{x'_{i}\}_{1\le i\le r}$ of $X'$ and chains of Lagrangian subsheaves of $\cQ$
\begin{equation}
\begin{tikzcd}
\cL_{0} \ar[r, dashed, "{f_{1}}"]   & \cL_{1} \ar[r,dashed, "{f_{2}}"]  & \cdots \ar[r, dashed, "{f_{r}}"]  & \cL_{r}
\end{tikzcd}
\end{equation}
where the dashed arrow $f_{i}$ are modifications at $x'_{i}\cup \s(x'_{i})$, similar to those in Definition \ref{defn: Hk U(n)}. There is a natural map $\cZ^{r}_{\cE}(a)\to \Hk^{r}_{\Lagr(\cQ)}$ sending a point $(\{x'_{i}\}, \{t_{i}: \cE\to \cF_{i}\})$ of $\cZ^{r}_{\cE}(a)$ to the collection of (necessarily Lagrangian) subsheaves $\ov \cF_{i}=\coker(t_{i})\subset\cQ=\s^{*}\cE^{\vee}/\cE$. This map fits into a Cartesian diagram
\begin{equation}\label{ZLag}
\xymatrix{\cZ^{r}_{\cE}(a)\ar[r]\ar[d] & \Hk^{r}_{\Lagr(\cQ)} \ar[d] \\
\Lagr(\cQ)\ar[r]^-{(\Id,\Frob)} & \Lagr(\cQ)\times \Lagr(\cQ)}
\end{equation}
Now both $\Lagr(\cQ)$ and $\Hk^{r}_{\Lagr(\cQ)}$ are proper schemes over $k$, hence the same is true for  $\cZ^{r}_{\cE}(a)$.  The diagram also makes it clear that $\cZ^{r}_{\cE}(a)$ only depends on $(\cQ,\ov a)$.

\subsection{Smoothness}

\begin{lemma}\label{l:Md smooth}  The map 
\begin{equation}
\b_{d}: \cM_{d}\to \Coh_{d}(X')\times \Bun_{U(n)}
\end{equation}
sending $(\cE\xr{t}\cF,h)$ to $(\coker(t), (\cF,h))$ is smooth of relative dimension $dn$. In particular, $\cM_{d}$ is smooth of pure dimension $dn+n^{2}(g-1)$. 
\end{lemma} 
\begin{proof}
Consider the stack $\cM'_{d}$ classifying $(\cT,\cF,h, s)$ where $\cT\in \Coh_{d}(X')$, $(\cF,h)\in \Bun_{U(n)}$ and an $\cO_{X'}$-linear map $s:\cF\to \cT$. Let $\b'_{d}: \cM'_{d}\to \Coh_d(X') \times \Bun_{U(n)}$ be the natural map. Due to the vanishing of $\Ext^{1}(\cF,\cT)$, $\b'_{d}$ exhibits $\cM'_{d}$ as a vector bundle over $\Coh_d(X') \times \Bun_{U(n)}$ of rank $dn=\dim\Hom_{X'}(\cF,\cT)$. Now $\cM_{d}$ is the open substack of $\cM'_{d}$ where $s$ is surjective. Therefore $\b_{d}$ is also smooth of relative dimension $dn$. 
\end{proof}

\begin{prop}\label{prop:Ad Herm smooth}
The map $g \co \Cal{A}_d \rightarrow \Herm_{2d}$ is smooth. 
\end{prop}
%

%
%
%

\begin{proof}
We have a map $s^{\Lagr}_{2d} \co \Lagr_{2d}\to X'_{d}$ sending $(\cQ,h_{\cQ},\cL)$ to the divisor of $\cL$. Recall also the map $s^{\Herm}_{2d}: \Herm_{2d}\to X_{d}$ sending $(\cQ,h_{\cQ})$ to the descent of the divisor of $\cQ$ to $X$.

Recall in \S\ref{ss:HSpr stalk} we introduced the open subset 
\[
(X_d')^{\dm} = \{D' \subset X' \co D' \cap \sigma(D') = \vn\} \subset X_d'.
\]
Let $\Lagr_{2d}^{\dm} \subset \Lagr_{2d}$ and $\Cal{M}_d^{\dm} \subset \Cal{M}_d$ be the preimages of $(X'_{d})^{\dm}$ under the maps $s^{\Lagr}_{2d}$ and $s^{\Lagr}_{2d}\c g_{\cM}$. 


We claim that both squares in the diagram 
\[
\begin{tikzcd}
\Cal{M}_d^{\dm} \ar[r] \ar[d, "f_d"] & \Lagr_{2d}^{\dm} \ar[d] \ar[r] & (X'_d)^{\dm}\ar[d] \\
\Cal{A}_d \ar[r,"g"] & \Herm_{2d} \ar[r] & X_d
\end{tikzcd}
\]
are Cartesian. The left square is Cartesian by definition. Now we show that the right square is Cartesian. Let $(\cQ ,h_{\cQ})\in \Herm_{2d}$,  $D' \in (X_d')^{\dm}$ lying over $D=s^{\Herm}_{2d}(\cQ,h_{\cQ})$. Since  $D'\cap\sigma(D')=\vn$, there is a unique Lagrangian subsheaf $\cL \subset \cQ$ supported on the support of $D'$, namely $\cL=\cQ|_{\supp D'}$. This gives the unique point $(\cQ,h_{\cQ}, \cL)\in \Lagr^{\dm}_{2d}$ mapping to $(\cQ,h_{\cQ})\in \Herm_{2d}$ and $D'\in (X'_{d})^{\dm}$. 

Note that the map $(X'_d)^{\dm} \rightarrow X_d$ is faithfully flat: it is clearly surjective, and the map $\nu_{d}: X'_d \rightarrow X_d$ is a finite morphism between smooth schemes, hence flat. We will show that $\Cal{M}_d^{\dm} \rightarrow \Lagr_d^{\dm}$ is smooth. By fppf descent it then follows that $\Cal{A}_d \rightarrow \Herm_d$ is also smooth. 

Recall from \S\ref{ss:HSpr stalk} that $\e'_{d}: \Lagr_{2d}\to \Coh_{d}(X')$ (recording only $\cL$) restricts to an isomorphism $\Lagr_{2d}^{\dm}\isom \Coh_{d}(X')^{\dm}:=\Coh_d(X')|_{(X_d')^{\dm}}$. Therefore it suffices to show that the composition $\cM_{d}^{\dm}\xr{g_{\cM}}\Lagr_{2d}^{\dm}\xr{\e'_{d}} \Coh_{d}(X')^{\dm}$ is smooth. This follows from the smoothness of $\cM_{d}\to \Coh_{d}(X')$ proved in Lemma \ref{l:Md smooth}.
\end{proof}

\begin{cor}\label{c:small} 
The Hitchin fibration $f \co \cM_d \rightarrow \cA_d$ is small. The complex $Rf_{*}\Qlbar$ is a shifted perverse sheaf that is the middle extension from any dense open substack of $\cA_{d}$.
\end{cor}
\begin{proof}
By the smoothness of $g$ in Proposition \ref{prop:Ad Herm smooth} and the Cartesian diagram in Lemma \ref{lem: hitchin cartesian}, the smallness of $f$ follows from that of $\upsilon_{2d}: \Lagr_{2d}\to \Herm_{2d}$, which is proved in Lemma \ref{l:Lag small}.
\end{proof}


\subsection{Cycle class from Hitchin shtukas}\label{ss:ShtM cycle}

In this subsection we take $m=n$, so $\cM = \cM(n,n)$ and $\cA = \cA(n)$. Now consider the Hitchin shtukas for $\cM_{d}\subset \cM$. Let $N=\dim \cM_{d}$. By Corollary \ref{c:small}, $\dim \cA_{d}=N$. The Cartesian diagram \eqref{eq: shtuka as intersection} restricts to a Cartesian diagram
\begin{equation}\label{def ShtMd}
\xymatrix{\Sht_{\cM_{d}}^{r}\ar[d] \ar[r] & \Hk^{r}_{\cM_{d}}\ar[d]^{(\pr_{0},\pr_{r})}\\
\cM_{d}\ar[r]^-{(\Id, \Frob)} & \cM_{d}\times\cM_{d}}
\end{equation}
We would like to define a $0$-cycle class on $\Sht_{\cM_{d}}^{r}$ as the Gysin pullback of a cycle on $\Hk^{r}_{\cM_{d}}$ along the Frobenius graph of $\cM_{d}$. Although the virtual dimension of $\Hk^{r}_{\cM_{d}}$ is the same as $\dim\cM_{d}$, its actual dimension may be larger. For this reason we have to define in a roundabout way a virtual fundamental cycle on $\Hk^{r}_{\cM_{d}}$ of the virtual dimension by relating it to $\Hk^{1}_{\cM_{d}}$, which we show is smooth below. 

\begin{lemma}\label{l:HkMd sm} The stack $\Hk^{1}_{\cM_{d}}$ is smooth and equidimensional of the same dimension as $\cM_{d}$.
\end{lemma}
\begin{proof}
Let $(x',\Cal{F}_0 \dashrightarrow \Cal{F}_1) \in \Hk_{U(n)}^1$. Let $\cF^{\flat}:=\cF^{\flat}_{1/2}=\Cal{F}_0 \cap \Cal{F}_1$ as in Definition \ref{defn: Hk U(n)}. 
The generically compatible Hermitian structures on $\Cal{F}_0$ and $\Cal{F}_1$ equip this intersection with a Hermitian structure $h^{\flat}: \cF^{\flat} \inj \sigma^* (\Cal{F}^{\flat})^{\vee}$ whose cokernel has length $1$ at $x'$ and at $\sigma (x')$. We call such a bundle {\em almost Hermitian} with defect at the ordered pair of conjugate points $(x',\s(x'))$. Conversely, given $(\Cal{F}^{\flat},h^{\flat})$ almost Hermitian with defect at $(x', \s(x'))$, one can uniquely recover $\cF_{0}$ (resp. $\cF_{1}$) as the upper modification of $\cF^{\flat}$ at $x'$ (resp. $\s(x')$) inside $\s^{*}(\cF^{\flat})^{\vee}$.  

Let $\Bun^{\flat}_{U(n)}$ be the moduli stack parametrizing  $(x'\in X', \cF^{\flat}, h^{\flat})$ where $(\cF^{\flat}, h^{\flat})$ is an almost Hermitian bundle with defect at $(x',\sigma(x'))$. The discussion in the previous paragraph shows that there is an isomorphism $\Hk^{1}_{U(n)}\isom \Bun^{\flat}_{U(n)}$ over $X'$.  
Let $\cM^{\flat}_{d}$ be the moduli stack of $(x',\cE\xr{t}\cF^{\flat}, h^{\flat})$ where $(\cF^{\flat}, h^{\flat})$ is almost Hermitian with defect at $(x',\s(x'))$, $\cE$ is a vector bundle on $X'$ of rank $n$ and $\chi(X',\cE)=-d$, and $t$ is injective. Then we have an isomorphism $\Hk^{1}_{\cM_{d}}\cong \cM^{\flat}_{d}$.

We have a natural map 
\begin{equation}
\b^{\flat}_{d}: \cM^{\flat}_{d}\to \Coh_{d-1}(X')\times \Bun^{\flat}_{U(n)}
\end{equation}
sending $(x',\cE\xr{t}\cF^{\flat}, h^{\flat})$ to $(\coker(t), (x',\cF^{\flat}, h^{\flat}))$. The same argument as Lemma \ref{l:Md smooth} shows that $\b^{\flat}_{d}$ exhibits $\cM^{\flat}_{d}$ as an open substack in a vector bundle of rank $n(d-1)$ over $\Coh_{d-1}(X')\times \Bun^{\flat}_{U(n)}$. Now $\dim \Coh_{d-1}(X')=0$ and $\dim \Bun^{\flat}_{U(n)}=\dim \Hk^{1}_{U(n)}=\dim \Bun_{U(n)}+n$ by Lemma \ref{lem: shtuka geometry}(1). Therefore $\Hk^{1}_{\cM_{d}}\cong \cM^{\flat}_{d}$ is of pure dimension $\dim \Bun_{U(n)}+dn$, which is the same as $\dim \cM_{d}$ by Lemma \ref{l:Md smooth}.
\end{proof}


\begin{defn}\label{def:Phi} For any stack $S$ over $k$ we define a morphism
\begin{equation}
\Phi^{r}_{S}: S^{r+1}\to S^{2r+2}
\end{equation}
by the formula $\Phi^{r}_{S}(\xi_{0},\cdots, \xi_{r})=(\xi_{0},\xi_{1},\xi_{1},\xi_{2},\xi_{2},\cdots, \xi_{r-1}, \xi_r, \xi_r, \Frob(\xi_{0}))$. When $r$ is fixed in the context, we simply write $\Phi_{S}$.
\end{defn}

We rewrite $\Sht^{r}_{\cM_{d}}$ as the fiber product
\begin{equation}\label{rewrite ShtM}
\xymatrix{\Sht^{r}_{\cM_{d}} \ar[r]\ar[d] & (\Hk^{1}_{\cM_{d}})^{r}\ar[d]^{(\pr_{0},\pr_{1})^{r} \times \Delta} \times \cM_d\\
(\cM_{d})^{r+1}\ar[r]^-{\Phi^{r}_{\cM_{d}}} & (\cM_{d})^{2r+2}=(\cM_{d})^{2r}\times(\cM_{d})^{2}}
\end{equation}
Here the vertical map $(\pr_{0},\pr_{1})^{r}$ sends $(h_{1},\cdots, h_{r})\in (\Hk^{1}_{\cM_{d}})^{r} $ to 
$(\pr_{0}(h_{1}), \pr_{1}(h_{1}),\cdots, \pr_{0}(h_{r}),\pr_{1}(h_{r}))\in (\cM_{d})^{2r}$, while $\Delta$ is the diagonal map. 

\begin{defn}\label{def: hitchin shtukas class} We define a $0$-cycle classes $[\Sht^{r}_{\cM_{d}}]\in \Ch_{0}(\Sht^{r}_{\cM_{d}})$ as the image of the fundamental class of $(\Hk^{1}_{\cM_{d}})^{r} \times \cM_d$  (which is smooth of the same dimension as $(\cM_{d})^{r+1}$ by Lemma \ref{l:HkMd sm}) under the refined  Gysin map along $\Phi_{\cM_{d}}:(\cM_{d})^{r+1}\to (\cM_{d})^{2r+2}$ (which is defined since $\cM_{d}$ is smooth and equidimensional by Lemma \ref{l:Md smooth}; see \cite[Theorem 2.1.12(xi)]{Kr99}) 
\begin{equation}
[\Sht^{r}_{\cM_{d}}]:=(\Phi^{r}_{\cM_{d}})^{!}[(\Hk^{1}_{\cM_{d}})^{r} \times \cM_d]\in \Ch_{0}(\Sht^{r}_{\cM_{d}}).
\end{equation}
\end{defn}

\section{Special cycles of corank one}\label{ssec: geometric properties}

In this section we prove geometric properties of the special cycles $\cZ^{r}_{\cE}(a)$ when $m= \rank \cE = 1$ (where the number field analogues are called ``Kudla-Rapoport divisors''). In particular, we show that for $a\ne0$, $\cZ^{r}_{\cE}(a)$  are local complete intersections of dimension $r$ less than $\Sht^{r}_{U(n)}$.  When $a=0$,  we show that  $\cZ^{r}_{\cE}(a)^{*}$ has dimension at most $\dim\Sht^{r}_{U(n)}-r$. These geometric properties are proved by studying stratifications introduced and analyzed in \S \ref{ss:strat KR div} and \S \ref{ssec: m=1 a=0 case}.

\subsection{Stratification of the special cycle when $m=1$}\label{ss:strat KR div}
In this section we fix a line bundle $\cL$ on $X'_{k}$ and $a\in \cA_{\cL}(k)$, i.e.,  $a: \cL\to \s^{*}\cL^{\vee}$ is nonzero Hermitian. We now define a stratification of $\cZ:=\cZ_{\cL}^{r}(a)_{\ov k}$ and estimate the dimension of each stratum.

When $n=1$, we have $\dim \cZ  =0$. Indeed, Lemma \ref{l:KR cycle legs} implies that there is a $0$-dimensional closed subscheme $D \subset (X')^r$ such that the morphism $\cZ \rightarrow \Sht_{U(1)}^r$ takes values in the preimage of $D$ under the leg map $\Sht_{U(1)}^r \rightarrow (X')^r$. This preimage is $0$-dimensional by Lemma \ref{lem: shtuka geometry} (2). As $Z \rightarrow \Sht_{U(1)}^r$ is finite (by Proposition \ref{prop: Z finite}) with 0-dimensional image, we conclude that $\dim \cZ = 0$. 

Therefore in the rest of this subsection we will assume $n\ge2$. 
 
For each $(\cL\xr{t_{i}}\cF_{i})_{0\le i\le r}\in \cZ(\ov k)$ with legs $(x'_{i})_{1\le i\le r}\in X'(\ov k)^{r}$, let $D_{i}$ ($0\le i\le r$) be the divisor on $X'_{\ov k}$ such that $t_{i}: \cL(D_{i})\incl \cF_{i}$ is saturated. For each $1\le i\le r$, we have one of the  four cases:
\begin{enumerate}
\item[(0)] $D_{i}=D_{i-1}$;
\item[(+)] $D_{i}=D_{i-1}+\s(x'_{i})$;
\item[($-$)] $D_{i}=D_{i-1}-x'_{i}$;
\item[($\pm$)] $D_{i}=D_{i-1}-x'_{i}+\s(x'_{i})$.
\end{enumerate}
Since the composition $\cL\xr{t_{i}}\cF_{i}\xr{h}\s^{*}\cF_{i}^{\vee}\xr{\s^{*}t_{i}^{\vee}}\s^{*}\cL^{\vee}$ is equal to $a$, we see that $D_{i}+\s(D_{i})$ is a subdivisor of the divisor of $a$. Therefore $\nu(D_{i})\le D_{a}$ as divisors on $X(\ov{k})$. 

\subsubsection{Indexing set for strata} Consider the set $\frD$ of sequences of effective divisors $(D_{i})_{0\le i\le r}$ on $X'_{\ov{k}}$ satisfying
\begin{itemize}
\item $\nu(D_{i})\le D_{a}$ for all $0\le i\le r$.
\item For each $1\le i\le r$, the pair $(D_{i-1}, D_{i})$ falls into one of the cases $(0),(+),(-),(\pm)$ above for some $x'_{i}\in X'(\ov k)$.
\item $D_{r}={}^{\t}D_{0}$.
\end{itemize}
It is clear that $\frD$ is a finite set. This will be the index set for our stratification of $\cZ$. 

\subsubsection{Definition of strata} Fix $D_{\bu}=(D_{i})_{0\le i\le r}\in \frD$. Let $I_{0} :=\{1\le i\le r|D_{i}=D_{i-1}\}$. Similarly we define $I_{+}, I_{-}$ and $I_{\pm}$ as the set of those $i$ such that $(D_{i-1}, D_{i})$ falls into case $(+),(-)$ and $(\pm)$  respectively. Let $\cZ[D_{\bu}]$ be the substack of $\cZ_{\ol{k}}$ classifying 
\[
(\{x'_{i}\}\in X'^{r}; \{\cF_{i}\}\in \Hk^{r}_{U(n)}; \{\cL(D_{i})\xr{t'_{i}}\cF_{i}\}_{0\le i\le r})
\]
such that every $t'_{i}$ is saturated. Let 
\begin{equation}
\pi[D_{\bu}]: \cZ[D_{\bu}]\to  (X'_{\ol{k}})^{I_{0}}
\end{equation}
be the map recording those $x'_{i}$ for $i\in I_{0}$. Note that for $i\in I_{+}\cup I_{-}\cup I_{\pm}$, $x'_{i}$ is determined by $D_{\bu}$.


\begin{prop}\label{p:ShM dim} Let $n \geq 2$. \hfill
\begin{enumerate}	
\item The substacks $\cZ[D_{\bu}]$ for $D_{\bu}\in \frD$ give a partition of $\cZ$.
\item Each geometric fiber of $\pi[D_{\bu}]$ has dimension $\le (n-1)|I_{+}|+(n-2)|I_{0}|$.
\item We have $\dim \cZ[D_{\bu}]\le r(n-1)$. The equality can only be achieved when  $I_{0}=\{1,2,\cdots, r\}$, i.e., all $D_{i}$ are equal to the same divisor of $X'$, which is then necessarily defined over $k$. \footnote{In this case, $\cZ[D_{\bu}]$ can be identified with the open substack $\mathring{\cZ}^{r}_{\cL(D_{0})}(a')\subset \cZ^{r}_{\cL(D_{0})}(a')$ (where $a'$ is the map $\cL(D_{0}) \rightarrow \sigma^* \cL(D_{0})^{\vee}$ induced from $a$) defined by requiring all the maps $t'_{i}: \cL(D_{0})\to \cF_{i}$ be saturated.}
\end{enumerate}
\end{prop}
\begin{proof}
(1) Each geometric point of $z\in \cZ$ defines a (unique) point $D_{\bu}\in \frD$ by taking the zero divisor of $t_{i}$, and then $z\in \cZ[D_{\bu}]$ by definition. 

(2) Let $\cH[D_{\bu}]$ be the substack of the fiber of $(\Hk^{r}_{\cM})_{\ol{k}}$ over $(\cL,a)\in \cA_{d}(k)$ classifying data $(\{x'_{i}\},\cL\xr{t_{i}}\cF_{i})$ such that $t_{i}$ extends to a map $t'_{i}: \cL(D_{i})\to \cF_{i}$ which is saturated. Note that  for $i\in I_{+}\cup I_{-}\cup I_{\pm}$, $x'_{i}$ is determined by $D_{\bu}$. Let $\cM[D_{i}]$ be the substack of the fiber of $\cM(1,n)_{\ol{k}}$ over $(\cL,a)\in \cA(1,n)(k)$ classifying maps $t:\cL\to \cF$ that extend to a saturated map $t': \cL(D_{i})\to \cF$. Then we have a Cartesian diagram of stacks over $\ol{k}$
\begin{equation}\label{ZHM}
\xymatrix{\cZ[D_{\bu}]\ar[d]\ar[r] & \cH[D_{\bu}]\ar[d]^{(p_{0}, p_{r})}\\
\cM[D_{0}]\ar[r]^-{(\Id,\Frob)} & \cM[D_{0}]\times \cM[D_{r}].}
\end{equation}
Note since $D_{r}={}^{\t}D_{0}$, the Frobenius morphism sends $\cM[D_{0}]$ to $\cM[D_{r}]$. 

Let 
\begin{equation}
\Pi[D_{\bu}]: \cH[D_{\bu}]\to \cM[D_{r}]\times X'^{I_{0}}_{\ov k}
\end{equation}
be the projection $p_{r}$ and the map recording $x'_{i}$ for $i\in I^{0}$. 

\begin{claim} The map $\Pi[D_{\bu}]$ is smooth and representable of relative dimension $(n-1)|I_{+}|+(n-2)|I_{0}|$.
\end{claim}

Assuming the claim, we finish the proof of (2). Indeed, we may apply Lemma \ref{l:VL} below to the Cartesian diagram \eqref{ZHM} to conclude, except that $\cM[D_{0}]$ is generally not a scheme. To remedy, we may restrict to a finite type open substack  $\cM[D_{0}]^{\le P}\subset \cM[D_{0}]$ by bounding Harder-Narasimhan polygon of $(\cF,h)$, and impose level structures on the Hermitian bundle $(\cF,h)$ at a closed point $x\in |X|$ to arrive at a scheme $\wt\cM[D_{0}]^{\le P}$ which is a torsor over $\cM[D_{0}]^{\le P}$ by an algebraic group $H$. Truncating and imposing the same level structures to $\cH[D_{\bu}]$ gives a  scheme $\wt\cH[D_{\bu}]^{\le P}$ (with legs away from $x$) such that $\wt\cH[D_{\bu}]^{\le P}/H=\cH[D_{\bu}]^{\le P}$. Let $\wt\cZ[D_{\bu}]^{\le P}$ be defined by a Cartesian diagram similar to \eqref{ZHM}, with $\cH[D_{\bu}]$ replaced by $\wt\cH[D_{\bu}]^{\le P}$ and $\cM[D_{i}]$ replaced by $\wt\cM[D_{i}]^{\le P}$ for $i=0,r$. We apply  Lemma \ref{l:VL} to conclude that the fibers of $\wt \cZ[D_{\bu}]^{\le P}\to X'^{I_{0}}_{\ol{k}}$ have dimension $\le (n-1)|I_{+}|+(n-2)|I_{0}|$. Now $\wt \cZ[D_{\bu}]^{\le P}/H(k)\isom \cZ[D_{\bu}]^{\le P}$ and for varying $P$ and $x$, $\cZ[D_{\bu}]^{\le P}$ cover $\cZ[D_{\bu}]$, hence the same dimension estimate holds for $\cZ[D_{\bu}]$.
%

It remains to prove the claim.

For $r\ge j\ge 0$, let $\cH_{\ge j}$ be the moduli stack defined similarly to $\cH[D_{\bu}]$ but classifying saturated maps $\{t_{i}: \cL(D_{i}) \to \cF_{i}\}_{j\le i\le r}$ (lying over $a$) only for $i$ in the indicated range. We can factorize $\Pi[D_{\bu}]$ as
\begin{equation}
\Pi[D_{\bu}]: \cH[D_{\bu}]=\cH_{\ge0} \xr{ \Pi_{1}} \cH_{\ge1}\xr{\Pi_{2} }\cdots \xr{ \Pi_{r}}\cH_{\ge r}=\cM[D_{r}]\times X'^{I_{0}}_{\ov k}.
\end{equation}
The desired smoothness and relative dimension claims would follow from the following four statements:
\begin{enumerate}
\item[($H0$)] If $i\in I_{0}$, then $\Pi_{i}$ exhibits $\cH_{\ge i-1}$ as an open substack in a $\PP^{n-2}$-bundle over $\cH_{\ge i}$.
\item[($H+$)] If $i\in I_{+}$, then $\Pi_{i}$ exhibits $\cH_{\ge i-1}$ as an open substack in a $\PP^{n-1}$-bundle over $\cH_{\ge i}$.
\item[($H-$)] If $i\in I_{-}$, then $\Pi_{i}$ is an isomorphism.
\item[($H\pm$)] If $i\in I_{\pm}$, then $\Pi_{i}$ is an open immersion.
\end{enumerate}
We next establish each of these statements. 

\emph{Proof of }($H0$). When $i\in I_{0}$, $D_{i-1}=D_{i}$. We write the modification $\cF_{i-1}\dashrightarrow\cF_{i}$ as 
\begin{equation}\label{lower middle}
\begin{tikzcd}
\cF_{i-1} & \cF^{\flat}_{i-1/2}\ar[r, "\s(x'_{i})", hook] \ar[l, "{x'_{i}}"', hook']  & \cF_{i} 
\end{tikzcd}
\end{equation}
Here both arrows have cokernel of length one supported at the labelled points. Such modifications of $\cF_{i}$ are parametrized by a hyperplane $H$ in the fiber $\cF_{i}|_{\s(x'_{i})}$. 
The requirement that $t_{i}:\cL(D_{i})\to \cF_{i}$ should land in $\cF^{\flat}_{i-1/2}$ is equivalent to the (closed) condition that $H$ should contain the line given by the image of $\cL(D_{i})|_{\s(x'_{i})}$. This cuts out a $\PP^{n-2}$ in the space of hyperplanes $H \subset \cF_{i}|_{\s(x'_{i})}$. The further requirement that $t_{i-1}:\cL(D_{i})\to \cF^{\flat}_{i-1/2}\to \cF_{i-1}$ be saturated is an open condition. 

This argument globalizes in the evident way, exhibiting that $\Pi_i$ as an open substack in a $\PP^{n-2}$-bundle. This applies similarly for the analogous arguments below for the other cases.

\emph{Proof of} ($H+$). When $i\in I_{+}$, we have $D_{i-1}=D_{i}-\s(x'_{i})$. We write the modification of $\cF_{i}$ as in \eqref{lower middle}. This time the choice of the $\cF_{i}\leftarrow\cF^{\flat}_{i-1/2}$ is the open subset of those hyperplanes $H\in \PP(\cF_{i}|_{\s(x'_{i})})$ that do not contain the image of $\cL(D_{i})|_{\s(x'_{i})}$. The requirement that $\cL(D_{i-1})=\cL(D_{i}-\s(x'_{i}))\to \cF^{\flat}_{i-1/2}\to \cF_{i-1}$ be saturated at $x'_{i}$ imposes a further open condition. 

\emph{Proof of} ($H-$). When $i\in I_{-}$, we have $D_{i-1}=D_{i}+x'_{i}$. We write the modification as 
\begin{equation}\label{upper middle}
\begin{tikzcd}
\cF_{i-1}\ar[r, "\s(x'_{i})", hook] & \cF^{\sh}_{i-1/2} & \cF_{i}\ar[l, "x'_{i}"', hook' ]
\end{tikzcd} 
\end{equation}
where both arrows have cokernel of length one supported at the labelled points. Now $t_{i}: \cL(D_{i})\to \cF_{i}$ is required to extend to $\cL(D_{i}+x'_{i})\to \cF^{\sh}_{i-1/2}$. This determines the upper modification $\cF_{i}\to \cF^{\sh}_{i-1/2}$ uniquely, which in turn determines the lower modification $\cF^{\sh}_{i-1/2}\leftarrow \cF_{i-1}$ as well. We get a map $t'_{i-1}: \cL(D_{i-1})=\cL(D_{i}+x'_{i})\to \cF^{\sh}_{i-1/2}$. We claim that $t'_{i-1}$ automatically lands in $\cF_{i-1}$. Indeed, the claim is equivalent to saying that under the pairing between $\cF_{i}|_{x'_{i}}$ and $\cF_{i}|_{\s(x'_{i})}$, the images of $t_{i}(x'_{i})$ and $t_{i}(\s(x'_{i}))$ pair to zero. The latter statement is equivalent to saying that the induced Hermitian map $a'_{i}: \cL(D_{i})\to \s^{*}(\cL(D_{i}))^{\vee}$ vanishes at $x'_{i}$. But we know that the divisor of $a'_{i}$ is $\nu^{*}D_{a}-D_{i}-\s(D_{i})$. Since $D_{i-1}=D_{i}+x'_{i}$ satisfies $\nu(D_{i-1})\le D_{a}$ by assumption, we see that $\nu^{*}D_{a}-D_{i}-\s(D_{i})\ge x'_{i}+\s(x'_{i})$, hence $a'_{i}$ is guaranteed to vanish at $x'_{i}$ and $\s(x'_{i})$. This shows that there is a unique lifting of any point of $\cH_{\ge i}$ to $\cH_{\ge i-1}$, hence $\Pi_{i}$ is an isomorphism.  

\emph{Proof of} ($H\pm$). When $i\in I_{\pm}$, we have $D_{i}+x'_{i}=D_{i-1}+\s(x'_{i})$. We write the modification as in \eqref{upper middle}. As in the case ($H-$), the requirement that $t_{i}: \cL(D_{i})\to \cF_{i}$ should extend to $\cL(D_{i}+x'_{i})\to \cF^{\sh}_{i-1/2}$ determines the modification. Then we automatically get a map $t_{i-1}: \cL(D_{i-1})=\cL(D_{i}+x'_{i}-\s(x'_{i}))\to \cF_{i-1}$; the requirement that $t_{i-1}$ be saturated is an open condition. Therefore $\Pi_{i}$ is an open immersion in this case.

(3) By (2) we have 
\[
\dim \cZ[D_{\bu}]\le (n-1)|I_{+}|+(n-2)|I_{0}|+|I_{0}|=(n-1)(|I_{+}\cup I_{0}|)\le r(n-1).
\]
Equality holds only if $I_{-}$ and $I_{\pm}$ are empty. However, by degree reasons we have $|I_{+}|=|I_{-}|$, so in the equality case we must have $I_{+}=\vn$ as well. We conclude that equality can only be achieved if $I_{0}=\{1,2,\cdots, r\}$; in other words, all $D_{i}$ must be the same. In particular, since $D_r = \ft D_0$, this forces $D_0$ to be defined over $k$. 
\end{proof}

The Lemma below, a slight variant of \cite[Lemma 2.13]{Laff18}, was used above.

\begin{lemma}[Variant of {\cite[Lemma 2.13]{Laff18}}]\label{l:VL}
Let $W,Z, T$ be schemes of finite type over $\ov{k}$. Let $Z^{(1)}$ be the Frobenius twist of $Z$ (i.e., the pullback of $Z$ under the $q$-Frobenius $\Spec \ov{k}\to \Spec \ov{k}$). Let $\wt h = (h_1, h_T) \co W \rightarrow Z^{(1)} \times T$ be smooth of relative dimension $d$, and $h_{0}\co W \rightarrow Z$ be an arbitrary map. Define $V$ as the fibered product
\[
\begin{tikzcd}
V \ar[r] \ar[d] & W \ar[d, "{(h_{0},h_{1})}"] \\
Z \ar[r, "{(\Id, \Frob)}"] & Z \times Z^{(1)}
\end{tikzcd}
\]
Then  each fiber of the composition map $V \rightarrow W\xr{h_{T}}T$ has dimension $\le d$. 
\end{lemma}

\begin{proof} Restricting $W$ over a point $t\in T(\ov{k})$, we reduce to the case $T$ itself is the point $\Spec\ov{k}$. We may assume $Z=\Spec R$ where $R=\ov k[x_{1},\ldots, x_{l}]/I$. Let $R^{(1)}=\ov k\ot_{\ov k}R\cong \ov k[\xi_{1},\ldots, \xi_{l}]/I^{(1)}$ be the base change of $R$ under $\Frob_{q}$, where $\xi_{i}=1\ot x_{i}$ . Since $h_{1}: W\to Z^{(1)}$ is smooth of relative dimension $d$,  by Zariski localizing we may assume $W=\ov k[\xi_{1},\ldots, \xi_{l}, y_{1},\ldots, y_{m+d}]/(I^{(1)}, r_{1},\ldots, r_{m})$,  with $(\pderiv{r_i}{y_j})_{j=1}^m$ having rank $m$ ($r_{i}\in \ov k[\xi_{1},\ldots, \xi_{l}, y_{1},\ldots, y_{m+d}]$). Under $h_{0}: W\to Z$, the coordinates $x_{i}$ of $Z$ pullback to functions $\ov f_{i}$ on $W$, $1\le i\le l$. We lift $f_{i}$ to polynomials $f_{i}\in \ov k[\xi_{1},\ldots, \xi_{l}, y_{1},\ldots, y_{m+d}]$. 

By definition, $V$ has the form 
\begin{eqnarray*}
V & \cong&  \Spec\dfrac{ \ov k[\xi_{1},\ldots, \xi_{l}, y_1, \ldots, y_m, y_{m+1} \ldots, y_{m+d}  ]}
{(I^{(1)}(\xi), g_{1},\ldots, g_{l}, r_1, \ldots, r_m)}
\end{eqnarray*}
where $g_{i}=\xi_{i}-f_{i}^{q}$. In particular,  $V$ is a closed subscheme of 
\begin{equation}
U:=\Spec\Big( \ov k[\xi_{1},\ldots, \xi_{l}, y_1, \ldots, y_m, y_{m+1} \ldots, y_{m+d}  ]/ (g_{1},\ldots, g_{l}, r_1, \ldots, r_m)\Big).
\end{equation}
The Jacobian matrix for the defining equations of $U$ has the form 
\[
\begin{pmatrix}
\boxed{\pderiv{g_i}{\xi_j}}_{j=1}^l  & \boxed{\pderiv{g_i}{y_j}}_{j=1}^m & \boxed{\pderiv{g_i}{y_j}}_{j=m+1}^{m+d}  \\
\boxed{\pderiv{r_i}{\xi_j}}_{j=1}^l  & \boxed{\pderiv{r_i}{y_j}}_{j=1}^m & \boxed{\pderiv{r_i}{y_j}}_{j=m+1}^{m+d} 
\end{pmatrix} = \begin{pmatrix} \Id_l & 0 &  0 \\ * & \text{invertible}_m & *  \end{pmatrix},
\]
which evidently has rank $l+m$. Hence $U$ is smooth of dimension $d$. Since $V\incl U$, $\dim V\le d$.
\end{proof}

\subsection{The case $m=1$ and $a=0$}\label{ssec: m=1 a=0 case}
We keep the notations from \S\ref{ss:strat KR div}. In this subsection we extend the discussion in \S\ref{ss:strat KR div} to the case $a=0$. Fix a line bundle $\cL\in \Pic(X')$. Recall from Definition  \ref{def:Z(0)} that $\cZ_{\cL}^{r}(0)^{\circ}$ is the moduli stack classifying hermitian shtukas $(\{x'_{i}\},\{\cF_{i}\})$ together with compatible maps $\{\cL\xr{t_{i}}\cF_{i}\})$ with $t_{i}$ injective (fiberwise over the test scheme $S$) and the image of $t_{i}$ being isotropic.  In this subsection, let
\begin{equation}
\cZ:=\cZ_{\cL}^{r}(0)^{\circ}_{\ov k}.
\end{equation}
If $n=1$ then $\cZ_{\cL}^{r}(0)^{\circ}=\vn$. We always assume $n\ge2$ below. 

\subsubsection{Indexing set for strata} Let $I_{0}\sqcup I_{+}\sqcup I_{-}$ be a partition of $\{1,2,\cdots, r\}$ such that $|I_{+}|=|I_{-}|$.  We denote this partition simply by $I_{\bu}$. For any $N\in \Z_{\ge0}$, define $\frD(N;I_{\bu})$ to be the moduli space of sequences of effective divisors $(D_{i})_{0\le i\le r}$ on $X'_{\ov k}$ such that
\begin{enumerate}
\item $\deg(D_{0})\le N$.
\item For $?=0,+$ or $-$, and $i\in I_{?}$, the pair $(D_{i-1},D_{i})$ belongs to the corresponding Case (?) listed in the beginning of \S\ref{ss:strat KR div} (for some $x'_{i}\in X'_{\ov k}$ in the case $?=+$ or $-$).
\item $D_{r}={}^{\t}D_{0}$. 
\end{enumerate}
We have a map recording the points $x'_{i}$ for $i\in I_{+}\cup I_{-}$:
\begin{equation}
(\pi_{+}, \pi_{-}): \frD(N;I_{\bu})\to (X'_{\ol{k}})^{I_{+}}\times (X'_{\ol{k}})^{I_{-}}.
\end{equation}

\begin{lemma}\label{l:dim D}
The map $\pi_{+}: \frD(N;I_{\bu})\to (X'_{\ol{k}})^{I_{+}}$ is quasi-finite.
\end{lemma}
\begin{proof}
For a fixed geometric point $(x'_{i})_{i\in I_{+}} \in (X'_{\ol{k}})^{I_{+}}(\ol{k})$, its fiber in $\frD(N;I_{\bu})$ consists of $(D_{0}, \{x'_{i}\}_{i\in I_{-}})$ such that $\deg D_{0}\le N$ and 
\begin{equation}\label{eqn D0}
D_{0}+\sum_{i\in I_{+}} \s(x'_{i})={}^{\t}D_{0}+\sum_{i\in I_{-}} x'_{i}.
\end{equation}
Let $v\in |X'|$ be a closed point that intersects $\supp D_{0}$. If $\deg(v)>N$, then $D_{0}$ cannot contain all geometric points over $v$ and hence there exists a geometric point $y|v$ such that $y\in \supp {}^{\t}D_{0}$ but $y\notin \supp D_{0}$. By \eqref{eqn D0}, $y=\s(x'_{i})$ for some $i\in I_{+}$.  Therefore points in $D_{0}$ are either over closed points of degree $\le N$, or in the Galois orbit of $\s(x'_{i})$ for some $i\in I_{+}$. This leaves finitely many possibilities for $D_{0}$, hence for $\{x'_{i}\}_{i\in I_{-}}$ as well.
\end{proof}


\subsubsection{Definition of strata} For a partition $I=(I_{0},I_{+},I_{-})$ of $\{1,2,\cdots, r\}$, define $\cZ[N;I_{\bu}]$ to be the stack classifying 
\[
(\{D_{i}\}_{0\le i\le r},( \{x'_{i}\}_{1\le i\le r}, \{\cF_{i}\}_{0\le i\le r})\in \Hk^{r}_{U(n)}, \{\cL\xr{t_{i}}\cF_{i}\}_{0\le i\le r})
\]
such that $\{D_{i}\}\in \frD(N;I_{\bu})$ with image $\{x'_{i}\}_{i\in I_{?}}$ under $\pi_{?}$ ($?=+,-$), and $t_{i}$ extends to a saturated embedding $\cL(D_{i})\inj \cF_{i}$. We have a map
\begin{equation}
\pi[N;I_{\bu}]: \cZ[N;I_{\bu}]\to (X'_{\ol{k}})^{I_{0}}\times \frD(N;I_{\bu}).
\end{equation}
The following is the analog of Proposition \ref{p:ShM dim} when $a=0$.

\begin{prop}\label{p:dim ZL0} Let $n \geq 2$. \hfill 
\begin{enumerate}
\item For varying $N\in \Z_{\ge0}$ and partitions $I_{\bu}$ of $\{1,2,\cdots, r\}$ such that $|I_{+}|=|I_{-}|$, the substacks  $\cZ[N;I_{\bu}]$ give a partition of $\cZ$. 
\item The fibers of the map $\pi[N;I_{\bu}]$ have dimension $\le(n-1)|I_{+}|+(n-2)|I_{0}|$.
\item We have $\dim \cZ[N;I_{\bu}]\le r(n-1)$. Moreover, when $n\ge3$, the equality  can only be achieved when $I_{0}=\{1,2,\cdots, r\}$, i.e., all $D_{i}$ are equal to the same divisor of $X'$ defined over $k$.
\end{enumerate}
\end{prop}
\begin{proof}
(1) is similar to Proposition \ref{p:ShM dim}(1), except we have to argue that the strata are only non-empty for $|I_+|=|I_-|$, and that Case $(\pm)$ cannot appear for points in $\cZ=\cZ^{r}_{\cL}(0)^{\circ}$. The first statement follows from the assumption that $D_r = \ft D_0$ has the same degree as $D_0$. For the second statement, suppose Case $(\pm)$ happens for the modification 
\[
\cF_{i-1}\leftarrow \cF_{i-1/2}^{\flat}\to \cF_{i},
\]
and let $H\subset \cF_{i-1,x'_{i}}$ be the hyperplane that is the image of $\cF_{i-1/2}^{\flat}$. Then $H^{\bot}\subset \cF_{i-1,\s(x'_{i})}$ is the line along which the upper modification $\cF_{i-1/2}^{\flat}\inj \cF_{i}$ is performed. Let $\ell_{x'_{i}}$ (resp. $\ell_{\s(x'_{i})}$) be the image of $\cL(D_{i-1})\to \cF_{i-1}$ at $x'_{i}$ (resp. at $\s(x'_{i})$). Since the image of $\cL(D_{i-1})$ is isotropic (because $a=0$), $(\ell_{x'_{i}},\ell_{\s(x'_{i})})=0$ under the pairing between $\cF_{i-1,x'_{i}}$ and $\cF_{i-1,\s(x'_{i})}$. The condition $D_{i}+x'_{i}=D_{i-1}+\s(x'_{i})$ happens only if $\ell_{x'_{i}}\not\subset H$ and $\ell_{\s(x'_{i})}=H^{\bot}$. This contradicts the fact that $(\ell_{x'_{i}},\ell_{\s(x'_{i})})=0$.

(2) is proved in the same way as Proposition \ref{p:ShM dim}(2). 

(3) Applying (2) and Lemma \ref{l:dim D} we get
\begin{eqnarray*}
\dim\cZ[N;I_{\bu}]&=&  (n-1)|I_{+}|+(n-2)|I_{0}|+|I_{0}|+\dim \frD(N;I_{\bu})\\
&\le & (n-1)|I_{+}|+(n-1)|I_{0}|+|I_{+}|=(n-1)|I_{0}|+n|I_{+}|.
\end{eqnarray*}
Since $n\ge2$,  we have $n\le2(n-1)$, therefore the above is $\le (n-1)(|I_{0}|+2|I_{+}|)=r(n-1)$. When $n\ge3$, we have strict inequalities $n<2(n-1)$, so equality can only be achieved when $I_{+}$, hence $I_{-}$, are all empty.
\end{proof}


\section{Comparison of two cycle classes}\label{s:compare cycles}
The goal in this section is to show the following theorem.

\begin{thm}\label{th:compare KR} Let $\cE\in \Bun_{\GL_{n}'}(k)$ and $a\in \cA_{\cE}(k)$. Let $s:\cE'=\op_{i=1}^{n}\cL_{i}\inj \cE$ be a good framing of $(\cE,a)$ in the sense of Definition \ref{def: auxiliary E'}. Let $a':\cE'\to\s^{*}(\cE')^{\vee}$ be the Hermitian map induced from $a$. Then we have an equality in the Chow group
\begin{equation}\label{KR vs ShtM} 
\z^{r}_{\cL_{1},\cdots, \cL_{n}}(a')|_{\cZ^{r}_{\cE}(a)}=[\Sht^{r}_{\cM(n,n)}]|_{\cZ^{r}_{\cE}(a)}\in \Ch_{0}(\cZ^{r}_{\cE}(a)).
\end{equation}
Here the restriction on the RHS is via \eqref{eq: hitchin shtuka decomposition E and A inj}, noting that $\cZ_{\cE}^r(a)^{\circ}=\cZ_{\cE}^r(a)$, and the class $[\Sht^{r}_{\cM(n,n)}]$ is as in Definition \ref{def: hitchin shtukas class}.

In particular, the cycle class $[\cZ^{r}_{\cE}(a)]$ as in Definition-Proposition \ref{defn:general KR} is well-defined (i.e., independent of the choice of a good framing). 
\end{thm}

Below we consider the case where $X'$ is geometrically connected. At the end of this section (\S\ref{ss:compare cycle split}) we comment on how to modify the argument in the case $X'=X\coprod X$ or $X'=X_{k'}$, where $k'/k$ is the quadratic extension.

\subsection{First reductions}
For a vector bundle $\cE$ on $X'$ let $\mu_{\min}(\cE)\in \Q$ be the smallest slope that appears in the Harder-Narasimhan filtration of $\cE$.  For $\cE$ of rank $n$ and $a\in \cA_{\cE}(k)$, a good framing $s:\op_{i=1}^{n}\cL_{i}\inj \cE$ for $(\cE,a)$ is called {\em very good} if it satisfies the additional condition
\begin{enumerate}
\item[(3)] $\mu_{\min}(\cE)>\max\{\deg\cL_{i}+2g'-1\}_{1\le i\le n}$.
\end{enumerate}

Most of the work in this section will be devoted to proving the slightly weaker statement below. 
\begin{thm}\label{th:compare KR weak} Suppose $X'$ is connected. Then the identity \eqref{KR vs ShtM} holds if $s:\cE'=\op_{i=1}^{n}\cL_{i}\inj \cE$ is a very good framing of $(\cE,a)$.
\end{thm}

\begin{lemma}\label{l:weak implies strong} Theorem \ref{th:compare KR weak} implies Theorem \ref{th:compare KR}.
\end{lemma}
\begin{proof} Choose effective divisors $D_{i}$ on $X'$ ($1\le i\le n$) such that $\nu(D_{1}+\cdots+D_{n})$ is multiplicity-free and disjoint from $\nu^{-1}(D_{a})$. Let $\cL'_{i}=\cL_{i}(-D_{i})$. When the $D_{i}$'s have sufficiently large degree,  the resulting map
\begin{equation}
s': \op_{i=1}^{n}\cL'_{i}\inj \op_{i=1}^{n}\cL_{i}\inj \cE
\end{equation}
is a very good framing. Let $a''$ be the induced Hermitian map $\op\cL'_{i} \rightarrow \sigma^* (\op\cL'_{i})^{\vee}$. By Theorem \ref{th:compare KR weak} we have
\begin{equation}
\z^{r}_{\cL'_{1},\cdots, \cL'_{n}}(a'')|_{\cZ^{r}_{\cE}(a)}=[\Sht^{r}_{\cM(n,n)}]|_{\cZ^{r}_{\cE}(a)}.
\end{equation}
Therefore, to prove \eqref{KR vs ShtM} it suffices to show 
\begin{equation}
\z^{r}_{\cL_{1},\cdots, \cL_{n}}(a')|_{\cZ^{r}_{\cE}(a)}=\z^{r}_{\cL'_{1},\cdots, \cL'_{n}}(a'')|_{\cZ^{r}_{\cE}(a)}\in\Ch_{0}(\cZ^{r}_{\cE}(a)).
\end{equation}
Let $U'$ be the complement of $\cup_{i=1}^{n}\supp(D_{i}+\s D_{i})$ in $X'$. By construction, $U'$ contains $\nu^{-1}(D_{a})$, therefore $\cZ_{\cE}^{r}(a)|_{U'^{r}}=\cZ_{\cE}^{r}(a)$ by Lemma \ref{l:KR cycle legs}. Let $\z_{i}=[\cZ^{r}_{\cL_{i}}(a'_{ii})^{\circ}]|_{U'^{r}}\in \Ch_{r(n-1)}(\cZ^{r}_{\cL_{i}}(a'_{ii})^{\circ}|_{U'^{r}})$. Similarly define $\z'_{i}$ using $(\cL'_{i}, a''_{ii})$. Then it suffices to show the equality
\begin{equation}\label{zz}
(\z_{1}\cdot \z_{2}\cdot \dots \cdot\z_{n})|_{\cZ^{r}_{\cE}(a)}=(\z'_{1}\cdot \z'_{2}\cdot \dots \cdot\z'_{n})|_{\cZ^{r}_{\cE}(a)}\in \Ch_{0}(\cZ^{r}_{\cE}(a)),
\end{equation}
where the intersection products are taken over $\Sht^{r}_{U(n)}|_{U'^{r}}$. Applying Lemma \ref{l:change E}  to each injection $\cL'_{i}\inj \cL_{i}$, we see that $\cZ^{r}_{\cL_{i}}(a'_{ii})|_{U'^{r}}\inj \cZ^{r}_{\cL'_{i}}(a''_{ii})|_{U'^{r}}$ is open and closed. Therefore $\cZ^{r}_{\cL_{i}}(a'_{ii})^{\circ}|_{U'^{r}}\inj \cZ^{r}_{\cL'_{i}}(a''_{ii})^{\circ}|_{U'^{r}}$ is open. This shows that the fundamental class $\z_{i}$ is the open restriction of $\z'_{i}$ to $\cZ^{r}_{\cL_{i}}(a'_{ii})^{\circ}|_{U'^{r}}$. The equality \eqref{zz} then follows.
\end{proof}

\subsection{Auxiliary moduli spaces}
Let $\un d=(d_{i})_{1\le i\le n}\in \Z^{n}_{\ge0}$ and $e\in \Z_{\ge0}$. Write $d=\sum d_{i}$. 

Recall that $\cM_{e}\subset \cM(n,n)$ is the open-closed substack where $\chi(X',\cE)=-e$. Let $\cM_{\un d}$ be the moduli stack classifying $(\{\cL_{i}\}_{1\le i\le n}, (\cF,h),\{t'_{i}:\cL_{i}\to \cF\}_{1\le i\le n})$ where
\begin{itemize}
\item $\cL_{i}$ is a line bunde on $X'$ with $\chi(X',\cL_{i})=-d_{i}$ for $1\le i\le n$;
\item $(\cF,h)\in \Bun_{U(n)}$ satisfying
\begin{equation}\label{slope F}
\mu_{\min}(\cF)>\max\{-d_{i}+3g'-2\}_{1\le i\le n}.
\end{equation}
\item For each $1\le i\le n$, $t'_{i}: \cL_{i}\to \cF$ is an injective map (fiberwise over the test scheme).
\end{itemize}
We define $\Hk^{r}_{\cM_{\un d}}$ to be the moduli stack classifying
\begin{equation*}
(\{\cL_{i}\}_{1\le i\le n}, (\{x'_{i}\}_{1\le i\le r}, \{(\cF_{j},h_{j})\}_{0\le j\le r})\in \Hk^{r}_{U(n)},\{t'_{ij}:\cL_{i}\to \cF_{j}\})
\end{equation*}
where  $\cL_{i}$ are the same  as in $\cM_{\un d}$, each $\cF_{j}$ satisfies the analogue of \eqref{slope F} with $\cF$ replaced by $\cF_{j}$, and fiberwise injective maps $t'_{ij}:\cL_{i}\to \cF_{j}$ are compatible with the isomorphisms between $\cF_{j-1}$ and $\cF_{j}$ away from $x'_{i}$.

\begin{lemma}\label{l:low deg M sm} The stacks $\cM_{\un d}$ and $\Hk^{1}_{\cM_{\un d}}$ are smooth stacks of pure dimension $dn - (n^2-2n)(g-1)$.
\end{lemma}
\begin{proof}
We first prove the statement for $\cM_{\un d}$. 
Consider the map $\g_{\un d}: \cM_{\un d}\to \prod_{i=1}^{n}\Pic^{-d_{i}+g'-1}_{X'}\times \Bun_{U(n)}$ sending $(\{\cL_{i}\},(\cF,h),\{t_i\})$ to $(\{\cL_{i}\},(\cF,h))$.  For $(\cF,h)\in \Bun_{U(n)}$ and $\cL_{i}\in \Pic^{-d_{i}+g'-1}_{X'}$, the condition $\mu_{\min}(\cF)\ge\max\{-d_{i}+3g'-2\}_{1\le i\le n}=\max\{\deg\cL_{i}+2g'-1\}_{1\le i\le n}$ guarantees that $\Ext^{1}(\cL_{i}, \cF)=\Hom(\cF,\cL_{i}\ot\om_{X'})^{\vee}=0$. Noting that $\deg \cF = n(g'-1)$, the Riemann-Roch formula implies that $\g_{\un d}$ exhibits $\cM_{\un d}$ as an open substack of a vector bundle of rank $\dim\Hom(\op\cL_{i}, \cF)=-n\sum_{i}\deg\cL_{i}=dn-n^{2}(g'-1)$ over the base. In particular, $\cM_{\un d}$ is smooth and equidimensional. Since $\dim \Pic^{-d_{i}+g'-1}_{X'}=g'-1$ and $\dim \Bun_{U(n)}=n^{2}(g-1)$, we conclude that 
\begin{align*}
\dim\cM_{\un d} & =dn-n^{2}(g'-1) + n (g'-1) + n^{2}(g-1) = dn- (n^2-2n)(g-1).
\end{align*}

The argument for $\Hk^{1}_{\cM_{\un d}}$ is similar. The natural map 
$$
\Hk^{1}_{\cM_{\un d}}\to \prod_{i=1}^{n}\Pic^{-d_{i}+g'-1}_{X'}\times\Hk^{1}_{U(n)}
$$
exhibits $\Hk^{1}_{\cM_{\un d}}$ as an open substack of a vector bundle of rank $\dim\Hom(\op\cL_{i}, \cF_{0}\cap \cF_{1})=dn-n$ (using that $\deg(\cF_{0}\cap\cF_{1})=n(g'-1)-1$) over the base. Here we need the stronger inequality $\mu_{\min}(\cF_{0})>\max\{-d_{i}+3g'-2\}$ to guarantee $\mu_{\min}(\cF_{0}\cap\cF_{1})\ge\max\{-d_{i}+3g'-2\}$. In particular, $\Hk^{1}_{\cM_{\un d}}$ is smooth and equidimensional, and 
\begin{align*}
\dim \Hk^{1}_{\cM_{\un d}} & =dn-n -n^{2}(g'-1)+ n(g'-1)+\dim \Hk^{1}_{U(n)}  \\
&= dn-n -2n^{2}(g-1)+ 2n(g-1) + n + n^2(g-1) =\dim\cM_{\un d},
\end{align*}
as desired. 
\end{proof}

Let $\cM_{\un d,e}$ be the  moduli stack of $(\{\cL_{i}\}_{1\le i\le n}, \cE, \cF,h,s,t)$ where 
\begin{itemize}
\item $\cL_{i}\in \Pic_{X'}$ satisfies $\chi(X',\cL_{i})=-d_{i}$ for $i=1,\cdots, n$; 
\item $\cE\in \Bun_{\GL_n'}$ satisfies $\chi(X',\cE)=-e$;
\item $(\cF,h)\in \Bun_{U(n)}$; 
\item $t: \cE\to \cF$ is an injective map;
\item $s: \op_{i=1}^{n}\cL_{i}\to\cE$ is a very good framing for $(\cE,a)$, where $a=\s^{*}t^{\vee}\c h\c t$ is the induced Hermitian map on $\cE$.
\end{itemize}
Note that being a very good framing requires $-d_{i}<\mu_{\min}(\cE)-(3g'-2)$ for all $i$, which imposes an open condition on $\cE$.  We view $\cM_{\un d,e}$ as a correspondence
\begin{equation}\label{corr M}
\xymatrix{ & \cM_{\un d, e}\ar[dl]_{w_{1}}\ar[dr]^{w_{2}}\\
\cM_{\un d} & &\cM_{e}\,.}
\end{equation}
Here $w_{1}$ records $\op_{i=1}^{n}\cL_{i}\xr{t\circ s}\cF$ and $w_{2}$ records $\cE\xr{t}\cF$.

We denote the Hitchin bases for $\cM_{\un d}, \cM_{e}$ and $\cM_{\un d,e}$ by $\cA_{\un d}$,  $\cA_{e}$ and $\cA_{\un d,e}$ respectively. Here $\cA_{\un d}$ parametrizes $(\{\cL_{i}\}_{1\le i\le n},a'=(a'_{ij}))$ (where $a'_{ij}: \cL_{j}\to \s^{*}\cL_{i}^{\vee}$) such that $a':\op\cL_{i}\to \s^{*}(\op\cL_{i})^{\vee}$ is an injective Hermitian map.  The base $\cA_{e}$ classifies $(\cE, a)$ with $a$ an injective Hermitian map. The base $\cA_{\un d,e}$ is the moduli stack of $(\{\cL_{i}\}_{1\le i\le n}, \cE, s,a)$ where  $(\cE,a)\in \cA_{e}$, $s: \op_{i=1}^{n}\cL_{i}\to\cE$ is a very good framing of $(\cE,a)$ and $\cL_{i}\in \Pic_{X'}$ with $\chi(X',\cL_{i})=-d_{i}$.  We view $\cA_{\un d,e}$ as a correspondence 
\begin{equation}\label{corr A}
\xymatrix{ & \cA_{\un d, e}\ar[dl]_{u_{1}}\ar[dr]^{u_{2}}\\
\cA_{\un d} & &\cA_{e}\,.}
\end{equation}
We have Hitchin maps
\begin{eqnarray*}
f_{\un d}: &\cM_{\un d}\to \cA_{\un d},\\
f_{e}: &\cM_{e}\to \cA_{e},\\
f_{\un d, e}: &\cM_{\un d, e}\to \cA_{\un d,e}.
\end{eqnarray*} 
These maps together give a map of correspondences \eqref{corr M} to \eqref{corr A}. Note $f_{\un d}$ is not necessarily proper because we have imposed an open condition on the minimal slope of $\cF$.

Similarly we define the Hecke version $\Hk^{r}_{\cM_{\un d ,e}}$ of $\cM_{\un d,e}$ as the moduli stack of $(\{x'_{i}\}_{1\le i\le r}, \{\cE\to\cF_{i}\}_{0\le i\le r})\in \Hk^{r}_{\cM_{e}}$ together with a very good framing $s:\op \cL_{i}\inj \cE$ for $(\cE,a)$ with $\chi(X',\cL_{i})=-d_{i}$. Again we view $\Hk^{r}_{\cM_{\un d,e}}$ as a correspondence
\begin{equation*}
\xymatrix{ & \Hk^{r}_{\cM_{\un d,e}}\ar[dl]_{h_{1}}\ar[dr]^{h_{2}}\\
\Hk^{r}_{\cM_{\un d}} & &\Hk^{r}_{\cM_{e}}\,.}
\end{equation*}

\begin{lemma}\label{l:uwh et}
The maps $w_{1},u_1$ and $h_{1}$ are \'etale.
\end{lemma}
\begin{proof}
We first prove that $u_{1}$ is \'etale.  Let $(X_{d-e}\times X_{e})^{\hs}$ be the open subscheme of divisors $(D_{1},D_{2})\in X_{d-e}\times X_{e}$ such that $D_{1}$ is multiplicity-free and disjoint from $D_{2}$; let $(X'_{d-e}\times X_{e})^{\hs}$ be the preimage of $(X_{d-e}\times X_{e})^{\hs}$ in $X'_{d-e}\times X_{e}$. We have a map $\a: X'_{d-e}\times X_{e}\to X_{d}$ sending $(D'_{1}, D_{2})$ to $\nu(D_{1})+D_{2}$. Let $\a^{\hs}$ be the restriction of $\a$ to  $(X'_{d-e}\times X_{e})^{\hs}$. By factorizing $\a^{\hs}$ as the composition 
\begin{equation}
 (X'_{d-e}\times X_{e})^{\hs} \xr{\nu_{d-e}\times \Id}  (X_{d-e}\times X_{e})^{\hs}\xr{\mrm{add}}X_{d} 
 \end{equation}
we see that $\a^{\hs}$ is \'etale.  From the definition we have a map
\begin{equation}\label{Ade}
j=(u_{1}, j'_{d-e}, j_{e}): \cA_{\un d,e}\to \cA_{\un d}\times_{X_{d}}(X'_{d-e}\times X_{e})^{\hs}.
\end{equation}
where $j'_{d-e}: \cA_{\un d,e}\to X'_{d-e}$ sends $(\{\cL_{i}\}, \op \cL_{i}\xr{s}\cE,a)$ to $\Div(s)$ (the divisor of $\det(s)$) and $j_{e}: \cA_{\un d,e}\to X_{e}$ sends it to $D_{a}$  (see Definition \ref{def:Da}).  The map $\cA_{\un d}\to X_{d}$ used in the fiber product records the divisor $D_{a'}$ of the Hermitian map $a'$ on $\op \cL_{i}$. We claim that $j$ is an open immersion. Indeed, given $(\{\cL_{i}\},a')\in \cA_{\un d}$ and $(D'_{1}, D_{2})\in (X'_{d-e}\times X_{e})^{\hs}$ such that $\nu(D_{1})+D_{2}=D_{a'}$, by the disjointness of $D'_{1},\s(D'_{1})$ and $\nu^{-1}(D_{2})$, there is one and only one coherent sheaf $\cE$ such that $\op\cL_{i}\subset \cE\subset \s^{*}(\op\cL_{i})^{\vee}$, $\cE/\op\cL_{i}$ is supported on $D'_{1}$, and $\s^{*}(\op\cL_{i})^{\vee}/\cE$ is supported away from $D'_{1}$. This would give a very good framing of $\cE$ if the open condition $\mu_{\min}(\cE)>\max\{-d_{i}+3g'-2\}_{1\le i\le n}$ is satisfied. This shows that $j$ is an open immersion. Since $\a^{\hs}$ is \'etale, we conclude that $\cA_{\un d,e}$ is \'etale over $\cA_{\un d}$.

To show $w_{1}$ is \'etale, we observe that $\cM_{\un d,e}\cong \cM_{\un d}\times_{\cA_{\un d}}\cA_{\un d,e}$. Since $u_{1}$ is \'etale, so is $w_{1}$. 



Finally,  $\Hk^{r}_{\cM_{\un d,e}}$ is the open substack of $\Hk^{r}_{\cM_{\un d}}\times_{\cA_{\un d}}\cA_{\un d,e}$  where the legs avoid $\Div(s)$. Since $u_{1}$ is \'etale,  so is $h_{1}$.
\end{proof}


\subsection{Auxiliary Hitchin shtukas}
We define $\Sht^{r}_{\cM_{\un d}}$ and $\Sht^{r}_{\cM_{\un d,e}}$ as the fiber product
\begin{equation}\label{defn ShtMde}
\xymatrix{ \Sht^{r}_{\cM_{\un d}}\ar[d] \ar[r] & \Hk^{r}_{\cM_{\un d}}\ar[d]^{(\pr_{0},\pr_{r})}& \Sht^{r}_{\cM_{\un d,e}}\ar[d] \ar[r] & \Hk^{r}_{\cM_{\un d,e}}\ar[d]^{(\pr_{0},\pr_{r})}\\
\cM_{\un d}\ar[r]^-{(\Id, \Frob)} & \cM_{\un d}\times \cM_{\un d} & \cM_{\un d,e}\ar[r]^-{(\Id, \Frob)} & \cM_{\un d,e}\times \cM_{\un d,e}}
\end{equation}

The maps $w_{i}$ and $h_{i}$ induce maps
\begin{equation}\label{corr ShtM}
\xymatrix{& \Sht^{r}_{\cM_{\un d,e}}\ar[dl]_{u_{1}}\ar[dr]^{u_{2}}\\
\Sht^{r}_{\cM_{\un d}} & & \Sht^{r}_{\cM_{e}}}
\end{equation}
The stack $\Sht^{r}_{\cM_{\un d,e}}$ decomposes into a disjoint union of open-closed substacks indexed by $(\{\cL_{i}\}, \op\cL_{i}\xr{s} \cE, a)\in \cA_{\un d,e}(k)$ 
\begin{equation}
\Sht^{r}_{\cM_{\un d,e}}=\coprod_{(\op\cL_{i}\to \cE, a)\in \cA_{\un d,e}(k)}\cZ^{r}_{\cE}(a).
\end{equation}
Correspondingly, the diagram \eqref{corr ShtM} decomposes into the disjoint union indexed by $\cA_{\un d,e}(k)$ of diagrams of the form
\begin{equation}
\xymatrix{ & \cZ^{r}_{\cE}(a)\ar@{^{(}->}[dl]_{u_{1}}\ar@{=}[dr]^{u_{2}}\\
\cZ^{r}_{\op\cL_{i}}(a')^{\hs} & & \cZ^{r}_{\cE}(a) }
\end{equation}
Here $\cZ^{r}_{\op\cL_{i}}(a')^{\hs}\subset \cZ^{r}_{\op\cL_{i}}(a')$ (where $a'$ is the Hermitian map on $\op\cL_{i}$ induced from $a$) is cut by the open condition $\mu_{\min}(\cF_{j})>\max\{-d_{i}+3g'-2\}$ for all $0\le j\le r$. From this description  and Corollary \ref{c:ZE open closed}, we see that:

\begin{lemma}\label{l:ShtM isom}
The map $u_{1}$ (resp. $u_{2}$), when restricted to each connected component of $\Sht^{r}_{\cM_{\un d,e}}$, is an isomorphism onto a connected component of $\Sht^{r}_{\cM_{\un d}}$ (resp. $\Sht^{r}_{\cM_{e}}$).
\end{lemma}

\subsection{Zero cycles on auxiliary Hitchin shtukas}
Similar to the definition of $[\Sht^{r}_{\cM_{e}}]$ given in \S\ref{ss:ShtM cycle}, we define $0$-cycles supported on $\Sht^{r}_{\cM_{\un d,e}}$ and $\Sht^{r}_{\cM_{\un d}}$ as follows.

We rewrite $\Sht^{r}_{\cM_{\un d}}$ as the fiber product
\begin{equation}\label{rewrite ShtM und}
\xymatrix{\Sht^{r}_{\cM_{\un d}}\ar[d]\ar[r] & (\Hk^{1}_{\cM_{\un d}})^{r}\ar[d]^{(\pr_{0},\pr_{1})^{r} \times \Delta} \times \cM_{\un d}\\
(\cM_{\un d})^{r+1}\ar[r]^{\Phi_{\cM_{\un d}}} & (\cM_{\un d})^{2r+2}
}
\end{equation}
Here $\Phi_{\cM_{\un d}}=\Phi^{r}_{\cM_{\un d}}$ and the vertical maps are defined as in Definition \ref{def:Phi}. By the smoothness of $\Hk^{1}_{\cM_{\un d}}$ and $\cM_{\un d}$ proved in Lemma \ref{l:low deg M sm} and the dimension calculation there, we define
\begin{equation}
[\Sht^{r}_{\cM_{\un d}}]:=\Phi_{\cM_{\un d}}^{!}[(\Hk^{1}_{\cM_{\un d}})^{r} \times \cM_{\un d}]\in\Ch_{0}(\Sht^{r}_{\cM_{\un d}}).
\end{equation}

Similarly, using the Cartesian diagram
\begin{equation}\label{rewrite ShtMde}
\xymatrix{\Sht^{r}_{\cM_{\un d, e}}\ar[d]\ar[r] & (\Hk^{1}_{\cM_{\un d, e}})^{r}\ar[d]^{(\pr_{0},\pr_{1})^{r} \times \Delta} \times \cM_{\un d, e}\\
(\cM_{\un d ,e})^{r+1}\ar[r]^{\Phi_{\cM_{\un d,e}}} & (\cM_{\un d, e})^{2r+2}
}
\end{equation}
and the smoothness and dimension calculations of $\Hk^{1}_{\cM_{\un d,e}}$ and $\cM_{\un d, e}$ (which follow from Lemma \ref{l:uwh et} and Lemma \ref{l:low deg M sm}), we define
\begin{equation}
[\Sht^{r}_{\cM_{\un d, e}}]:=\Phi_{\cM_{\un d,e}}^{!}[(\Hk^{1}_{\cM_{\un d, e}})^{r} \times \cM_{\un d, e}]\in\Ch_{0}(\Sht^{r}_{\cM_{\un d ,e}}).
\end{equation}

\begin{lemma}\label{l:u1}
We have $u_{1}^{*}[\Sht^{r}_{\cM_{\un d}}]=[\Sht^{r}_{\cM_{\un d,e}}]\in \Ch_{0}(\Sht^{r}_{\cM_{\un d,e}})
$.
\end{lemma}
\begin{proof}
This is because the maps $w_{1}$ and $h_{1}$ are both \'etale by Lemma \ref{l:uwh et}.
\end{proof}

\begin{lemma}\label{l:u2}
We have $u_{2}^{*}[\Sht^{r}_{\cM_{e}}]=[\Sht^{r}_{\cM_{\un d,e}}]\in \Ch_{0}(\Sht^{r}_{\cM_{\un d,e}})
$.
\end{lemma}
\begin{proof}
The diagram \eqref{rewrite ShtMde} is obtained from \eqref{rewrite ShtM} (for $\cM_{d}$ replaced by $\cM_{e}$) by base changing termwise along the map of the following two Cartesian diagrams induced by $u_{2}: \cA_{\un d,e}\to \cA_{e}$:
\begin{equation}
\xymatrix{\cA_{\un d,e}(k)\ar[r]\ar[d] & (\cA_{\un d,e})^{r+1}\ar[d]^{\Delta^{r+1}}="a" & \cA_{e}(k)\ar[r] \ar[d]_{}="b" & (\cA_{e})^{r+1}\ar[d]^{\Delta^{r+1}}\\
(\cA_{\un d,e})^{r+1}\ar[r]^{\Phi_{\cA_{\un d,e}}} & (\cA_{\un d,e})^{2r+2} & (\cA_{e})^{r+1}\ar[r]^{\Phi_{\cA_{e}}} & (\cA_{e})^{2r+2}
}
\end{equation}
Note that $u_{2}: \cA_{\un d,e}\to \cA_{e}$ is smooth since it exhibits $\cA_{\un d,e}$ as an open substack of a vector bundle over $\cA_{e}$ (using the condition $\mu_{\min}(\cE)>\max\{-d_{i}+3g'-2\}$). We conclude by applying Proposition \ref{p:base change} below.
\end{proof}

\subsubsection{Compatibility of cycle classes under Hitchin base change} To state the next result, we need some notations. Suppose we are given:
\begin{itemize}
\item stacks $S,M$ and $H$ that are locally of finite type over $k$ and can be stratified  into locally closed substacks that are global quotient stacks;
\item the stack $M$ is smooth of pure dimension $N$ with a map $f:M\to S$;
\item a map $\wt h: H  \to S^r \times_{\D^r, S^{2r}}M^{2r}$ (the fiber product uses the $r$-fold product of the diagonal  $\D^{r}: S^r \to S^{2r}$). 
\end{itemize}
Define $h \co H  \rightarrow M^{2r}$ as $\wt{h}$ followed by the second projection. Form the Cartesian square
\begin{equation}\label{ShtH}
\xymatrix{\Sht_{H}\ar[d]\ar[r] &  H \times M \ar[d]^{h \times \Delta}\\
M^{r+1}\ar[r]^{\Phi_{M}} & M^{2r+2}}
\end{equation}
Let $u: S'\to S$ be a smooth representable morphism of pure relative dimension $D$. Let $M'=M\times_{S}S'$, $H'=H\times_{S^{r}}S'^{r}$ with natural maps $\wt h': H'\to S'^{r}\times_{\D^{r}, S'^{2r}}M'^{2r}$. Let $h':H'\to M'^{2r}$ be the resulting map. Let $u_{M}: M'\to M$ and $u_{H}: H'\to H$ be the natural maps. Form the Cartesian square
\begin{equation}\label{Sht'H}
\xymatrix{\Sht'_{H}\ar[d]\ar[r] &  H'\ar[d]^{ h' \times \Delta} \times M' \\
M'^{r+1}\ar[r]^{\Phi_{M'}} & M'^{2r+2}}
\end{equation}
Since $S^{r+1}\times_{\Phi_{S}, S^{2r+2}, \D^{r+1}}S^{r+1}=S(k)$,  $\Sht_{H}$ decomposes as 
\begin{equation}
\Sht_{H}=\coprod_{s\in S(k)}\Sht_{H}(s).
\end{equation}
Similarly $\Sht'_{H}$ decomposes into the disjoint union of $\Sht_{H}'(s')$ indexed by $s'\in S'(k)$.
Then the natural map $u_{\Sht}: \Sht'_{H}\to \Sht_{H}$ lifts to an isomorphism
\begin{equation}
\Sht'_{H}=\Sht_{H}\times_{S(k)}S'(k)=\coprod_{s'\in S'(k)}\Sht_{H}(u(s')).
\end{equation}

\begin{prop}\label{p:base change}  Let $\z\in \Ch_{*}(H)$ and $[M]$ be the fundamental class of $M$. Then we have
\begin{equation}
u_{\Sht}^{*}\Phi_{M}^{!}(\z \times [M])=\Phi_{M'}^{!}u_{H}^{*}(\z \times [M])\in\Ch_{*-rN}(\Sht'_{H}).
\end{equation}
\end{prop}
\begin{proof}
Consider first the diagram where all squares are Cartesian
\begin{equation}
\xymatrix{\Sht_{H}\times_{S^{r}}S'^{r}\ar[r]\ar[d]^{v} & H' \times M' \ar[r]\ar[d]^{u_{H} \times u_M} & S'^{r+1}\ar[d]^{u^{r+1}}\\
\Sht_{H}\ar[r]\ar[d] & H \times M \ar[r]\ar[d]^{h \times \Delta} & S^{r+1}\\
M^{r+1}\ar[r]^-{\Phi_{M}} & M^{2r+2}
}
\end{equation}
Here the top vertical arrows are smooth and representable. By the compatibility of  Gysin map with flat pullback \cite[Theorem 2.1.12(ix)]{Kr99}, we have
\begin{equation}\label{vPhiu}
v^{*}\Phi_{M}^{!}(\z \times [M])=\Phi_{M}^{!}u_{H}^{*}(\z \times [M])\in \Ch_{*-rN+(r+1)D}(\Sht_{H}\times_{S^{r+1}}S'^{r+1}).
\end{equation}
Here we recall that $D$ is the relative dimension of $u$.
We have
\begin{equation}\label{ShtHSS'}
\Sht_{H}\times_{S^{r+1}}S'^{r+1}=\coprod_{s\in S(k)}\Sht_{H}(s)\times (S'_{s})^{r+1}
\end{equation}
where $S'_{s}=u^{-1}(s)$, which is a smooth scheme over $k$. Factorize $u_{\Sht}$ as the composition
\begin{equation}
\begin{tikzcd}
\Sht'_{H} \ar[r, "{i}"] & \Sht_{H}\times_{S^{r}}(S'^{r}) \ar[r, "v"]  & \Sht_{H}.
\end{tikzcd}
\end{equation}
From \eqref{ShtHSS'} we see that $i$ is a regular embedding of codimension $(r+1)D$. Now $v$ and $u_{\Sht}$ are both smooth. Applying \cite[Theorem 2.1.12(ix)]{Kr99} we have $u^{*}_{\Sht}(-)=i^{!}v^{*}$ as maps $\Ch_{*}(\Sht_{H})\to \Ch_*(\Sht'_{H})$. Therefore
\begin{equation}\label{iuPhi}
u_{\Sht}^{*}\Phi_{M}^{!}(\z \times [M])=i^{!}v^{*}\Phi_{M}^{!}(\z \times [M])\in \Ch_{*-rN}(\Sht_{H}').
\end{equation}	
On the other hand, consider the following diagram where all squares are Cartesian
\begin{equation}
\xymatrix{\Sht'_{H} \ar[r]^{i}\ar[d] & \Sht_{H}\times_{S^{r+1}}S'^{r+1}\ar[r]\ar[d] & H' \times M' \ar[d]^{ h' \times \Delta}\\
M'^{r+1}\ar[r]^-{\Phi_{1}} &  M^{r+1}\times_{\Phi_{S}\c f^{r+1}, S^{2r+2}}S'^{2r+2}\ar[r]^-{\Phi_{2}}\ar[d] & M'^{2r+2}\ar[d]^{u^{2r+2}_{M}}\\
& M^{r+1}\ar[r]^{\Phi_{M}} & M^{2r+2}
}
\end{equation}
Here $\Phi_{1}$ is the base change of $\Phi_{S'}$ and $\Phi_{2}$ is the base change of $\Phi_{M}$. The outer square of the top rows give \eqref{Sht'H}. By the transitivity of Gysin maps, we have
\begin{equation}
\Phi_{M'}^{!}u_{H}^{*}(\z\times [M])=\Phi^{!}_{1}\Phi^{!}_{2}u_{H}^{*}(\z \times [M]).
\end{equation}
Since $u^{2r+2}_{M}$  is smooth representable, we have $\Phi_{2}^{!}u_{H}^{*}(\z \times [M])=\Phi_{M}^{!}u_{H}^{*}(\z \times [M])$. Hence \begin{equation}
\Phi_{M'}^{!}u_{H}^{*}(\z\times [M])=\Phi^{!}_{1}\Phi_{M}^{!}u_{H}^{*}(\z\times [M]).
\end{equation}
Since both $\Phi_{1}$ and $i$ are regular embeddings of the same codimension, we have $\Phi^{!}_{1}(-)=i^{!}(-)$ as maps $\Ch_{*}(\Sht_{H}\times_{S^{r+1}}S'^{r+1})\to \Ch_{*-(r+1)D}(\Sht'_{H})$, by \cite[Theorem 2.1.12(xi)]{Kr99} and \cite[Theorem 6.2(c)]{Ful98}. Therefore
\begin{equation}\label{iPhiu}
\Phi_{M'}^{!}u_{H}^{*}(\z\times [M])=i^{!}\Phi_{M}^{!}u_{H}^{*}(\z\times [M])\in \Ch_{*-rN}(\Sht'_{H}).
\end{equation}
Combining \eqref{vPhiu}, \eqref{iuPhi} and \eqref{iPhiu} we conclude
\begin{equation}
u_{\Sht}^{*}\Phi_{M}^{!}(\z\times [M])=i^{!}v^{*}\Phi_{M}^{!}(\z\times [M])=i^{!}\Phi_{M}^{!}u_{H}^{*}(\z\times [M])=\Phi_{M'}^{!}u_{H}^{*}(\z\times [M])\in \Ch_{*-rN}(\Sht'_{H}).
\end{equation}
\end{proof}

\begin{lemma}\label{l:M und KR} Let  $(\{\cL_{i}\}_{1\le i\le n}, a')\in \cA_{\un d}(k)$ and $\cE':=\op_{i=1}^{n}\cL_{i}$.  Then we have an equality
\begin{equation}
[\Sht^{r}_{\cM_{\un d}}]|_{\cZ^r_{\cE'}(a')}= \z^r_{\cL_{1},\cdots, \cL_{n}}(a')\in \Ch_{0}(\cZ^r_{\cE'}(a')).
\end{equation}
For the definition of $\z^r_{\cL_{1},\cdots, \cL_{n}}(a')$ see \S\ref{ssec: KR cycle classes}. 
\end{lemma}
\begin{proof}
We will apply the Octahedron Lemma \cite[Theorem A.10]{YZ} to a diagram of moduli stacks in our setting. Since the Octahedron Lemma requires certain stacks in question to be Deligne-Mumford, we need to rigidify our moduli stacks to satisfy these requirements. This is a minor technical issue which we encourage the reader to ignore: it is simply because the Octahedron Lemma in \cite{YZ} is not stated and proved in the most general form.

Let $v\in |X'|$. Let $P_{v}$ be the moduli space (a scheme!) of line bundles on $X'$ together with a trivialization of their fibers over $v$. Let $\G_{v}=\Res^{k_{v}}_{k}\G_{m}$. Then $P_{v}\to\Pic_{X'}$ is a $\G_{v}$-torsor.

Now for each moduli stacks $\cM_{\un d}, \cA_{\un d}$ and $\Hk^{r}_{\dot\cM_{\un d}}$ that involve an $n$-tuple of line bundles $\{\cL_{i}\}$, we write  $\dot\cM_{\un d}, \dot\cA_{\un d}$ and $\dot\Hk^{r}_{\dot\cM_{\un d}}$ to mean their rigidified versions where $\cL_{i}\in\Pic_{X'}$ is  replaced by $\dot\cL_{i}\in P_{v}$. Note that we do not impose any compatibility condition between the rigidifcation on $\cL_{i}$ and the rest of the structures classified by these moduli stacks. Define $\Sht^{r}_{\dot\cM_{\un d}}$ using the dotted version of the left one of the Cartesian diagrams in  \eqref{defn ShtMde}.

Note that $\dot\cM_{\un d}, \dot\cA_{\un d}$ and $\dot\Hk^{r}_{\dot\cM_{\un d}}$ are now schemes, and they are $\G_{v}^{n}$-torsors over their undotted counterparts. The dotted version of Lemma \ref{l:low deg M sm} remains valid if we add $n\deg(v)=\dim\G_{v}^{n}$ to the dimensions.  Also, $\Sht^{r}_{\dot\cM_{\un d}}\simeq \Sht^{r}_{\cM_{\un d}}\times_{\cA_{\un d}(k)}\dot\cA_{\un d}(k)$.  Since $\dot\cA_{\un d}(k)\to\cA_{\un d}(k)$ is surjective, to prove the Lemma, it suffices to prove its dotted version: for any $(\{\dot\cL_{i}\}, a')\in \dot\cA_{\un d}(k)$, writing $\cE':=\op_{i=1}^{n}\cL_{i}$, then there is an open and closed embedding $\cZ^r_{\cE'}(a')\incl \Sht^{r}_{\dot\cM_{\un d}}$ (using the rigifications $\dot\cL_{i}$ of $\cL_{i}$); then we shall prove
\begin{equation}\label{dot cycle eqn}
[\Sht^{r}_{\dot\cM_{\un d}}]|_{\cZ^r_{\cE'}(a')}=\z^r_{\cL_{1},\cdots, \cL_{n}}(a')\in \Ch_{0}(\cZ^r_{\cE'}(a')).
\end{equation}


For $i=1,\cdots, n$,  let $\cN_{d_{i}}$ be the open substack of $\cM(1,n)$ consisting of points $(\cL\inj \cF,h)$ where $\chi(X',\cL)=-d_{i}$ and $\mu_{\min}(\cF)>-d_{i}+3g'-2$. Similarly define $\Hk^{1}_{\cN_{d_{i}}}$ and $\Sht^{r}_{\cN_{d_{i}}}$; these are open substacks of $\Hk^1_{\cM(1,n)}$ and $\Sht^r_{\cM(1,n)}$ respectively. 

Let $\dot\cN_{d_{i}}, \Hk^{1}_{\dot\cN_{d_{i}}}$ be the rigidified versions of $\cN_{d_{i}}$ and $\Hk^{1}_{\cN_{d_{i}}}$ where $\cL\in \Pic_{X'}$ is replaced by $\dot\cL\in P_{v}$. Let $\om_{i}: \dot\cN_{d_{i}}\to \Bun_{U(n)}$ and $\wt\om_{i}: \Hk^{1}_{\dot\cN_{d_{i}}}\to \Hk^{1}_{U(n)}$ be the forgetful maps.

We shall apply the Octahedron Lemma \cite[Theorem A.10]{YZ} to the following diagram:
\begin{equation}\label{9 term}
\xymatrix{ (\Hk^{1}_{U(n)})^{r} \times \Bun_{U(n)}\ar[r]^-{\Delta}\ar[d]^{(\pr_{0},\pr_{1})^{r} \times \Delta} & \prod_{i=1}^{n} ((\Hk^{1}_{U(n)})^r \times \Bun_{U(n)} ) \ar[d]^{\prod ((\pr_{0},\pr_{1})^{r} \times \Delta)} &  \prod_{i=1}^{n} ((\Hk^{1}_{\dot\cN_{d_{i}}})^r \times \dot\cN_{d_{i}}) \ar[l]_-{\prod (\wt\om_{i}^{r} \times \omega_i)}\ar[d]^{\prod ((\pr_{0},\pr_{1})^{r} \times \Delta )}\\
\Bun_{U(n)}^{2r+2} \ar[r]^-{\Delta} & \prod_{i=1}^{n}\Bun_{U(n)}^{2r+2} & \prod_{i=1}^{n}\dot\cN_{d_{i}}^{2r+2}\ar[l]_-{\prod \om_{i}^{2r+2}}\\
\Bun_{U(n)}^{r+1} \ar[r]^-{\Delta}\ar[u]_{\Phi_{\Bun_{U(n)}}} & \prod_{i=1}^{n}\Bun_{U(n)}^{r+1}\ar[u]_{\prod\Phi_{\Bun_{U(n)}}}  & \prod_{i=1}^{n}\dot\cN_{d_{i}}^{r+1}\ar[u]_{\prod\Phi_{\dot\cN_{d_{i}}}}\ar[l]_-{\prod\om_{i}^{r+1}} 
}
\end{equation}
The fiber products of the three columns are
\begin{equation}\label{col Sht}
\xymatrix{\Sht^{r}_{U(n)}\ar[r]^-{\Delta} & \prod_{i=1}^{n}\Sht^{r}_{U(n)}  & \prod_{i=1}^{n}\Sht^{r}_{\dot\cN_{d_{i}}}\ar[l]}
\end{equation}
where
\begin{equation}\label{ShtN decomp}
\Sht^{r}_{\dot\cN_{d_{i}}}=\coprod_{\substack{\dot\cL_{i}\in P_{v}(k), \\ \chi(X',\cL_{i})=-d_{i}}}\left(\cZ^{r}_{\cL_{i}}(0)^{*}\coprod\left(\coprod_{ a'_{ii}\in \cA_{\cL_{i}}(k)}\cZ^{r}_{\cL_{i}}(a'_{ii})\right)\right).
\end{equation}

Let $\cMa_{\un d}$ be the moduli stack of $(\{\dot\cL_{i}\}, \op_{i=1}^{n}\cL_{i}\xr{t'} \cF,h)$ defined similarly as $\dot\cM_{\un d}$ but without the condition that $t'$ be injective, only that $t'_{i}=t|_{\cL_{i}}$ be injective. Then $\dot\cM_{\un d}\incl \cMa_{\un d}$ is open. Similarly define $\Hk^{1}_{\cMa_{\un d}}$ and $\Sht^{r}_{\cMa_{\un d}}$. Note that $\cMa_{\un d}$ is exactly the fiber product of
\begin{equation}\label{Mbar d}
\xymatrix{\Bun_{U(n)}\ar[r]^-{\Delta} & \prod_{i=1}^{n}\Bun_{U(n)} & \prod_{i=1}^{n}\dot\cN_{d_{i}}\ar[l]_-{\prod\om_{i}}.}
\end{equation}
Similar remarks apply to $\Hk^{1}_{\cMa_{\un d}}$.  Therefore the fiber products of the three rows of \eqref{9 term} are
\begin{equation}\label{row Hk}
\xymatrix{ (\Hk^{1}_{\cMa_{\un d}})^{r} \times \cMa_{\un d} \ar[d]^{(\pr_{0},\pr_{1})^{r} \times \Delta} \\
\cMa_{\un d}^{2r+2} \\
\cMa_{\un d}^{r+1}\ar[u]_{\Phi_{\cMa_{\un d}}}}
\end{equation}
The common fiber product of \eqref{col Sht} and \eqref{row Hk} is $\Sht^{r}_{\cMa_{\un d}}$, which decomposes as a disjoint union over the groupoid $\cB(k)$ of $(\dot\cL_{i},a'_{ii})_{1\le i\le n}$ where $\dot\cL_{i}\in P_{v}(k)$ with $\chi(X',\cL_{i})=-d_{i}$ and $a'_{ii}:\cL_{i}\to\s^{*}\cL_{i}^{\vee}$ injective Hermitian. For a point ${(\dot\cL_{i}, a'_{ii})_{1\le i\le n}}$ of $\cB(k)$, we have 
\begin{equation}
\Sht^{r}_{\cMa_{\un d}}|_{(\dot\cL_{i}, a'_{ii})_{1\le i\le n}} = \cZ^{r}_{\cL_{1},\cdots,\cL_{n}}(a'_{11},\cdots, a'_{nn})^{\circ},
\end{equation}
where $\Sht^{r}_{\cMa_{\un d}}|_{(\dot\cL_{i}, a'_{ii})_{1\le i\le n}}$ means the pullback of $
\Sht^{r}_{\cMa_{\un d}}$ to $\Spec k$ along the corresponding $\Spec k \rightarrow \cB(k)$. 

We check that the assumptions for applying the Octahedron Lemma are satisfied (the numbering below refers to that in \cite[Theorem A.10]{YZ}). 
\begin{enumerate}
\item All members in the diagram \eqref{9 term} are smooth and equidimensional. This is clear for $\Bun_{U(n)}$ and $\Hk^1_{U(n)}$. The same argument as in Lemma \ref{l:low deg M sm} proves that $\dot\cN_{d_{i}}$ and $\Hk^{1}_{\dot\cN_{d_{i}}}$ are smooth of pure dimension $d_{i}n+(n^2-2n+2)(g-1)+\deg(v)$.

\item We check that, in forming the fiber products of the middle and bottom rows and the left and middle columns,  the intersections are proper intersections with smooth equidimensional outcomes with the expected dimension. Here we use Lemma \ref{lem: smooth Bun_G} and \ref{lem: shtuka geometry} to argue for the left and middle columns. For the rows, the same argument as in Lemma \ref{l:low deg M sm} proves that $\cMa_{\un d}$ and $\Hk^{1}_{\cMa_{\un d}}$ are smooth of the same dimension as $\dot\cM_{\un d}$, which is $dn-(n^2-2n)(g-1)+n\deg(v)$. This is the virtual dimension for $\cMa_{\un d}$ as the fiber product of \eqref{Mbar d}, since 
\begin{align*}
\sum_{i=1}^{n}\dim \cN_{d_{i}}-(n-1)\dim \Bun_{U(n)}  &=\sum_{i=1}^{n} \left(d_{i}n+(n^{2}-2n+2)(g-1)+\deg(v)\right)-n^{2}(n-1)(g-1) \\
&=dn-(n^{2}-2n)(g-1)+n\deg(v).
\end{align*}

\item We check the fiber products of the top row and right column of \eqref{9 term} satisfy the conditions for  \cite[A.2.10]{YZ}. The fiber product of the top row is also a proper intersection: this follows from the same calculation as for the middle and bottom rows. The fiber product of the right column is also a proper intersection: this uses the decomposition \eqref{ShtN decomp} and the calculation of the dimension of $\cZ^{r}_{\cL_{i}}(0)^{*}$ in Proposition \ref{p:dim ZL0} and the dimension of $\cZ^{r}_{\cL_{i}}(a'_{ii})^{*}$ in Proposition \ref{p:ShM dim}.

The only issue is that $(\Hk^{1}_{\cMa_{\un d}})^{r+1}$ may not be a Deligne-Mumford stack, which was part of the requirement of \cite[A.2.10]{YZ}. However, we argue that this is not really an issue. The proof of the Octahedron Lemma allows the following flexibility: since eventually we only care about the $0$-cycles restricted to $\Sht^{r}_{\dot\cM_{\un d}}$, in the middle steps of forming the fiber products, we may restrict to open substacks as long as the final fiber product contains $\Sht^{r}_{\dot\cM_{\un d}}$ and only need check the relevant requirements there. Now in \eqref{row Hk} we may restrict to the open substack $(\Hk^{1}_{\dot\cM_{\un d}})^{r+1}\subset (\Hk^{1}_{\cMa_{\un d}})^{r+1}$, which is  a scheme.

\item The same remark as above shows that it suffices to check that the fiber squares obtained from \eqref{col Sht} and \eqref{row Hk}, after replacing $\cMa_{\un d}$ by $\dot\cM_{\un d}$,  each satisfy the condition \cite[A.2.8]{YZ}. Therefore it suffices to check
\begin{itemize}
\item $\Sht^{r}_{\dot\cM_{\un d}}$ admits a finite flat presentation in the sense of \cite[Definition A.1]{YZ}. This is true because $\Sht^{r}_{\dot\cM_{\un d}}$ is a scheme.
\item The diagonal map $\Delta: \Sht^{r}_{U(n)}\incl \prod_{i=1}^{n}\Sht^{r}_{U(n)}$ is a regular local immersion. This is true because $\Sht^{r}_{U(n)}$ is a smooth Deligne-Mumford stack.
\item The map $\Phi_{\dot\cM_{\un d}}: \dot\cM_{\un d}^{r+1}\to \dot\cM_{\un d}^{2r+2}$ is a regular local immersion. This is true because $\dot\cM_{\un d}$ is a smooth equidimensional scheme by the dotted version of Lemma \ref{l:low deg M sm}. 
\end{itemize}
\end{enumerate}
The conclusion of the (variant of) Octahedron Lemma says that the following two elements in $\Ch_{0}(\Sht^{r}_{\cMa_{\un d}})$
\begin{equation}
\Phi_{\cMa_{\un d}}^{!}\Delta_{(\Hk^{1}_{U(n)})^{r}}^{!}[\prod_{i=1}^{n}(\Hk^{1}_{\dot\cN_{d_{i}}})^{r} \times \dot\cN_{d_{i}}] \quad \mbox{and}\quad \Delta_{\Sht^{r}_{U(n)}}^{!}(\prod\Phi_{\dot\cN_{d_{i}}})^{!}[\prod_{i=1}^{n}(\Hk^{1}_{\dot\cN_{d_{i}}})^{r} \times \dot\cN_{d_{i}}]
\end{equation}
become the same when restricted to $\Sht^{r}_{\dot\cM_{\un d}}$. Further restricting to $\cZ_{\cE'}^r(a')$ we get the desired identity \eqref{dot cycle eqn}. 
\end{proof}

\begin{proof}[Proof of Theorem \ref{th:compare KR weak}] Restricting the equality in Lemma \ref{l:M und KR} to $\cZ_{\cE}^r(a)$, which is open and closed in $\cZ_{\cE'}^r(a')$ by Corollary \ref{c:ZE open closed}, we get
\begin{equation}\label{ZL}
\z^r_{\cL_{1},\cdots, \cL_{n}}(a')|_{\cZ^r_{\cE}(a)}=[\Sht^{r}_{\cM_{\un d}}]|_{\cZ_{\cE}^r(a)}.
\end{equation}
For fixed $(\op\cL_{i}\inj \cE,a)\in \cA_{\un d,e}(k)$, $\cZ^{r}_{\cE}(a)$ can be viewed as a finite \'{e}tale cover of an open-closed substack in $\Sht^{r}_{\cM_{\un d}},\Sht^{r}_{\cM_{\un d,e}}$ and $\Sht^{r}_{\cM_{e}}$ by Lemma \ref{l:ShtM isom}. By Lemma \ref{l:u1} and Lemma \ref{l:u2} we have
\begin{equation}
[\Sht^{r}_{\cM_{\un d}}]|_{\cZ^r_{\cE}(a)}=[\Sht^{r}_{\cM_{\un d,e}}]|_{\cZ^r_{\cE}(a)}=[\Sht^{r}_{\cM_{e}}]|_{\cZ^{r}_{\cE}(a)}.
\end{equation}
Combining this with \eqref{ZL} proves the theorem.
\end{proof}

\subsection{Proof of Theorem \ref{th:compare KR} for $X'=X\coprod X$ or $X_{k'}$}\label{ss:compare cycle split} Here $k'/k$ is the quadratic extension. 

In the case $X'=X\coprod X$,  we have $\Bun_{U(n)}\cong \Bun_{\GL_{n}}$. We shall identify a Hermitian bundle $\cF$ on $X'$ with a pair of vector bundles $(\cF_{1},\cF_{2})$ equipped with an isomorphism $\cF_{2}\cong \cF_{1}^{\vee}$, each living on one copy of $X$.  A vector bundle $\cE$ on $X'$ of rank $n$ corresponds to two rank $n$ vector bundles $(\cE_{1},\cE_{2})$, each living on one copy of $X$. Now $\cA_{\cE}(k)$ is the set of injective maps $a: \cE_{1}\to \cE^{\vee}_{2}$.  A good framing $s:\op_{i=1}^{n}\cL_{i}\inj \cE$ for $(\cE,a)$ now consists of line bundles $\cL_{i}=(\cL_{i,1}, \cL_{i,2})$ ($1\le i\le n$) satisfying the same conditions in Definition \ref{def: auxiliary E'}; it is called {\em very good} if it satisfies the additional conditions
\begin{enumerate}
\item[($3_{1}$)] $\mu_{\min}(\cE_{1})>\max\{\deg\cL_{i,1}+2g-1\}_{1\le i\le n}$, and
\item[($3_{2}$)] $\mu_{\min}(\cE_{2})>\max\{\deg\cL_{i,2}+2g-1\}_{1\le i\le n}$.
\end{enumerate}
The same argument of Lemma \ref{l:weak implies strong} shows that it suffices to prove the analogue of Theorem \ref{th:compare KR weak}, i.e., prove Theorem \ref{th:compare KR} for very good framings.

In both the $X'=X\coprod X$ and $X'=X_{k'}$ case, we need to modify the definitions of $\cM_{\un d}$ and $\cM_{\un d,e}$ as follows. In the definition of $\cM_{\un d,e}$, we use the notion of very good framing just defined over geometric fibers of $X'_{S}\to S$ (which are of the form $X_{\ov s}\coprod X_{\ov s}$).  In the definition of $\cM_{\un d}$, we change the inequality \eqref{slope F} to two inequalities over the geometric fibers of $X'_{S}\to S$
\begin{eqnarray*}
\mu_{\min}(\cF_{1})>\max\{\deg\cL_{i,1}+2g-1\}_{1\le i\le n},\\
\mu_{\min}(\cF_{2})>\max\{\deg\cL_{i,2}+2g-1\}_{1\le i\le n}.
\end{eqnarray*}
The same inequalities should be imposed in the definition of $\cN_{d_{i}}$ that appear in the proof of Lemma \ref{l:M und KR}. With these changes, the argument for proving Theorem \ref{th:compare KR weak} goes through.

\section{Local intersection number and trace formula}\label{s:int trace}

\subsection{Local nature of the intersection problem}Recall from Proposition \ref{prop: KR cycle proper} that $\cZ^{r}_{\cE}(a)$ only depends on the Hermitian torsion sheaf $\cQ=\coker(a)$. In this subsection we show that the $0$-cycle class $[\cZ^{r}_{\cE}(a)]$ also only depends on $\cQ$. 

Recall the stacks $\Herm_{2d}=\Herm_{2d}(X'/X)$ and $\Lagr_{2d}$ from \S \ref{sec: Herm}. We have a self-correspondence $\Hk^{r}_{\Lagr_{2d}}$ of $\Lagr_{2d}$ over $\Herm_{2d}$: it classifies 
$(\cQ,h,\{\cL_{i}\}_{0\le i\le r})$ where $(\cQ,h)\in \Herm_{2d}$, $\cL_{i}\subset \cQ$ are Lagrangian subsheaves such that $\cL_{i}/(\cL_{i}\cap \cL_{i-1})$ has length one for $1\le i\le r$. 
Define the local version $\Sht^{r}_{\Lagr_{2d}}$ of $\Sht^{r}_{\cM_{d}}$ by the Cartesian diagram
\begin{equation}\label{Sht Lag}
\xymatrix{\Sht^{r}_{\Lagr_{2d}}\ar[d] \ar[r] & \Hk^{r}_{\Lagr_{2d}}\ar[d]^{(\pr_{0},\pr_{r})}\\
\Lagr_{2d}\ar[r]^-{(\Id,\Frob)} & \Lagr_{2d}\times\Lagr_{2d}}
\end{equation}
We have a decomposition into open-closed substacks
\begin{equation}
\Sht^{r}_{\Lagr_{2d}}=\coprod_{(\cQ,h)\in\Herm_{2d}(k)} \cZ^{r}_{\cQ}.
\end{equation}

\begin{lemma}
The stack $\Hk^{1}_{\Lagr_{2d}}$ is smooth of dimension zero.
\end{lemma}
\begin{proof}
We may identify $\Hk^{1}_{\Lagr_{2d}}$ with the moduli stack of $(0\subset \cL'\subset \cL\subset \cQ, h)$ where $(\cL\subset\cQ,h)\in\Lagr_{2d}$ and $\cL/\cL'$ has length one. Under the local chart for $\Herm_{2d}$ described in Lemma \ref{l:Herm local}, $\Hk^{1}_{\Lagr_{2d}}$ becomes $[\frp/P]$, where $P\subset \Og_{2d}=\Og(V)$ is the parabolic subalgebra stabilizing a pair of subspaces  $L'\subset L$ with $L$ Lagrangian and $\dim L'=d-1$. This local description implies that $\Hk^{1}_{\Lagr_{2d}}$ is smooth of dimension zero.
\end{proof}

Rewriting $\Sht^{r}_{\Lagr_{2d}}$ as the fibered product (cf. Definition \ref{def:Phi})
\begin{equation}\label{Sht Lagr}
\xymatrix{\Sht^{r}_{\Lagr_{2d}}\ar[r]\ar[d] & (\Hk^{1}_{\Lagr_{2d}})^{r} \times \Lagr_{2d} \ar[d]^{(\pr_{0},\pr_{1})^{r} \times \Delta}\\
(\Lagr_{2d})^{r+1}\ar[r]^{\Phi_{\Lagr_{2d}}} & (\Lagr_{2d})^{2r+2}}
\end{equation}
we define a $0$-cycle class
\begin{equation}
[\Sht^{r}_{\Lagr_{2d}}]:=\Phi_{\Lagr_{2d}}^{!}[(\Hk^{1}_{\Lagr_{2d}})^{r} \times \Lagr_{2d}]\in \Ch_{0}(\Sht^{r}_{\Lagr_{2d}}).
\end{equation}
Restricting to $\cZ^{r}_{\cQ}$ we get
\begin{equation}
[\cZ^{r}_{\cQ}]:=[\Sht^{r}_{\Lagr_{2d}}]|_{\cZ^{r}_{\cQ}}\in \Ch_{0}(\cZ^{r}_{\cQ}).
\end{equation}

Recall the maps $g_{\cM}: \cM_{d}\to \Lagr_{2d}$ and $g:\cA_{d}\to \Herm_{2d}$ defined in \S\ref{ss:vb to tor}. We also have a map  $g_{\Hk}: \Hk^{r}_{\cM_{d}}\to \Hk^{r}_{\Lagr_{2d}}$ sending $(\{x'_{i}\} , \{\cE\xr{t_{i}}\cF_{i}\})$ to $(\cQ=\coker(a), h_{\cQ}, \{\coker(t_{i})\})$ ($a$ is the induced Hermitian map on $\cE$).   The maps $g_{\Hk}, g_{\cM}$ and $g$ exhibit the diagram \eqref{def ShtMd} as the pullback of the diagram \eqref{Sht Lag} via the base change $g: \cA_{d}\to\Herm_{2d}$.  In particular we have a natural map
\begin{equation}
g_{\Sht}: \Sht^{r}_{\cM_{d}}\to \Sht^{r}_{\Lagr_{2d}}.
\end{equation}
For fixed $(\cE,a)\in \cA_{d}(k)$ with image $\cQ=\coker(a)\in \Herm_{2d}(k)$, $g_{\Sht}$ restricts to an isomorphism to the open-closed subschemes
\begin{equation}\label{gSht}
g_{\Sht}|_{\cZ^{r}_{\cE}(a)}: \cZ^{r}_{\cE}(a)\isom \cZ^{r}_{\cQ}.
\end{equation}

\begin{prop} We have an equality
\begin{equation}\label{ShtM ShtLag}
[\Sht^{r}_{\cM_{d}}]=g_{\Sht}^{*}[\Sht^{r}_{\Lagr_{2d}}]\in \Ch_{0}(\Sht^{r}_{\cM_{d}}).
\end{equation}
\end{prop} 
\begin{proof}
Apply Proposition \ref{p:base change} to the diagram \eqref{Sht Lagr}, the fundamental class $\z=[(\Hk^{1}_{\Lagr_{2d}})^{r} \times \Lagr_{2d}]$ and the base change map $u=g: \cA_{d}\to \Herm_{2d}$. By Proposition \ref{prop:Ad Herm smooth}, $g$ is smooth. We then have
\begin{equation}
\Phi_{\cM_{d}}^{!}g_{\Hk}^{*}[(\Hk^{1}_{\Lagr_{2d}})^{r} \times \Lagr_{2d}]=g^{*}_{\Sht}\Phi_{\Lagr_{2d}}^{!}[(\Hk^{1}_{\Lagr_{2d}})^{r} \times \Lagr_{2d}]\in \Ch_{0}(\Sht^{r}_{\cM_{d}}).
\end{equation}
Since $g_{\Hk}$ is smooth,  $g_{\Hk}^{*}[(\Hk^{1}_{\Lagr_{2d}})^{r} \times \Lagr_{2d}]=[(\Hk^{1}_{\cM_{d}})^{r} \times \Hk^{1}_{\cM_{d}}]$. The above equality then becomes \eqref{ShtM ShtLag}.
\end{proof}

Combined with Theorem \ref{th:compare KR}, we get a local description of the cycle class $[\cZ^{r}_{\cE}(a)]$:
\begin{cor}\label{c:Zloc}For any $(\cE,a)\in \cA_{d}(k)$ with image $\cQ=\coker(a)\in \Herm_{2d}(k)$, $[\cZ^{r}_{\cE}(a)]$ is the same as $[\cZ^{r}_{\cQ}]$ under the isomorphism \eqref{gSht}. In particular,
\begin{equation}
\deg[\cZ^{r}_{\cE}(a)]=\deg[\cZ^{r}_{\cQ}].
\end{equation}
\end{cor}

\subsection{Sheaves on $\Herm_{2d}$}
To describe the direct image complex $Rf_{*}\Qlbar$ on $\cA_{d}$, by the Cartesian diagram \eqref{eq: hitchin cartesian}, we first need to understand $R(\upsilon_{2d})_*\Qlbar$ on $\Herm_{2d}$.

\begin{lemma}\label{l:Lag push}
The perverse sheaf $R(\upsilon_{2d})_{*}\Qlbar$ on $\Herm_{2d}$ is canonically isomorphic to $(\Spr^{\Herm}_{2d})^{S_{d}}$ (see Proposition \ref{p:SH Spr}(2)). Here the $S_{d}$-action on $\Spr^{\Herm}_{2d}$ is the restriction of the Springer $W_{d}$-action.
\end{lemma}
\begin{proof} We have a Cartesian diagram
\begin{equation}
\xymatrix{\wt\Herm_{2d}\ar[d]^{\l_{2d}} \ar[r]^{\wt\e_{d}} & \wt\Coh_{d}(X') \ar[d]^{\pi^{\Coh}_{X',d}} \\
\Lagr_{2d} \ar[r]^{\e'_{d}} & \Coh_{d}(X')
}
\end{equation}
where $\e'_{d}$ sends $(\cQ,h_{\cQ},\cL)$ to $\cL$ and $\wt\e_{d}$ sends $(\cQ_{1}\subset \cdots\subset\cQ_{d}\subset \cdots\subset \cQ,h)$ to $(\cQ_{1}\subset \cdots\subset\cQ_{d})$. By proper base change
\begin{equation}
R\l_{2d*}\Qlbar\cong \e'^{*}_{d}\Spr_{d, X'}.
\end{equation}
In particular, $R\l_{2d*}\Qlbar$ carries an action of $S_{d}$. Moreover, the induced $S_{d}$-action  on  $\Spr^{\Herm}_{2d}\cong R(\upsilon_{2d})_{*}R\l_{2d*}\Qlbar$ is  the restriction  of the Springer $W_{d}$-action to $S_{d}$:  this can be easily checked over $\Herm^{\c}_{2d}$, and then the statement holds over $\Herm_{2d}$ since $\Spr^{\Herm}_{2d}$ is the middle extension from its restriction to $\Herm^{\c}_{2d}$ by Proposition \ref{p:SH Spr}(2). Since  $(R\l_{2d*}\Qlbar)^{S_{d}}\cong \e'^{*}_{d}(\Spr_{d,X'})^{S_{d}}\cong \Qlbar$, we conclude that $(\Spr^{\Herm}_{2d})^{S_{d}}\cong R(\upsilon_{2d})_{*}(R\l_{2d*}\Qlbar)^{S_{d}}\cong R(\upsilon_{2d})_{*}\Qlbar$, as desired.
\end{proof}

It is an elementary exercise to see that $\Ind_{S_{d}}^{W_{d}}\one$ decomposes into irreducible representations
\begin{equation}\label{rho i}
\Ind_{S_{d}}^{W_{d}}\one=\bigoplus_{i=0}^{d}\r_{i}
\end{equation}
where $\r_{i}=\Ind^{W_{d}}_{(\Z/2\Z)^{d}\rtimes (S_{i}\times S_{d-i})}(\chi_{i}\bt \one)$, and $\chi_{i}: (\Z/2\Z)^{d}\rtimes (S_{i}\times S_{d-i})\to \{\pm1\}$ is the character that is nontrivial on the first $i$ factors of $(\Z/2\Z)^{i}$, trivial on the rest and trivial on $S_{i}\times S_{d-i}$.  The decomposition also shows up in \cite[\S8.1.1]{YZ}.

Recall the notation $\Spr^{\Herm}_{2d}[\r]$ from Definition \ref{d:Spr isotypic}.
\begin{cor}\label{c:Lag decomp} There is a canonical decomposition 
\begin{equation}\label{decomp ph}
R(\upsilon_{2d})_{*}\Qlbar\cong\bigoplus_{i=0}^{d}\Spr^{\Herm}_{2d}[\r_{i}].
\end{equation}
\end{cor}
\begin{proof}
By Lemma \ref{l:Lag push} and Frobenius reciprocity,  we have
\begin{equation}
R(\upsilon_{2d})_{*}\Qlbar\cong \Hom_{W_{d}}(\Ind_{S_{d}}^{W_{d}}\one, \Spr^{\Herm}_{2d})=\Spr^{\Herm}_{2d}[\Ind_{S_{d}}^{W_{d}}\one].
\end{equation}
The desired decomposition then follows from \eqref{rho i}.
\end{proof}

\begin{defn}\label{def:Kd}
Define the graded perverse sheaf on $\Herm_{2d}(X'/X)$ 
\begin{equation}
\cK_{d}^{\mrm{Int}}(T) := \bigoplus_{i=0}^d \Spr^{\Herm}_{2d}[\r_{i}] T^{i}.
\end{equation}
\end{defn}

The fundamental class of the self-correspondence $\Hk^{1}_{\Lagr_{2d}}$ of $\Lagr_{2d}$ is viewed as a cohomological correspondence of the constant sheaf on $\Lagr_{2d}$ with itself. It induces an endomorphism (see notation from \cite[A.4.1]{YZ})
\begin{equation}
(\upsilon_{2d})_{!}[\Hk^{1}_{\Lagr_{2d}}]:R(\upsilon_{2d})_{*}\Qlbar\to R(\upsilon_{2d})_{*}\Qlbar.
\end{equation}

\begin{prop}\label{p:action Hk1} The action of $(\upsilon_{2d})_{!}[\Hk^{1}_{\Lagr_{2d}}]$ on $R(\upsilon_{2d})_{*}\Qlbar$ preserves the decomposition \eqref{decomp ph},  and it acts on $\Spr^{\Herm}_{2d}[\r_{i}]$ by multiplication by $(d-2i)$.
\end{prop}
\begin{proof}By Proposition \ref{p:SH Spr}(2),  $\Spr^{\Herm}_{2d}[\r_{i}]$ is the middle extension from its restriction to $\Herm^{\c}_{2d}$, it suffices to prove the same statement on $\Herm_{2d}^{\c}$. Let $I_d' \subset X_d' \times X'$ be the universal divisor $\{(D \in X_d', y \in X') \co y \in D\}$. The maps $\pr(D,y) = D$ and $q(D,y) = D-y+\sigma(y)$ define the incidence correspondence as in \cite[proof of Proposition 8.3]{YZ},
\begin{equation}\label{eq: incidence correspondence}
\begin{tikzcd}
& I_d' \ar[dl, "\pr"'] \ar[dr, "q"] \\
X_d'  & & X_d' 
\end{tikzcd}
\end{equation}
Now over $\Herm_{2d}^{\c}$, the map $\Herm_{2d}^{\c}\to X^{\c}_{d}$ is smooth, and $\Hk^{1}_{\Lagr_{2d}}$ is the pullback of \eqref{eq: incidence correspondence}. We therefore reduce to checking the statement for the action of $[I'_{d}]$ on the direct image sheaf of $\nu_{d}: X'_{d}\to X_{d}$, which is done in the proof of \cite[Proposition 8.3]{YZ}.
\end{proof}

\subsection{Lefschetz trace formula}
We shall give a slight generalization of the Lefschetz trace formula \cite[Proposition A.12]{YZ} expressing the intersection number of a cycle with the graph of Frobenius as a trace. Instead of the graph of Frobenius, we need to intersect along $\Phi_{M}:M^{r+1}\to M^{2r+2}$. Consider the following situation:
\begin{itemize} 
\item Let $S$ be an algebraic stack locally of finite type over $k=\F_q$. Assume $S$ can be stratified by locally closed substacks that are global quotients.
\item Let $M$ be a smooth equidimensional stack over $k=\F_q$ of dimension $N$ with a proper representable map $f: M\to S$. 
\item For $1\le i\le r$,  let $(\pr^{i}_{0},\pr^{i}_{1}): C_{i}\to M\times_{S}M$ be a self-correspondence of $M$ over $S$. Assume $\pr^{i}_{0}$ is proper and representable.
\end{itemize}
Form the Cartesian diagram
\begin{equation}\label{ShtC}
\xymatrix{\Sht_{C} \ar[r]\ar[d] & (\prod_{i=1}^{r}C_{i} ) \times M\ar[d]^{(\pr^{i}_{0}, \pr^{i}_{1})_{1\le i\le r}}\\
M^{r+1}\ar[r]^{\Phi_{M}} & M^{2r+2}\,.}
\end{equation}
Then $\Sht_{C}$ decomposes as
\begin{equation}
\Sht_{C}=\coprod_{s\in S(k)}\Sht_{C}(s)
\end{equation}
where $\Sht_C(s)$ is the fibered product of $\Sht_C \rightarrow S(k)$ against $s \in S(k)$. For $s \in S(k)$, we write $\Sht_{C}|_s$ for the fibered product 
\[
\begin{tikzcd}
\Sht_C|_s \ar[r] \ar[d] & \Sht_C \ar[d] \\
\Spec k \ar[r, "s"] & S(k)
\end{tikzcd}
\]
Then $\Sht_C|_s \rightarrow \Sht_C(s)$ is a torsor for $\Aut(s)$, and in particular a finite \'{e}tale cover. 

Suppose we are given cycle classes
\begin{equation}
\z_{i}\in \Ch_{N}(C_{i}), \quad 1\le i\le r.
\end{equation}
The cycle class  $\cla(\z_{i})\in \hBM{2N}{C_{i},\Qlbar(-N)}$ is viewed as a cohomological correspondence between the constant sheaf on $M$ and itself. Therefore it induces an endomorphism of $Rf_{!}\Qlbar$ which we denote by $f_{!}\cla(\z_{i})$.

\begin{prop}\label{p:Lef} For each $s\in S(k)$ we have
\begin{equation*}
\deg((\Phi_{M}^{!}(\z_{1}\times\cdots\times\z_{r} \times [M] ))|_{\Sht_{C}|_s})=\Tr(f_{!}\cla(\z_{1})\c\cdots\c f_{!}\cla(\z_{r})\c\Frob_{s}, (Rf_{!}\Qlbar)_{s}).
\end{equation*}
\end{prop}
\begin{proof} We first prove the formula when $S$ is a scheme of finite type.  In this case $M,C_{i}$ are also schemes of finite type over $k$. Let $C=C_{1}\times_{M}\times\cdots\times_{M}C_{r}$ be the composition correspondence, with maps $\pr_{i}: C\to M$ for $0\le i\le r$.  Consider the diagram where all squares are Cartesian
\begin{equation}
\xymatrix{\Sht_{C}\ar[d] \ar[r] & C \ar[r]\ar[d]^{(\pr_{0},\cdots, \pr_{r}, \pr_r)} & (\prod_{i=1}^{r}C_{i} ) \times M\ar[d]^{\prod (\pr^{i}_{0}, \pr^{i}_{1}) \times \Delta}\\
M^{r+1} \ar[d]_{\pr_{0}} \ar[r]^{\Phi_{1}} & M^{r+2}\ar[d]^{(\pr_{0}, \pr_{r+1})} \ar[r]^{\Phi_{2}} & M^{2r+2}\\
M\ar[r]_-{(\Id,\Frob_{M})} & M\times M}
\end{equation}
Here
\begin{align*}
\Phi_{1}(\xi_{0},\cdots,  \xi_{r-1}, \xi_r) & =(\xi_{0},\cdots, \xi_r, \Frob_{M}(\xi_{0})),\\
\Phi_{2}(\xi_{0},\cdots, \xi_r,\xi_{r+1}) & =(\xi_{0},\xi_{1},\xi_{1},\cdots,\xi_{r},\xi_{r}, \xi_{r+1}).
\end{align*}
We have $\Phi_{M}=\Phi_{2}\c\Phi_{1}$. Let $\z=\Phi_{2}^{!}(\z_{1}\times\z_{2}\times\cdots\times\z_{r} \times [M])\in \Ch_{N}(C)$. 

On the one hand, by the transitivity of the Gysin maps, 
\begin{equation}\label{Phizz}
\Phi^{!}_{M}(\z_{1}\times\cdots\times\z_{r} \times [M] )=\Phi_{1}^{!}(\z)=(\Id,\Frob_{M})^{!}(\z).
\end{equation}
Applying the Lefschetz trace formula \cite[Proposition  A.12]{YZ}, we get
\begin{equation}\label{deg tr z}
\deg((\Id,\Frob_{M})^{!}(\z)|_{\Sht_{C}|_s})=\Tr(f_{!}\cla(\z)\c \Frob, (Rf_{!}\Qlbar)_{s}).
\end{equation}
One the other hand, by a diagram chase, we see that $\cla(\z)$ is the composition of the cohomological correspondences $\cla(\z_{i})$ ($1\le i\le r$), hence $f_{!}\cla(\z)\in \End(Rf_{!}\Qlbar)$ is the composition of $f_{!}\cla(\z_{1})\c f_{!}\cla(\z_{2})\c\cdots\c f_{!}\cla(\z_{r})$. Combining this fact with \eqref{Phizz} and \eqref{deg tr z} we get the desired formula.

Now consider the general case where $S$ is a stack locally of finite type over $k$ and we aim to prove the formula for $s\in S(k)$. We claim that there exists a scheme $S'$ of finite type over $k$ and a smooth map $u: S'\to S$ such that $u(S'(k))$ contains $s$.  Indeed, pick any smooth map $u_{1}: S_{1}\to S$ with $S_{1}$ a scheme of finite type over $k$, such that  $s$ is contained in the image of $u$. Let $s_{1}\in S_{1}(\F_{q^{m}})$ be a point that maps to $s$. Let $(S_{1}/S)^{m}$ be the $m$-fold fibered product of $S_1$ over $S$, based changed to $\ol{k}$. We equip  $(S_{1}/S)^{m}$ with the $\Frob$-descent datum given by $(x_{1},\cdots, x_{m})\mapsto (\Frob(x_{m}), \Frob(x_{1}),\cdots, \Frob(x_{m-1}))$. This gives a descent of $(S_{1}/S)^{m}$ to a scheme $S'$ over $k$ equipped with a map $u:S'\to S$ which is still smooth since $u_{1}$ is. Now $s_{1}$ gives rise to a $k$-point $s'=(s_{1},\Frob(s_{1}),\cdots, \Frob^{m-1}(s_{1}))\in S'(k)$ such that $u(s')=s$. 

Let $M'=M\times_{S}S'$, $C'_{i}=C_{i}\times_{S}S'$ and let $u_{C_{i}}:C'_{i}\to C_{i}$ be the projection. Define $\Sht'_{C}$ using the analog of the diagram \eqref{ShtC} with $M$ and $C_{i}$ replaced by $M'$ and $C'_{i}$. Then $\Sht'_{C}\cong \Sht_{C}\times_{S(k)}S'(k)$. For $s'\in S'(k)$ such that $u(s')=s$, we get an isomorphism $\Sht'_{C}|_{s'}\isom \Sht_{C}|_{s}$. Let $\z'_{i}=u_{C_{i}}^{*}\z_{i}$. Now we apply Proposition \ref{p:base change} to the diagram \eqref{ShtC} along the base change map $u:S'\to S$ to get
\begin{equation}
u_{\Sht}^{*}\Phi_{M}^{!}(\z_{1}\times\cdots\times\z_{r} \times [M])\cong \Phi_{M'}^{!}(\z'_{1}\times\cdots\times\z'_{r} \times [M])\in \Ch_{0}(\Sht'_{C}).
\end{equation}
Restricting to $\Sht'_{C}|_{s'}\cong \Sht_{C}|_s$ and taking degrees we get
\begin{equation}\label{degSS'}
\deg(\Phi_{M}^{!}(\z_{1}\times\cdots\times\z_{r}  \times [M] )|_{\Sht_{C}|_s})=\deg(\Phi_{M'}^{!}(\z'_{1}\times\cdots\times\z'_{r}  \times [M])|_{\Sht'_{C}|_{s'}}).
\end{equation}
On the other hand, letting $f':M'\to S'$, by smooth base change we have
\begin{equation}\label{trSS'}
\Tr(f_{!}\cla(\z_{1})\c\cdots\c f_{!}\cla(\z_{r})\c\Frob_{s}, (Rf_{!}\Qlbar)_{s}) = \Tr(f'_{!}\cla(\z'_{1})\c\cdots\c f'_{!}\cla(\z'_{r})\c\Frob_{s'}, (Rf'_{!}\Qlbar)_{s'}).
\end{equation}
Since the right sides of \eqref{degSS'} and \eqref{trSS'} are equal by the scheme case that is already proven, the left sides of \eqref{degSS'} and \eqref{trSS'} are also equal, proving the proposition in general.
\end{proof}

Recall the graded perverse sheaf $\cK^{\mrm{Int}}_{d}(T)$ on $\Herm_{2d}$ from Definition \ref{def:Kd}. 

\begin{cor}\label{cor: trace of correspondence} Let $(\cQ,h_{\cQ})\in \Herm_{2d}(k)$. We have
\begin{eqnarray*}
\deg[\cZ^{r}_{\cQ}]&=&\Tr( [\Hk^{1}_{\Lagr_{2d}}]^r\c\Frob, (R(\upsilon_{2d})_{*}\Qlbar)_{\cQ}) = \sum_{i=0}^d (d-2i)^r \Tr(\Frob, \Spr^{\Herm}_{2d}[\r_{i}]_{\cQ})  \\
&=& \frac{1}{(\log q)^r} \left(\frac{d}{ds}\right)^r\Big|_{s=0} \left(q^{ds}\Tr(\Frob, \cK_d^{\mrm{Int}}(q^{-2s})_{\cQ})\right).
\end{eqnarray*}
\end{cor}
\begin{proof}
The first equality is an application of Proposition \ref{p:Lef} to the case $S=\Herm_{2d}$, $M=\Lagr_{2d}$, $C_{i}=\Hk^{1}_{\Lagr_{2d}}$ and $\z_{i}=[\Hk^{1}_{\Lagr_{2d}}]$.  The second equality follows from Proposition \ref{p:action Hk1}. The third one is a direct calculation.
\end{proof}

Combining Corollary \ref{cor: trace of correspondence} with Corollary \ref{c:Zloc} we get:
\begin{cor}\label{cor: KR degree as trace of Frob} Let $(\cE,a)\in \cA_{d}(k)$ with image $(\cQ,h_{\cQ})\in \Herm_{2d}(k)$. Then we have 
\[
\deg [\Cal{Z}^{r}_{\Cal{E}}(a)] = \frac{1}{(\log q)^r} \left(\frac{d}{ds}\right)^r\Big|_{s=0} \left(q^{ds}\Tr(\Frob, \cK_d^{\mrm{Int}}(q^{-2s})_{\cQ})\right).
\]
\end{cor}

\subsection{Symmetry} This subsection is not used in the proof of the main theorem. The graded perverse sheaf $\cK^{\Int}_{d}(T)$ has a palindromic symmetry that we spell out. First, the \'etale double covering $\nu:X'\to X$ gives a local system $\y_{X'/X}$ on $X$ with monodromy in $\pm1$.  It induces a local system $\y_{d}$ on $X_{d}$ with monodromy in $\pm1$: its stalk at a divisor $x_{1}+\cdots+x_{d}\in X_{d}(\ov{k})$ is $\ot_{i=1}^{d}(\y_{X'/X})_{x_{i}}$. Let 
\begin{equation}
\y^{\Herm}_{2d}:=s^{\Herm*}_{2d}\y_{d},
\end{equation}
where $s^{\Herm}_{2d}: \Herm_{2d}\to X_{d}$ is the support map. This is a rank-one local system on $\Herm_{2d}$ with monodromy in $\pm1$.

\begin{lemma} We have a canonical isomorphism of perverse sheaves on $\Herm_{2d}$:
\begin{equation}
\Spr^{\Herm}_{2d}[\r_{d}]\cong\y^{\Herm}_{2d}.
\end{equation}
\end{lemma}
\begin{proof}
By Proposition \ref{p:SH Spr}(2), $\Spr^{\Herm}_{2d}[\r_{d}]$ is the middle extension of its restriction to the open dense substack $\Herm_{2d}^{\c}$ (preimage of $X^{\c}_{d}$). The same is true for $\y_{d}$  because it is a local system and $\Herm_{2d}$ is smooth. Therefore it suffices to check their equality over $\Herm_{2d}^{\c}$, over which both are obtained by pushing out the $W_{d}$-torsor $(X'^{d})^{\c}\to X_{d}^{\c}$ along the character $\chi_{d}: W_{d}\to \{\pm1\}$.
\end{proof}

\begin{lemma} There is an isomorphism of graded perverse sheaves on $\Herm_{2d}$
\begin{equation}
T^{d}\cK_d^{\Int}(T^{-1})\cong \y^{\Herm}_{2d}\ot \cK^{\Int}_{d}(T).
\end{equation}
\end{lemma}
\begin{proof}
The equality amounts to 
\begin{equation}
\Spr^{\Herm}_{2d}[\r_{d-i}]\cong \y^{\Herm}_{2d}\ot\Spr^{\Herm}_{2d}[\r_{i}].
\end{equation}
Both sides are middle extensions from $\Herm^{\c}_{2d}$ by Proposition \ref{p:SH Spr}(2), over which they correspond to representations $\r_{d-i}$ and $\chi_{d}\ot \r_{i}$ of $W_{d}$. By definition, 
\begin{equation}
\chi_{d}\ot\r_{i}\cong \chi_{d}\ot\Ind^{W_{d}}_{W_{i}\times W_{d-i}}(\chi_{i}\bt\one).
\end{equation}
Inserting $\chi_{d}|_{W_{i}\times W_{d-i}}\cong \chi_{i}\bt \chi_{d-i}$ to the right side above gives 
\begin{equation}
\Ind^{W_{d}}_{W_{i}\times W_{d-i}}(\chi_{d}|_{W_{i}\times W_{d-i}}\ot(\chi_{i}\bt\one))\cong \Ind^{W_{d}}_{W_{i}\times W_{d-i}}(\one\bt\chi_{d-i})\cong \r_{d-i}.
\end{equation}
\end{proof}

\begin{lemma}
If $(\cQ,h_{\cQ})\in \Herm_{2d}(k)$ is the image of some $(\cE,a)\in \cA_{d}(k)$, then $\Tr(\Frob, \y^{\Herm}_{2d}|_{\cQ})=1$.
\end{lemma}
\begin{proof} If $X'/X$ is split, then $\y^{\Herm}_{2d}$ is trivial, and there is nothing to prove. Below we assume $X'/X$ is nonsplit. The local system $\y_{d}$ on $X_{d}$ is pulled back from a local system $\y_{\Pic}$ on $\Pic_{X}$ via the Abel-Jacobi map $\mrm{AJ}_{d}:X_{d}\to \Pic^{d}_{X}\subset \Pic_{X}$. The Frobenius trace function of $\y_{\Pic}$ is the id\`ele class character 
\[
\y_{F'/F}: F^{\times}\bs\BA^{\times}_{F}/\wh\cO^{\times}=\Pic_{X}(k)\to\{\pm1\}
\]
trivial on the image of $\Nm_{X'/X}: \Pic_{X'}(k)\to \Pic_{X}(k)$. Denote by $\det_{X}(\cQ)$ the image of $\cQ$ under $\Herm_{2d}\to X_{d}\xr{\mrm{AJ}_{d}} \Pic^{d}_{X}$. We have $\y^{\Herm}_{2d}|_{\cQ}\cong \y_{\Pic}|_{\det_{X}(\cQ)}$ as $\Frob$-modules. Now $(\cQ,h_{\cQ})$ comes from $(\cE,a)$, which implies
\begin{equation*}
{\det}_{X}(\cQ)\cong \Nm_{X'/X}(\det\cE)^{-1}\ot \om_{X}^{\ot n}.
\end{equation*}
By \cite[p.291, Theorem 13]{Weil}, $\om_{X}$ is a square in $\Pic_{X}(k)$, hence $\y_{F'/F}(\om_{X}^{\ot n})=1$. Since $\y_{F'/F}(\Nm_{X'/X}(\det\cE))=1$, we see that $\y_{F'/F}(\det_{X}(\cQ))=1$, hence $\Tr(\Frob, \y^{\Herm}_{2d}|_{\cQ})=1$.
\end{proof}

\begin{cor}\label{cor:FE Int}
Let $(\cE,a)\in \cA_{d}(k)$ with image  $(\cQ,h_{\cQ})\in \Herm_{2d}(k)$. Then $s\mapsto q^{ds}\Tr(\Frob, \cK_d^{\mrm{Int}}(q^{-2s})_{\cQ})$ is an even function in $s$. In particular, its odd order derivatives at $s=0$ vanish.
\end{cor}
By Corollary \ref{cor: KR degree as trace of Frob}, this implies $\deg[\cZ^{r}_{\cE}(a)]=0$ for $r$ odd. However, we know from Lemma \ref{lem: shtuka empty parity} that $\Sht^{r}_{U(n)}=\vn$ when $r$ is odd, which implies $\cZ^{r}_{\cE}(a)=\vn$.

\part{The comparison}

\section{Matching of sheaves}

\subsection{Recap} 
Let 
\begin{equation}\label{eq: normalized Eisenstein coefficients}
\wt{E}_{a}(m(\cE),s,\Phi)   = E_a(m(\cE), s, \Phi) \cdot \chi(\det(\cE))^{-1} q^{\deg(\cE)( s-\frac{n}{2})+\frac{1}{2}n^2\deg(\omega_X) } \cdot  \sL_n(s) = \Den(q^{-2s}, (\cE,a))
\end{equation}
where the notation is as in Theorem \ref{th:Eis Den}, being a renormalization of the $a^{\mrm{th}}$ Fourier coefficient of $E_a(m(\cE), s, \Phi)$.

We emphasize that, in keeping with \S \ref{ssec: notation}, $X$ is proper and $\nu \co X' \rightarrow X$ is a finite \'{e}tale double cover (possibly trivial). 

\begin{thm}\label{thm: main} Keep the notations above. Let $(\cE,a)\in \cA_{d}(k)$. Then we have 
\begin{equation}\label{eq: SWF}
\deg [\cZ_{\cE}^r(a)] = \frac{1}{(\log q)^r}  \left(\frac{d}{ds}\right)^r\Big|_{s=0} \left( q^{ds} \wt{E}_{a}(m(\cE),s,\Phi) \right). 
\end{equation}
\end{thm}

In the previous parts, we have found sheaves on $\cA_d$ which correspond to the two sides of \eqref{eq: SWF}, in the sense of the function-sheaf dictionary. Let us summarize the situation. 

 On the analytic side, we proved a formula expressing the non-singular Fourier coefficient of the Siegel--Eisenstein series in terms of the Frobenius trace of a graded virtual perverse sheaf $\cK_{d}^{\Eis}(T)$ on $\Herm_{2d}(X'/X)$.

\begin{thm}[Combination of Theorems \ref{th:Eis Den} and \ref{th:Den Wd}] Let $(\cE,a)\in \cA_{d}(k)$. Then we have 
\begin{equation}\label{eq: comparison 1}
\wt{E}_{a}(m(\cE),s,\Phi) =  \Tr(\Frob,\cK_{d}^{\Eis}(q^{-2s})_{\cQ}).
\end{equation}
\end{thm}

On the geometric side,  in Corollary \ref{cor: KR degree as trace of Frob}, we found a formula expressing the degree of the special $0$-cycle in terms of $r^{\mrm{th}}$ derivative of the Frobenius trace of another graded perverse sheaf $\cK^{\Int}_{d}(T)$ on $\Herm_{2d}(X'/X)$, repeated below: 
\begin{equation}\label{eq: comparison 2}
\deg [\Cal{Z}^{r}_{\Cal{E}}(a)] = \frac{1}{(\log q)^r} \left(\frac{d}{ds}\right)^r\Big|_{s=0} \left(q^{ds}\Tr(\Frob, \cK_d^{\mrm{Int}}(q^{-2s})_{\cQ})\right).
\end{equation}

\subsection{Proof of the main theorem} Comparing \eqref{eq: comparison 1} and \eqref{eq: comparison 2}, we see that in order to prove Theorem \ref{thm: main}, it remains to match the graded sheaves $\cK_d^{\mrm{Int}}(T)$ and $\cK_{d}^{\Eis}(T)$ on $\Herm_{2d}(X'/X)$.

\begin{prop}
We have $ \Cal{K}^{\Int}_{d}(T) \cong \cK_{d}^{\Eis}(T)$ as graded perverse sheaves  on $\Herm_{2d}(X'/X)$.
\end{prop}
\begin{proof} Both sides can be written as $\Spr^{\Herm}_{2d}[\r]$ for some graded virtual representation $\r$ of $W_{d}$. By definition (Definition \ref{def:Kd}), the sheaf $\cK^{\Int}_{d}(T)$ corresponds to 
\[
\r_{\cK_{d}^{\Int}}(T)=\sum_{i=0}^{d} \Ind_{W_i \times W_{d-i}}^{W_d} (\chi_i \boxtimes 1) T^{i}.
\]
We calculate the (a priori virtual) representation of $W_d$ which corresponds under Springer theory to the $\cK_{d}^{\Eis}(T)$ from Definition \ref{def: K_Eis}. The operation $R\oll{f}_{i,!}R\orr{f}^{*}_{i}$ corresponds under Springer theory to $\Ind_{W_{d-i} \times S_i}^{W_d}$ and $\frP_{d-i}(T)$ corresponds to $\sum_{j=0}^{d-i} (-1)^j \Ind_{W_j \times W_{d-i-j}}^{W_{d-i}} (\sgn_j \boxtimes \mbf{1}) T^j$ (cf. Definition \ref{d:Spr isotypic}  for $\HSpr_{d}$). Hence $
\cK_{d}^{\Eis}(T)$ corresponds to 
\[
\sum_{i=0}^d \Ind_{W_{d-i} \times S_i}^{W_d} \left( \sum_{j=0}^{d-i}(-1)^j  \Ind_{W_j \times W_{d-i-j}}^{W_{d-i}}(\sgn_j \boxtimes \one_{d-i-j}) \boxtimes \one_i T^{i+j} \right) .
\]
After re-indexing, we have
\[
\r_{\cK_{d}^{\Eis}}(T)=\sum_{i=0}^{d} \sum_{j=0}^{i}(-1)^j\Ind_{S_{i-j} \times W_j \times W_{d-i}}^{W_d} (\one \bt \ol{\sgn}_{j} \bt \one ) T^{i}.
\]
The desired statement then follows from the Lemma below (whose notation has been re-indexed) by comparing each coefficient.
\end{proof}

\begin{lemma} We have the identity of virtual representations of $W_d$: 
\[
\chi_d = \sum_{j=0}^{d} (-1)^j \Ind_{S_{d-j} \times W_j}^{W_d}(1 \boxtimes \ol{\sgn}_{j}).
\]
\end{lemma}
\begin{proof}
We will prove this by comparing traces of an arbitrary element $g \in W_d$.  For $g \in W_d$, 
\begin{equation}\label{eq: trace of induced}
\Tr(g, \Ind_{S_{d-j} \times W_j}^{W_d} (1 \times \ol{\sgn}_j))=\sum_{\substack{w\in W_{d}/(S_{d-j} \times W_j) \\ w^{-1}gw\in S_{d-j} \times W_j} } \ol{\sgn}_j(g'').
\end{equation}
Here, when $w^{-1}gw\in S_{d-j} \times W_j$, we write $w^{-1}gw=(g',g'')$ for $g'\in S_{d-j}$ and $g''\in W_{j}$. 

Identify $W_d$ with the group of permutations of $\{\pm 1, \ldots, \pm d\}$ that commute with the involution $\s$ exchanging $j \leftrightarrow -j$ for all $1 \leq j \leq d$. The subgroup $S_{d-j}\times W_{j}$ is the stabilizer of $\{1,2,\cdots, d-j\}$. Therefore the coset space $W_{d}/(S_{d-j} \times W_j)$ is in natural bijection with subsets $J\subset \{\pm 1, \ldots, \pm d\}$ such that $|J|=d-j$ and $J\cap (-J)=\vn$. Let $\frJ_{g}$ be the set of $J\subset \{\pm 1, \ldots, \pm d\}$ such that $|J|=d-j$, $J\cap (-J)=\vn$ and $gJ=J$.  Let $g''_{J}$ be the permutation of $g$ on $\{\pm 1, \ldots, \pm d\}\bs (J\cup (-J))$. Combining this with \eqref{eq: trace of induced}, we obtain 
\begin{equation}
\sum_{j=0}^{d}(-1)^{j}\Tr(g, \Ind_{S_{d-j} \times W_j}^{W_d} (1 \times \ol{\sgn}_j))=\sum_{J\in \frJ_{g}}(-1)^{d-|J|}\ov\sgn(g''_{J}).
\end{equation}




For any $g \in W_d$,  the cycle decomposition of $g$ can be grouped into a decomposition $g = g_1 \ldots g_r$ (unique up to reordering) where $g_i$ is one of the two forms:
\begin{itemize}
\item (positive bicycle) $g_i$ is a product of two disjoint cycles $c_{i}\s(c_{i})$ (in particular, no two elements appearing in $c_{i}$ are negatives of each other). 
\item (negative cycle) $g_i$ is a single cycle invariant under the involution $\s$. 
\end{itemize}
Let $C^{+}_{g}$ be the set of cycles of $g$ that are part of a positive bicycle (i.e., $C_g^+$ contains both $c_{i}$ and $\s(c_{i})$ for each positive bicycle $g_{i}$). For any $x\in W_{d}$ we denote by $\un x\subset \{\pm1,\cdots, \pm d\}$ the set of elements that are not fixed by $x$. For a cycle $c$ we let $|c|$ be its length. From this description we see that $J\in \frJ_{g}$ if and only if $J$ is a union of $\un c$ for a subset of cycles  $c\in C^{+}_{g}$. In other words, consider the set $\frI_{g}$ of subsets $I\subset C^{+}_{g}$ such that $\un I$ is disjoint from $\s(\un I)$. Then we have a bijection $\frI_{g}\isom \frJ_{g}$ sending $I\in \frI_{g}$ to $J:=\cup_{c\in I}\un c$.

For $I\in \frI_{g}$, let $g'_{I}$ be the product of  $g_{i}$ such that $g_{i}$ contains a cycle in common with $I$; let $g''_{I}$ be the product of the remaining $g_{i}$'s. The above discussion allows us to rewrite
\begin{equation}
\sum_{j=0}^{d}(-1)^{j}\Tr(g, \Ind_{S_{d-j} \times W_j}^{W_d} (1 \times \ol{\sgn}_j))=\sum_{I\in \frI_{g}}(-1)^{\sum_{c\in C^{+}_{g}}|c|}\ov\sgn(g''_{I}).
\end{equation}

This sum factorizes as a product over the $g_i$ with individual factors as follows: 
\begin{itemize}
\item For a positive bicycle $g_i = c_i \s(c_{i})$, the local factor is the sum of three contributions, corresponding to whether $c_{i}\in I$, $\s(c_{i})\in I$ or neither  $c_{i}$ nor $\s(c_{i})$ is in $I$. The first two cases each contribute $1$. The last case leads to a contribution of $(-1)^{|c_i|} \sgn(c_i)=-1$. The total contribution of the factor corresponding to a positive bicycle $g_{i}$ is therefore $1+1+(-1)=1$.

\item For each negative cycle $g_i$, since it always appears in $g''_{I}$, its contribution is  $(-1)^{|g_{i}|/2}\ol{\sgn}(g_i)$. Let  $\ov g_{i}$ be the image of $g_{i}$ in $S_{d}$, which is a cycle of length $|g_{i}|/2$. Then $(-1)^{|g_{i}|/2}\ol{\sgn}(g_i)=(-1)^{|\ov g_{i}|}\sgn(\ov g_{i})=-1$. Therefore the contribution of the factor corresponding to a negative cycle $g_{i}$ is $-1$.
\end{itemize}
Summarizing, we have found 
\begin{equation}\label{Trg factor}
\sum_{j=0}^{d} (-1)^{j}\Tr(g, \Ind_{S_{d-j} \times W_j}^{W_d} (1 \times \ol{\sgn}_j)) = \left( \prod_{g_i \text{ positive }} 1 \right) \cdot  \left( \prod_{g_i \text{ negative} }(-1) \right).
\end{equation}

On the other hand,  we have
\begin{equation}
\chi_{d}(g_{i})=\begin{cases} 1, & g_{i} \mbox{ is positive;} \\ -1, &g_{i} \mbox{ is negative.} \end{cases}
\end{equation}
Indeed, if $g_{i}=c_{i}\s(c_{i})$ is positive, then we have $\chi_{d}(c_{i})=\chi_{d}(\s(c_{i}))=1$ because both $c_{i}$ and $\s(c_{i})$ can be conjugated into $S_{d}$. If $g_{i}$ is negative, then up to conjugacy we may assume $g_{i}$ is the cyclic permutation $(1,2,\cdots, m, -1,\cdots, -m)$ for some $1\le m\le d$. Then $g_{i}=(1,-1)(1,2,\cdots, m) (-1,-2,\cdots, -m)$, from which we see $\chi_{d}(g_{i})=-1$.

We conclude that the right side of \eqref{Trg factor} is $\prod_{}\chi_{d}(g_{i})=\chi_{d}(g)$. This completes the proof.

%

\end{proof}

\subsection{The split case $X'=X\coprod X$}
We make our result more explicit in the split case $X'=X\coprod X$.
  
  On the analytic side in \S\ref{ss:Eis}, the group $H_n=\GL_{2n,F}$ and $P_n$ is the standard parabolic corresponding to the partition $(n,n)$, with Levi $M_n\simeq \GL_{n,F}\times \GL_{n,F}$. We then have the degenerate principal series
  $$
I_n(s)=\Ind_{P_n(\BA)}^{H_n(\BA)}( |\cdot|_{F}^{s+n/2}\times |\cdot|_{F}^{-s-n/2}),\quad s\in \mathbb{C}.
$$
Let $\cE=(\cE_1,\cE_2)\in \Bun_{M_n}(k)\simeq\Bun_{\GL_n}(k)\times \Bun_{\GL_n}(k)$, and let $a:\cE_1\to\cE_2^\vee$ be an injective map of $\cO_X$-modules. Then by \S\ref{ss:FC} the Siegel--Eisenstein series has a well-defined $a^{\mrm{th}}$ Fourier coefficient $E_{a}(m(\cE),s,\Phi)$ at $(\cE_1,\cE_2)$. By Theorem \ref{th:Eis Den}  and \ref{th:Eis Den torsion} we have
  \begin{eqnarray*}
E_{a}(m(\cE),s,\Phi)=q^{-(\deg(\cE_1)+\deg(\cE_2))( s-n/2)-\frac{1}{2}n^2\deg\omega_X }\sL_n(s) ^{-1} \Den(q^{-2s},  \cE_2^\vee/\cE_1),
\end{eqnarray*}
where $\sL_n(s)=\prod_{i=1}^n\zeta_F(i+2s)$ and, for a torsion $\cO_X$-module $\cQ$, the density polynomial is given by
$$
 \Den(T,  \cQ)=\sum_{0 \subset \cI_1 \subset \cI_2\subset \cQ} T^{\dim_{k}\cI_1+\dim_k \cQ/\cI_2}\prod_{v\in |X|}\fm_{v}(t_v(\cI_2/\cI_1);T^{\deg(v)}).
$$
Here see \eqref{def:maT} for $\fm_{v}(t_v;T)$.
The normalized Fourier coefficient \eqref{eq: normalized Eisenstein coefficients} is
\[
\wt{E}_{a}(m(\cE),s,\Phi)  =  \Den(q^{-2s},  \cE_2^\vee/\cE_1).
\]

Next we come to the geometric side. We have a natural partition
$$
(X')^r=\coprod_{\mu\in \{\pm 1\}^r} X^r.
$$
The moduli of hermitian shtukas $\Sht_{U(n)}^r$ defined in \S\ref{sec: unitary shtukas} is then partitioned into 
$$
\Sht_{U(n)}^r=\coprod_{\mu\in \{\pm 1\}^r} \Sht_{U(n)}^{\mu},
$$
and there is a natural isomorphism
$$
 \Sht_{U(n)}^{\mu}\simeq \Sht_{\GL_n}^{\mu}.
$$
Here we recall that $ \Sht_{\GL_n}^{\mu}$ is the moduli of shtukas for $\GL_{n}$, cf. \cite{YZ}, whose
$S$-points are given by the groupoid of the following data: 
\begin{enumerate}
\item $x_i \in X(S)$ for $i = 1, \ldots, r$.
\item  $\Cal{F}_0, \ldots, \Cal{F}_n\in \Bun_{\GL_n} ( S)$.
\item  An elementary modification $f_i \co \Cal{F}_{i-1} \dashrightarrow
 \Cal{F}_i$ at the graph of $x_i$, which is of upper of length 1 if $\mu_i=+1$ and of lower of length 1 if $\mu_i=-1$.

\item An isomorphism $\varphi \co \Cal{F}_r \cong \ft \Cal{F}_0$. 
\end{enumerate}
In particular, $\Sht_{\GL_n}^{\mu}$ is empty unless $\sum_{i=1}^r \mu_i=0$.

For the special cycle $\cZ_\cE^r $ (cf. Definition \ref{def: Z}) associated to $\cE=(\cE_1,\cE_2)$ above, we have a partition
$$
\cZ_\cE^r =\coprod_{\mu\in \{\pm 1\}^r} \cZ^\mu_\cE,
$$
where an object in $\cZ^\mu_\cE(S)$ is an object as above in $ \Sht_{\GL_n}^{\mu}(S)$ together with maps
\begin{equation}\label{eq: t1t2}
\Cal{E} _1\boxtimes \cO_S \xrightarrow{t_i^{(1)}} \Cal{F}_i  \xrightarrow{t_i^{(2)}}  \Cal{E}_2^{\vee} \boxtimes \cO_S,\quad i = 1, \ldots, r,
\end{equation}
 such that the diagram commutes
\[
\begin{tikzcd}
\Cal{E}_1 \boxtimes \cO_S \ar[d, "t_0^{(1)}"] \ar[r, equals] &  \Cal{E}_1 \boxtimes \cO_S \ar[r, equals]  \ar[d, "t_1^{(1)}"]& \ldots \ar[d] \ar[r, equals] & \Cal{E}_1\boxtimes \cO_S  \ar[r, "\sim"] \ar[d, "t_{r}^{(1)}"] & \ft (\Cal{E}_1 \boxtimes \cO_S) \ar[d, "\ft t_0^{(1)}"] \\
\Cal{F}_0  \ar[d, "t_0^{(2)}"] \ar[r, dashed, "f_0"] & \Cal{F}_1 \ar[r, dashed, "f_1"]   \ar[d, "t_1^{(2)}"]&  \ldots \ar[d] \ar[r, dashed, "f_r"] &  \Cal{F}_{r}\ar[d, "t_{r}^{(2)}"]  \ar[r, "\sim"] & \ft \Cal{F}_0  \ar[d, "\ft t_{0}^{(2)}"]\\
\Cal{E}_2^\vee \boxtimes \cO_S \ar[r, equals] &  \Cal{E}_2^\vee  \boxtimes \cO_S \ar[r, equals] & \ldots  \ar[r, equals] & \Cal{E}_2^\vee  \boxtimes \cO_S\ar[r, "\sim"] & \ft ( \Cal{E}_2^\vee \boxtimes \cO_S)
\end{tikzcd}
\]
Let $a:\cE_1\to\cE_2^\vee$ be a map of $\cO_X$-modules.  Then $\cZ^\mu_\cE(a)$ is the open-closed subscheme of $\cZ^\mu_\cE$ such that the common composition \eqref{eq: t1t2} is equal to $a \boxtimes \Id_{\cO_S} $.

For an injective $a:\cE_1\to\cE_2^\vee$, our \S\ref{s:int prob} shows that $\cZ_{\cE}^\mu(a)$ is proper over $\Spec k$ and  defines a class $[\cZ_{\cE}^\mu(a)]\in \Ch_0( \cZ_{\cE}^\mu(a))$ for each $\mu\in \{\pm 1\}^r$.  Then our main Theorem asserts 
\begin{equation}\label{eq: SWF spl}
\sum_{\mu\in \{\pm 1\}^r} \deg [\cZ_{\cE}^\mu(a)] = \frac{1}{(\log q)^r}  \left(\frac{d}{ds}\right)^r\Big|_{s=0}  \left( q^{ds} \wt{E}_{a}(m(\cE),s,\Phi) \right),
\end{equation}
where $d=-(\chi(X,\cE_1)+\chi(X,\cE_2))$. We remark that $\deg [\cZ_{\cE}^\mu(a)]$ is not independent of $\mu \in \{\pm 1\}^r$, even if we restrict our attention to those $\mu$ with $\sum_{i=1}^r \mu_i = 0$.

\bibliographystyle{amsalpha}
\bibliography{Bibliography}

\end{document}